\pdfoutput=1

\documentclass[12pt]{amsart}
\usepackage[margin=1in]{geometry} 

\usepackage{amsmath,amssymb,amsthm,bm,colonequals,graphicx,mathrsfs,mathtools,microtype,stmaryrd}
\usepackage[shortlabels]{enumitem}
\usepackage[colorinlistoftodos]{todonotes}

\numberwithin{equation}{section}

\renewcommand{\Re}{\operatorname{Re}}
\renewcommand{\Im}{\operatorname{Im}}

\theoremstyle{plain}

\newtheorem{theorem}{Theorem}[section] 
\newtheorem*{theorem*}{Theorem}

\newtheorem{lemma}[theorem]{Lemma} 

\newtheorem{proposition}[theorem]{Proposition} 

\newtheorem{proposition-definition}[theorem]{Proposition-Definition} 

\newtheorem{corollary}[theorem]{Corollary} 

\newtheorem{conjecture}[theorem]{Conjecture} 

\theoremstyle{definition}

\newtheorem{definition}[theorem]{Definition} 

\theoremstyle{remark}

\newtheorem{remark}[theorem]{Remark} 

\newcommand{\Aff}{\mathbb{A}} 

\newcommand{\bd}{\mathbf}

\newcommand{\CC}{\mathbb{C}}
\newcommand{\EE}{\mathbb{E}}
\newcommand{\FF}{\mathbb{F}}

\newcommand{\PP}{\mathbb{P}}
\newcommand{\QQ}{\mathbb{Q}}
\newcommand{\RR}{\mathbb{R}}

\newcommand{\ZZ}{\mathbb{Z}}

\newcommand{\lcm}{\operatorname{lcm}}
\newcommand{\sgn}{\operatorname{sgn}}

\newcommand{\conj}{\overline}

\newcommand{\vol}{\operatorname{vol}}

\newcommand{\abs}[1]{\lvert #1 \rvert}
\newcommand{\card}[1]{\lvert #1 \rvert}
\newcommand{\norm}[1]{\lVert #1 \rVert}
\newcommand{\floor}[1]{\lfloor #1 \rfloor}
\newcommand{\ceil}[1]{\lceil #1 \rceil}

\newcommand{\eps}{\epsilon}

\newcommand{\rad}{\operatorname{rad}}

\newcommand{\Supp}{\operatorname{Supp}}

\newcommand{\disc}{\operatorname{disc}}
\newcommand{\Aut}{\operatorname{Aut}}
\newcommand{\Gal}{\operatorname{Gal}}

\newcommand{\map}{\operatorname}
\newcommand{\mscr}{\mathscr}
\newcommand{\mcal}{\mathcal}
\newcommand{\mf}{\mathfrak}

\newcommand{\ol}{\overline}

\newcommand{\texpdf}{\texorpdfstring}
\newcommand{\GL}{\operatorname{GL}}

\newcommand{\defeq}{\colonequals}
\newcommand{\eqdef}{\equalscolon}
\newcommand{\maps}{\colon}

\newcommand{\belongs}{\subseteq}

\newcommand{\set}[1]{\{#1\}}

\newcommand{\grad}{\nabla}

\usepackage[pagebackref=true]{hyperref}
\usepackage[alphabetic,initials]{amsrefs}

\begin{document}

\title{Sums of cubes and the Ratios Conjectures}

\date{}
\author{Victor Y. Wang}
\address{Department of Mathematics, Princeton University, Princeton, NJ, USA}
\address{Courant Institute of Mathematical Sciences, New York University, New York, NY, USA}
\email{vywang@alum.mit.edu}


\begin{abstract}
Works of Hooley and Heath-Brown imply a near-optimal bound on the number $N$ of integral solutions to $x_1^3+\dots+x_6^3 = 0$ in expanding regions,
conditional on automorphy and GRH for certain Hasse--Weil $L$-functions;
for regions of diameter $X\ge 1$, the bound takes the form $N\le C(\varepsilon) X^{3+\varepsilon}$ ($\varepsilon>0$).
We attribute the $\varepsilon$ to several subtly interacting proof factors;
we then remove the $\varepsilon$ assuming some standard number-theoretic hypotheses, mainly featuring the Ratios and Square-free Sieve Conjectures.
In fact,
our softest hypotheses imply
conjectures of Hooley and Manin on $N$,
and show that almost all integers $a\not\equiv \pm 4 \bmod{9}$ are sums of three cubes.
Our fullest hypotheses are capable of proving power-saving asymptotics for $N$,
and producing almost all primes $p\not\equiv \pm 4 \bmod{9}$.
\end{abstract}

\maketitle

\setcounter{tocdepth}{1}
\tableofcontents
\setcounter{tocdepth}{3}

\section{Introduction}
\label{SEC:intro}

Let $F_0=F_0(x,y,z) \defeq x^3+y^3+z^3$.
For each $a\in \ZZ$, the cubic surface $F_0=a$
has a fairly rich set of \emph{rational} points \cite{segre1943note}.
On the other hand, Mordell has suggested that producing large, general \emph{integer} solutions to $F_0=a$ for $a=3$ (or for any other fixed $a\in \ZZ$) could be as hard as ``finding when an assigned sequence, e.g.~$123456789$, occurs in the decimal expansion of $\pi$'' \cite{mordell1953integer}*{p.~505}.
The recent work \cite{booker2021question} of Booker and Sutherland resolves Mordell's specific question for $a=3$, but the spirit of Mordell's suggestion certainly remains.

Heath-Brown has conjectured that $F_0=a$ should have infinitely many solutions $(x,y,z)\in \ZZ^3$ for any fixed $a\not\equiv \pm 4 \bmod{9}$ (see \cite{heath1992density}*{p.~623}).
To represent all $a\not\equiv \pm 4 \bmod{9}$ even once, one must allow both positive and negative values of $x,y,z$.
The set $F_0(\ZZ_{\ge 0}^3)$ has upper density $\leq \Gamma(4/3)^3/6 = 0.1186788\ldots$ in $\ZZ_{\geq 0}$ \cite{davenport1939waring};
Deshouillers, Hennecart, and Landreau have given a model and evidence suggesting a precise density of $0.0999425\ldots$ \cite{deshouillers2006density}.
Wooley has shown, unconditionally, that $F_0(\ZZ_{\ge 0}^3)$ contains $\gg A^{0.91709477}$ integers $a\in [0,A]$ for reals $A\ge 1$ \cites{wooley1995breaking,wooley2000sums,wooley2015sums}.
We now recall a result of
Hooley \cites{hooley1986HasseWeil,hooley_greaves_harman_huxley_1997} and Heath-Brown \cite{heath1998circle},
and state our main result building on it; we then give further details and background.

\begin{theorem*}
[Hooley; Heath-Brown]
For certain Hasse--Weil $L$-functions,
assume automorphy and GRH.
Then $\gg_\eps A^{1-\eps}$ integers $a\in [0,A]$ lie in $F_0(\ZZ_{\ge 0}^3)$,
for any $\eps>0$.
\end{theorem*}

\begin{theorem*}
For certain Hasse--Weil $L$-functions,
assume automorphy, GRH, and the Ratios Conjectures.
For a certain polynomial $\Delta$, assume the Square-free Sieve Conjecture.
Then $\gg A$ integers $a\in [0,A]$ lie in $F_0(\ZZ_{\ge 0}^3)$,
and $100\%$ of integers $a\not\equiv\pm4\bmod{9}$ lie in $F_0(\ZZ^3)$.
\end{theorem*}

Both results require estimating sums that roughly take the following form:
\begin{equation}
\label{EXPR:rough-picture-of-delta-method-analysis}
\sum_{\bm{c}\in \ZZ^6} \int_{t\in \RR} \frac{dt}{\prod_{p\nmid \Delta(\bm{c})} (\textnormal{local $L$-factors})}
\cdot (\textnormal{real analysis})
\cdot \prod_{p\mid \Delta(\bm{c})} (\textnormal{geometry/$\FF_p$ + analysis/$\ZZ_p$}).
\end{equation}
The utility of GRH in this context has been highlighted by Bombieri; see \cite{bombieri2006riemann}*{p.~111}.
Also, in a function-field setting,
\cite{glas2022question} has made \cite{heath1998circle} unconditional,
and one could likely simplify our present hypotheses accordingly (since GRH and the Square-free Sieve Conjecture are known over function fields).
We only focus on $\QQ$ for practical reasons.

For a typical $a$, the integer solutions to $F_0 = a$ are expected to be at least exponentially sparse, if they in fact exist.
Heath-Brown's conjecture would imply that the only obstructions to solubility for $F_0 = a$, for any $a$, are local.
The naive local-to-global analog of Heath-Brown's conjecture for $5x^3+12y^3+9z^3$ is known to fail
(see e.g.~\cite{ghosh2017integral}*{p.~691, footnote~3}),
due to Brauer--Manin obstructions that do not apply to $F_0$ \cite{colliot2012groupe}*{p.~1304}.

We restrict ourselves to a statistical analysis of $F_0 = a$ over $a\in \ZZ$, for $(x,y,z)\in \ZZ^3$ lying in carefully chosen regions, conducted using second moments and the variance framework of \cites{ghosh2017integral,diaconu2019admissible,wang2022thesis,wang2023prime}.
This connects naturally to difficult open questions in $6$ variables, e.g.~\cite{hooley1986some}*{Conjecture~2},
which lie beyond the square-root barrier in the classical Hardy--Littlewood circle method.
We will attack these questions under standard number-theoretic hypotheses, primarily regarding $L$-function statistics of Random Matrix Theory (RMT) type.
Our work opens with the delta method of \cites{duke1993bounds,heath1996new} (a clean modern form of the Kloosterman method of \cite{kloosterman1926representation}, more precise than the upper-bound variant used in \cites{hooley1986HasseWeil,hooley_greaves_harman_huxley_1997}), whose harmonic analysis in principle allows for cancellation over the difficult classical minor arcs.
We prove three levels of results, under three levels of hypotheses (the first two levels being relatively soft and qualitative; see Conjectures~\ref{CNJ:(R2')} and~\ref{CNJ:(RA1o)}).
We first recall a general weighted version of Hooley's conjecture.

Fix a cubic form $F(\bm{x})=F(x_1,\dots,x_m)\in \ZZ[x_1,\dots,x_m]$ in $m\geq 4$ variables with nonzero discriminant.
Let $V$ be the hypersurface $F=0$ in $\PP^{m-1}_\QQ$.
Let $\Upsilon$ be the set of $\floor{m/2}$-dimensional vector spaces $L\belongs \QQ^m$ over $\QQ$ on which $F$ vanishes.
On a first reading, we suggest assuming that $m=6$ and $F$ is diagonal,
though we will often work generally.

Given a real $X\geq 1$, and a function $w\in C^\infty_c(\RR^m)$, let
\begin{equation}
\label{EQN:define-N_w(X),N(X)}
N_{F,w}(X) \defeq \sum_{\bm{x}\in \ZZ^m:\, F(\bm{x})=0} w(\bm{x}/X),
\quad
N_F(X) \defeq \sum_{\bm{x}\in [-X,X]^m:\, F(\bm{x})=0} 1;
\end{equation}
and if $m=6$ (our main case of interest, in which $\floor{m/2}=3$), let
\begin{equation}
\label{EQN:define-HLH-error-E_w(X)}
E_{F,w}(X)\defeq N_{F,w}(X)
- \mf{S}_F\cdot \sigma_{\infty,F,w}\cdot X^3
- \sum_{L\in \Upsilon} \sum_{\bm{x}\in L\cap \ZZ^6} w(\bm{x}/X),
\end{equation}
where $\mf{S}_F\ll_F 1$ is the familiar singular series defined in \cite{wang2023_isolating_special_solutions}*{\S6}, and where
\begin{equation}
\label{EQ:define-real-density}
\sigma_{\infty,F,w}
\defeq \lim_{\eps\to 0}{(2\eps)^{-1} \int_{\abs{F(\bm{x})}\le \eps} d\bm{x}\, w(\bm{x})}
\ll_{F,w} 1.
\end{equation}
One could attribute to Hooley \cite{hooley1986some}*{Conjecture~2}, Manin (see e.g.~\cite{franke1989rational}), Vaughan--Wooley \cite{vaughan1995certain}*{Appendix}, Peyre \cite{peyre1995hauteurs}, et al.~the following conjecture:
\begin{equation}
\label{EQN:soft-HLH-general-homogeneous-weight}
\lim_{X\to \infty} X^{-3} E_{F,w}(X) = 0.
\end{equation}
(The original \cite{hooley1986some}*{Conjecture~2 for $l=3$} would follow from \eqref{EQN:soft-HLH-general-homogeneous-weight} for $F = x_1^3+\dots+x_6^3$, applied to a suitable sequence of weights $w$.)
See \cite{ding2020variance} for another related problem.

Unconditionally, $N_F(X)\ll_\eps X^{7/2}/(\log{X})^{5/2-\eps}$ for $X\geq 2$ \cite{vaughan2020some}, when $m=6$ and $F$ is diagonal.
Under standard hypotheses on the varieties $V_{\bm{c}}\belongs \PP^{m-1}_\QQ$ cut out by $F(\bm{x}) = \bm{c}\cdot\bm{x} = 0$,
one can prove the near-optimal Theorem~\ref{THM:Hoo-HB-6-var-cubic-conditional-3+eps-bound} for the same $F$.
(Here $\bm{c}\cdot \bm{x} \defeq \sum_{1\le i\le m} c_ix_i$.)
Let $\Delta(\bm{c})\in \ZZ[\bm{c}]=\ZZ[c_1,\dots,c_m]$ be
the discriminant polynomial defined in \S\ref{SEC:background-on-discriminants-and-the-delta-method}.
Let
\begin{equation}
\label{EQN:define-S_0,S_1-for-Delta-vanishing-and-nonzero-loci}
    \mcal{S}_0\defeq \set{\bm{c}\in \ZZ^m: \Delta(\bm{c})=0},
    \quad
    \mcal{S}_1\defeq \set{\bm{c}\in \ZZ^m: \Delta(\bm{c})\neq 0}.
\end{equation}
For each $\bm{c}\in \mcal{S}_1$, one can package local data on $V_{\bm{c}}$
into a Hasse--Weil $L$-function $L(s,V_{\bm{c}})$, defined in \S\ref{SEC:background-on-individual-L-functions} along with the rest of the list \eqref{full-list-of-Hasse-Weil-L-functions}.

\begin{theorem}
[\cites{hooley1986HasseWeil,hooley_greaves_harman_huxley_1997}; \cite{heath1998circle}]
\label{THM:Hoo-HB-6-var-cubic-conditional-3+eps-bound}
Assume $m=6$ and $F$ is diagonal.
For each $\bm{c}\in \mcal{S}_1$,
assume Conjecture~\ref{CNJ:(HW2)} for $L(s,V_{\bm{c}})$.
Then $N_F(X)\ll_\eps X^{3+\eps}$ for all $\eps>0$.
\end{theorem}

Conjecture~\ref{CNJ:(HW2)} concerns automorphy and the Grand Riemann Hypothesis (GRH).
Strictly speaking, Hooley and Heath-Brown assume (in their ``Hypothesis~HW'') GRH plus certain Selberg-type axioms, e.g.~analyticity, but it is natural to assume automorphy in place of such axioms.
A nice reference bridging these two perspectives is \cite{farmer2019analytic}.

The proof of Theorem~\ref{THM:Hoo-HB-6-var-cubic-conditional-3+eps-bound} involves GRH on $1/L(s,V_{\bm{c}})$, surprisingly,
coupled with subtle algebro-geometric factors
of a different nature.
See \S\ref{SEC:background-on-discriminants-and-the-delta-method} for details.
The $L$-function ingredients extend directly to general $F$,
but some other ingredients have yet to be generalized.

GRH gives a pointwise upper bound on $1/L(s,V_{\bm{c}})$ for $\Re(s)>1/2$, sufficient for Theorem~\ref{THM:Hoo-HB-6-var-cubic-conditional-3+eps-bound}.
Going past the critical line (in mean value over $\bm{c}$) turns out to require much new work, both with $L$-functions and with other factors.
In terms of $L$-functions, we mainly use $L(s,V_{\bm{c}})$ for $\Re(s)\geq 1/2-\delta$, as well as $L(s,V_{\bm{c}},\bigwedge^2)$, $\zeta(s)$, and $L(s,V)$ for $\Re(s)\geq 1-\delta$.
It will be convenient to assume Conjecture~\ref{CNJ:(HW2)} in full,
but average versions might also suffice.

\begin{conjecture}
[HW2]
\label{CNJ:(HW2)}
Let $\bm{c}\in \mcal{S}_1$.
Let $L(s)$ be one of the Hasse--Weil $L$-functions
\begin{equation}
\label{full-list-of-Hasse-Weil-L-functions}
    L(s, V_{\bm{c}}),\;
    L(s, V_{\bm{c}}, {\textstyle \bigotimes^2}),\;
    L(s, V_{\bm{c}}, \map{Sym}^2),\;
    L(s, V_{\bm{c}}, {\textstyle \bigwedge^2}),\;
    \zeta(s),\;
    L(s, V).
\end{equation}
Let $L_v(s)$ be the local factors (including $L_\infty(s)$, the gamma factor) associated to $L$.
\begin{enumerate}
    \item There exists an integer $d\ge 1$, and an isobaric automorphic representation $\Pi$ of $\GL_d(\bd{A}_\QQ)$, such that $L_v(s) = L_v(s,\Pi)$ at all places $v\leq \infty$.
    
    \item $L(s,\Pi)$ has no zeros in the half-plane $\Re(s)>1/2$.
\end{enumerate}
\end{conjecture}

For the background needed to interpret Conjecture~\ref{CNJ:(HW2)}~(HW2), see \S\ref{SEC:background-on-individual-L-functions}.

\begin{theorem}
\label{THM:level-1-positive-density}
    
    
    
Suppose $F=x_1^3+\dots+x_6^3$.
Assume Conjectures~\ref{CNJ:(HW2)}, \ref{CNJ:(R2')}, and~\ref{CNJ:(SFSCp)}.
Then
\begin{equation}
\label{INEQ:main-level-1-goal-eps-free-point-count}
    N_F(X)\ll X^3.
\end{equation}
Let $S\belongs \ZZ_{\ge 0}$.
If $S$ has positive lower density in $\ZZ_{\ge 0}$,
then so does $F_0(S^3)$.
\end{theorem}


We use mean-value RMT-type predictions derived from the Moments and Ratios Conjectures of \cites{conrey2005integral,conrey2007applications,conrey2008autocorrelation}.
(See also \cites{diaconu2003multiple,vcech2022ratios}, and references within, for another important perspective on such conjectures.)
The Hasse--Weil $L$-functions $L(s,V_{\bm{c}})$ form a geometric family, in the sense of \cite{sarnak2016families}.
For each $\bm{c}\in \mcal{S}_1$, define the Dirichlet series
\begin{equation}
\label{EQN:define-Phi_1}
    \Phi^{\bm{c},1}(s)
    \defeq \zeta(2s)^{-1} L(s+1/2, V)^{-1} L(s,V_{\bm{c}})^{-1}
    = 1/{\zeta(2s) L(s+1/2, V) L(s, V_{\bm{c}})}.
\end{equation}
Note that $\zeta(2s)$ and $L(s+1/2, V)$ are independent of $\bm{c}$.



\begin{conjecture}
[R2']
\label{CNJ:(R2')}
Suppose $2\mid m$.
Let $f\maps \CC\to \CC$ be entire, with $f(s)\ll_{f,b} (1+\abs{\Im(s)})^{-b}$ on the strip $0\le \Re(s)\le 2$ for all $b\in \ZZ_{\ge 1}$.
Let $Z,N\in \RR_{\ge 1}$ with $N\le Z^3$.
If $\sigma_0\in (1, 2)$, then
\begin{equation}
\label{INEQ:R2'-goal}
\sum_{\bm{c}\in \mcal{S}_1\cap [-Z,Z]^m}
\left\lvert\int_{(\sigma_0)} ds\, \Phi^{\bm{c},1}(s) \cdot f(s)N^s\right\rvert^2
\ll_{F} Z^m N
\sup_{0\le \sigma\le 2} \int_{\RR} dt\, (1+\abs{t})^2 \abs{f(\sigma+it)}^2.
\end{equation}
\end{conjecture}

The contour $(\sigma_0)$ in \eqref{INEQ:R2'-goal} runs from $s=\sigma_0-i\infty$ to $s=\sigma_0+i\infty$.
The left-hand side of \eqref{INEQ:R2'-goal} is independent of $\sigma_0\in (1,2)$ unconditionally, or further under (HW2).
Other versions of Conjecture~\ref{CNJ:(R2')} might also suffice for our purposes;
for instance, an $\ell^{1+\delta}$ analog of \eqref{INEQ:R2'-goal}
(with some nontrivial adjustments)
might suffice,
the precise norm of $f$ on the right-hand side of \eqref{INEQ:R2'-goal} is not very important,
and one might not need to allow such general $f$.

There are no log factors on the right-hand side of \eqref{INEQ:R2'-goal};
the factor $\zeta(2s)^{-1}L(s+1/2,V)^{-1}$ in \eqref{EQN:define-Phi_1},
and the integral over $(\sigma_0)$ in \eqref{INEQ:R2'-goal},
play a mollifying role.
The statement (R2') can be derived from (HW2) and the Ratios Conjecture~\ref{CNJ:(R2o)}~(R2$o$); see Proposition~\ref{PROP:R2-implies-R2'}.
However,
there may well be another route to (R2') not passing through (R2$o$).
The statement (R2') essentially concerns cancellation in the coefficients of $\Phi^{\bm{c},1}(s)$ over moduli $n$ in a dyadic range; see Proposition~\ref{PROP:R2'-implies-R2'E}.
A similar log-free cancellation statement, \cite{li2022moments}*{(1.3)}, has recently played a crucial role in another context.
Furthermore, over function fields (or under (HW2) over $\QQ$), it should already be possible to obtain partial results towards (R2'), using ideas of \cites{soundararajan2009moments,harper2013sharp,bui2021ratios,florea2021negative,bui2023negative} (after Cauchy--Schwarz over $s$).

For our main results, we also need the Square-free Sieve Conjecture (cf.~\cite{miller2004one}*{p.~956} and \cites{granville1998abc,poonen2003squarefree,bhargava2014geometric}) for the polynomial $\Delta$, restricted to a certain range $1\le P\le Z^{3/2}$.
This hypothesis concerns ``unlikely divisors'' of the outputs of $\Delta$.
Such hypotheses can be made unconditional over function fields; see e.g.~\cite{poonen2003squarefree}*{Lemma~7.1}.

\begin{conjecture}
[SFSC$_{p,3}$]
\label{CNJ:(SFSCp)}
There exists $\eta_0=\eta_0(\Delta)\in \RR_{>0}$ such that
if $Z, P \in \RR_{\ge 1}$ and $P\leq Z^{3/2}$, then
$\#\set{\bm{c} \in \ZZ^m \cap [-Z,Z]^m
:\exists\;\textnormal{a prime $p \in [P,2P)$ with $p^2 \mid \Delta(\bm{c})$}}
\ll_{\Delta} Z^m P^{-\eta_0}$.
\end{conjecture}

The Square-free Sieve Conjecture (SFSC$_{p,3}$) is used to confront some novel algebro-geometric issues in our work.
We use it to prove, for diagonal $F$,
a geometric relative (Conjecture~\ref{CNJ:(B3G)}) of the automorphic Sarnak--Xue Density Hypothesis (\cite{sarnak1991bounds}*{Conjecture~1}, concerning the extent to which a naive generalization of the Ramanujan Conjecture can fail).
This involves sieve-theoretic ideas
weighted by somewhat dangerous factors.
Even though we do not currently see how to eliminate the use of (SFSC$_{p,3}$) over $\QQ$, it is fortunate to be able to reduce Conjecture~\ref{CNJ:(B3G)} to (SFSC$_{p,3}$) when $F$ is diagonal.
One can also reduce (SFSC$_{p,3}$) for diagonal $F$ to the case $m=4$ (but not to $m=2$, it seems);
see Proposition~\ref{PROP:SFSC-for-m=4-implies-SFSC-for-diagonal}.

After Theorem~\ref{THM:level-1-positive-density}, our next main result is the following:

\begin{theorem}
\label{THM:level-2-almost-all-integers}
Suppose $m=6$ and $F$ is diagonal.
Assume Conjectures~\ref{CNJ:(HW2)}, \ref{CNJ:(R2')}, \ref{CNJ:(SFSCp)}, and~\ref{CNJ:(RA1o)}.
Then \eqref{EQN:soft-HLH-general-homogeneous-weight} holds for all functions $w\in C^\infty_c(\RR^m)$ for which
\begin{equation}
\label{COND:clean-weight-condition-in-diagonal-case}
\ol{\set{\bm{x}\in \RR^m: w(\bm{x})\ne 0}}
\belongs \set{\bm{x}\in \RR^m: x_1\cdots x_m\neq 0}.
\end{equation}
Therefore, the Hasse principle holds for $V$.
Furthermore, if $F=x_1^3+\dots+x_6^3$, then
$100\%$ of integers $a\not\equiv\pm4\bmod{9}$ lie in $F_0(\ZZ^3)$.
\end{theorem}



For another conditional approach to the Hasse principle for $V$ when $F$ is diagonal, see \cite{swinnerton2001solubility}.
Over function fields of characteristic $\geq 7$, the Hasse principle for $V$ is already known in general when $m=6$ \cite{tian2017hasse}.
But our approach has quantitative advantages, which become qualitative when applied to sums of three cubes.

We expect that the condition \eqref{COND:clean-weight-condition-in-diagonal-case} could be removed with enough work.
When $F=x_1^3+\dots+x_6^3$, it is in fact possible to do this for free: Theorem~\ref{THM:level-2-almost-all-integers} has the following corollary.
(A similar but messier statement is possible for arbitrary diagonal $F$ when $m=6$.)

\begin{corollary}
\label{COR:remove-Hessian-assumption-for-free-in-soft-Fermat-case}
Suppose $F=x_1^3+\dots+x_6^3$.
Assume Conjectures~\ref{CNJ:(HW2)}, \ref{CNJ:(R2')}, \ref{CNJ:(SFSCp)}, and~\ref{CNJ:(RA1o)}.
Then \eqref{EQN:soft-HLH-general-homogeneous-weight} holds for all $w\in C^\infty_c(\RR^m)$.
Consequently, \cite{hooley1986some}*{Conjecture~2 for $l=3$} holds.
\end{corollary}









Theorem~\ref{THM:level-2-almost-all-integers} makes use of a first-moment estimate for the quantity \eqref{EQN:define-Phi_1} over $\mcal{S}_1$,
and over some mildly localized pieces of $\mcal{S}_1$ (``adelic perturbations'' of $\mcal{S}_1$ restricted by a parameter $M$).
For $\bm{b}=(b_1,\dots,b_m)\in \ZZ^m$ and $M\in \RR_{>0}$, consider the box
\begin{equation}
\label{EQN:define-box-B_M(b)}
\mcal{B}_M(\bm{b})
\defeq [b_1/M, (b_1+1)/M) \times \dots \times [b_m/M, (b_m+1)/M) \belongs \RR^m.
\end{equation}
For $Z\in \RR_{>0}$, let $Z\cdot \mcal{B}_M(\bm{b})$ denote the dilate $\set{Z\bm{r}: \bm{r}\in \mcal{B}_M(\bm{b})}\belongs \RR^m$.
For each $Z\in \RR_{\ge 2}$, let
\begin{equation}
\label{EQN:sigma(Z)}
\sigma(Z)
\defeq 1/2+1/\log{Z}.
\end{equation}

\begin{conjecture}
[RA1$o$]
\label{CNJ:(RA1o)}
Suppose $2\mid m$, and assume Conjecture~\ref{CNJ:(HW2)}.
Let $M\in \RR_{\ge 1}$;
let $n_0\in \ZZ\cap [1,M]$ and $\bm{a}, \bm{b}\in \ZZ^m\cap [-M,M]^m$.
Let $A_{F,1}^{\bm{a},n_0}(s)$ be defined as in \S\ref{SUBSUBSEC:deriving-(RA1)} (in terms of $F$, $\bm{a}$, $n_0$).
If $Z\in \RR_{\ge 2}$ and $t\in [-\log{Z}, \log{Z}]$, then for $s = \sigma(Z)+it$, we have
\begin{equation}
\label{EQN:soft-RA1-goal}
\sum_{\substack{\bm{c}\in \mcal{S}_1\cap Z\cdot \mcal{B}_M(\bm{b}): \\ \bm{c}\equiv\bm{a}\bmod{n_0}}}
\Phi^{\bm{c},1}(s)
= \sum_{\substack{\bm{c}\in \mcal{S}_1\cap Z\cdot \mcal{B}_M(\bm{b}): \\ \bm{c}\equiv\bm{a}\bmod{n_0}}}
(1 + o_{F, M; Z\to \infty}(1)) \cdot A_{F,1}^{\bm{a},n_0}(s).
\end{equation}
\end{conjecture}


Here $A_{F,1}^{\bm{a},n_0}(s)$ is a Dirichlet series with an Euler product, absolutely convergent on the half-plane $\Re(s)>1/3$.
The terms $o_{F, M; Z\to \infty}(1)$ are required to tend to $0$ as $Z\to \infty$ (when $F$, $M$ are fixed).
See Proposition~\ref{PROP:RA1o-implies-RA1o'E} for the main use of Conjecture~\ref{CNJ:(RA1o)}.



The Ratios Conjectures include (RA1$o$), even with a power saving;
see \S\ref{SUBSUBSEC:deriving-(RA1)} for details.
The choice \eqref{EQN:sigma(Z)} is permissible according to \cite{conrey2007applications}*{(2.11b)}.
The essential feature of \eqref{EQN:sigma(Z)} for us is that $(\sigma(Z)-\frac12)\cdot \log{Z}$ is positive and independent of $Z$ (but its precise constant value is not important).
We could get away with a larger choice of $\sigma(Z)$ if we assumed a correspondingly stronger error term in \eqref{EQN:soft-RA1-goal};
but the present formulation of (RA1$o$) is clean, and easy to compare with other literature.\footnote{Note that if $s = \sigma(Z)$,
then in \eqref{EQN:define-Phi_1}, we have $\zeta(2s) \asymp \log{Z}$ and $L(s+1/2,V) \asymp (\log{Z})^{r_F}$, where $r_F\in \ZZ_{\ge 0}$ is given explicitly for diagonal $F$ by \cite{wang2022thesis}*{Lemma~8.6.7}.
The moments of $1/L$ we consider (over our orthogonal family of $L$-functions $L(s,V_{\bm{c}}$)) are thus analogous to central moments over symplectic families.}
The paper \cite{florea2021negative} seems to come close to proving (RA1$o$) for a different family of $L$-functions, over a function field.

We believe (RA1$o$), like (R2'), represents a tantalizing research direction.
There is another direction worth mentioning.
In light of the log-free square-root cancellation in $\ell^2$ conjectured in \eqref{INEQ:R2'-goal}, one may hope that ``better than square-root cancellation'' occurs in $\ell^{2-\eps}$, by analogy with \cite{harper2023typical}*{(1.2)} (an attractive conjecture based on random multiplicative functions and multiplicative chaos; see e.g.~\cites{gorodetsky2021magic,harper2023typical} for details and references).
If true, this would provide additional cancellation in Proposition~\ref{PROP:R2'E-implies-R2'E'} (one of the key ingredients for Theorem~\ref{THM:level-1-positive-density}), and thus provide an alternative approach to Theorem~\ref{THM:level-2-almost-all-integers} (but not to Theorem~\ref{THM:level-3-power-saving}).

Our final main result, Theorem~\ref{THM:level-3-power-saving}, goes beyond Theorem~\ref{THM:level-2-almost-all-integers}.

\begin{theorem}
\label{THM:level-3-power-saving}
Suppose $m=6$ and $F$ is diagonal.
Assume Conjectures~\ref{CNJ:(HW2)}, \ref{CNJ:(SFSCp)} with $\eta_0$, \ref{CNJ:(RA1delta)} with $\eta_1$, and~\ref{CNJ:(EKL)} with $H$.
Then there exists a real $\delta=\delta(\eta_0,\eta_1,\deg{H})>0$ such that
$E_{F,w}(X)\ll_{F,w} X^{3-\delta}$
holds for all functions $w\in C^\infty_c(\RR^m)$ satisfying \eqref{COND:clean-weight-condition-in-diagonal-case}.
\end{theorem}

Our methods would allow one to prove a version of Theorem~\ref{THM:level-3-power-saving} uniform over small archimedean and non-archimedean perturbations to $E_{F,w}(X)$.
By \cite{wang2023prime}*{Theorem~1.2}, one could then show under the hypotheses of Theorem~\ref{THM:level-3-power-saving} that if $F=x_1^3+\dots+x_6^3$, then $100\%$ of primes $p\not\equiv\pm4\bmod{9}$ lie in $F_0(\ZZ^3)$.
One would also be able to give a power-saving analog of Corollary~\ref{COR:remove-Hessian-assumption-for-free-in-soft-Fermat-case}.
But to give full details would obscure
our exposition.

The power saving in Theorem~\ref{THM:level-3-power-saving} is small and complicated.
(Egregious ``exponent divisions'' occur in Lemma~\ref{LEM:(EKL')} and Proposition~\ref{PROP:SFSCp-implies-SFSCq}, due to our use of \eqref{INEQ:multivarite-zero-density-mod-q}; and similarly in \eqref{INEQ:lambda-decay-of-S_1(C,lambda)}, due to \eqref{INEQ:multivarite-near-zero-real-density}.)
It would be very interesting to understand the limits of what one can hope for.

\begin{conjecture}
[RA1$\delta$]
\label{CNJ:(RA1delta)}
Suppose $2\mid m$, and assume Conjecture~\ref{CNJ:(HW2)}.
Let $Z\in \RR_{\ge 2}$
and $M\in [1, Z^{\eta_1}]$.
Let $n_0\in \ZZ\cap [1,M]$ and $\bm{a}, \bm{b}\in \ZZ^m\cap [-M,M]^m$.
There exists a real $\eta_1=\eta_1(F)>0$, depending only on $F$,
such that if $t\in [-M,M]$ and $s = \sigma(Z)+it$, then
\begin{equation}
\label{EQN:power-saving-RA1-goal}
\sum_{\substack{\bm{c}\in \mcal{S}_1\cap Z\cdot \mcal{B}_M(\bm{b}): \\ \bm{c}\equiv\bm{a}\bmod{n_0}}}
\Phi^{\bm{c},1}(s)
= \sum_{\substack{\bm{c}\in \mcal{S}_1\cap Z\cdot \mcal{B}_M(\bm{b}): \\ \bm{c}\equiv\bm{a}\bmod{n_0}}}
(1 + O_F(Z^{-\eta_1})) \cdot A_{F,1}^{\bm{a},n_0}(s).
\end{equation}
\end{conjecture}


For Theorem~\ref{THM:level-3-power-saving}, certain degenerate residue classes play a larger role in local calculations than for Theorem~\ref{THM:level-2-almost-all-integers}.
To pacify these residue classes, we need effective control on the variation of an individual local factor $L_p(s, V_{\bm{c}})$ over $\bm{c}\in \mcal{S}_1$.
We work $p$-adically for convenience, using Remark~\ref{RMK:extended-local-L-factor-definitions} to define $L_p(s, V_{\bm{c}})$ for each $\bm{c}\in \ZZ_p^m$ with $\Delta(\bm{c})\ne 0$.

\begin{conjecture}
[EKL]
\label{CNJ:(EKL)}
There exists a \emph{nonzero} homogeneous polynomial $H\in \ZZ[\bm{c}]$, with $H/\Delta\in \ZZ[\bm{c}]$, such that for all primes $p$ and tuples $\bm{a}, \bm{b}\in \ZZ_p^m$ with $H(\bm{b}) \ne 0$ and $\bm{a} \equiv \bm{b} \bmod{p H(\bm{b})}$, we have $H(\bm{a})\ne 0$ and $L_p(s, V_{\bm{a}}) = L_p(s, V_{\bm{b}})$.
\end{conjecture}

Conjecture~\ref{CNJ:(EKL)} is an effective Krasner-type statement for $L_p(s, V_{\bm{c}})$.
A soft
version
follows from \cite{kisin1999local}, and suffices for Theorem~\ref{THM:level-2-almost-all-integers} but not for Theorem~\ref{THM:level-3-power-saving}.
When $m=4$, it should be possible to prove Conjecture~\ref{CNJ:(EKL)} (with $H\in \QQ^\times\cdot \Delta$) using a minimal model for the Jacobian of $V_{\bm{c}}$.
In general, one might hope to take $H$ to be a power of $\Delta$ (or perhaps $\Delta$ itself).

\subsection{Proof overview}

\S\ref{SEC:background-on-discriminants-and-the-delta-method} gives background on discriminants and the delta method.
The delta method (see \eqref{EQN:delta-method}) connects the point count $N_{F,w}(X)$ (from \eqref{EQN:define-N_w(X),N(X)}) to the local behavior of the intersections $F(\bm{x}) = \bm{c}\cdot\bm{x} = 0$ over $\FF_p$, $\ZZ_p$, $\RR$, and other rings, as $\bm{c}\in \ZZ^m$ and $p$ vary.
Cf.~\eqref{EXPR:rough-picture-of-delta-method-analysis}.
We highlight several distinct sources of epsilon in the Hooley--Heath-Brown Theorem~\ref{THM:Hoo-HB-6-var-cubic-conditional-3+eps-bound},
and state a result from \cite{wang2023_isolating_special_solutions} (over $\mcal{S}_0$) addressing one such source.

\S\ref{SEC:background-on-individual-L-functions} provides background on Hasse--Weil $L$-functions and automorphic $L$-functions.

\S\ref{SEC:local-control-on-polynomials-and-L-functions} gives some local control on polynomials and $L$-factors (based in part on \cite{kisin1999local}); we need this for some local estimates and calculations.

\S\ref{SEC:separation-lemmas} gives a useful ``smooth framework'' for dyadic decomposition and separation of variables.
This lets us break certain key sums throughout the paper into more manageable pieces.

\S\ref{SEC:statistics-of-L-function-families} derives some $L$-function statistics over $\bm{c}\in \mcal{S}_1$,
after first doing local calculations (in the spirit of the Deligne--Katz equidistribution theorem) connected to RMT Symmetry Types (cf.~\cite{sarnak2016families}*{Universality Conjecture}).
In particular, we state, and build on, some cases of the Ratios Conjectures for $L(s,V_{\bm{c}})$.
Importantly here, the Ratios Recipe (see \S\ref{SUBSEC:applying-CFZ-ratios-conjecture-recipe}) can only apply once we restrict to $\bm{c}\in \mcal{S}_1$.
The recipe, naively extended to all $\bm{c}\in \ZZ^m$, would give false results (failing to detect the special subvarieties $L\in \Upsilon$ on $F=0$ isolated in \cite{wang2023_isolating_special_solutions}).

\S\ref{SEC:adapting-L-statistics-to-the-delta-method} begins to connect the ``pure'' $L$-function statistics from \S\ref{SEC:statistics-of-L-function-families} to the delta method.
We approximate certain Dirichlet series ``past'' the critical line, in a reasonably simple and uniform way over $\bm{c}\in \mcal{S}_1$,
and control the resulting ``approximation errors'' on average.
Handling these ``errors'' demands careful use of H\"{o}lder and other ideas.
For example, by algorithmic tree-like means, we construct in \S\ref{SUBSEC:handling-variation-of-error-factors} a small exceptional set away from which one may apply the ``pure'' Conjectures~\ref{CNJ:(RA1o)} and~\ref{CNJ:(RA1delta)}; it is also here that Conjecture~\ref{CNJ:(EKL)} plays some key role.

The ``mollified'' series $\Phi^{\bm{c},1}(s)$ from \eqref{EQN:define-Phi_1} not only makes the formulas in \S\ref{SEC:statistics-of-L-function-families} nicer, but (as we will explain in \S\ref{SEC:adapting-L-statistics-to-the-delta-method}) also holds significance in \eqref{EQN:delta-method};
this double significance, though innocent at first glance, secretly reflects a randomness property (connected to Deligne--Katz) stemming from the fact that $\deg{F}\ge 3$.
Throughout the paper,
we take much advantage of the structure of \eqref{EQN:define-Phi_1}; this is essential in Conjecture~\ref{CNJ:(R2')} and Lemma~\ref{LEM:basic-main-term-bound-for-RA1'E'}, for instance.

\S\ref{SEC:new-bounds-on-integral-J} proves new integral bounds sensitive to some real geometry (involving discriminants).
Our approach shares some important features with \cite{huang2020density} (a beautiful recent paper on \emph{approximate} integral points).
In addition, we have several parameters of interest, and must obtain genuinely multivariate decay.
Keeping track of uniformity is tricky.

\S\ref{SEC:new-bounds-on-bad-sums-S}, like \S\ref{SEC:new-bounds-on-integral-J}, proves some new ``discriminating'' pointwise estimates, but on complete exponential sums instead of oscillatory integrals.
We then apply these to formulate and address (under Conjecture~\ref{CNJ:(SFSCp)}) a geometric analog of the Sarnak--Xue Density Hypothesis.
Finite-field geometry (see \cite{wang2023dichotomous}) and the geometric sieve (see \cites{ekedahl1991infinite,bhargava2014geometric}) both play an important role here,
as does a nice result of Bus\'{e} and Jouanolou on discriminants (see Theorem~\ref{THM:Buse-Jouanolou-saturated-discriminant-gradient-square-decomposition}).
It is crucial throughout \S\ref{SEC:new-bounds-on-bad-sums-S} that $\deg{F}\ge 3$ (see e.g.~Remark~\ref{RMK:failure-of-prime-boundedness-criteria-for-quadratic-F}); this reflects a ``randomness'' not present for quadrics.

\S\ref{SEC:delta-endgame} ties everything together to prove our main results.
We also isolate ``axioms'' that---if true---would allow for non-diagonal $F$; see Theorems~\ref{THM:axiomatic-level-1-endgame-theorem}, \ref{THM:axiomatic-level-2-endgame-theorem}, and~\ref{THM:axiomatic-level-3-endgame-theorem}.
Here
we only consider $\bm{c}\in \mcal{S}_1$;
there are also separate issues for $\bm{c}\in \mcal{S}_0$ (see \cite{wang2023_isolating_special_solutions}*{Remark~1.6}).



For Theorem~\ref{THM:level-1-positive-density}, see \S\ref{SUBSEC:sharp-upper-bounds-via-R2'}.
Here we use an ``entirely positive'' H\"{o}lder argument over $\bm{c}\in \mcal{S}_1$:
we do not detect any cancellation over $\bm{c}$ that would go beyond
a log-free ``Mertens-type
heuristic on average''
(cf.~\cite{ng2004distribution}*{Theorem~1(iii)})
over $\bm{c}$.
Despite the ``decoupling'' convenience and power of H\"{o}lder,
we must therefore be careful in \S\ref{SUBSEC:sharp-upper-bounds-via-R2'} to obtain \emph{$\eps$-free} bounds.
(The structure of \S\ref{SUBSEC:sharp-upper-bounds-via-R2'} is inspired by our work with the large sieve in \cite{wang2023_large_sieve_diagonal_cubic_forms}.)

The proof of Theorem~\ref{THM:level-1-positive-density} tells us
(conditionally)
that \emph{even if one takes absolute values over $\bm{c}\in \mcal{S}_1$} in \eqref{EQN:delta-method},
the $X^\eps$ allowance of \cites{hooley_greaves_harman_huxley_1997,heath1998circle} is unnecessary; see \eqref{INEQ:absolute-log-free-bound-in-delta-method} in Theorem~\ref{THM:axiomatic-level-1-endgame-theorem}.
Theorem~\ref{THM:level-1-positive-density} also highlights
a nontrivial use of
the \emph{log-free} order of magnitude in Conjecture~\ref{CNJ:(R2')},
a robust qualitative prediction
that (if true) could perhaps be explained in other ways
(not just following the rather arithmetic Ratios Recipe).

For Theorem~\ref{THM:level-2-almost-all-integers}, Corollary~\ref{COR:remove-Hessian-assumption-for-free-in-soft-Fermat-case}, and Theorem~\ref{THM:level-3-power-saving}, see \S\ref{SUBSEC:cancellation-via-RA1}.
In most ranges, the ``entirely positive'' moment estimates of \S\ref{SUBSEC:sharp-upper-bounds-via-R2'} still suffice.
But this time, in a few key ranges,
we identify cancellation over $\bm{c}\in \mcal{S}_1$ (via \S\ref{SEC:adapting-L-statistics-to-the-delta-method}).
One critical step here is a reduction, via \S\ref{SEC:new-bounds-on-integral-J}, to large moduli in \eqref{EQN:delta-method} (over $\bm{c}\in \mcal{S}_1$), over which certain ``mollified'' RMT-type main terms vanish (cf.~Lemma~\ref{LEM:basic-main-term-bound-for-RA1'E'}).
There is also an alternative approach to cancellation (which we do not pursue): instead of Lemma~\ref{LEM:basic-main-term-bound-for-RA1'E'}, we could use the fact that $\int_{\RR^m} d\bm{c}\, J_{\bm{c},X}(n) = 0$ (provided $w(\bm{0})=0$), where $J_{\bm{c},X}(n)$ is defined as in \eqref{EQN:define-normalized-S-tilde-and-J}; cf.~\cite{wang2021_HLH_vs_RMT}*{Observation~10.7}.

\subsection{Conventions}
\label{SUBSEC:conventions}

We write $f\ll_S g$, or $g\gg_S f$, to mean $\abs{f} \leq Cg$ for some $C = C(S)>0$.
We let $O_S(g)$ denote a quantity that is $\ll_S g$.
We write $f\asymp_S g$ if $f\ll_S g\ll_S f$.
We let $o_{S; X\to \infty}(g)$ denote a quantity $f$ such that for every $\eps>0$, there exists $X_0=X_0(\eps,S)>0$ such that $\abs{f}\le \eps g$ holds for all $X\ge X_0$.
When making estimates, we think of $m$, $F$, $w$ as fixed, but may still occasionally write $\ll_F$ (or similar) for emphasis.

We frequently use indicator notation, letting $\bm{1}_E\defeq 1$ if $E$ holds, and $\bm{1}_E\defeq 0$ otherwise.
For any nonempty set $S$ with an obvious measure (e.g.~the counting measure on a finite set, or the usual Haar measure on $\ZZ_p^m$), we let $\EE_{b\in S}[f(b)]$ denote the average of $f(b)$ over $b\in S$.

We let $\ZZ_{\ge 0}\defeq \set{a\in \ZZ: a\ge 0}$, and similarly define sets like $\ZZ_{\ne 0}$, $\RR_{>1}$, $\RR_{\ge 2}$.
For $c\in \ZZ_{\ne 0}$, we let $v_p(c)$ denote the $p$-adic valuation of $c$.
For $n\in \ZZ_{\ge 1}$, we let $\varphi(n)$ denote the totient function,
$\omega(n)$ the number of distinct prime factors of $n$, and $\rad(n)$ the radical of $n$.

We let $C^\infty_c(\RR^s)$ (resp.~$C^\infty_c(\RR^s)\otimes \CC$) denote the set of smooth compactly supported functions $\RR^s\to \RR$ (resp.~$\RR^s\to \CC$).
For any function $f=f(\bm{u})$, we let $\Supp{f}$ denote the closure of $\set{\bm{u}: f(\bm{u})\neq 0}$ in the domain of $f$; so for instance, the left-hand side of \eqref{COND:clean-weight-condition-in-diagonal-case} equals $\Supp{w}$.

We let $e(t)\defeq e^{2\pi it}$, and $e_r(t)\defeq e(t/r)$.
In integrals, we use notation analogous to summation notation.
For instance, we write $\int_X dx\,f(x)$ to mean $\int_X f(x)\,dx$ (in conventional notation), and we then write $\int_{X\times Y} dx\,dy\,f(x,y)$ to mean $\int_X dx\,(\int_Y dy\,f(x,y))$.

We need concise notation for $L^p$-norms and $\ell^p$-norms.
If $f$ is a quantity depending on a scalar or vector variable $t$ (and possibly also on other variables),
we write $\norm{f}_{L^p_t(S)}\defeq (\int_{t\in S} dt\, \abs{f}^p)^{1/p}$ or $\norm{f}_{\ell^p_t(S)}\defeq (\sum_{t\in S} \abs{f}^p)^{1/p}$ to denote the $p$-norm of $f$ over $t\in S$, according as $t$ is a continuous or discrete variable, respectively.
If the variable $t$
is clear from context, we may omit it.

We let $\partial_{u}\defeq \partial/\partial u$ for $u\in \RR$.
When doing calculus in $d$ dimensions, a \emph{multi-index} is a tuple of $d$ nonnegative integers (where $d$ will always be clear from context).
Given a multi-index $\bm{\alpha}\ge 0$, we let $\abs{\bm{\alpha}}$ denote the sum of the coordinates of $\bm{\alpha}$.
For a vector $\bm{u}\in \RR^s$, we let $\norm{\bm{u}}\defeq \max_{i}(\abs{u_i})$, and write $d\bm{u}\defeq du_1\cdots du_s$ and $\partial_{\bm{u}}^{\bm{\alpha}}\defeq \partial_{u_1}^{\alpha_1} \cdots \partial_{u_s}^{\alpha_s}$.
For example, if $f=f(\bm{u}) \in C^\infty_c(\RR^s)\otimes \CC$ and $k\in \ZZ_{\ge 0}$, then under our conventions,
\begin{equation*}
\max_{\abs{\bm{\alpha}}\le k}
\max_{\bm{u}\in \RR^s}{\abs{\partial_{\bm{u}}^{\bm{\alpha}}{f}}}
= \max_{\alpha_1,\dots,\alpha_s\ge 0:\, \alpha_1+\dots+\alpha_s\leq k}\,
\max_{u_1,\dots, u_s\in \RR}{\abs{\partial_{u_1}^{\alpha_1} \cdots \partial_{u_s}^{\alpha_s} f}}.
\end{equation*}

For (repeated) later use,
we fix a function $\nu_0\in C^\infty_c(\RR^m)$, supported on $[-2,2]^m$, with $0\le \nu_0\le 1$ everywhere, $\nu_0=1$ on $[-1/2,1/2]^m$, and $\int_{\RR^m} d\bm{x}\,\nu_0(\bm{x}) = 1$.
We also fix a radial\footnote{in the Euclidean sense (so that the value of $\nu_1(\bm{x})$ depends only on $x_1^2+\dots+x_m^2$)} function $\nu_1\in C^\infty_c(\RR^m)$, supported on the annulus $m\le x_1^2+\dots+x_m^2\le m^2$ (so that $\nu_1\vert_{[-1,1]^m} = 0$ and $\Supp{\nu_1}\belongs [-m,m]^m$), such that $\int_{\lambda>0} \frac{d\lambda}{\lambda}\,\nu_1(\bm{x}/\lambda) = 1$ for all $\bm{x}\in \RR^m\setminus \set{\bm{0}}$.

\section{Background on discriminants and delta}
\label{SEC:background-on-discriminants-and-the-delta-method}

Let $\disc(F)\neq 0$ be the discriminant of $F$.
Let $m_\ast\defeq m-3$.
Define the bi-homogeneous polynomial expression $\disc(F, \bm{c})$ as in \cite{wang2023dichotomous}*{\S3}, so that if $\bm{c}\in \CC^m$, then $\disc(F, \bm{c})\neq 0$ if and only if the variety $F(\bm{x}) = \bm{c}\cdot\bm{x} = 0$ in $\PP^{m-1}_\CC$ is smooth of dimension $m_\ast$.
Let
\begin{equation}
\label{EQN:define-Delta-as-product-of-disc}
    \Delta(\bm{c}) \defeq \disc(F) \cdot \disc(F, \bm{c});
\end{equation}
then for any $\bm{c}\in \ZZ^m$ and prime $p$ with $p\nmid \Delta(\bm{c})$, the varieties $F(\bm{x})=0$ and $F(\bm{x})=\bm{c}\cdot\bm{x}=0$ in $\PP^{m-1}_{\FF_p}$ are both smooth.
Now recall $\mcal{S}_0$, $\mcal{S}_1$ from \eqref{EQN:define-S_0,S_1-for-Delta-vanishing-and-nonzero-loci}.
For each $\bm{c}\in \mcal{S}_1$, let
\begin{equation}
\label{EQN:define-moduli-sets-N^c,N_c}
\mcal{N}^{\bm{c}} \defeq \set{n\geq 1: p\mid n\Rightarrow p\nmid \Delta(\bm{c})},
\quad
\mcal{N}_{\bm{c}} \defeq \set{n\geq 1: p\mid n\Rightarrow p\mid \Delta(\bm{c})}.
\end{equation}

\begin{lemma}
\label{LEM:count-R-N_c-infty-divisors}
Let $N, R\geq 1$ be integers.
Then $\card{\set{N\le n<2N: n\mid R^\infty}} \ll_\eps (R N)^\eps$, where $n\mid R^\infty$ means $\map{rad}(n)\mid R$.
In particular, if $\bm{c}\in \mcal{S}_1$, then
$\card{\set{N\le n<2N: n\in \mcal{N}_{\bm{c}}}}\ll_\eps \norm{\bm{c}}^\eps N^\eps$.
\end{lemma}

\begin{proof}
See e.g.~\cite{heath1998circle}*{antepenultimate display of p.~683}.
\end{proof}

As usual, we write $\grad{F} = (\partial F/\partial x_1,\dots,\partial F/\partial x_m)$.
It is known that the equation $\Delta(\bm{c}) = 0$ defines the projective dual variety of $V$ (see e.g.~\cite{wang2023dichotomous}*{Proposition~4.4}).
So
\begin{equation}
\label{INEQ:Disc-grad-composite-divisibility-by-F}
\Delta(\grad{F}(\bm{x}))/F(\bm{x}) \in \QQ[\bm{x}].
\end{equation}
Also, $\Delta(\bm{c})$ is irreducible in $\ol{\QQ}[c_1,\dots,c_m]$,
and has total degree $\deg{\Delta} = 3\cdot 2^{m-2}$ in $c_1,\dots,c_m$.
In particular, $\deg \Delta \ge 12$ (since $m\ge 4$), so by \cite{castryck2020dimension}*{Theorem~1}, we have
\begin{equation}
\label{INEQ:dimension-growth-bound-on-S_0}
    \card{\mcal{S}_0 \cap [-Z,Z]^m} \ll_m Z^{m-2}
\end{equation}
for all reals $Z\ge 1$.
(For diagonal $F$, more is known, e.g.~$\card{\mcal{S}_0 \cap [-Z,Z]^m} \ll_{m,\eps} Z^{m/2+\eps}$.)

One can express $\disc(F, \bm{c})$ in terms of the discriminant of a cubic form in $m-1$ variables:
\begin{equation}
\label{EQN:disc(F,c)-conversion}
    \disc(F,\bm{c})
    = \pm \disc(F(x_1,\dots,x_{m-1},-(c_1x_1+\dots+c_{m-1}x_{m-1})/c_m))\cdot c_m^{3\cdot 2^{m-2}}
\end{equation}
\cite{wang2023dichotomous}*{Proposition~3.2}.
This lets us bring into play a nice result of \cite{buse2014discriminant} (though the weaker classical result \cite{hooley1988nonary}*{(84) on p.~62} would also suffice for most of our needs):

\begin{theorem}
[\cite{buse2014discriminant}*{Corollary~4.30}]
\label{THM:Buse-Jouanolou-saturated-discriminant-gradient-square-decomposition}
Let $R$ be a ring and let $f\in R[x_1,\dots,x_{m-1}]$ be a homogeneous polynomial.
If $1\le a\le m-1$, then there exists $N\ge 0$ such that $x_a^N \disc(f)$ lies in the homogeneous ideal of $R[x_1,\dots,x_{m-1}]$ generated by
\begin{equation*}
f,\; \set{(\partial_{x_i}{f}) (\partial_{x_j}{f}): 1\leq i,j\leq m-1}.
\end{equation*}
\end{theorem}

Via \eqref{EQN:disc(F,c)-conversion}, Theorem~\ref{THM:Buse-Jouanolou-saturated-discriminant-gradient-square-decomposition} has the following corollary, useful in \S\ref{SEC:new-bounds-on-bad-sums-S}:

\begin{corollary}
\label{COR:Delta(c)-in-double-saturation}
If $1\le a,b\le m$, then there exists $N\ge 0$ such that $x_a^N c_b^N \Delta(\bm{c})$ lies in the homogeneous ideal of $\ZZ[x_1,\dots,x_m,c_1,\dots,c_m]$ generated by
\begin{equation*}
F(\bm{x}),\; \bm{c}\cdot \bm{x},\;
\set{(c_k\cdot \partial_{x_i}{F} - c_i\cdot \partial_{x_k}{F})(c_k\cdot \partial_{x_j}{F} - c_j\cdot \partial_{x_k}{F}): 1\le i,j,k\le m}.
\end{equation*}
\end{corollary}

\begin{proof}
Assume $b=m$.
Then use \eqref{EQN:define-Delta-as-product-of-disc}, \eqref{EQN:disc(F,c)-conversion}, Theorem~\ref{THM:Buse-Jouanolou-saturated-discriminant-gradient-square-decomposition} with $R = \ZZ[c_1/c_m,\dots,c_{m-1}/c_m]$,
and the congruence $-(c_1x_1+\dots+c_{m-1}x_{m-1})/c_m\equiv x_m\bmod{(\bm{c}/c_m)\cdot \bm{x}}$ in $R[x_1,\dots,x_m]$.
\end{proof}

Fix $w\in C^\infty_c(\RR^m)$ satisfying \eqref{COND:clean-weight-condition-in-diagonal-case} if $F$ is diagonal, or satisfying
\begin{equation}
\label{COND:clean-weight-condition-in-general}
    \Supp{w}\belongs \set{\bm{x}\in \RR^m: \det(\map{Hess}{F}(\bm{x}))\neq 0}
\end{equation}
in general.\footnote{For convenience, we will maintain this hypothesis for the rest of the paper, except when specified otherwise.}
Recall $N_{F,w}(X)$ from \eqref{EQN:define-N_w(X),N(X)}.
The delta method of \cites{duke1993bounds,heath1996new} allows one to express $N_{F,w}(X)$, up to a negligible error, as a sum over $\bm{c}\in \ZZ^m$ of ``adelic'' data.
We use the precise setup from \cite{wang2023_isolating_special_solutions}*{\S1} (based on \cites{duke1993bounds,heath1996new,heath1998circle}).

Let $h\maps (0,\infty) \times \RR \to \RR$ be the smooth function given by \cite{heath1998circle}*{(2.3)}; we will need the full definition later, in \S\ref{SEC:new-bounds-on-integral-J}.
For each real $X\ge 1$, let $Y\defeq X^{3/2}$.
Let
\begin{align}
I_{\bm{c}}(n)
&\defeq \int_{\bm{x}\in \RR^m}d\bm{x}\,w(\bm{x}/X)
h(n/Y,F(\bm{x})/Y^2)
e(-\bm{c}\cdot\bm{x}/n),
\label{EQN:define-I_c(n)} \\
S_{\bm{c}}(n)
&\defeq \sum_{1\leq a\leq n:\, \gcd(a,n)=1}\,
\sum_{1\leq \bm{x}\leq n} e_n(aF(\bm{x}) + \bm{c}\cdot\bm{x}),
\label{EQN:define-S_c(n)} \\
S^\natural_{\bm{c}}(n) &\defeq n^{-(m+1)/2} S_{\bm{c}}(n),
\quad J_{\bm{c},X}(n) \defeq X^{-m}I_{\bm{c}}(n).
\label{EQN:define-normalized-S-tilde-and-J}
\end{align}
A simple rearrangement of \cite{wang2023_isolating_special_solutions}*{(1.3)} gives (for all $A>0$)
\begin{equation}
\label{EQN:delta-method}
(1+O_A(Y^{-A}))\cdot N_{F,w}(X)/X^{m-3}
= \sum_{n\geq 1} \sum_{\bm{c}\in \ZZ^m} n^{(1-m)/2} S^\natural_{\bm{c}}(n) J_{\bm{c},X}(n).
\end{equation}
The infinite sums in \eqref{EQN:delta-method} are essentially finite, by the following standard result:

\begin{proposition}
\label{PROP:basic-integral-facts}
For some constant $A_0=A_0(F,w)>0$, we have $J_{\bm{c},X}(n) = 0$ for all $n\ge A_0 Y$.
Also, $J_{\bm{c},X}(n)\ll_{\eps,A} X^{-A}$ holds whenever $\eps,A>0$ and $\norm{\bm{c}}\ge X^{1/2+\eps}$.
\end{proposition}

\begin{proof}
See e.g.~\cite{wang2023_isolating_special_solutions}*{Proposition~5.1}.
\end{proof}


We now recall some background (cf.~\cite{hooley1986HasseWeil}*{\S\S5--6}) on the sums $S_{\bm{c}}(n)$.
It is known that $S_{\bm{c}}(n)$, and thus $S^\natural_{\bm{c}}(n)$ too, is multiplicative in $n$.
The Dirichlet series
\begin{equation}
\label{EQN:define-Phi}
    \Phi(\bm{c},s)\defeq \sum_{n\ge 1} n^{-s} S^\natural_{\bm{c}}(n)
\end{equation}
thus has an Euler product.
At prime powers, $S_{\bm{c}}$ is related to certain point counts in projective space.
Given a prime power $q$,
let $\mcal{V}(\FF_q)$ be the set of $\FF_q$-points on the variety $F(\bm{x})=0$ in $\PP^{m-1}_{\FF_q}$,
let $\mcal{V}_{\bm{c}}(\FF_q)$ be the set of $\FF_q$-points on the variety $F(\bm{x})=\bm{c}\cdot\bm{x}=0$ in $\PP^{m-1}_{\FF_q}$,
and let
\begin{equation*}
    E_F(q) \defeq \card{\mcal{V}(\FF_q)} - \card{\PP^{1+m_\ast}(\FF_q)},
    \quad E_{\bm{c}}(q) \defeq \card{\mcal{V}_{\bm{c}}(\FF_q)} - \card{\PP^{m_\ast}(\FF_q)},
\end{equation*}
where $\card{\PP^d(\FF_q)} = (q^{d+1}-1)/(q-1)$.
Then let
\begin{equation}
\label{EQN:define-normalized-point-count-errors-E_F,E_c}
    E^\natural_F(q)\defeq q^{-(1+m_\ast)/2} E_F(q),
    \quad E^\natural_{\bm{c}}(q)\defeq q^{-m_\ast/2} E_{\bm{c}}(q).
\end{equation}
If $p\nmid\bm{c}$ (e.g.~if $p\nmid \Delta(\bm{c})$), then by \cite{wang2022thesis}*{Proposition~3.2.4}, we have
\begin{equation}
\label{EQN:rewrite-S_c(p)-via-E_c}
    S^\natural_{\bm{c}}(p) = E^\natural_{\bm{c}}(p) - p^{-1/2}E^\natural_F(p).
\end{equation}
If $p\nmid \Delta(\bm{c})$ and $l\geq 2$, then by \cite{wang2022thesis}*{Proposition~3.2.6}, we have
\begin{equation}
\label{EQN:S_c(p^l)-vanishing}
    S^\natural_{\bm{c}}(p^l)=0.
\end{equation}

Recall $\mcal{S}_0$, $\mcal{S}_1$ from \eqref{EQN:define-S_0,S_1-for-Delta-vanishing-and-nonzero-loci}.
For any set $\mcal{S}\belongs \ZZ^m$, let
\begin{equation}
\label{EQN:define-sum-Sigma(X,S)}
    \Sigma^\natural(X, \mcal{S})
    \defeq \sum_{n\geq 1} \sum_{\bm{c}\in \mcal{S}} n^{(1-m)/2} S^\natural_{\bm{c}}(n) J_{\bm{c},X}(n),
    \quad \Sigma(X, \mcal{S})\defeq X^{m-3}\Sigma^\natural(X, \mcal{S}).
\end{equation}
If $\bm{c}\in \mcal{S}_0$, then $\Phi(\bm{c},s)$ can resemble $\zeta(s-1/2)$ (in some sense), leading to the following result:

\begin{theorem}
[\cite{wang2023_isolating_special_solutions}]
\label{THM:Sigma(X,S_0)-theorem-black-box}
Suppose $m=6$ and $F$ is diagonal.
Then
\begin{equation}
\label{EQN:disc-locus-S_0-evaluation}
    \Sigma(X, \mcal{S}_0)
    = O_{F,w,\eps}(X^{2.75+\eps})
    + \mf{S}_F\cdot \sigma_{\infty,F,w}\cdot X^3
    + \sum_{L\in \Upsilon} \sum_{\bm{x}\in L\cap \ZZ^6} w(\bm{x}/X),
\end{equation}
unconditionally.
In particular, $\Sigma(X, \mcal{S}_0) \ll_{F,w} X^3$.
\end{theorem}

\begin{proof}
This follows from \cite{wang2023_isolating_special_solutions}*{Corollary~1.2 and (1.7)}, even if we relax the condition \eqref{COND:clean-weight-condition-in-general} to $\bm{0}\notin \Supp{w}$.
The earlier paper \cite{heath1998circle} had proved
$\Sigma(X, \mcal{S}_0) \ll_{F,w,\eps} X^{3+\eps}$.
\end{proof}



Since $N_{F,w}(X) = \Sigma(X, \mcal{S}_0) + \Sigma(X, \mcal{S}_1) + O_A(X^{-A})$,
we may thus concentrate on $\Sigma(X, \mcal{S}_1)$.
The rest of \S\ref{SEC:background-on-discriminants-and-the-delta-method} provides some technical context for our work (in comparison with Theorem~\ref{THM:Hoo-HB-6-var-cubic-conditional-3+eps-bound} due to \cites{hooley_greaves_harman_huxley_1997,heath1998circle}), but is not logically necessary for the paper.

If $\bm{c}\in \mcal{S}_1$, then by \eqref{EQN:rewrite-S_c(p)-via-E_c}, \eqref{EQN:S_c(p^l)-vanishing}, and \eqref{EQN:compute-lambda_V_c-p}, one might expect $\Phi(\bm{c},s)$ to resemble $L(s,V_{\bm{c}})^{(-1)^{m_\ast}}$,
up to a factor absolutely convergent for $\Re(s)>1/2$;
cf.~\cite{wang2023_large_sieve_diagonal_cubic_forms}*{(2.4)}.
With this intuition in mind, let us now recall how one can prove (as is key for Theorem~\ref{THM:Hoo-HB-6-var-cubic-conditional-3+eps-bound})
\begin{equation}
\label{INEQ:near-optimal-diagonal-GRH-bound-over-smooth-locus-S_1}
\Sigma(X, \mcal{S}_1) \ll_\eps X^{3(m-2)/4+\eps},
\quad\textnormal{or equivalently }\Sigma^\natural(X, \mcal{S}_1) \ll_\eps X^{(6-m)/4+\eps},
\end{equation}
under Conjecture~\ref{CNJ:(HW2)}~(HW2) for $L(s,V_{\bm{c}})$, when $m$ is even and $F$ is diagonal.
We find it illuminating to work in this generality, but key here is that $3(m-2)/4 = 3$ when $m=6$.

\begin{proof}
[Conditional proof sketch for \eqref{INEQ:near-optimal-diagonal-GRH-bound-over-smooth-locus-S_1}]
For this proof sketch only, let $Z\defeq X^{1/2+\eps_0}$.
Let
\begin{equation}
\label{EQN:define-Psi_1,Psi_2}
    \Psi^{\bm{c},1}(s) \defeq 1/L(s,V_{\bm{c}}),
    \quad \Psi^{\bm{c},2}(s) \defeq \Phi(\bm{c},s) L(s,V_{\bm{c}}),
\end{equation}
for $\bm{c}\in \mcal{S}_1$; then $\Phi = \Psi^{\bm{c},1}\Psi^{\bm{c},2}$.
Let $\mu_{\bm{c}}(n)$, $a^{!}_{\bm{c}}(n)$ be the $n$th coefficients of the Dirichlet series $\Psi^{\bm{c},1}$, $\Psi^{\bm{c},2}$, respectively.
Then the following hold (see e.g.~\cite{wang2023_large_sieve_diagonal_cubic_forms}*{Proposition~4.12}):
\begin{enumerate}[{label=[B\arabic*']}]
    \item For $N\leq Z^3$,
    we have $\sum_{\bm{c}\in \mcal{S}_1\cap [-Z,Z]^m} \sum_{N\le n<2N} \abs{a^{!}_{\bm{c}}(n)} \ll_\eps Z^{m+\eps} N^{1/2}$.
    
    
    \item For $N\leq Z^3$,
    we have $\sum_{\bm{c}\in \mcal{S}_1\cap [-Z,Z]^m} (\sum_{N\le n<2N} \abs{a^{!}_{\bm{c}}(n)})^2 \ll_\eps Z^{m+\eps} N$.
    
\end{enumerate}
The arguments in \cites{hooley1986HasseWeil,hooley_greaves_harman_huxley_1997} and \cite{heath1998circle}
can then loosely be interpreted as
\begin{enumerate}[label=(H\arabic*)]
    \item using partial summation over $n\in [N, 2N)$ to ``factor out'' $J_{\bm{c},X}(n)$ from the sum over $n$ in \eqref{EQN:define-sum-Sigma(X,S)} (for $\mcal{S}=\mcal{S}_1$),
    and then bounding the $J$-contribution in $\ell^\infty_n([N, 2N))$;
    
    \item expanding $S^\natural_{\bm{c}} = \mu_{\bm{c}}\ast a^{!}_{\bm{c}}$ using $\Phi = \Psi^{\bm{c},1}\Psi^{\bm{c},2}$;
    
    \item using GRH to bound the $\Psi^{\bm{c},1}$-contribution in $\ell^\infty_{\bm{c}}(\mcal{S}_1\cap [-Z,Z]^m)$;
    and
    
    \item using [B1'] afterwards, to bound the $\Psi^{\bm{c},2}$-contribution in $\ell^1_{\bm{c}}(\mcal{S}_1\cap [-Z,Z]^m)$.
\end{enumerate}
Upon dyadic summation over $1\ll N\ll Y$, one gets $\Sigma^\natural(X, \mcal{S}_1) \ll_{\eps_0} X^{(6-m)/4+O(\eps_0)}$.

In place of GRH, one could use an elementary $\ell^2$ statement in the spirit of a large sieve inequality.
One would then use [B2'] instead of [B1'].
(See \cite{wang2023_large_sieve_diagonal_cubic_forms}.)
\end{proof}

We now diagnose
(and sketch ``cures'' for)
the key ``sources of $\eps$'' above:
\begin{enumerate}
    \item In (H3),
    the \emph{pointwise} GRH bound $\sum_{N\le n<2N} \mu_{\bm{c}}(n)\ll_\eps \norm{\bm{c}}^\eps N^{1/2+\eps}$.
    (\emph{Cure}: RMT-type predictions such as Conjectures~\ref{CNJ:(R2')} and~\ref{CNJ:(RA1o)}.)
    
    \item In (H4),
    the bound [B1'].
    (Or [B2'],
    for the argument of \cite{wang2023_large_sieve_diagonal_cubic_forms}.)
    In fact, upon closer inspection,
    the proofs of [B1'] and [B2'] each have \emph{two} sources of $\eps$:
    \begin{enumerate}
        \item \emph{Good} prime factors $p\nmid \Delta(\bm{c})$ of $n$,
        via the ``first-order error'' present in $\Psi^{\bm{c},2}$.
        (\emph{Cure}: Replacing $\Psi^{\bm{c},1}$ with a ``better approximation'' of $\Phi$.)
        
        \item \emph{Bad} prime factors $p\mid \Delta(\bm{c})$ of $n$,
        via the \emph{sometimes large} failure of square-root cancellation
        in \emph{individual} sums of the form $S_{\bm{c}}(p^l)$.
        (\emph{Cure}: New and old pointwise bounds on $\abs{S_{\bm{c}}(p^l)}$,
        and the \emph{average}-type Conjecture~\ref{CNJ:(SFSCp)}.)
    \end{enumerate}
    
    \item The fact that over $1\ll N\ll Y$,
    each dyadic range $N\le n<2N$ contributes \emph{roughly equally} to the final bound \eqref{INEQ:near-optimal-diagonal-GRH-bound-over-smooth-locus-S_1};
    cf.~\cite{wang2023_large_sieve_diagonal_cubic_forms}*{Remark~5.3}.
    If unaddressed,
    then the following sources of $\eps$ would arise
    (in our work):
    \begin{enumerate}
        \item A fatal $\log{X}$ factor in our proof of Theorem~\ref{THM:level-1-positive-density},
        via summation over $N$.
        
        \item A large contribution from \emph{relatively small} $n$
        in our proof of Theorems~\ref{THM:level-2-almost-all-integers} and~\ref{THM:level-3-power-saving},
        when handling variation of $J_{\bm{c},X}(n)$ over $\bm{c}$---a step needed
        when applying Conjecture~\ref{CNJ:(RA1o)} or~\ref{CNJ:(RA1delta)} to beat pointwise GRH.
    \end{enumerate}
    (\emph{Cure}: New integral bounds
    that decay,
    as $n\to0$,
    fairly uniformly over $\bm{c}$.)
    
    \item The lack of ``$\eps$-care''
    in bounds and ``decay cutoffs'' for integrals.
    The best recorded integral estimates (valid at least for some $w$) seem to be
    \begin{equation}
    \label{INEQ:Hooley-best-Hessian-free-integral-bound-at-least-for-textbook-weights?}
    J_{\bm{c},X}(n),\,
    n\cdot \partial{J_{\bm{c},X}(n)}/\partial{n}
    \ll_A (1+\norm{\bm{c}}/X^{1/2})^{-A} \cdot f(2+X\norm{\bm{c}}/n),
    \end{equation}
    with $f(r) = r^{1-m/2} (\log{r})^m$, from \cite{hooley2014octonary}*{p.~252, (31)}.
    Summing over $n/X^{1-\eps} \le \norm{\bm{c}} \le X^{1/2}$ carefully,
    or summing over $X^{1/2} \le \norm{\bm{c}} \le X^{1/2+\eps}$ carelessly,
    incurs $\eps$-losses.
    (\emph{Cure}: Summing carefully over $\bm{c}$,
    and using the ``cure to (3)'' over small $n$.)
\end{enumerate}




\section{Background on individual \texpdf{$L$}{L}-functions}
\label{SEC:background-on-individual-L-functions}

\subsection{Geometric background}

We need to give a precise meaning to Conjecture~\ref{CNJ:(HW2)}.
We first define the necessary Hasse--Weil $L$-functions and their local factors, following \cite{serre1969facteurs}.
(Another option, not pursued here, would be to follow \cite{taylor2004galois}.)
This is technical, but allows us to capture ``variation in $p$'' in a representation-theoretic framework.
At most primes, the data captured is very concrete; see e.g.~\eqref{EQN:compute-lambda_V_c-p}.

For any perfect field $K$, let $G_K\defeq \Gal(\ol{K}/K)$.
Let $\Gamma_\RR(s)\defeq \pi^{-s/2}\Gamma(s/2)$ and $\Gamma_\CC(s)\defeq (2\pi)^{-s}\Gamma(s)$.
The case of the Riemann zeta function $\zeta(s)$ in \eqref{full-list-of-Hasse-Weil-L-functions} is familiar (with $\zeta_p(s) = (1-p^{-s})^{-1}$ and $\zeta_\infty(s) = \Gamma_\RR(s)$), so we focus on the other cases.
Let $\bm{c}\in \mcal{S}_1$.
Let $\ell$ be a prime, and (viewing $V$, $V_{\bm{c}}$ as subvarieties of $\PP^{m-1}$) consider the $\ell$-adic Galois representations
\begin{equation*}
    \rho_V \maps G_\QQ \to
    \frac{H^{1+m_\ast}(V \times_\QQ \ol{\QQ}, \QQ_\ell)}
    {H^{1+m_\ast}(\PP^{m-1}_{\ol{\QQ}}, \QQ_\ell)},
    \quad
    \rho_{V_{\bm{c}}} \maps G_\QQ \to
    \frac{H^{m_\ast}(V_{\bm{c}} \times_\QQ \ol{\QQ}, \QQ_\ell)}
    {H^{m_\ast}(\PP^{m-1}_{\ol{\QQ}}, \QQ_\ell)}.
\end{equation*}
It is known that
the representations $\rho_V$, $\rho_{V_{\bm{c}}}$, $\bigotimes^2 \rho_{V_{\bm{c}}}$, $\map{Sym}^2 \rho_{V_{\bm{c}}}$, $\bigwedge^2 \rho_{V_{\bm{c}}}$
of $G_\QQ$
have dimensions depending only on $m$, and are pure of weight $(1+m_\ast)/2$, $m_\ast/2$, $m_\ast$, $m_\ast$, $m_\ast$, respectively.

Let $\rho\maps G_\QQ\to M$ be one of these five representations, and let $d_\rho$, $w_\rho$ be the dimension and weight of $\rho$, respectively.
Define $\Gamma_\rho(s)$ using Hodge theory, following \cite{serre1969facteurs}*{\S3.2}
(after passing from $M$ to a singular cohomology group $M'/\CC$ independent of $\ell$);
then for certain integers $h_\rho(+), h_\rho(-), h_\rho(a,b) \ge 0$ with
$h_\rho(+)+h_\rho(-)+2\sum_{0\le a<b:\, a+b=w_\rho} h_\rho(a,b) = d_\rho$,
we have
\begin{equation*}
    \Gamma_\rho(s)
    = \Gamma_\RR(s - w_\rho/2)^{h_\rho(+)} \Gamma_\RR(s - w_\rho/2 + 1)^{h_\rho(-)}
    \prod_{0\le a<b:\, a+b=w_\rho} \Gamma_\CC(s - a)^{h_\rho(a,b)}.
\end{equation*}
Let $L_\infty(s, \rho) \defeq \Gamma_\rho(s + w_\rho/2)$;
then $L_\infty(s, \rho)$ is holomorphic on the half-plane $\Re(s)>0$.
Note that if we let $\bm{c}$, $\rho$ vary,
the number of possible functions $L_\infty(s, \rho)$ could be is $\ll_m 1$.

For each prime $p\neq \ell$, we may restrict $\rho$ to $G_{\QQ_p}$; let $P_p(T)\defeq \det(1-\pi_p T)\in \QQ_\ell[T]$ denote the reverse characteristic polynomial of geometric Frobenius on $M^{I_{\QQ_p}}$ (the inertia invariants of $M$), following \cite{serre1969facteurs}*{\S2.2}.
Clearly $\deg P_p \le d_\rho$, with equality if and only if $\rho$ is unramified at $p$.
Write $P_p(T) = \prod_{1\le j\le \deg P_p} (1-\alpha_{\rho,j}(p)T)$, and let
\begin{equation}
\label{EQN:define-normalized-eigenvalues-alpha-and-L-factor}
    \tilde{\alpha}_{\rho,j}(p)\defeq p^{-w_\rho/2} \alpha_{\rho,j}(p),
    \quad
    L_p(s, \rho)
    \defeq \prod_{1\le j\le \deg P_p} (1 - \tilde{\alpha}_{\rho,j}(p) p^{-s})^{-1}
    = P_p(p^{-s-w_\rho/2})^{-1}.
\end{equation}
Because $V$, $V_{\bm{c}}$ are complete intersections in $\PP^{m-1}_\QQ$,
it is now known\footnote{thanks to \cite{laskar2017local}*{Corollary~1.2 and its proof} (cf.~\cite{saito2003weight}*{Corollary~0.6}), which builds on progress of \cite{scholze2012perfectoid} on the weight-monodromy conjecture} that the polynomial $P_p(T)$ lies in $\QQ[T]$ and is independent of $\ell$ (so that $\set{\alpha_{\rho,j}(p)}_j$ is a multiset of algebraic numbers),
and furthermore (for any embedding of $\alpha_{\rho,j}(p)$ into the complex numbers) we have
\begin{equation}
\label{geometric-Ramanujan-bound-via-WMC}
\abs{\alpha_{\rho,j}(p)}\le p^{w_\rho/2},
\quad
\abs{\tilde{\alpha}_{\rho,j}(p)}\le 1,
\end{equation}
for all $p$, $j$.
One might also be able to directly (without automorphy) define a conductor and root number for $\rho$ independent of $\ell$ (following \cite{serre1969facteurs} and \cite{deligne1969constantes}), but we need not do so.

Let $L_\infty(s, V)\defeq L_\infty(s, \rho_V)$ for any $\ell$.
Given a prime $p$, let $L_p(s, V)\defeq L_p(s, \rho_V)$ for any $\ell\ne p$, and let $\tilde{\alpha}_{V,j}(p)\defeq \tilde{\alpha}_{\rho_V,j}(p)$ and $d_{V,p}\defeq \deg P_p$.
Let
\begin{equation*}
    L(s, V) \defeq \prod_{p<\infty} L_p(s, V)
    \eqdef \sum_{n\ge 1} \lambda^\natural_V(n) n^{-s},
\end{equation*}
so that $\lambda^\natural_V(p) = \sum_{1\le j\le d_{V,p}} \tilde{\alpha}_{V,j}(p)$ for all $p$.
Let $\deg L(s, V)\defeq \max_p d_{V,p} = d_{\rho_V}$.
Make analogous definitions for $L(s, V_{\bm{c}})$ and its tensor squares in \eqref{full-list-of-Hasse-Weil-L-functions}, in terms of $\rho_{V_{\bm{c}}}$ and its tensor squares.
If $p\nmid \Delta(\bm{c})$ and $d = \deg L(s, V_{\bm{c}})$, then by \eqref{EQN:define-normalized-point-count-errors-E_F,E_c}, \eqref{EQN:define-normalized-eigenvalues-alpha-and-L-factor}, smooth proper base change, and the Grothendieck--Lefschetz trace formula, we have (for instance)
\begin{align}
    \lambda^\natural_V(p)
    &= \sum_{1\le j\le \deg L(s,V)} \tilde{\alpha}_{V,j}(p)
    = (-1)^{1+m_\ast} E^\natural_F(p), \label{EQN:compute-lambda_V-p} \\
    \lambda^\natural_{V_{\bm{c}}}(p)
    &= \sum_{1\le j\le d} \tilde{\alpha}_{V_{\bm{c}},j}(p)
    = (-1)^{m_\ast} E^\natural_{\bm{c}}(p), \label{EQN:compute-lambda_V_c-p} \\
    \lambda^\natural_{V_{\bm{c}}, \map{Sym}^2}(p)
    &= \sum_{1\le i\le j\le d} \tilde{\alpha}_{V_{\bm{c}},i}(p) \tilde{\alpha}_{V_{\bm{c}},j}(p)
    = \lambda^\natural_{V_{\bm{c}}}(p^2), \label{EQN:compute-symmetric-square-lambda_V_c-p} \\
    \lambda^\natural_{V_{\bm{c}}, \bigwedge^2}(p)
    &= \sum_{1\le i<j\le d} \tilde{\alpha}_{V_{\bm{c}},i}(p) \tilde{\alpha}_{V_{\bm{c}},j}(p)
    = \lambda^\natural_{V_{\bm{c}}}(p)^2 - \lambda^\natural_{V_{\bm{c}}}(p^2), \label{EQN:compute-exterior-square-lambda_V_c-p} \\
    (-1)^{m_\ast} E^\natural_{\bm{c}}(p^2)
    &= \sum_{1\le j\le d} \tilde{\alpha}_{V_{\bm{c}},j}(p)^2
    = E^\natural_{\bm{c}}(p)^2 - 2\lambda^\natural_{V_{\bm{c}}, \bigwedge^2}(p), \label{EQN:compute-E_c-p^2}
\end{align}
where for all $j$ we have (by the Weil conjectures, since $p\nmid \Delta(\bm{c})$)
\begin{equation}
\label{EQN:smooth-case-Deligne-purity}
    \abs{\tilde{\alpha}_{V,j}(p)} = 1, \quad \abs{\tilde{\alpha}_{V_{\bm{c}},j}(p)} = 1.
\end{equation}

\begin{remark}
\label{RMK:extended-local-L-factor-definitions}
By working locally from the beginning,
one can define (for any integer $n\ge 1$) the quantities $\prod_{p\mid n} L_p(s,V_{\bm{c}})$ and $\lambda^\natural_{V_{\bm{c}}}(n)$
on all of $\set{\bm{c}\in\prod_{p\mid n}\ZZ_p^m:\Delta(\bm{c})\neq0}$.
Extended definitions like this will be convenient
for local calculations in \S\ref{SUBSEC:computing-local-averages} and \S\ref{SUBSEC:handling-variation-of-error-factors}.
\end{remark}

\subsection{Automorphic background}

We need some background on automorphic representations $\Pi$ of $\GL_d(\bd{A}_\QQ)$ for $d\ge 1$.
We will only work with
\emph{cuspidal} $\Pi$'s,
or more generally,
\emph{isobaric} $\Pi$'s.
These $\Pi$'s have
well-defined $L$-functions $L(s,\Pi)$,
and good formal properties
(due to Rankin, Selberg,
Langlands,
Godement, Jacquet, Shalika,
and others):
\begin{enumerate}
    \item If $\Pi$ is cuspidal,
    then $L(s,\Pi)$ is \emph{primitive} in the sense of
    \cite{farmer2019analytic}*{Lemma~2.4},
    and has certain familiar analytic properties \cite{farmer2019analytic}*{Theorem~3.1}.
    
    \item For each isobaric $\Pi$,
    there is a \emph{unique}
    multiset $\set{\Pi_1,\dots,\Pi_r}$,
    consisting of \emph{cuspidals},
    such that $L(s,\Pi)=L(s,\Pi_1)\cdots L(s,\Pi_r)$.
    We call the $\Pi_i$'s \emph{cuspidal constituents} of $\Pi$.
    
    \item \emph{Strong multiplicity one}:
    If $\Pi$, $\Pi'$ are isobaric,
    and $L_p(s,\Pi)=L_p(s,\Pi')$ for all but finitely many primes $p$,
    then $L(s,\Pi)=L(s,\Pi')$.
\end{enumerate}

Conjecture~\ref{CNJ:(HW2)} has a host of standard consequences (which may be treated as a black box).

\begin{proposition}
\label{PROP:HW2-consequences}
Let $\bm{c}\in \mcal{S}_1$.
Let $L(s)$ be one of the Hasse--Weil $L$-functions in \eqref{full-list-of-Hasse-Weil-L-functions}.
Assume Conjecture~\ref{CNJ:(HW2)} for $L(s)$ holds for some $d$, $\Pi$.
Then the following hold:
\begin{enumerate}
    \item $d = \deg L(s)$; in particular, $d \ll_m 1$.
    \label{ITEM:automorphic-degree-equals-Galois-degree}
    
    \item The conductor $q(\Pi)\in \ZZ_{\ge 1}$ of $\Pi$ satisfies $q(\Pi)\mid \Delta(\bm{c})^{O_m(1)}$.
    \label{ITEM:automorphic-conductor-divides-discriminant-power}
    
    \item Each cuspidal constituent of $\Pi$ has unitary, finite-order central character.
    \label{ITEM:unitary-finite-order-central-characters}
    
    \item $L(s)$ is holomorphic on $\CC$, except possibly for poles at $s = 1$ corresponding to trivial constituents of $\Pi$.
    \label{ITEM:holomorphic-except-trivial-constituents}
    
    \item $(s-1)^d L(s)$ is an entire function of order $1$.
    \label{ITEM:entire-function-of-order-1}
    
    \item $L(s)$ has a standard functional equation (with critical line $\Re(s) = 1/2$), involving $q(\Pi)$ and some root number $\varepsilon(\Pi)\in \set{z\in \CC: \abs{z}=1}$.
    \label{ITEM:functional-equation-critical-line-and-root-number}
    
    \item $L(s)$ has real coefficients, $\Pi$ is self-dual, and $\varepsilon(\Pi) \in \set{-1, +1}$.
    \label{ITEM:real-coefficients-self-dual-and-root-number}
    
    \item Let $\psi$ be a cuspidal or isobaric constituent of $\Pi$.
    Then $1/L(s,\psi) \ll_{m, \eps} \norm{\bm{c}}^\eps (1+\abs{s})^\eps$ for $\Re(s)\ge 1/2+\eps$.
    If $\mu_\psi(n)$ denotes the $n$th coefficient of the Dirichlet series $1/L(s,\psi)$, then $\sum_{1\le n\le N} \mu_\psi(n) n^{-it} \ll_{m, \eps} \norm{\bm{c}}^\eps (1+\abs{t})^\eps N^{1/2+\eps}$ for all $t\in \RR$ and $N\in \RR_{>0}$.
    \label{ITEM:1/L-bound-from-GRH}
\end{enumerate}
\end{proposition}

\begin{proof}
If $p$ is a prime, then $L_p(s,\Pi)$ has degree $\le d$, with equality if and only if $p\nmid q(\Pi)$; cf.~\cite{farmer2019analytic}*{Axiom~3(b) and (3.3)}.

(\ref{ITEM:automorphic-degree-equals-Galois-degree}):
Compare the degrees of $L_p(s)$, $L_p(s, \Pi)$ at a prime $p\nmid q(\Pi)\Delta(\bm{c})$.





(\ref{ITEM:automorphic-conductor-divides-discriminant-power}):
For some $\nu\in \ZZ$, the local Dirichlet polynomial $L_p(s-\nu/2)$ has rational coefficients for all primes $p$.
So by \cite{shin2014fields}*{(3.2)}, there exists $u\in \RR$ such that for all $p\nmid q(\Pi)$, the representation $\Pi_p \otimes (\abs{\cdot}^u\circ \det)$ of $\GL_d(\QQ_p)$ has field of rationality $\QQ$, in the sense of \cite{shin2014fields}*{Definition~2.2}.
So by strong multiplicity one, the field of rationality of $\Pi \otimes (\abs{\cdot}^u\circ \det)$ is $\QQ$.
Now consider any $p\mid q(\Pi)$.
Then $L_p(s, \Pi)$ has degree $<d$, so $L_p(s)$ has degree $<d$, whence $p\mid \Delta(\bm{c})$.
By \cite{shin2014fields}*{Lemmas~3.11 and~3.13}, then, $v_p(q(\Pi))\ll_d 1$.
So $q(\Pi)\mid \Delta(\bm{c})^{O_d(1)}\mid \Delta(\bm{c})^{O_m(1)}$.

(\ref{ITEM:unitary-finite-order-central-characters}):
Let $\Pi'$ be a cuspidal constituent of $\Pi$, so $\Pi'$ is a cuspidal automorphic representation of $\GL_{d'}(\bd{A}_\QQ)$ for some $d'\ge 1$.
Let $\omega'\maps \QQ^\times \backslash \bd{A}_\QQ^\times \to \CC^\times$ be the central character of $\Pi'$.
Then $\omega'$ corresponds to a classical character $n\mapsto \abs{n}^z \chi(n)$ on $\ZZ$, where $z\in \CC$ and $\chi$ is a Dirichlet character of conductor dividing $q(\Pi)$, such that $\omega'_p(p) = p^z \chi(p)$ for all primes $p\nmid q(\Pi)$.
But at each prime $p\nmid q(\Pi)$, if we write $L_p(s,\Pi)=\prod_{1\leq j\leq d'}(1-\alpha_{\Pi,j}(p)p^{-s})^{-1}$, then $\omega'_p(p) = \prod_{1\leq j\leq d'}\alpha_{\Pi,j}(p)$; cf.~\cite{farmer2019analytic}*{(3.3) and its proof}.
By \eqref{EQN:smooth-case-Deligne-purity} and the algebraicity of the eigenvalues $\tilde{\alpha}_j(p)$, it follows that for infinitely many primes $p$, we have $\abs{p^z} = 1$ and $p^z\in \ol{\QQ}$.
So $\Re(z)=0$, i.e.~$\omega'$ is unitary; and then $\Im(z)=0$ by the six exponentials theorem (cf.~\cite{farmer2019analytic}*{proof of Lemma~4.9}), so $\omega'$ has finite order.

(\ref{ITEM:holomorphic-except-trivial-constituents}), (\ref{ITEM:entire-function-of-order-1}), (\ref{ITEM:functional-equation-critical-line-and-root-number}):
Use (\ref{ITEM:unitary-finite-order-central-characters}) and results of Godement and Jacquet; cf.~\cite{farmer2019analytic}*{Theorem~3.1}.

(\ref{ITEM:real-coefficients-self-dual-and-root-number}):
For some $\nu\in \ZZ$, the coefficients of $L(s-\nu/2)$ are all rational.
Hence $L(s)$ has real coefficients.
So by (\ref{ITEM:unitary-finite-order-central-characters}) and strong multiplicity one, we have $L(s,\Pi)=L(s,\Pi^\vee)$ and thus $\Pi$ is self-dual.
The functional equation from (\ref{ITEM:functional-equation-critical-line-and-root-number}) then implies $\varepsilon(\Pi)\in \RR$, so $\varepsilon(\Pi) = \pm 1$.

(\ref{ITEM:1/L-bound-from-GRH}):
By \eqref{geometric-Ramanujan-bound-via-WMC}, (\ref{ITEM:entire-function-of-order-1}), (\ref{ITEM:functional-equation-critical-line-and-root-number}), and Conjecture~\ref{CNJ:(HW2)}(2), we have $1/L(s,\psi) \ll_{d, \eps} q(\psi)^\eps (1+\abs{s})^\eps$ for $\Re(s)\ge 1/2+\eps$; see e.g.~\cite{iwaniec2004analytic}*{Theorem~5.19 and the ensuing paragraph}.
But $q(\psi)\mid q(\Pi)\mid \Delta(\bm{c})^{O_m(1)}$, so the desired bound on $1/L(s,\psi)$ follows.
One can then prove $\sum_{1\le n\le N} \mu_\psi(n) n^{-it} \ll_{m, \eps} \norm{\bm{c}}^\eps (1+\abs{t})^\eps N^{1/2+\eps}$ using Perron's formula (\cite{iwaniec2004analytic}*{Proposition~5.54}) and contour integration; cf.~\cite{hooley1986HasseWeil}*{p.~75, proof of Lemma~10}.
\end{proof}

\section{Local control on polynomials and \texpdf{$L$}{L}-functions}
\label{SEC:local-control-on-polynomials-and-L-functions}

We first recall some standard bounds on the local near-zero loci of a fixed polynomial; cf.~\cite{serre1981quelques}*{p.~146, (57)} and \cite{ganzburg2001polynomial}.
Given $f\in \ZZ[x]$, let
$N(f;q)$ be the number of solutions $x\bmod{q}$ to $f(x)\equiv 0\bmod{q}$,
and let $\mu_\RR(f;\lambda)$ be the Lebesgue measure of the set $\set{x\in [-1,1]: \abs{f(x)}\le \lambda}$.
Similarly, for $P\in \ZZ[y_1,\dots,y_n]$, define
\begin{equation*}
\begin{split}
    N(P;q) &\defeq \#{\set{(y_1,\dots,y_n)\in (\ZZ/q\ZZ)^n: P(y_1,\dots,y_n)\equiv 0\bmod{q}}}, \\
    \mu_\RR(P;\lambda) &\defeq \vol{\set{(y_1,\dots,y_n)\in [-1,1]^n: \abs{P(y_1,\dots,y_n)}\le \lambda}}.
\end{split}
\end{equation*}



\begin{proposition}
Suppose $f\in \ZZ[x]$ has leading term $ax^d$ with $a\neq 0$ and $d\geq 1$.
Then
\begin{equation}
\label{INEQ:univariate-zero-density-mod-q}
    N(f;q) \ll_d \abs{a}^{1/d} q^{1 - 1/d},
\end{equation}
uniformly over integers $q\geq 1$.
Also, uniformly over reals $\lambda > 0$, we have
\begin{equation}
\label{INEQ:univariate-near-zero-real-density}
    \mu_\RR(f;\lambda) \ll_d \abs{a}^{-1/d} \lambda^{1/d}.
\end{equation}
\end{proposition}

\begin{proof}
The bound \eqref{INEQ:univariate-near-zero-real-density} goes back to P\'{o}lya (see e.g.~\cite{ganzburg2001polynomial}*{Theorem~1.1}).
Now let $h$ denote the greatest common divisor of $q$ and the coefficients of $f$.
The bound \eqref{INEQ:univariate-zero-density-mod-q} follows from \cite{konjagin1979number} if $h=1$, and then in general from the inequality $N(f;q) \le h\cdot N(f/h;q/h)$.
\end{proof}





\begin{corollary}
Fix a nonconstant polynomial $P\in \ZZ[y_1,\dots,y_n]$, where $n\ge 1$.
Then
\begin{equation}
\label{INEQ:multivarite-zero-density-mod-q}
    N(P;q)\ll_P q^{n - 1/\deg{P}},
\end{equation}
uniformly over integers $q\geq 1$.
Also, uniformly over reals $\lambda > 0$, we have
\begin{equation}
\label{INEQ:multivarite-near-zero-real-density}
    \mu_\RR(P;\lambda) \ll_P \lambda^{1/\deg{P}}.
\end{equation}
\end{corollary}

\begin{proof}
Let $d = \deg{P}$.
Given a matrix $A\in \GL_n(\QQ)$ with integral entries, let $Q = (\det{A})^d\cdot P\circ A^{-1} \in \ZZ[y_1,\dots,y_n]$.
Via the $\ZZ$-linear map $\bm{y}\mapsto A\bm{y}$,
we have $N(P;q)\ll_A N(Q;q)$ and $\mu_\RR(P;\lambda)\ll_A \mu_\RR(Q;\lambda)$.
By \cite{eisenbud1995commutative}*{p.~283, proof of Lemma~13.2.c}, we may choose $A$ so that $Q$ has $y_1$-leading term $a_1 y_1^d$, with $a_1\ne 0$.
Fix $y_2,\dots,y_n$; then $\#\set{y_1: Q(\bm{y})\equiv 0} \ll_{d, a_1} q^{1 - 1/d}$ by \eqref{INEQ:univariate-zero-density-mod-q}.
Summing over $y_2,\dots,y_n$ gives \eqref{INEQ:multivarite-zero-density-mod-q}.
Similarly, \eqref{INEQ:univariate-near-zero-real-density} implies \eqref{INEQ:multivarite-near-zero-real-density}.
\end{proof}

We now turn to local $L$-factors.
(For $\bm{c}\in \ZZ_p^m$ with $\Delta(\bm{c})\ne 0$, we define $L_p(s,V_{\bm{c}})$ using Remark~\ref{RMK:extended-local-L-factor-definitions}.)

\begin{proposition}
[\cite{kisin1999local}]
\label{PROP:ineffective-Krasner-lemma-of-Kisin}
Fix a prime $p$ and a tuple $\bm{b}\in \ZZ_p^m$ with $\Delta(\bm{b})\ne 0$.
Then there exists an integer $l\ge 0$, depending only on $p$ and $\bm{b}$, such that for all tuples $\bm{a}\in \ZZ_p^m$ with $\bm{a}\equiv \bm{b}\bmod{p^{1+l}}$, we have $\Delta(\bm{a})\ne 0$ and $L_p(s, V_{\bm{a}}) = L_p(s, V_{\bm{b}})$.
\end{proposition}

\begin{proof}
Let $S$ be the open subscheme $\Delta(\bm{c})\ne 0$ of $\Aff^m_{\QQ_p}$.
Let $W$ be the closed subscheme $F(\bm{x})=\bm{c}\cdot\bm{x}=0$ of $\PP^{m-1}_S$.
Let $\ell\ne p$ be a prime.
The maps $f_1\maps W\to S$ and $f_2\maps \PP^{m-1}_S\to S$ induce local systems $\mcal{L}_j \defeq R^{m_\ast}(f_j)_\ast(\ZZ_\ell)$ on $S$.
By \cite{kisin1999local}*{Theorem~5.1, case~(2), and its proof}, there exists a $p$-adic neighborhood $U$ of $\bm{b}$ in $S$ such that the Galois representations $\rho_{j,\bm{a}}\maps G_{\QQ_p}\to \Aut(\mcal{L}_{j,\bm{a}})$ for $\bm{a}\in U(\QQ_p)$ all factor through $\pi_1(U)$ in an appropriate sense.
So the isomorphism class of the representation $G_{\QQ_p} \to \Aut(\QQ_\ell\otimes (\mcal{L}_{1,\bm{a}}/\mcal{L}_{2,\bm{a}}))$ is constant over $\bm{a}\in U(\QQ_p)$, and thus $L_p(s, V_{\bm{a}})$ is too.
\end{proof}

\begin{lemma}
\label{LEM:density-form-of-ineffective-Krasner-lemma}
Fix a prime $p$.
For each $\bm{b}\in \ZZ_p^m$ with $\Delta(\bm{b})\ne 0$, let $l(p,\bm{b})\ge 0$ be the smallest integer for which the conclusion of Proposition~\ref{PROP:ineffective-Krasner-lemma-of-Kisin} holds.
Then for each integer $l\ge 0$, the set
\begin{equation}
\label{EXPR:p-adic-set-of-b-with-l(p,b)-at-most-l}
    \set{\bm{b}\in \ZZ_p^m: \Delta(\bm{b})\ne 0,\; l(p,\bm{b})\le l}
\end{equation}
is invariant under translation by any element of $p^{l+1}\ZZ_p^m$.
Furthermore, the measure of \eqref{EXPR:p-adic-set-of-b-with-l(p,b)-at-most-l} tends to $1$ as $l\to \infty$.
(Here we use the usual Haar measure on $\ZZ_p^m$.)
\end{lemma}

\begin{proof}
Suppose $\bm{b}, \bm{c}\in \ZZ_p^m$ with $\Delta(\bm{b})\ne 0$ and $\bm{c}\equiv \bm{b}\bmod{p^{1+l(p,\bm{b})}}$.
Then $\Delta(\bm{c})\ne 0$ and $L_p(s, V_{\bm{c}}) = L_p(s, V_{\bm{b}})$, and thus $l(p,\bm{c})\le l(p,\bm{b})$ (because for every $\bm{a}\equiv \bm{c}\bmod{p^{1+l(p,\bm{b})}}$, we have $\bm{a}\equiv \bm{b}\bmod{p^{1+l(p,\bm{b})}}$ and thus $\Delta(\bm{a})\ne 0$ and $L_p(s, V_{\bm{a}}) = L_p(s, V_{\bm{b}}) = L_p(s, V_{\bm{c}})$).
But then $\bm{b}\equiv \bm{c}\bmod{p^{1+l(p,\bm{c})}}$, so a similar argument gives $l(p,\bm{b})\le l(p,\bm{c})$, whence
\begin{equation}
\label{EQN:local-constancy-of-l(p,b)}
l(p,\bm{b}) = l(p,\bm{c}).
\end{equation}
Therefore, if $\bm{b}$ lies in \eqref{EXPR:p-adic-set-of-b-with-l(p,b)-at-most-l} for some $l\ge 0$, then \eqref{EXPR:p-adic-set-of-b-with-l(p,b)-at-most-l} indeed contains the set $\bm{b}+p^{l+1}\ZZ_p^m$.

Next,
let $A\in \ZZ_{\ge 0}$.
The set $S_A = \set{\bm{c}\in \ZZ_p^m: v_p(\Delta(\bm{c}))\le A}$ is closed in $\ZZ_p^m$, and thus compact.
The function $\bm{b}\mapsto l(p,\bm{b})$ on $S_A$ is locally constant (by \eqref{EQN:local-constancy-of-l(p,b)}), and thus has a (finite) maximum value.
Therefore, for all $l$ sufficiently large in terms of $A$, the set \eqref{EXPR:p-adic-set-of-b-with-l(p,b)-at-most-l} contains $S_A$.
But by \eqref{INEQ:multivarite-zero-density-mod-q} (or by \cite{serre1981quelques}*{p.~146, Corollaire}), the measure of $S_A$ tends to $1$ as $A\to \infty$.
\end{proof}

\begin{remark}
\label{RMK:EKL-implication-on-l}
If Conjecture~\ref{CNJ:(EKL)} holds, then $l(p,\bm{b})\le v_p(H(\bm{b}))$ (whenever $H(\bm{b})\ne 0$).
\end{remark}

\section{General separation technique}
\label{SEC:separation-lemmas}

At several points in the paper, we need to understand quantities that vary with $\bm{c}$ and $n$.
A key tool we use for this is smooth dyadic decomposition (minimizing convergence issues) followed by separation of variables (via Mellin inversion); see Lemma~\ref{LEM:dyadic-partial-Mellin-summation}.

Let $d^\times{r}\defeq dr/r$ and $\partial_{\log{r}}\defeq r\cdot \partial_r$ for $r\in \RR_{>0}$.
For any $k\in \ZZ_{\ge 1}$ and $\bm{r}\in \RR_{>0}^k$, let $d^\times{\bm{r}} = d^\times{r_1}\cdots d^\times{r_k}$ and $\partial_{\log{\bm{r}}}^{\bm{\alpha}}\defeq \partial_{\log{r_1}}^{\alpha_1}\cdots \partial_{\log{r_k}}^{\alpha_k}$.
Given $g\in C^\infty_c(\RR_{>0}^k) \otimes \CC$ and $\bm{s}\in \CC^k$, let
\begin{equation}
\label{EQN:define-multivariate-Mellin-transform}
g^\vee(\bm{s}) \defeq \int_{\bm{r}>0} d^\times\bm{r}\, g(\bm{r})\bm{r}^{\bm{s}}
= \int_{r_1,\dots,r_k>0} d^\times{r_1}\cdots d^\times{r_k}\, g(r_1,\dots,r_k) r_1^{s_1}\cdots r_k^{s_k},
\end{equation}
so that Mellin inversion (see e.g.~\cite{iwaniec2004analytic}*{p.~90, (4.106)}) gives (for all $\bm{\sigma}\in \RR^k$)
\begin{equation}
\label{EQN:Mellin-inversion-compact-support-case}
g(\bm{r}) = (2\pi)^{-k} \int_{\bm{t}\in \RR^k} d\bm{t}\,
g^\vee(\bm{\sigma}+i\bm{t})\cdot \bm{r}^{-\bm{\sigma}-i\bm{t}}.
\end{equation}

\begin{proposition}
[Standard Mellin bound]
\label{PROP:standard-general-Mellin-bound}
Fix a compact set $I\belongs\RR_{>0}$.
Let $k\in \ZZ_{\ge 1}$ and $(\bm{M}, \bm{s})\in \RR_{>0}^k \times \CC^k$.
Let $g\in C^\infty_c(\RR_{>0}^k) \otimes \CC$ with $\Supp{g}\belongs \prod_{1\le i\le k} (M_i\cdot I)$.
Then
\begin{equation}
\label{INEQ:standard-Mellin-bound-goal}
g^\vee(\bm{s}) \ll_{k,b} O_I(1)^k
\cdot \frac{\norm{g \circ \exp}_{b,\infty}}{(1+\norm{\bm{s}})^b}
\cdot \prod_{1\le i\le k} (O_I(1)^{\abs{\Re(s_i)}} M_i^{\Re(s_i)}),
\end{equation}
for all $b\in \ZZ_{\ge 0}$.
Here $\norm{g \circ \exp}_{b,\infty}\defeq \sum_{\abs{\bm{\alpha}}\leq b}
\sup_{\bm{r}\in \RR_{>0}^k}{\abs{\partial_{\log\bm{r}}^{\bm{\alpha}} g(\bm{r})}}$.
\end{proposition}

\begin{proof}
By \eqref{EQN:define-multivariate-Mellin-transform}, and the scale invariance of $d^\times{r_i}$, $\partial_{\log{r_i}}$ for each $i$,
we can reduce to the case where $M_1=\cdots=M_k=1$.
Now let $\vol(\log{I})\defeq\int_{r>0} d^\times{r}\,\bm{1}_{r\in I}<\infty$.
If $\norm{\bm{s}}\leq 1$, we may assume $b=0$.
Now in general, choose $i$ with $\abs{s_i} = \max(\abs{s_1},\dots,\abs{s_d})$.
Then on the right-hand side of \eqref{EQN:define-multivariate-Mellin-transform}, integrate by parts $b$ times in $\log{r_i}$,
to rewrite $g^\vee(\bm{s})$ as
\begin{equation*}
(-1)^b \int_{\bm{r}>0} d^\times{\bm{r}}\, \frac{\bm{r}^{\bm{s}}}{s_i^b} \partial_{\log{r_i}}^b g(\bm{r})
\ll \vol(\log{I})^k \frac{\norm{g \circ \exp}_{b,\infty}}{\abs{s_i}^b} \prod_{1\le j\le k} O_{\Supp{I}}(1)^{\abs{\Re(s_j)}}.
\end{equation*}
This suffices for \eqref{INEQ:standard-Mellin-bound-goal}.
\end{proof}

Fix a function $\nu_2\in C^\infty_c(\RR_{>0})$, supported on $[1, 2]$, with
\begin{equation}
\label{EQN:nu_2-squared-integrates-to-1}
    \int_{r>0} d^\times{r}\,\nu_2(r)^2 = 1.
\end{equation}

\begin{lemma}
[``Dyadic partial Mellin summation'']
\label{LEM:dyadic-partial-Mellin-summation}
Let $k\in \ZZ_{\ge 1}$.
Let $a\maps \ZZ_{\ge 1}^k \to \CC$ be a function.
Let $f\maps \RR_{>0}^k \to \CC$ be a smooth function supported on $\prod_{1\le j\le k} r_j \le A$ for some real $A\ge 1$.
Let $\nu=\nu_2$; let $\nu(\bm{r}/\bm{N}) \defeq \prod_{1\le j\le k} \nu(r_j/N_j)$ and $g_{\bm{N}}(\bm{r}) \defeq f(\bm{r}) \nu(\bm{r}/\bm{N})$.
Then
\begin{equation}
\label{EQN:dyadic-partial-Mellin-summation-goal}
    \sum_{\bm{n}\ge 1} a(\bm{n}) f(\bm{n})
    = (2\pi)^{-k} \int_{\bm{N}\in [1/2, \infty)^k} d^\times\bm{N}
    \int_{\bm{t}\in \RR^k} d\bm{t}\, g_{\bm{N}}^\vee(i\bm{t})
    \sum_{\bm{n}\ge 1} \nu(\bm{n}/\bm{N}) a(\bm{n}) \bm{n}^{i\bm{t}}.
\end{equation}
\end{lemma}

\begin{proof}
For all $\bm{n}\ge 1$, we have
$f(\bm{n}) = \int_{\bm{N}\ge 1/2} d^\times{\bm{N}}\, \nu(\bm{n}/\bm{N}) g_{\bm{N}}(\bm{n})$
by \eqref{EQN:nu_2-squared-integrates-to-1} (since $\nu(\bm{n}/\bm{N}) = 0$ for $\bm{N}<1/2$).
Also, $g_{\bm{N}}(\bm{r}) \in C^\infty_c(\RR_{>0}^k) \otimes \CC$ for each $\bm{N}\ge 1/2$.
So by \eqref{EQN:Mellin-inversion-compact-support-case}, we get
\begin{equation}
\label{EQN:pre-partial-Mellin-sumation}
    \sum_{\bm{n}\ge 1} a(\bm{n}) f(\bm{n})
    = \sum_{\bm{n}\ge 1} a(\bm{n})
    \int_{\bm{N}\ge 1/2} d^\times\bm{N}\, \nu(\bm{n}/\bm{N})
    (2\pi)^{-k} \int_{\bm{t}\in \RR^k} d\bm{t}\, g_{\bm{N}}^\vee(i\bm{t}) \bm{n}^{i\bm{t}}.
\end{equation}

If $\bm{r} \in \Supp{g_{\bm{N}}}$, then $\prod_j r_j\le A$ and $\bm{r} \in \prod_j [N_j, 2N_j]$, so $\prod_j N_j \le A$.
Therefore, there exist compact sets $\mcal{K}_1, \mcal{K}_2 \belongs \RR_{>0}^k$ such that if $(\bm{n}, \bm{N})\notin \mcal{K}_1 \times \mcal{K}_2$ in \eqref{EQN:pre-partial-Mellin-sumation}, then $\nu(\bm{n}/\bm{N}) \cdot g_{\bm{N}}^\vee(i\bm{t}) = 0$.
Furthermore, there exists a compact set $\mcal{K}_3 \belongs \RR_{>0}^k$ such that if $\bm{N}\in \mcal{K}_2$, then $\Supp{g_{\bm{N}}} \belongs \mcal{K}_3$.
So Proposition~\ref{PROP:standard-general-Mellin-bound} gives $g_{\bm{N}}^\vee(i\bm{t}) \ll_{f, \nu, b} (1+\norm{\bm{t}})^{-b}$ for all $b\ge 0$, uniformly over $\bm{N}\in \mcal{K}_2$.
Thus
$\sum_{\bm{n}\ge 1} \int_{\bm{N}\ge 1/2} d^\times\bm{N} \int_{\bm{t}\in \RR^k} d\bm{t}\,
\abs{a(\bm{n})} \cdot \abs{\nu(\bm{n}/\bm{N})} \cdot \abs{g_{\bm{N}}^\vee(i\bm{t})}$
is $\ll_{a, f, \nu, b} \sum_{\bm{n}\in \mcal{K}_1\cap \ZZ}
\int_{\bm{N}\in \mcal{K}_2} d^\times\bm{N}
< \infty$,
provided $b > k$.
Now \eqref{EQN:dyadic-partial-Mellin-summation-goal} follows from \eqref{EQN:pre-partial-Mellin-sumation} by Fubini.
\end{proof}

For later use, we now do a general dyadic calculation.

\begin{lemma}
\label{LEM:general-dyadic-sum-split-into-2-geometric-series}
Let $q,a,b\in \RR$ with $a\le q<b$.
Let $r_1,r_2,\tau \in \RR_{>0}$.
Then
\begin{equation*}
\sum_{r_1\le r\le r_2:\, \log_2(r)\in \ZZ}
\frac{r^q}{r^a + \tau^{b-a} r^b},
\quad
\int_{r_1\le r\le r_2} d^\times{r}\, \frac{r^q}{r^a + \tau^{b-a} r^b},
\end{equation*}
are both $\ll_{q,a,b} \min(\tau^{-1}, r_2)^{q-a} \cdot (\log(1+r_2/r_1))^{\bm{1}_{q=a}}$.
\end{lemma}

\begin{proof}
We may assume $r_1\le r_2$, or else the sum and integral both vanish.
By subtracting $q$, $a$, $b$ by $a$, we may also assume $a=0$.
If $r_2\le \tau^{-1}$, then the result follows from the bound $1+\tau^b r^b\ge 1$ and a geometric series (or corresponding integral).
If $r_2>\tau^{-1}$, then separately considering $r\le \tau^{-1}$ and $r>\tau^{-1}$ leads to the result.
\end{proof}

\section{Statistics of families of \texpdf{$L$}{L}-functions}
\label{SEC:statistics-of-L-function-families}

Throughout \S\ref{SEC:statistics-of-L-function-families}, assume $m$ is even.
Recall the eigenvalue and coefficient notation from \S\ref{SEC:background-on-individual-L-functions}.

\subsection{Computing local averages}
\label{SUBSEC:computing-local-averages}

For convenience, let $\tilde{\alpha}_{\bm{c},j}(p)\defeq \tilde{\alpha}_{V_{\bm{c}},j}(p)$ and $\lambda^\natural_{\bm{c}}(n)\defeq \lambda^\natural_{V_{\bm{c}}}(n)$.
Let $\mu_{\bm{c}}(n)$ be the $n$th coefficient of the Dirichlet series $1/L(s,V_{\bm{c}})$.
We have
\begin{equation}
\label{EQN:inverse-local-L-factor-1/L_p(s,V_c)}
1/L_p(s,V_{\bm{c}}) = {\textstyle \prod_{j} (1-\tilde{\alpha}_{\bm{c},j}(p)p^{-s})}
\end{equation}
by \eqref{EQN:define-normalized-eigenvalues-alpha-and-L-factor}.
So if $p\nmid \Delta(\bm{c})$, then \eqref{EQN:compute-lambda_V_c-p}, \eqref{EQN:compute-exterior-square-lambda_V_c-p}, \eqref{EQN:compute-E_c-p^2}, and $2\nmid m_\ast$ imply
\begin{equation}
\label{EQN:compute-1/L-Mobius-mu_c-p,p^2-coefficients}
\mu_{\bm{c}}(p)
= -\lambda^\natural_{\bm{c}}(p)
= E^\natural_{\bm{c}}(p),
\quad
\mu_{\bm{c}}(p^2)
= \lambda^\natural_{V_{\bm{c}},\bigwedge^2}(p)
= \tfrac12(E^\natural_{\bm{c}}(p)^2+E^\natural_{\bm{c}}(p^2)).
\end{equation}

The local statistics of $\lambda^\natural_{\bm{c}}(n)$, $\mu_{\bm{c}}(n)$ over $\bm{c}$
play a basic role in the global statistics of $L(s,V_{\bm{c}})$ over $\bm{c}$.
To prove that certain averages \emph{exist},
we will use \eqref{geometric-Ramanujan-bound-via-WMC} and Lemma~\ref{LEM:density-form-of-ineffective-Krasner-lemma}.
But to \emph{estimate} said averages, we will take a point-counting approach (though one could use monodromy groups instead; see e.g.~\cite{sarnak2016families}*{\S2.11}).
The result is Proposition~\ref{PROP:(LocAvSp)} below.



Let $\mcal{B}\belongs \RR^m$ be a region of the form $I_1 \times \cdots \times I_m$, where $I_1,\dots,I_m\belongs \RR$ are compact intervals of positive length.
Let $\bm{a}\in \ZZ^m$ and $n_0,n,n_1,n_2\in \ZZ_{\ge 1}$.
Let $\EE^{\bm{a},n_0}_{\bm{c}\in S}[f]$ be the average of $f$ over $\set{\bm{c}\in S: \bm{c}\equiv\bm{a}\bmod{n_0}}$ (assuming this set is nonempty).
Let $\EE^{\bm{a},n_0}_{1\le \bm{c}\le n}[f] \defeq \EE^{\bm{a},n_0}_{\bm{c}\in \set{1,2,\dots,n}^m}[f]$.

\begin{proposition}
[LocAv]
\label{PROP:(LocAvSp)}
The following two limits exist, and are independent of $\mcal{B}$:
\begin{equation*}\begin{split}
\bar{\mu}_{F,1}^{\bm{a},n_0}(n)
&\defeq \lim_{Z\to \infty}
\EE^{\bm{a},n_0}_{\bm{c}\in \mcal{S}_1\cap Z\cdot \mcal{B}}
[\mu_{\bm{c}}(n)], \\
\bar{\mu}_{F,2}^{\bm{a},n_0}(n_1,n_2)
&\defeq \lim_{Z\to \infty}
\EE^{\bm{a},n_0}_{\bm{c}\in \mcal{S}_1\cap Z\cdot \mcal{B}}
[\mu_{\bm{c}}(n_1)\mu_{\bm{c}}(n_2)].
\end{split}\end{equation*}
The quantity $\bar{\mu}_{F,1}^{\bm{a},n_0}(n)$ is multiplicative in $n$:
if $\gcd(n,n')=1$, then
\begin{equation}
\label{EQN:1-point-mu-average-multiplicative}
\bar{\mu}_{F,1}^{\bm{a},n_0}(n)\bar{\mu}_{F,1}^{\bm{a},n_0}(n')
= \bar{\mu}_{F,1}^{\bm{a},n_0}(nn').
\end{equation}
The quantity $\bar{\mu}_{F,2}^{\bm{a},n_0}(n_1,n_2)$ is multiplicative in $(n_1,n_2)$:
if $\gcd(n_1n_2,n'_1n'_2)=1$, then
\begin{equation}
\label{EQN:2-point-mu-correlation-multiplicative}
\bar{\mu}_{F,2}^{\bm{a},n_0}(n_1,n_2)\bar{\mu}_{F,2}^{\bm{a},n_0}(n'_1,n'_2)
= \bar{\mu}_{F,2}^{\bm{a},n_0}(n_1n'_1,n_2n'_2).
\end{equation}

Now let $p$ be a prime, and let $l,l_1,l_2\ge 0$ be integers.
Then (uniformly over $p,l,l_1,l_2,\bm{a},n_0$)
\begin{equation}
\label{INEQ:GRC-bound-for-averages-mu-bar}
\bar{\mu}_{F,1}^{\bm{a},n_0}(p^l) \ll_\eps p^{l\eps},
\quad \bar{\mu}_{F,2}^{\bm{a},n_0}(p^{l_1},p^{l_2}) \ll_\eps p^{(l_1+l_2)\eps}.
\end{equation}
Furthermore, if $p\nmid n_0$, then
\begin{gather}
\bar{\mu}_{F,1}^{\bm{a},n_0}(p) = \lambda^\natural_V(p)p^{-1/2}+O(p^{-1}),
\quad
\bar{\mu}_{F,1}^{\bm{a},n_0}(p^2) = 1+O(p^{-1}), \label{EQN:LocAv1} \\
\bar{\mu}_{F,2}^{\bm{a},n_0}(p^l,1) = \bar{\mu}_{F,2}^{\bm{a},n_0}(1,p^l)
= \bar{\mu}_{F,1}^{\bm{a},n_0}(p^l),
\quad
\bar{\mu}_{F,2}^{\bm{a},n_0}(p,p) = 1+O(p^{-1}). \label{EQN:LocAv2}
\end{gather}
\end{proposition}

\begin{proof}

For convenience, define $\mu_{\bm{c}}(n) \bm{1}_{\Delta(\bm{c})\ne 0}$ to be $\mu_{\bm{c}}(n)$ if $\Delta(\bm{c})\ne 0$, and $0$ if $\Delta(\bm{c}) = 0$.
(Here we allow $\bm{c}\in \ZZ^m$, or more generally, $\bm{c}\in \prod_{p\mid n} \ZZ_p^m$.)
By \eqref{geometric-Ramanujan-bound-via-WMC} and \eqref{EQN:inverse-local-L-factor-1/L_p(s,V_c)}, we have
\begin{equation}
\label{INEQ:applied-GRC-to-mu_c(n)}
\mu_{\bm{c}}(p^l) \bm{1}_{\Delta(\bm{c})\ne 0}\ll_m 1,
\quad
\mu_{\bm{c}}(n) \bm{1}_{\Delta(\bm{c})\ne 0}\ll_\eps n^\eps.
\end{equation}
(In contrast, for $\lambda^\natural_{\bm{c}}$,
we have
$\lambda^\natural_{\bm{c}}(p^l) \bm{1}_{\Delta(\bm{c})\ne 0}\ll (l+1)^{O_m(1)} \ll_{m,\eps} p^{l\eps}$ and $\lambda^\natural_{\bm{c}}(n) \bm{1}_{\Delta(\bm{c})\ne 0}\ll_\eps n^\eps$.)

Since $\mu_{\bm{c}}(1) \bm{1}_{\Delta(\bm{c})\ne 0} = \bm{1}_{\Delta(\bm{c})\ne 0}$, we have $\bar{\mu}_{F,1}^{\bm{a},n_0}(1) = \bar{\mu}_{F,2}^{\bm{a},n_0}(1,1) = 1$ (since $\card{\mcal{S}_0\cap Z\cdot \mcal{B}} = o_{\mcal{B};Z\to\infty}(Z^m)$ by \eqref{INEQ:dimension-growth-bound-on-S_0}).
On the other hand, if $n\ge 2$ and $k\ge 1$, then by Lemma~\ref{LEM:density-form-of-ineffective-Krasner-lemma}, the quantity $\mu_{\bm{c}}(n) \bm{1}_{\Delta(\bm{c})\ne 0}$ depends only on $\bm{c}\bmod{n^k}$, unless $\bm{c}$ lies in one of $o_{n;k\to\infty}(n^{km})$ exceptional residue classes of $\ZZ^m$ (or $\prod_{p\mid n} \ZZ_p^m$) modulo $n^k$.
This, together with \eqref{INEQ:applied-GRC-to-mu_c(n)} and the Chinese remainder theorem, implies that (for any $n\ge 1$) the three quantities
\begin{equation*}
\bar{\mu}_{F,1}^{\bm{a},n_0}(n),
\quad \lim_{k\to\infty} \EE^{\bm{a},n_0}_{1\le \bm{c}\le n^k}
[\mu_{\bm{c}}(n) \bm{1}_{\Delta(\bm{c})\ne 0}],
\quad \prod_{p\mid n} \EE^{\bm{a},n_0}_{\bm{c}\in \ZZ_p^m}
[\mu_{\bm{c}}(p^{v_p(n)}) \bm{1}_{\Delta(\bm{c})\ne 0}]
\end{equation*}
all exist and equal one another.
Similarly, the following exist and equal one another:
\begin{equation*}
\bar{\mu}_{F,2}^{\bm{a},n_0}(n_1,n_2),
\quad \lim_{k\to\infty} \EE^{\bm{a},n_0}_{1\le \bm{c}\le (n_1n_2)^k}
[\mu_{\bm{c}}(n_1)\mu_{\bm{c}}(n_2) \bm{1}_{\Delta(\bm{c})\ne 0}],
\quad \prod_{p\mid n_1n_2} \EE^{\bm{a},n_0}_{\bm{c}\in \ZZ_p^m}
[\mu_{\bm{c}}(p^{v_p(n_1)})\mu_{\bm{c}}(p^{v_p(n_2)}) \bm{1}_{\Delta(\bm{c})\ne 0}].
\end{equation*}
This establishes the required existence, independence, and multiplicativity of limits.


We now turn to the required estimates.
First, \eqref{INEQ:GRC-bound-for-averages-mu-bar} follows from \eqref{INEQ:applied-GRC-to-mu_c(n)}.
Now assume $p\nmid n_0$; we must prove \eqref{EQN:LocAv1} and \eqref{EQN:LocAv2}.
But $p\nmid n_0$ implies $\set{\bm{c}\in \ZZ_p^m: \bm{c}\equiv \bm{a}\bmod{n_0}} = \ZZ_p^m$, so $\EE^{\bm{a},n_0}_{\bm{c}\in \ZZ_p^m}[f] = \EE_{\bm{c}\in \ZZ_p^m}[f]$ for any quantity $f$.
So by our $p$-adic interpretations (from the previous paragraph) of 
$\bar{\mu}_{F,1}^{\bm{a},n_0}(p^l)$, $\bar{\mu}_{F,1}^{\bm{a},n_0}(p^{l_1},p^{l_2})$, it remains to prove the following:
\begin{enumerate}
    \item $\EE_{\bm{c}\in\ZZ_p^m}[\mu_{\bm{c}}(p) \bm{1}_{\Delta(\bm{c})\ne 0}]
    = \lambda^\natural_V(p)p^{-1/2}
    + O(p^{-1})$.
    
    \item $\EE_{\bm{c}\in\ZZ_p^m}[\mu_{\bm{c}}(p^2) \bm{1}_{\Delta(\bm{c})\ne 0}]
    = 1 + O(p^{-1})$.
    
    \item $\EE_{\bm{c}\in\ZZ_p^m}[\mu_{\bm{c}}(p)^2 \bm{1}_{\Delta(\bm{c})\ne 0}]
    = 1 + O(p^{-1})$.
\end{enumerate}

We now prove (1)--(3).
By \cite{wang2023dichotomous}*{Corollary~1.7} (and our definition \eqref{EQN:define-normalized-point-count-errors-E_F,E_c}),
we have
\begin{align}
    \EE_{\bm{c}\in\FF_p^m}[E^\natural_{\bm{c}}(p)
    \bm{1}_{p\nmid \Delta(\bm{c})}]
    &= \lambda^\natural_V(p)p^{-1/2}
    + O(p^{-1}), \label{EQN:E_c(p)-average} \\
    \EE_{\bm{c}\in\FF_p^m}[E^\natural_{\bm{c}}(p^2)
    \bm{1}_{p\nmid \Delta(\bm{c})}]
    &= 1 + O(p^{-1}), \label{EQN:E_c(p^2)-average} \\
    \EE_{\bm{c}\in\FF_p^m}[E^\natural_{\bm{c}}(p)^2
    \bm{1}_{p\nmid \Delta(\bm{c})}]
    &= 1 + O(p^{-1}). \label{EQN:E_c(p)^2-average}
\end{align}
But for each $f\in \set{\mu_{\bm{c}}(p), \mu_{\bm{c}}(p^2), \mu_{\bm{c}}(p)^2}$,
we have $\EE_{\bm{c}\in\ZZ_p^m}[f \bm{1}_{\Delta(\bm{c})\ne 0}] = \EE_{\bm{c}\in\ZZ_p^m}[f \bm{1}_{p\nmid \Delta(\bm{c})}] + O(p^{-1})$, by \eqref{INEQ:applied-GRC-to-mu_c(n)} and Lang--Weil (for $\Delta(\bm{c}) \equiv 0\bmod{p}$).
After rewriting $f \bm{1}_{p\nmid \Delta(\bm{c})}$ using \eqref{EQN:compute-1/L-Mobius-mu_c-p,p^2-coefficients},
we conclude that \eqref{EQN:E_c(p)-average} implies (1),
that \eqref{EQN:E_c(p)^2-average} implies (3),
and that \eqref{EQN:E_c(p^2)-average}--\eqref{EQN:E_c(p)^2-average} imply (2).
\end{proof}

For the rest of \S\ref{SEC:statistics-of-L-function-families}, assume that $m$ is even and that Conjecture~\ref{CNJ:(HW2)} holds.

\subsection{The Sarnak--Shin--Templier framework}

Using Conjecture~\ref{CNJ:(HW2)} and
\eqref{EQN:E_c(p)-average}--\eqref{EQN:E_c(p)^2-average},
we can obtain useful statistical information on $L(s, V_{\bm{c}})$ over $\bm{c}\in \mcal{S}_1$.

Let $\Pi_{\bm{c}}$ be an isobaric automorphic representation over $\QQ$ corresponding to $L(s, V_{\bm{c}})$ in Conjecture~\ref{CNJ:(HW2)}.
Let $\mcal{S}_2$ be the set of $\bm{c}\in \mcal{S}_1$ for which $\Pi_{\bm{c}}$ is cuspidal, self-dual, and symplectic, in the sense of \cite{sarnak2016families}*{p.~533}.
For each $\bm{c}\in \mcal{S}_2$, the $L$-function $L(s,V_{\bm{c}},\bigwedge^2)$ has a pole at $s=1$, whence there exists an isobaric automorphic representation $\phi_{\bm{c},2}$ over $\QQ$ with
\begin{equation}
\label{EQN:factor-exterior-square-L-function}
    L(s,V_{\bm{c}},{\textstyle\bigwedge^2}) = \zeta(s) L(s,\phi_{\bm{c},2}).
\end{equation}

\begin{proposition}
\label{PROP:SST-framework}
Let $Z\in \RR_{\ge 1}$.
Then $\card{(\mcal{S}_1\setminus \mcal{S}_2)\cap [-Z,Z]^m} \ll_{m, \eps} Z^{m - 1/2 + \eps}$.
\end{proposition}

\begin{proof}
We want to show that the family $\bm{c}\mapsto\Pi_{\bm{c}}$ indexed by $\bm{c}\in \mcal{S}_1$ is \emph{essentially}
cuspidal, self-dual, and symplectic
(in the sense of \cite{sarnak2016families}*{p.~538, (i)--(iii)}),
with a power-saving exceptional set.
We follow the GRH strategy suggested in \cite{sarnak2016families}.

Let $\nu_0$ be as in \S\ref{SUBSEC:conventions}.
Fix $\eps\in (0, \frac12)$ and let $P = Z^{1-\eps}$.

By Proposition~\ref{PROP:HW2-consequences}(\ref{ITEM:real-coefficients-self-dual-and-root-number}), each $\Pi_{\bm{c}}$ is self-dual.
Let $\mcal{S}_{1.5}$ be the set of $\bm{c}\in \mcal{S}_1$ for which $\Pi_{\bm{c}}$ is cuspidal.
Then $L(s,V_{\bm{c}},\bigotimes^2)$ has a pole at $s=1$ of order exactly $1$ if $\bm{c}\in \mcal{S}_{1.5}$, and at least $2$ if $\bm{c}\in \mcal{S}_1\setminus \mcal{S}_{1.5}$; this follows from the theory of unramified Rankin--Selberg $L$-functions (cf.~\cite{farmer2019analytic}*{proof of Lemma~2.3}).
A calculation with $\frac{L'}{L}(s,V_{\bm{c}},\bigotimes^2)$ (using \eqref{geometric-Ramanujan-bound-via-WMC}, GRH, and \cite{iwaniec2004analytic}*{\S5.6's Exercise~6 and \S5.7's Theorem~5.15}) then yields
\begin{equation}
\label{INEQ:tensor-square-polar-lower-bound}
\sum_{\bm{c}\in \mcal{S}_1} \nu_0(\bm{c}/Z)
\sum_{p\leq P:\,p\nmid \Delta(\bm{c})} (\log{p})
\cdot (\lambda^\natural_{V_{\bm{c}},\bigotimes^2}(p) - 1)
\geq P \sum_{\bm{c}\in \mcal{S}_1} \nu_0(\bm{c}/Z)
(\bm{1}_{\bm{c}\notin \mcal{S}_{1.5}} + O_\eps(P^{-1/2+\eps})).
\end{equation}


On the other hand, $L(s,V_{\bm{c}},\map{Sym}^2)$ has a pole at $s=1$ if $\bm{c}\in \mcal{S}_{1.5}\setminus \mcal{S}_2$; this follows from \cite{sarnak2016families}*{p.~533}.
So a calculation with $\frac{L'}{L}(s,V_{\bm{c}},\map{Sym}^2)$ gives
\begin{equation}
\label{INEQ:symmetric-square-polar-lower-bound}
\sum_{\bm{c}\in \mcal{S}_1} \nu_0(\bm{c}/Z)
\sum_{p\leq P:\,p\nmid \Delta(\bm{c})} (\log{p})
\cdot \lambda^\natural_{V_{\bm{c}},\map{Sym}^2}(p)
\ge P \sum_{\bm{c}\in \mcal{S}_1} \nu_0(\bm{c}/Z)
(\bm{1}_{\bm{c}\in \mcal{S}_{1.5}\setminus \mcal{S}_2} + O_\eps(P^{-1/2+\eps})).
\end{equation}

But for all primes $p$ and tuples $\bm{c}\in \ZZ^m$ with $p\nmid \Delta(\bm{c})$, we have $\bm{c}\in \mcal{S}_1$, and we can use \eqref{EQN:compute-lambda_V_c-p}, \eqref{EQN:compute-symmetric-square-lambda_V_c-p}, \eqref{EQN:compute-exterior-square-lambda_V_c-p} to write $\lambda^\natural_{V_{\bm{c}},\bigotimes^2}(p) = \lambda^\natural_{\bm{c}}(p)^2 = E^\natural_{\bm{c}}(p)^2$ and $\lambda^\natural_{V_{\bm{c}},\map{Sym}^2}(p) = \frac12(E^\natural_{\bm{c}}(p)^2 - E^\natural_{\bm{c}}(p^2))$ (since $2\nmid m_\ast$).
So the left-hand side of \eqref{INEQ:tensor-square-polar-lower-bound} equals
\begin{equation*}
\sum_{p\leq P} (\log{p})
\sum_{\bm{c}\in \ZZ^m} \nu_0(\bm{c}/Z)
(E^\natural_{\bm{c}}(p)^2 - 1) \bm{1}_{p\nmid \Delta(\bm{c})}
= \sum_{p\leq P} (\log{p}) Z^m (0 + O(p^{-1}))
\ll_\eps Z^m P^\eps,
\end{equation*}
by Poisson summation and \eqref{EQN:E_c(p)^2-average}.
Similarly, by \eqref{EQN:E_c(p)^2-average} and \eqref{EQN:E_c(p^2)-average}, the left-hand side of \eqref{INEQ:symmetric-square-polar-lower-bound} is $\ll_\eps Z^m P^\eps$.
Thus from \eqref{INEQ:tensor-square-polar-lower-bound} we get $\card{(\mcal{S}_1\setminus \mcal{S}_{1.5})\cap [-Z/2,Z/2]^m} \ll_{m, \eps} Z^m P^{-1/2+\eps}$,
and from \eqref{INEQ:symmetric-square-polar-lower-bound} we get $\card{(\mcal{S}_{1.5}\setminus \mcal{S}_2)\cap [-Z/2,Z/2]^m} \ll_{m, \eps} Z^m P^{-1/2+\eps}$.
Now let $\eps\to 0$.
\end{proof}

\subsection{The Ratios Recipe}
\label{SUBSEC:applying-CFZ-ratios-conjecture-recipe}

By Proposition~\ref{PROP:(LocAvSp)} and \cite{sarnak2016families}*{pp.~534--535, Geometric Families and Remark~(i)},
the recipe \cite{conrey2008autocorrelation}*{\S5.1} makes sense for the family $\bm{c}\mapsto\Pi_{\bm{c}}$.
We will soon derive Conjecture~\ref{CNJ:(RA1o)} accordingly, along with the following:

\begin{conjecture}
[R2$o$]
\label{CNJ:(R2o)}
Let $A_{F,2}(s_1,s_2)$ be defined as in \S\ref{SUBSUBSEC:deriving-(R1)--(R2)} (in terms of $F$).
For each real $Z\ge 2$, let $\sigma(Z)$ be as in \eqref{EQN:sigma(Z)}, and write $s_j=\sigma(Z)+it_j$.
Then there exists a real $\hbar>0$ such that uniformly over reals $Z\ge 2$ and $t_1, t_2\in [-Z^\hbar, Z^\hbar]$, we have
\begin{equation}
\label{EQN:soft-R2-goal}
\sum_{\bm{c}\in \mcal{S}_1\cap [-Z,Z]^m}
\Phi^{\bm{c},1}(s_1) \Phi^{\bm{c},1}(s_2)
= \sum_{\bm{c}\in \mcal{S}_1\cap [-Z,Z]^m}
(1 + o_{Z\to \infty}(1)) \cdot A_{F,2}(s_1,s_2) \zeta(s_1+s_2).
\end{equation}
\end{conjecture}

Here $A_{F,2}(s_1,s_2)$ is an Euler product absolutely convergent for $\Re(s_1), \Re(s_2) > 1/3$.

Before proceeding, we make some remarks on our specific Ratios Conjectures.


\begin{remark}
The $L$-functions $L(s,V_{\bm{c}})$ are not all primitive,
as the recipe in \cites{conrey2005integral,conrey2008autocorrelation} requires.
But by Proposition~\ref{PROP:SST-framework} and GRH, there is no real difference (in Proposition~\ref{PROP:(LocAvSp)} and in our Ratios Conjectures) between $\mcal{S}_1$ and $\mcal{S}_2$ (on which each $L(s,V_{\bm{c}})$ is primitive).
\end{remark}

\begin{remark}
We do not order our families by conductor
(or by discriminant, for that matter).
We are indexing by different level sets,
as is natural for families like ours;
cf.~\cite{sarnak2016families}*{p.~535, Remark~(i); and p.~560, second paragraph after (25)}.
\end{remark}

\begin{remark}
In Conjectures~\ref{CNJ:(RA1o)}, \ref{CNJ:(RA1delta)}, and~\ref{CNJ:(R2o)},
we restrict $t$
to (comfortably) respect the constraint \cite{conrey2007applications}*{(2.11c)} on ``vertical shifts''.
But it would be reasonable to allow $t\in[-Z^A,Z^A]$ for arbitrarily large $A>0$;
cf.~\cite{bettin2020averages}*{p.~4, the sentence before Conjecture~2}.
\end{remark}

\subsubsection{Deriving (RA1)}
\label{SUBSUBSEC:deriving-(RA1)}

To derive Conjecture~\ref{CNJ:(RA1o)},
first use \eqref{EQN:define-Phi_1} to write $\Phi^{\bm{c},1}(s)$ in terms of $1/L(s,V_{\bm{c}})$,
and then replace each term
$L(s,V_{\bm{c}})^{-1} = \sum_{n\ge 1} \mu_{\bm{c}}(n)n^{-s}$
on the left-hand side of \eqref{EQN:soft-RA1-goal} with
its ``naive expected value over $\set{\bm{c}\in \mcal{S}_1\cap Z\cdot \mcal{B}_M(\bm{b}): \bm{c}\equiv\bm{a}\bmod{n_0}}$ as $Z\to \infty$'' (computed using Proposition~\ref{PROP:(LocAvSp)}), i.e.~the Dirichlet series
\begin{equation}
\label{EXPR:1/L-average-series-over-c}
\sum_{n\ge 1} \bar{\mu}_{F,1}^{\bm{a},n_0}(n) n^{-s}.
\end{equation}
It turns out that the series \eqref{EXPR:1/L-average-series-over-c} behaves much like $\zeta(2s) L(s+1/2,V)$, as we now explain.

Define $A_{F,1}^{\bm{a},n_0}(s)$ to be the product of \eqref{EXPR:1/L-average-series-over-c} and $\zeta(2s)^{-1} L(s+1/2,V)^{-1}$, so that \eqref{EXPR:1/L-average-series-over-c} factors as $A_{F,1}^{\bm{a},n_0}(s) \zeta(2s) L(s+1/2,V)$.
Let $\bar{a}_{F,1}^{\bm{a},n_0}(n)$ be the $n$th coefficient of the Dirichlet series $A_{F,1}^{\bm{a},n_0}(s)$.
Then $\bar{a}_{F,1}^{\bm{a},n_0}(n)$ is multiplicative in $n$ by \eqref{EQN:1-point-mu-average-multiplicative};
and by \eqref{INEQ:GRC-bound-for-averages-mu-bar}, \eqref{EQN:LocAv1} we have
\begin{equation}
\label{INEQ:a_F,1(n)-coefficient-bound}
    \bar{a}_{F,1}^{\bm{a},n_0}(n)\ll_\eps n^\eps,
    \quad
    \bar{a}_{F,1}^{\bm{a},n_0}(p) \bm{1}_{p\nmid n_0} \ll p^{-1},
    \quad
    \bar{a}_{F,1}^{\bm{a},n_0}(p^2) \bm{1}_{p\nmid n_0} \ll p^{-1/2}
\end{equation}
(because $\bar{a}_{F,1}^{\bm{a},n_0}(p) = \bar{\mu}_{F,1}^{\bm{a},n_0}(p) - \lambda^\natural_V(p)p^{-1/2}$
and $\bar{a}_{F,1}^{\bm{a},n_0}(p^2) = \bar{\mu}_{F,1}^{\bm{a},n_0}(p^2) - \bar{\mu}_{F,1}^{\bm{a},n_0}(p) \lambda^\natural_V(p)p^{-1/2} - 1 + O(p^{-1})$).
So $A_{F,1}^{\bm{a},n_0}(s)$ has an Euler product, and satisfies
\begin{equation}
\label{INEQ:easy-absolute-A_F,1-bound}
\abs{A_{F,1}^{\bm{a},n_0}(s)}
\le \sum_{n\ge 1} \frac{\abs{\bar{a}_{F,1}^{\bm{a},n_0}(n)}}{n^{1/3+\eps}}
\le \prod_{p} \biggl(1 + O(\bm{1}_{p\mid n_0}) + \frac{O(p^{-1})}{p^{1/3+\eps}}
+ \frac{O(p^{-1/2})}{p^{2/3+2\eps}} + \frac{O_\eps(p^{\eps})}{p^{1+3\eps}}\biggr)
\ll_\eps n_0^\eps
\end{equation}
for $\Re(s)\ge 1/3+\eps$ (for any $\eps>0$).
On the other hand, $\zeta(2s) L(s+1/2,V)$ has a pole of order $\ge 1$ at $s=1/2$, and thus is more dominant than $A_{F,1}^{\bm{a},n_0}(s)$ in \eqref{EXPR:1/L-average-series-over-c}.

In view of the above, the Ratios Recipe \cite{conrey2008autocorrelation}*{\S5.1} produces the conjecture
\begin{equation*}
\sum_{\substack{\bm{c}\in \mcal{S}_1\cap Z\cdot \mcal{B}_M(\bm{b}): \\ \bm{c}\equiv\bm{a}\bmod{n_0}}}
L(s,V_{\bm{c}})^{-1}
= \sum_{\substack{\bm{c}\in \mcal{S}_1\cap Z\cdot \mcal{B}_M(\bm{b}): \\ \bm{c}\equiv\bm{a}\bmod{n_0}}}
(1 + o_{M; Z\to \infty}(1)) \cdot A_{F,1}^{\bm{a},n_0}(s) \zeta(2s) L(s+1/2,V)
\end{equation*}
(for $n_0$, $t$, $s$ as in Conjecture~\ref{CNJ:(RA1o)}).
This rearranges (upon division by $\zeta(2s) L(s+1/2,V)$) to \eqref{EQN:soft-RA1-goal}, giving Conjecture~\ref{CNJ:(RA1o)}.
Furthermore, the fullest Ratios Conjectures include a power-saving error term, leading naturally to Conjecture~\ref{CNJ:(RA1delta)}.

Here $\zeta(2s)$, $L(s+1/2,V)$ are called \emph{polar factors} (or \emph{polar terms}).
As we will see shortly, an additional polar factor, $\zeta(s_1+s_2)$, arises in (R2).

\begin{remark}
The Ratios Recipe involves the approximate functional equation for $L$
when there are $L$'s in the numerator,
but not when there are only $L$'s in the denominator.
\end{remark}

\subsubsection{Deriving (R2)}
\label{SUBSUBSEC:deriving-(R1)--(R2)}

For Conjecture~\ref{CNJ:(R2o)},
use \eqref{EQN:define-Phi_1} to write $\Phi^{\bm{c},1}(s_j)$ in terms of $1/L(s_j,V_{\bm{c}})$,
and then replace each term
$L(s_1,V_{\bm{c}})^{-1}L(s_2,V_{\bm{c}})^{-1}
= \sum_{n_1,n_2\ge 1} \mu_{\bm{c}}(n_1)\mu_{\bm{c}}(n_2) n_1^{-s_1}n_2^{-s_2}$
on the left-hand side of \eqref{EQN:soft-R2-goal} with
its ``naive expected value over $\bm{c}\in \mcal{S}_1\cap [-Z,Z]^m$ as $Z\to \infty$'' (computed using Proposition~\ref{PROP:(LocAvSp)}), i.e.~the series
\begin{equation*}
\sum_{n_1,n_2\ge 1} \bar{\mu}_{F,2}^{\bm{0},1}(n_1,n_2) n_1^{-s_1}n_2^{-s_2},
\end{equation*}
which by \eqref{EQN:2-point-mu-correlation-multiplicative} factors as
$A_{F,2}(s_1,s_2)\zeta(s_1+s_2)\prod_{1\le j\le 2}(\zeta(2s_j)L(s_j+1/2,V))$ for some Euler product $A_{F,2}(s_1,s_2)$.
If $\Re(s_1),\Re(s_2)\ge 1/3+\eps$, then by \eqref{INEQ:GRC-bound-for-averages-mu-bar} and \eqref{EQN:LocAv2}, we have
\begin{equation}
\label{INEQ:easy-absolute-A_F,2-bound}
\abs{A_{F,2}(s_1,s_2)}
\le \sum_{n_1,n_2\ge 1} \abs{\bar{a}_{F,2}(n_1,n_2)} \, (n_1n_2)^{-1/3-\eps}
\ll_\eps 1.
\end{equation}
(The justification is similar to that for $A_{F,1}^{\bm{a},n_0}(s)$.
Note in particular that if $\bar{a}_{F,2}(n_1,n_2)$ is the $(n_1,n_2)$th coefficient of the double Dirichlet series $A_{F,2}(s_1,s_2)$, then $\bar{a}_{F,2}(p^l,1) = \bar{a}_{F,2}(1,p^l) = \bar{a}_{F,1}^{\bm{0},1}(p^l)$ and $\bar{a}_{F,2}(p,p) = \bar{\mu}_{F,2}^{\bm{0},1}(p,p) - 2\bar{\mu}_{F,1}^{\bm{0},1}(p)\lambda^\natural_V(p)p^{-1/2} - 1 + \lambda^\natural_V(p)^2p^{-1}$.)

Division by $\prod_{1\le j\le 2}(\zeta(2s_j)L(s_j+1/2,V))$ (cf.~\S\ref{SUBSUBSEC:deriving-(RA1)}) leads to Conjecture~\ref{CNJ:(R2o)}.

\subsection{From (R2) to (R2') and (R2'E)}
\label{SUBSEC:deriving-(R2')-and-(R2'E)}

Write $A_{F,2}(\bm{s})=A_{F,2}(s_1,s_2)$.
Conjecture~\ref{CNJ:(R2o)} implies, uniformly over $Z\ge 2$ and $\bm{t}\in \RR^2$, that for $s_j = \sigma(Z)+it_j$, we have (for all $\eps>0$)
\begin{equation}
\label{EQN:soft-R2-height-1-implication}
\sum_{\bm{c}\in \mcal{S}_1\cap [-Z,Z]^m}
\Phi^{\bm{c},1}(s_1) \Phi^{\bm{c},1}(s_2)
= o_{\eps;Z\to \infty}(Z^m (1+\norm{\bm{t}})^\eps)
+ \sum_{\bm{c}\in \mcal{S}_1\cap [-Z,Z]^m} A_{F,2}(\bm{s}) \zeta(s_1+s_2);
\end{equation}
this follows for $\norm{\bm{t}}\le Z^\hbar$ by Conjecture~\ref{CNJ:(R2o)}, and for $\norm{\bm{t}} > Z^\hbar$ by GRH (see Proposition~\ref{PROP:HW2-consequences}(\ref{ITEM:1/L-bound-from-GRH})).
Using \eqref{EQN:soft-R2-height-1-implication}, we proceed to derive Conjecture~\ref{CNJ:(R2')}~(R2').

\begin{proposition}
\label{PROP:R2-implies-R2'}
Assume Conjectures~\ref{CNJ:(HW2)} and~\ref{CNJ:(R2o)}.
Then Conjecture~\ref{CNJ:(R2')} holds.
\end{proposition}

\begin{proof}
Let $\conj{f}(s)\defeq \conj{f(\conj{s})}$.
Let $Z, N \ge 1$ with $N\le Z^3$.
The left-hand side of \eqref{INEQ:R2'-goal} is independent of $\sigma_0>1/2$, since each integrand is holomorphic on $\Re(s)>1/2$ by GRH.
Now shift contours to $\Re(s) = \sigma(Z)$,
and expand squares using self-duality of the $L$-functions in \eqref{EQN:define-Phi_1} (see Proposition~\ref{PROP:HW2-consequences}(\ref{ITEM:real-coefficients-self-dual-and-root-number})),
to equate the left-hand side of \eqref{INEQ:R2'-goal} with the quantity
\begin{equation*}
\Sigma_0 \defeq
\sum_{\bm{c}\in \mcal{S}_1\cap [-Z,Z]^m}
\int_{\sigma(Z)-i\infty}^{\sigma(Z)+i\infty} ds_1
\int_{\sigma(Z)+i\infty}^{\sigma(Z)-i\infty} ds_2\,
\Phi^{\bm{c},1}(s_1) \Phi^{\bm{c},1}(s_2) \cdot f(s_1)\conj{f}(s_2)N^{s_1+s_2}.
\end{equation*}
After switching the order of $\bm{c}$ and $\bm{s}$ in $\Sigma_0$ (using Fubini),
and plugging in the estimate \eqref{EQN:soft-R2-height-1-implication} and the bound $1+\norm{\bm{t}}\le (1+\abs{t_1})(1+\abs{t_2})$ for each $\bm{s}$, we get the estimate
\begin{equation}
\label{EQN:decompose-Sigma_0-in-R2'-under-R2}
    \Sigma_0 = o_{\eps;Z\to\infty}(Z^m \Sigma_1^2) + \Sigma_2,
\end{equation}
where $\Sigma_1 \defeq \int_{t\in \RR} dt\, (1+\abs{t})^\eps \cdot \abs{f(\sigma(Z)+it)} N^{\sigma(Z)}$ and
\begin{equation*}
\Sigma_2 \defeq \sum_{\bm{c}\in \mcal{S}_1\cap [-Z,Z]^m}
\int_{\sigma(Z)-i\infty}^{\sigma(Z)+i\infty} ds_1
\int_{\sigma(Z)+i\infty}^{\sigma(Z)-i\infty} ds_2\,
\zeta(s_1+s_2) A_{F,2}(\bm{s}) \cdot f(s_1)\conj{f}(s_2)N^{s_1+s_2}.
\end{equation*}
Here $N\le Z^{O(1)}$, so $\Sigma_1\ll N^{1/2} \sup_{0\le \sigma\le 2} \norm{(1+\abs{t})^\eps f(\sigma+it)}_{L^1_t(\RR)}$.
And by Cauchy--Schwarz,
\begin{equation}
\label{INEQ:norm-Cauchy-for-R2'}
    \norm{(1+\abs{t})^\eps f(\sigma+it)}_{L^1_t(\RR)}
    \le \norm{(1+\abs{t})^{-1/2-\eps}}_{L^2_t(\RR)}
    \cdot \norm{(1+\abs{t})^{1/2+2\eps} f(\sigma+it)}_{L^2_t(\RR)}.
\end{equation}

Now let $\delta=1/20$, and assume $Z$ is large enough that $1-\sigma(Z)-\delta\ge 1/3+\delta$.
Shifting $s_2$ (in $\Sigma_2$) from $\Re(s_2) = \sigma(Z)$ to $\Re(s_2) = 1-\sigma(Z)-\delta$ yields $\Sigma_2 = \Sigma_3+O_{F,\eps}(Z^m \Sigma_4)$, where
\begin{equation*}
\Sigma_3 \defeq \sum_{\bm{c}\in \mcal{S}_1\cap [-Z,Z]^m}
\int_{\sigma(Z)-i\infty}^{\sigma(Z)+i\infty} ds_1\,
(-2\pi i) A_{F,2}(s_1, 1-s_1) \cdot f(s_1)\conj{f}(1-s_1) N
\end{equation*}
comes from the residue of $\zeta(s_1+s_2)$ at $s_2=1-s_1$,
and where on $\Re(s_2) = 1-\sigma(Z)-\delta$ we use \eqref{INEQ:easy-absolute-A_F,2-bound} and the consequence $\zeta(s_1+s_2)\ll_{\delta,\eps} (1+\abs{t_1+t_2})^{\eps}$ of RH to be able to take
\begin{equation*}
\Sigma_4 \defeq \int_{\bm{t}\in \RR^2} dt_1\,dt_2\,
(1+\abs{t_1+t_2})^\eps \abs{f(\sigma(Z)+it_1) f(1-\sigma(Z)-\delta-it_2)} N^{1-\delta}.
\end{equation*}

Shifting $s_1$ (in $\Sigma_3$) to $\Re(s_1) = 1/2$ yields
$\Sigma_3
\ll_F Z^m N \norm{f(1/2+it)}_{L^2_t(\RR)}^2$.
Also, $\Sigma_4\ll N^{1-\delta} \sup_{0\le \sigma\le 2} \norm{(1+\abs{t})^\eps f(\sigma+it)}_{L^1_t(\RR)}^2$, since $1+\abs{t_1+t_2}\le (1+\abs{t_1})(1+\abs{t_2})$.
But $\Sigma_0 \le O(Z^m \Sigma_1^2) + \Sigma_3 + O_{F,\eps}(Z^m \Sigma_4)$, by \eqref{EQN:decompose-Sigma_0-in-R2'-under-R2}.
By \eqref{INEQ:norm-Cauchy-for-R2'} with $\eps=1/4$, we get \eqref{INEQ:R2'-goal}.
\end{proof}



As a stepping stone from Conjecture~\ref{CNJ:(R2')}~(R2') to Conjecture~\ref{CNJ:(R2'E')}~(R2'E'),
we now state (R2'E).
Let $a_{\bm{c},1}(n)$ be the $n$th coefficient of the Dirichlet series $\Phi^{\bm{c},1}(s)$.

\begin{conjecture}
[R2'E]
\label{CNJ:(R2'E)}
Fix a function $D\in C^\infty_c(\RR_{>0})$.
Let $t_0\in \RR$.
Let $Z, N \in \RR_{>0}$ with $N\le Z^3$.
Then for some real $A_3>0$ depending only on $F$, we have
\begin{equation}
\label{INEQ:R2'E-goal}
\sum_{\bm{c}\in \mcal{S}_1\cap [-Z,Z]^m}\,
\Bigl\lvert{
\sum_{n\ge 1} D(n/N) \cdot n^{-it_0}a_{\bm{c},1}(n)
}\Bigr\rvert^2
\ll_{F, D} (1+\abs{t_0})^{A_3} Z^m N.
\end{equation}
\end{conjecture}

\begin{proposition}
\label{PROP:R2'-implies-R2'E}
Assume Conjecture~\ref{CNJ:(R2')}.
Then Conjecture~\ref{CNJ:(R2'E)} holds.
\end{proposition}

\begin{proof}
The case $N<1$ is trivial,
so assume $N\ge 1$.
Let $\sigma_0=1.5$.
By \eqref{EQN:Mellin-inversion-compact-support-case},
$D(n/N) = (2\pi)^{-1} \int_{t\in \RR} dt\, D^\vee(\sigma_0+it) (N/n)^{\sigma_0+it}$.
Apply Conjecture~\ref{CNJ:(R2')} with $f(s) = D^\vee(s-it_0)$, and bound $f$ using Proposition~\ref{PROP:standard-general-Mellin-bound}, to get \eqref{INEQ:R2'E-goal} with $A_3 = 2$.
\end{proof}




\subsection{From (RA1) to (RA1'E)}

We now build on Conjectures~\ref{CNJ:(RA1o)} and~\ref{CNJ:(RA1delta)}, introducing flexible weights over $\bm{c}$.
We need some terminology on \emph{residue classes} of $\ZZ^m$.

\begin{definition}
\label{DEFN:residue-class-modulus-and-average-coefficient-terminology}
If $\mcal{R}=\bm{a}+q\ZZ^m\belongs \ZZ^m$ (where $q\ge 1$),
let $q_{\mcal{R}}\defeq q$ be the \emph{modulus} of $\mcal{R}$,
and let $A_{F,1}^{\mcal{R}}(s)\defeq A_{F,1}^{\bm{a},q}(s)$ and $\bar{a}_{F,1}^{\mcal{R}}(n)\defeq \bar{a}_{F,1}^{\bm{a},q}(n)$ (where $A_{F,1}^{\bm{a},q}(s)$, $\bar{a}_{F,1}^{\bm{a},q}(n)$ are defined as in \S\ref{SUBSUBSEC:deriving-(RA1)}).
Given a nonempty set $\mscr{S} = \set{\mcal{R}}$ of residue classes $\mcal{R}\belongs \ZZ^m$, let $Q(\mscr{S})\defeq \max_{\mcal{R}\in \mscr{S}}(q_{\mcal{R}})$.
\end{definition}

Let $\mscr{P}=\set{\mcal{R}}$ be a partition of $\ZZ^m$ into finitely many residue classes $\mcal{R}\belongs \ZZ^m$.
(In \S\ref{SUBSEC:handling-variation-of-error-factors},
we will construct the partitions needed for our main results.)
Let $I\belongs \RR_{>0}$ be a compact set.
Let $\nu=\nu_{\bm{c}}(r)$ be a smooth function $\RR^m\times \RR \to \CC,\,(\bm{c}, r)\mapsto \nu_{\bm{c}}(r)$ supported on $[-1,1]^m \times I$.
Let
\begin{equation}
\label{EQN:define-norm-M_1,k(nu)}
\mcal{M}_{1, k}=\mcal{M}_{1, k}(\nu)
\defeq \sum_{\abs{\bm{\alpha}}\le 1}\, \sum_{0\le j\le k}\,
\sup_{(\tilde{\bm{c}}, r) \in \RR^m\times \RR_{>0}}
\left\lvert
\partial_{\tilde{\bm{c}}}^{\bm{\alpha}}
\partial_{\log r}^j{\nu_{\tilde{\bm{c}}}(r)}
\right\rvert.
\end{equation}

\begin{conjecture}
[RA1$o$'E]
\label{CNJ:(RA1o'E)}
Let $Z, N \ge 2Q(\mscr{P})$ be reals with $N\le Z^3$.
Then the quantity
\begin{equation}
\label{EXPR:RA1'E-main-quantity}
\sum_{\mcal{R}\in \mscr{P}}\,
\Bigl\lvert
\sum_{\bm{c}\in \mcal{S}_1\cap \mcal{R}}
\sum_{n\ge 1} \nu_{\bm{c}/Z}(n/N)
\cdot (a_{\bm{c},1}(n) - \bar{a}_{F,1}^{\mcal{R}}(n))
\Bigr\rvert
\end{equation}
is $\ll_{F,I} Z^m N^{1/2} \cdot o_{F,Q(\mscr{P});Z\to\infty}(1) \cdot \mcal{M}_{1, A_4}$,
for some real $A_4=A_4(F)>0$.
\end{conjecture}

\begin{proposition}
\label{PROP:RA1o-implies-RA1o'E}
Assume Conjectures~\ref{CNJ:(HW2)}, \ref{CNJ:(R2')}, and~\ref{CNJ:(RA1o)}.
Then Conjecture~\ref{CNJ:(RA1o'E)} holds.
\end{proposition}

\begin{proof}

(It would be nice to prove this only assuming \eqref{EQN:soft-RA1-goal} for $t\in [-M, M]$, say, but it will be convenient to assume \eqref{EQN:soft-RA1-goal} for all $t\in [-\log{Z}, \log{Z}]$, as in Conjecture~\ref{CNJ:(RA1o)}.)

By \eqref{EQN:Mellin-inversion-compact-support-case}, we have $\nu_{\bm{c}/Z}(n/N) = (2\pi i)^{-1} \int_{(\sigma(Z))} ds\, \nu_{\bm{c}/Z}^\vee(s) (N/n)^{s}$, so that
\begin{equation*}
    \sum_{n\ge 1} \nu_{\bm{c}/Z}(n/N) (a_{\bm{c},1}(n) - \bar{a}_{F,1}^{\mcal{R}}(n))
    = (2\pi i)^{-1} \int_{(\sigma(Z))} ds\, \nu_{\bm{c}/Z}^\vee(s) N^{s} (\Phi^{\bm{c},1}(s) - A_{F,1}^{\mcal{R}}(s)).
\end{equation*}
Let $M \ge Q(\mscr{P})$ be a real parameter; soon below, we will let $M$ tend slowly to infinity as $Z\to \infty$.
Recall $\mcal{B}_M(\bm{b})$ from \eqref{EQN:define-box-B_M(b)}.
Weight $\nu$ is supported on $[-1,1]^m \times I$, so the triangle inequality and the previous display imply that \eqref{EXPR:RA1'E-main-quantity} is at most
\begin{equation*}
\Sigma_5 \defeq \sum_{\mcal{R}\in \mscr{P}} \sum_{\bm{b}\in [-M,M]^m}\,
\Bigl\lvert
\sum_{\bm{c}\in \mcal{S}_1\cap \mcal{R}\cap Z\cdot \mcal{B}_M(\bm{b})}
\int_{(\sigma(Z))} ds\, \nu_{\bm{c}/Z}^\vee(s) N^{s} (\Phi^{\bm{c},1}(s) - A_{F,1}^{\mcal{R}}(s))
\Bigr\rvert.
\end{equation*}
Uniformly over $\bm{c}\in \RR^m$ and $s\in \CC$, Proposition~\ref{PROP:standard-general-Mellin-bound} and \eqref{EQN:define-norm-M_1,k(nu)} give (for all $b\in \ZZ_{\ge 0}$)
\begin{equation}
\label{INEQ:quotable-nu_c-Mellin-bound}
\nu_{\bm{c}/Z}^\vee(s) \ll_b O_I(1)^{1+\abs{\Re(s)}} \mcal{M}_{1,b}(\nu) (1+\abs{\Im(s)})^{-b}.
\end{equation}

For each pair $(\mcal{R},\bm{b})$, choose an element $\bm{c}(\mcal{R},\bm{b})$ of $\mcal{S}_1\cap \mcal{R}\cap Z\cdot \mcal{B}_M(\bm{b})$, if such an element exists.
Let $\Sigma_6$ be $\Sigma_5$ with $\nu_{\bm{c}/Z}^\vee(s) - \nu_{\bm{c}(\mcal{R},\bm{b})/Z}^\vee(s)$ in place of $\nu_{\bm{c}/Z}^\vee(s)$,
and let $\Sigma_7$ be $\Sigma_5$ with $\nu_{\bm{c}(\mcal{R},\bm{b})/Z}^\vee(s)$ in place of $\nu_{\bm{c}/Z}^\vee(s)$.
Clearly $\Sigma_5\le \Sigma_6+\Sigma_7$.

We need to split $\Sigma_7$ further according to the size of $t$,
with some analytic care (keeping in mind the entireness hypothesis on $f$ in Conjecture~\ref{CNJ:(R2')}).
Let $B=B_M(s) \defeq e^{(s/M)^2}$.
The following hold uniformly over $\sigma$ in any fixed finite interval:
\begin{enumerate}
    \item For all $t\in \RR$, we have $B(\sigma+it)\ll 1$.
    \label{ITEM:B_M(s)-uniformly-bounded}
    
    \item If $\abs{t}>\log{Z}$, then $B(\sigma+it) \ll_M Z^{-1}$.
    \label{ITEM:large-t-decay-of-B_M(s)}
    
    \item If $\abs{t}\le M^{1/2}$, then $B(\sigma+it) = 1 + O(M^{-1})$.
    \label{ITEM:small-t-estimate-for-B_M(s)}
\end{enumerate}
Let $\Sigma_7(A)$ be $\Sigma_5$ with $\nu_{\bm{c}(\mcal{R},\bm{b})/Z}^\vee(s)\cdot A(s)$ in place of $\nu_{\bm{c}/Z}^\vee(s)$, so that
\begin{equation}
\label{INEQ:Sigma_7-main-soft-RA1'E-integral-decomposition}
    \Sigma_7
    \le \Sigma_7(B\cdot \bm{1}_{\abs{t}\le \log{Z}})
    + \Sigma_7(B\cdot \bm{1}_{\abs{t}>\log{Z}})
    + \Sigma_7(1-B).
\end{equation}

We first bound $\Sigma_7(B\cdot \bm{1}_{\abs{t}\le \log{Z}})$ and $\Sigma_7(B\cdot \bm{1}_{\abs{t}>\log{Z}})$.
Plugging (\ref{ITEM:B_M(s)-uniformly-bounded}), \eqref{INEQ:quotable-nu_c-Mellin-bound} (with $b=2$), and Conjecture~\ref{CNJ:(RA1o)} (for $t\in [-\log{Z}, \log{Z}]$) into $\Sigma_7(B\cdot \bm{1}_{\abs{t}\le \log{Z}})$, we find that
\begin{equation*}
    \Sigma_7(B\cdot \bm{1}_{\abs{t}\le \log{Z}}) \le
    \sum_{\mcal{R}\in \mscr{P}} \sum_{\bm{b}\in [-M,M]^m}
    \sum_{\bm{c}\in \mcal{S}_1\cap \mcal{R}\cap Z\cdot \mcal{B}_M(\bm{b})} O_I(\mcal{M}_{1,2}(\nu)) N^{\sigma(Z)} o_{M; Z\to\infty}(1),
\end{equation*}
which is in turn $\ll_I Z^m N^{1/2} o_{M; Z\to\infty}(1) \mcal{M}_{1,2}(\nu)$.
Similarly, plugging (\ref{ITEM:large-t-decay-of-B_M(s)}), \eqref{INEQ:quotable-nu_c-Mellin-bound} (with $b=2$), GRH (see Proposition~\ref{PROP:HW2-consequences}(\ref{ITEM:1/L-bound-from-GRH})), and \eqref{INEQ:easy-absolute-A_F,1-bound} into $\Sigma_7(B\cdot \bm{1}_{\abs{t}>\log{Z}})$ reveals that
\begin{equation*}
    \Sigma_7(B\cdot \bm{1}_{\abs{t}>\log{Z}})
    \ll_{M, \eps} \sum_{\bm{c}\in \mcal{S}_1\cap [-Z,Z]^m} Z^{-1} O_I(\mcal{M}_{1,2}(\nu)) N^{\sigma(Z)} (Z^\eps + Q(\mscr{P})^\eps),
\end{equation*}
which (if $\eps=1/2$, say) is $\ll_I Z^m N^{1/2} o_{M; Z\to\infty}(1) \mcal{M}_{1,2}(\nu)$ (since $Q(\mscr{P})\le Z$).

By choosing $M=M(Z)\ge Q(\mscr{P})$ appropriately, we may ensure both (i) that $M\to \infty$ as $Z\to \infty$, and (ii) that the two ``$o_{M; Z\to\infty}(1)$'' terms in the previous paragraph are $o_{Z\to\infty}(1)$.

It remains to bound $\Sigma_6$ and $\Sigma_7(1-B)$.
To handle both at once, we need a Fourier analog of Lemma~\ref{LEM:dyadic-partial-Mellin-summation}.
Let $f\maps \CC\to \CC$ be one of the functions
$f_{\bm{c},6}\maps s\mapsto \nu_{\bm{c}/Z}^\vee(s) - \nu_{\bm{c}(\mcal{R},\bm{b})/Z}^\vee(s)$ (given $(\mcal{R},\bm{b},\bm{c})$ with $\bm{c}\in \mcal{S}_1\cap \mcal{R}\cap Z\cdot \mcal{B}_M(\bm{b})$),
$f_{\bm{c},7}\maps s\mapsto \nu_{\bm{c}/Z}^\vee(s) \cdot (1-B(s))$ (given $\bm{c}\in \mcal{S}_1$).
Let
\begin{equation}
\label{EQN:define-entire-function-gt0(t)}
    g_{t_0}(t) = f(it) \nu_2^\vee(-i(t-t_0)) / 2\pi.
\end{equation}
Note that $f$, $\nu_2^\vee$ are entire (and rapidly decaying in vertical strips), so $g_{t_0}$ is entire (and rapidly decaying in horizontal strips).
By Parseval's theorem
and \eqref{EQN:nu_2-squared-integrates-to-1}, we have
\begin{equation*}
    \int_{\RR} dt\, \nu_2^\vee(it) \nu_2^\vee(-it)
    = 2\pi \int_{r>0} d^\times{r}\, \nu_2(r)^2 = 2\pi.
\end{equation*}
So $f(it) = \int_{\RR} dt_0\, g_{t_0}(t) \nu_2^\vee(i(t-t_0))$ by \eqref{EQN:define-entire-function-gt0(t)}.
By Fourier inversion applied to $g_{t_0}$, we get
\begin{equation*}
    f(it)
    = \int_{\RR^2} dt_0\,dx\, \hat{g}_{t_0}(x) e(xt) \nu_2^\vee(i(t-t_0))
    = \int_{\RR} dt_0 \int_{y>0} d^\times{y}\, \hat{g}_{t_0}(\tfrac{\log{y}}{2\pi}) y^{it} \nu_2^\vee(i(t-t_0)).
\end{equation*}
By analytic continuation, it follows that for any $\sigma>1/2$, the quantity
\begin{equation}
\label{EXPR:general-integral-form-of-error-terms-in-RA1'E}
    \int_{(\sigma)} ds\, f(s) N^{s} (\Phi^{\bm{c},1}(s) - A_{F,1}^{\mcal{R}}(s))
\end{equation}
equals $\int_{\RR} dt_0 \int_{y>0} d^\times{y}\, \hat{g}_{t_0}(\tfrac{\log{y}}{2\pi}) \int_{(\sigma)} ds\, y^s \nu_2^\vee(s-it_0) N^s (\Phi^{\bm{c},1}(s) - A_{F,1}^{\mcal{R}}(s))$,
and thus has absolute value at most
\begin{equation}
\label{EXPR:key-upper-bound-in-RA1'E-to-remove-c-dependence}
\int_{\RR} dt_0 \int_{y>0} d^\times{y}\, \abs{\hat{g}_{t_0}(\tfrac{\log{y}}{2\pi})}
\Bigl\lvert{
\int_{(\sigma)} ds\, \nu_2^\vee(s-it_0) (Ny)^s (\Phi^{\bm{c},1}(s) - A_{F,1}^{\mcal{R}}(s))
}\Bigr\rvert.
\end{equation}
But for all $z\in \CC$, we have $\hat{g}_{t_0}(\tfrac{\log{y}}{2\pi}) = \int_{\RR} dt\, g_{t_0}(t) y^{-it} = \int_{\RR} dt\, g_{t_0}(t+z) y^{-i(t+z)}$ (by the definition of $\hat{g}_{t_0}$, if $z=0$; and then by analytic continuation, in general), and thus $\abs{\hat{g}_{t_0}(\tfrac{\log{y}}{2\pi})}\le y^{\Im(z)} \int_{\RR} dt\, \abs{g_{t_0}(t+z)}$.
Using \eqref{EQN:define-entire-function-gt0(t)}, Proposition~\ref{PROP:standard-general-Mellin-bound} (cf.~\eqref{INEQ:quotable-nu_c-Mellin-bound}), and the definitions of $\mcal{M}_{1,b}$ (see \eqref{EQN:define-norm-M_1,k(nu)}) to bound $\abs{g_{t_0}(t+z)}$ pointwise, we then get (by taking $z=\pm iu\in i\RR$)
\begin{equation}
\label{INEQ:analytic-Fourier-transform-estimate-for-g_t_0}
\hat{g}_{t_0}(\tfrac{\log{y}}{2\pi})
\ll_{u,I,b,\nu_2,B} M^{-1} \mcal{M}_{1, b}(\nu) (1+\abs{t_0})^{-b} (y^u+y^{-u})^{-1}
\end{equation}
for all reals $u\ge 0$ and integers $b\ge 0$.
(If $f=f_{\bm{c},6}$, the factor of $M^{-1}$ in \eqref{INEQ:analytic-Fourier-transform-estimate-for-g_t_0} arises from the bound $\norm{\bm{c} - \bm{c}(\mcal{R},\bm{b})} \ll Z/M$.
If $f=f_{\bm{c},7}$, the factor of $M^{-1}$ comes from (\ref{ITEM:small-t-estimate-for-B_M(s)}) over $\abs{t}\le M^{1/2}$, and from the decay of $g_{t_0}(t)$ in $t$, $t-t_0$ if $\abs{t}>M^{1/2}$.)

Note that for all $\sigma\ge 1/3+\eps$, we may apply \eqref{INEQ:easy-absolute-A_F,1-bound} and Proposition~\ref{PROP:standard-general-Mellin-bound} (after shifting contours to $\Re(s) = 1/3+\eps$) to get (provided $0<\eps\le 1/12$)
\begin{equation}
\label{INEQ:contour-shift-bound-for-A_F,1-integral}
\int_{(\sigma)} ds\, \nu_2^\vee(s-it_0) (Ny)^s A_{F,1}^{\mcal{R}}(s)
\ll_\eps \int_{t\in \RR} dt\, \frac{(Ny)^{1/3+\eps} Q(\mscr{P})^\eps}{(1+\abs{t-t_0})^2}
\ll (Ny+Q(\mscr{P}))^{1/2}.
\end{equation}
In view of the bound \eqref{EXPR:key-upper-bound-in-RA1'E-to-remove-c-dependence} for \eqref{EXPR:general-integral-form-of-error-terms-in-RA1'E}, we may now apply \eqref{INEQ:analytic-Fourier-transform-estimate-for-g_t_0}, \eqref{INEQ:contour-shift-bound-for-A_F,1-integral}, Conjecture~\ref{CNJ:(R2')} (with $Z + (Ny)^{1/3}$, $Ny$, $\nu_2^\vee(s-it_0)$ in place of $Z$, $N$, $f(s)$), and Cauchy--Schwarz to get
\begin{equation*}
\Sigma_6 \ll_I \int_{\RR} dt_0 \int_{y>0} d^\times{y}\,
\frac{M^{-1} \mcal{M}_{1, b}(\nu)}{(1+\abs{t_0})^b (y^u+y^{-u})}
(Z + (Ny)^{1/3})^m (Ny + Q(\mscr{P}))^{1/2} (1+\abs{t_0}),
\end{equation*}
which is $\ll_I M^{-1} \mcal{M}_{1,b}(\nu) Z^m N^{1/2}$ by Lemma~\ref{LEM:general-dyadic-sum-split-into-2-geometric-series} (provided $u\ge m/3+1$ and $b\ge 3$, and $N\ge Q(\mscr{P})$).
The same holds for $\Sigma_7(1-B)$.
Thus Conjecture~\ref{CNJ:(RA1o'E)} holds with $A_4 = 3$.
\end{proof}

\begin{proposition}
[RA1$\delta$'E]
\label{PROP:RA1delta-implies-RA1delta'E}
Assume Conjectures~\ref{CNJ:(HW2)} and~\ref{CNJ:(RA1delta)}.
Suppose $Z,N\ge 2Q(\mscr{P})$ with $N\le Z^3$ and $Q(\mscr{P})\le Z^{\eta_2}$.
Then the quantity \eqref{EXPR:RA1'E-main-quantity} is
$\ll_{F, I} Z^{m-\eta_2} N^{1/2} \mcal{M}_{1, A_5}$.
Here $\eta_2$, $A_5$ are positive reals depending only on $F$.
\end{proposition}

\begin{proof}
Mimic the proof of Proposition~\ref{PROP:RA1o-implies-RA1o'E},
but take $M = Z^{\eta_1}$,
replace Conjecture~\ref{CNJ:(RA1o)} with~\ref{CNJ:(RA1delta)},
replace every use of Conjecture~\ref{CNJ:(R2')} with GRH,
and replace \eqref{INEQ:Sigma_7-main-soft-RA1'E-integral-decomposition} with the bound $\Sigma_7\le \Sigma_7(\bm{1}_{\abs{t}\le Z^{\eta_1}}) + \Sigma_7(\bm{1}_{\abs{t}>Z^{\eta_1}})$.
The details simplify, since the desired final bound is not sensitive to losses of $Z^\eps$.
(We can take $\eta_2 = \eta_1-\eps$ and $A_5 = 3$.)
\end{proof}

When applying Propositions~\ref{PROP:RA1o-implies-RA1o'E} and~\ref{PROP:RA1delta-implies-RA1delta'E} (in \S\ref{SUBSEC:proving-(RA1'E')}), we need the following lemma:

\begin{lemma}
\label{LEM:basic-main-term-bound-for-RA1'E'}
Fix a real $\theta>0$.
Suppose $\mcal{S}\belongs \ZZ^m$ and $\card{\mcal{S}\cap [-C,C]^m} \ll C^\theta$ for all reals $C\ge 1$.
Let $Z, N \ge 1$ be reals.
Then
\begin{equation}
\label{INEQ:RA1'E-main-term-bound-goal}
\sum_{\mcal{R}\in \mscr{P}}\,
\Bigl\lvert
\sum_{\bm{c}\in \mcal{S}\cap \mcal{R}}
\sum_{n\ge 1} \nu_{\bm{c}/Z}(n/N)
\bar{a}_{F,1}^{\mcal{R}}(n)
\Bigr\rvert
\ll_{F,I,\theta,\eps} Z^\theta N^{1/3+\eps} Q(\mscr{P})^\eps \mcal{M}_{1, 0}.
\end{equation}
\end{lemma}

\begin{proof}
By the triangle inequality, the left-hand side of \eqref{INEQ:RA1'E-main-term-bound-goal} is at most
\begin{equation*}
    \sum_{\mcal{R}\in \mscr{P}}
    \sum_{\bm{c}\in \mcal{S}\cap \mcal{R}}
    \left(\sup_{r>0}{\abs{\nu_{\bm{c}/Z}(r)}}\right)
    \sum_{n\in N\cdot I} \abs{\bar{a}_{F,1}^{\mcal{R}}(n)}
    \le \biggl(\,\sup_{\mcal{R}\in \mscr{P}} \sum_{n\in N\cdot I} \abs{\bar{a}_{F,1}^{\mcal{R}}(n)}\biggr)
    \sum_{\bm{c}\in \mcal{S}} \sup_{r>0}{\abs{\nu_{\bm{c}/Z}(r)}},
\end{equation*}
since $\nu$ is supported on $\RR^m\times I$ and $\mscr{P}$ is a partition of $\ZZ^m$.
%
But $\sup_{\bm{c}\in \RR^m} \sup_{r>0}{\abs{\nu_{\bm{c}/Z}(r)}} \le \mcal{M}_{1, 0}(\nu)$ by \eqref{EQN:define-norm-M_1,k(nu)};
so $\sum_{\bm{c}\in \mcal{S}} \sup_{r>0}{\abs{\nu_{\bm{c}/Z}(r)}}\ll_\theta Z^\theta \mcal{M}_{1,0}(\nu)$ (since $\Supp{\nu}\belongs [-1,1]^m\times I$).
Also, $\sup_{\mcal{R}\in \mscr{P}} \sum_{n\in N\cdot I} \abs{\bar{a}_{F,1}^{\mcal{R}}(n)} \ll_{I,\eps} Q(\mscr{P})^\eps N^{1/3+\eps}$ by \eqref{INEQ:easy-absolute-A_F,1-bound}.
Multiplying gives \eqref{INEQ:RA1'E-main-term-bound-goal}.
\end{proof}

\subsection{Bounding exterior squares}
\label{SUBSEC:bounding-exterior-squares}

For each $\bm{c}\in \mcal{S}_1$, let $\mu_{\bm{c},2}(n)$ denote the $n$th coefficient of the Dirichlet series $\zeta(s) / L(s,V_{\bm{c}},\bigwedge^2)$.
Recall $\mcal{N}^{\bm{c}}$, $\mcal{N}_{\bm{c}}$ from \eqref{EQN:define-moduli-sets-N^c,N_c}.
Assume Conjecture~\ref{CNJ:(HW2)}.

\begin{proposition}
\label{PROP:good-moduli-restricted-wedge2E}
Let $A\ge 2$ be an even integer.
Let $Z, N, \eps \in \RR_{>0}$ with $N\le Z^{3/2}$.
Then 
\begin{equation}
\label{INEQ:good-moduli-restricted-wedge2E-goal}
\sum_{\bm{c}\in \mcal{S}_1\cap [-Z,Z]^m}\,
\Bigl\lvert{
\sum_{N\le n < 2N:\, n\in \mcal{N}^{\bm{c}}} \mu_{\bm{c},2}(n)
}\Bigr\rvert^A
\ll_{A, \eps} Z^m N^{A-1/3 + \eps}.
\end{equation}
\end{proposition}

\begin{proof}

For $N<2$, the left-hand side of \eqref{INEQ:good-moduli-restricted-wedge2E-goal} is $\ll_A Z^m$ by \eqref{geometric-Ramanujan-bound-via-WMC}.
Now assume $N\ge 2$.

Let $\nu_0$ be as in \S\ref{SUBSEC:conventions}.
Let $f(Z,N)$ denote the left-hand side of \eqref{INEQ:good-moduli-restricted-wedge2E-goal}.
Let $g(Z,N)$ denote $f(Z,N)$ with $\sum_{\bm{c}\in \mcal{S}_1} \nu_0(\bm{c}/Z)\cdots$ in place of $\sum_{\bm{c}\in \mcal{S}_1\cap [-Z,Z]^m}\cdots$.
Then by positivity,
\begin{equation}
\label{EQN:smoothing-above-and-below-for-restricted-wedge2E}
    g(Z/2,N) \le f(Z,N) \le g(2Z,N)
\end{equation}
for all $Z>0$.
But if $Z\ge (2N)^{A+1}$, then (if we let $\mu_{\bm{c},2}(\bm{n})\defeq \mu_{\bm{c},2}(n_1)\cdots \mu_{\bm{c},2}(n_A)$ and $P\defeq n_1\cdots n_A$, and note that $\bm{1}_{\gcd(P, \Delta(\bm{c})) = 1} \cdot \mu_{\bm{c},2}(\bm{n})$ is determined by $P$ and $\bm{c}\bmod{P}$)
\begin{equation*}
\begin{split}
    g(Z,N) &= \sum_{N\le n_1,\dots,n_A < 2N} \sum_{\bm{c}\in \ZZ^m} \nu_0(\bm{c}/Z)
    \bm{1}_{\bm{c}\in \mcal{S}_1} \bm{1}_{\gcd(P, \Delta(\bm{c})) = 1} \cdot \mu_{\bm{c},2}(\bm{n}) \\
    &= O_A(Z^m) + \sum_{N\le n_1,\dots,n_A < 2N} (Z/P)^m \sum_{\bm{c}\in (\ZZ/P\ZZ)^m} \bm{1}_{\gcd(P, \Delta(\bm{c})) = 1} \cdot \mu_{\bm{c},2}(\bm{n}),
\end{split}
\end{equation*}
by \eqref{geometric-Ramanujan-bound-via-WMC} (or \eqref{EQN:smooth-case-Deligne-purity}) and Poisson summation in residue classes modulo $P$ (cf.~the proof of Proposition~\ref{PROP:SST-framework});
note that $P\ge N^A > 1$, so the condition $\gcd(P, \Delta(\bm{c})) = 1$ automatically implies $\bm{c}\in \mcal{S}_1$.
If $\tilde{g}(Z,N)\defeq g(Z,N)/Z^m$, then for $Z\ge (2N)^{A+1}$ we conclude that
\begin{equation}
\label{EQN:large-Z-inflation-for-for-restricted-wedge2E}
    \tilde{g}(Z,N) = \tilde{g}((2N)^{A+1},N) + O_A(1).
\end{equation}

We now address $Z\le (2N)^{A+1}$.
Recall $\mcal{S}_2$ from Proposition~\ref{PROP:SST-framework}.
If $\bm{c}\in \mcal{S}_2$, then \eqref{EQN:factor-exterior-square-L-function} and Proposition~\ref{PROP:HW2-consequences}(\ref{ITEM:1/L-bound-from-GRH}) (applied to $L(s,\phi_{\bm{c},2})^{-1}$),
when combined with \eqref{geometric-Ramanujan-bound-via-WMC}
at primes $p\mid \Delta(\bm{c})$,
yield $\sum_{N\le n < 2N:\, n\in \mcal{N}^{\bm{c}}} \mu_{\bm{c},2}(n)\ll_\eps \norm{\bm{c}}^\eps N^{1/2+\eps}$.
If $\bm{c}\in \mcal{S}_1\setminus \mcal{S}_2$, we still have the bound $\sum_{N\le n < 2N:\, n\in \mcal{N}^{\bm{c}}} \mu_{\bm{c},2}(n)\ll_\eps N^{1+\eps}$ due to \eqref{geometric-Ramanujan-bound-via-WMC}.
Therefore, Proposition~\ref{PROP:SST-framework} yields
\begin{equation*}
    g(Z,N) \ll_{A,\eps} Z^{m+\eps} N^{A/2+\eps} + Z^{m-1/2+\eps} N^{A+\eps}
\end{equation*}
for all $Z>0$.
It follows that $\tilde{g}(Z,N) \ll_{A,\eps} N^{A/2+\eps} + N^{A-1/3+\eps}$ for $Z\ge N^{2/3}$ when $Z\le (2N)^{A+1}$, and thus (by \eqref{EQN:large-Z-inflation-for-for-restricted-wedge2E}) for all $Z\ge N^{2/3}$.
Now \eqref{INEQ:good-moduli-restricted-wedge2E-goal} follows from \eqref{EQN:smoothing-above-and-below-for-restricted-wedge2E}.
\end{proof}


To handle primes $p\mid \Delta(\bm{c})$, we prove the following (unconditional) result:

\begin{lemma}
\label{LEM:N_c-small-divisor-moment-bound}
If $Z, A, \eps\in \RR_{>0}$, then
$\sum_{\bm{c}\in \mcal{S}_1\cap [-Z,Z]^m}
\card{\mcal{N}_{\bm{c}} \cap [1,N]}^A
\ll_{A, \eps} Z^m N^\eps$.
\end{lemma}

\begin{proof}
Lemma~\ref{LEM:count-R-N_c-infty-divisors} immediately suffices if $Z<N^A$.
Now suppose $Z\ge N^A$.
By H\"{o}lder's inequality, we may assume $A\in \ZZ_{\ge 1}$.
The sum $\sum_{\bm{c}\in \mcal{S}_1\cap [-Z,Z]^m}
\card{\mcal{N}_{\bm{c}} \cap [1,N]}^A$ then equals
\begin{equation*}
    \Sigma_8 \defeq \sum_{\bm{c}\in \mcal{S}_1\cap [-Z,Z]^m}
    \sum_{u_1,\dots,u_A\le N:\,u_i\mid \Delta(\bm{c})^\infty} 1
    = \sum_{u_1,\dots,u_A\le N} \sum_{\bm{c}\in \mcal{S}_1\cap [-Z,Z]^m}
    \bm{1}_{\rad(u_1\cdots u_A)\mid \Delta(\bm{c})}.
\end{equation*}
In $\Sigma_8$ we have $u_1\cdots u_A\le N^A\le Z$, so by Lang--Weil and the Chinese remainder theorem,
\begin{equation*}
    \Sigma_8 \ll_\eps \sum_{u_1,\dots,u_A\le N} Z^m \rad(u_1\cdots u_A)^{\eps-1}
    \le \sum_{r\le N^A} Z^m r^{\eps-1} \sum_{u_1,\dots,u_A\le N:\, u_i\mid r^\infty} 1.
\end{equation*}
By Lemma~\ref{LEM:count-R-N_c-infty-divisors}, we conclude that $\Sigma_8 \ll_{A,\eps} \sum_{r\le N^A} Z^m r^{\eps-1} (Nr)^\eps \ll_\eps Z^m N^{(2A+1)\eps}$.
\end{proof}

In \S\ref{SEC:adapting-L-statistics-to-the-delta-method},
we need the following technical complement to Proposition~\ref{PROP:R2'-implies-R2'E}.
\begin{proposition}[$\bigwedge$2E]
\label{PROP:(bigwedge2E)}
Assume Conjecture~\ref{CNJ:(HW2)}.
Fix $A\in \RR_{>0}$ and $f(r)\in \set{\bm{1}_{r\le 1}} \cup C^\infty_c(\RR_{>0})$.
Let $Z, N \in \RR_{>0}$ with $N\le Z^{3/2}$.
Let $t_0\in \RR$.
For some $\eta_3=\eta_3(A)>0$, we have
\begin{equation}
\label{INEQ:general-weighted-wedge2E-goal}
\sum_{\bm{c}\in \mcal{S}_1\cap [-Z,Z]^m}\,
\Bigl\lvert{
\sum_{n\ge 1} f(n/N)n^{-it_0} \cdot \mu_{\bm{c},2}(n)
}\Bigr\rvert^A
\ll_{A,f} (1+\abs{t_0})^A Z^m N^{(1-\eta_3) A}.
\end{equation}
\end{proposition}

\begin{proof}
First suppose $(f,t_0) = (\bm{1}_{r\le 1},0)$ and $A\in \ZZ_{\ge 2}$.
Let $g(M) = \sum_{r\le M:\,r\in \mcal{N}^{\bm{c}}} \mu_{\bm{c},2}(r)$.
By multiplicativity, $\sum_{n\le N} \mu_{\bm{c},2}(n)
= \sum_{d\le N:\,d\in \mcal{N}_{\bm{c}}} \mu_{\bm{c},2}(d) g(N/d)$,
which is $\ll_\eps N^\eps \sum_{d\le N:\,d\in \mcal{N}_{\bm{c}}} \abs{g(N/d)}$ by \eqref{geometric-Ramanujan-bound-via-WMC}.
So by H\"{o}lder over $d$, the left-hand side of \eqref{INEQ:general-weighted-wedge2E-goal} is at most $O_{A,\eps}(N^\eps)$ times
\begin{equation*}
\Sigma_9 \defeq \sum_{\bm{c}\in \mcal{S}_1\cap [-Z,Z]^m}
\card{\mcal{N}_{\bm{c}} \cap [1,N]}^{A-1}
\sum_{d\le N:\,d\in \mcal{N}_{\bm{c}}} \abs{g(N/d)}^A.
\end{equation*}
Switching $\bm{c}$, $d$ yields $\Sigma_9 \le \sum_{d\le N}
\sum_{\bm{c}\in \mcal{S}_1\cap [-Z,Z]^m}
\card{\mcal{N}_{\bm{c}} \cap [1,N]}^{A-1} \abs{g(N/d)}^A$ (by positivity).
Then Cauchy--Schwarz over $\bm{c}$, followed by Lemma~\ref{LEM:N_c-small-divisor-moment-bound} and Proposition~\ref{PROP:good-moduli-restricted-wedge2E}, gives
\begin{equation*}
    \Sigma_9
    \ll_{A,\eps} \sum_{d\le N} (Z^m N^\eps)^{1/2} (Z^m (N/d)^{2A-1/3+\eps})^{1/2}
    \ll_{A,\eps} Z^m N^{A-1/6+\eps},
\end{equation*}
where we use $A-1/6+\eps/2\ge 1+\eps/2$ to evaluate the sum over $d$.
By H\"{o}lder, then, \eqref{INEQ:general-weighted-wedge2E-goal} holds with $\eta_3 = 1/7A$ if $A\ge 2$, and with $\eta_3 = 1/14$ if $0<A\le 2$.
The general case follows from partial summation and H\"{o}lder, since $\frac{\partial}{\partial n} D(n/N) \ll_D 1/N$ and $\frac{\partial}{\partial n} n^{it_0} \ll \abs{t_0}/n$.
\end{proof}

(With more work, one could relax the GRH assumption.
One could in fact unconditionally handle the case of very large $Z$,
using \eqref{EQN:E_c(p)-average}--\eqref{EQN:E_c(p)^2-average};
see \cite{wang2021_HLH_vs_RMT}*{\S7.6}.)

\section{Adapting \texpdf{$L$}{L}-function statistics to delta}
\label{SEC:adapting-L-statistics-to-the-delta-method}



\subsection{Factorization}

We need to mold the statistics from \S\ref{SEC:statistics-of-L-function-families} into a form friendlier for the delta method.
We first split the series $\Phi(\bm{c},s)$ (from \eqref{EQN:define-Phi}) into more manageable pieces.
Recall $\mcal{N}^{\bm{c}}$, $\mcal{N}_{\bm{c}}$ from \eqref{EQN:define-moduli-sets-N^c,N_c}.
Given $\bm{c}\in \mcal{S}_1$, consider the factorization $\Phi = \Phi^{\map{G}} \Phi^{\map{B}}$, where
\begin{align}
    \Phi^{\map{G}}(\bm{c},s)
    &\defeq \prod_{p\nmid \Delta(\bm{c})} \Phi_p(\bm{c},s)
    = \sum_{n\in \mcal{N}^{\bm{c}}} S^\natural_{\bm{c}}(n) n^{-s},
    \label{EQN:define-Phi^G} \\
    \Phi^{\map{B}}(\bm{c},s)
    &\defeq \prod_{p\mid \Delta(\bm{c})} \Phi_p(\bm{c},s)
    = \sum_{n\in \mcal{N}_{\bm{c}}} S^\natural_{\bm{c}}(n) n^{-s}.
    \label{EQN:define-Phi^B}
\end{align}
One can approximate $\Phi^{\map{G}}$ using Hasse--Weil $L$-functions.
It would be nice to also relate $\Phi_p$ for $p\mid \Delta(\bm{c})$ to $L$-functions, even in special cases like when $m=4$ and $v_p(\Delta(\bm{c})) = 1$.
For now, we study $\Phi^{\map{B}}$ by completely different means (see \S\ref{SEC:new-bounds-on-bad-sums-S}).
In \S\ref{SEC:adapting-L-statistics-to-the-delta-method}, we thus concentrate on $\Phi^{\map{G}}$.

For the rest of \S\ref{SEC:adapting-L-statistics-to-the-delta-method}, assume $2\mid m$.
We first factor $\Phi^{\map{G}}$ into three pieces: $\Phi^{\bm{c},1}$, $\Phi^{\bm{c},2}$, $\Phi^{\bm{c},3}$.

\begin{definition}
\label{DEFN:factor-Phi^G-into-Phi_1-Phi_2-Phi_3}
Let $\Phi^{\bm{c},1}(s)\defeq L(s,V_{\bm{c}})^{-1} L(1/2+s,V)^{-1} \zeta(2s)^{-1}$ as in \eqref{EQN:define-Phi_1}, and let
\begin{align}
\Phi^{\bm{c},2}(s)
&\defeq \zeta(2s)/L(2s,V_{\bm{c}},\textstyle{\bigwedge^2}),
\label{EQN:define-Phi_2} \\
\Phi^{\bm{c},3}(s)
&\defeq \Phi^{\map{G}}(\bm{c},s)
L(s,V_{\bm{c}})
L(1/2+s,V)
L(2s,V_{\bm{c}},\textstyle{\bigwedge^2}).
\label{EQN:define-Phi_3}
\end{align}
For each $j\in \set{1,2,3}$,
let $a_{\bm{c},j}(n)$ be the $n$th coefficient of the Dirichlet series $\Phi^{\bm{c},j}(s)$.
\end{definition}


The factors $\Phi^{\bm{c},2}$, $\Phi^{\bm{c},3}$ in $\Phi^{\map{G}} = \Phi^{\bm{c},1} \Phi^{\bm{c},2} \Phi^{\bm{c},3}$ hinder any attempt to apply statistics on $\Phi^{\bm{c},1}$ to $\Phi^{\map{G}}$.
Fortunately, $\Phi^{\bm{c},2}$, $\Phi^{\bm{c},3}$ turn out to behave as ``error factors'' on average.
Proposition~\ref{PROP:(bigwedge2E)} lets us handle large moduli in $\Phi^{\bm{c},2}$; note that $a_{\bm{c},2}(n)=0$ if $n$ is not a square, and
\begin{equation}
\label{EQN:relate-a_c,2-to-mu_c,2}
    a_{\bm{c},2}(n) = \mu_{\bm{c},2}(n^{1/2})
\end{equation}
otherwise (where $\mu_{\bm{c},2}$ is as in \S\ref{SUBSEC:bounding-exterior-squares}).
We now prove results to handle large moduli in $\Phi^{\bm{c},3}$.


The factor $\Phi^{\bm{c},3}$ measures the quality of $\Phi^{\bm{c},1} \Phi^{\bm{c},2}$ as an approximation to $\Phi^{\map{G}}$.
Recall the ``first-order approximation''
$\Phi(\bm{c},s) = \Psi^{\bm{c},1}(s) \Psi^{\bm{c},2}(s)$ given by $\Psi^{\bm{c},1}$, $\Psi^{\bm{c},2}$ from \eqref{EQN:define-Psi_1,Psi_2}; here $\Psi^{\bm{c},1}(s) = 1/L(s,V_{\bm{c}})$.
The ``first-order error'' $\Psi^{\bm{c},2}(s)$ is only expected to converge absolutely for $\Re(s)>1/2$.
As suggested in \S\ref{SEC:background-on-discriminants-and-the-delta-method}, this is a ``source of $\eps$'' in \eqref{INEQ:near-optimal-diagonal-GRH-bound-over-smooth-locus-S_1}.
On the other hand, the following result establishes absolute convergence for $\Phi^{\bm{c},3}(s)$ past the critical line $\Re(s)=1/2$.

\begin{proposition}
\label{PROP:approx-Phi-past-1/2}
Uniformly over $\bm{c}\in \mcal{S}_1$, primes $p$, and integers $l\ge 1$, we have
\begin{equation}
\label{INEQ:a_c,3(p^l)-estimate-goals}
a_{\bm{c},3}(p)\cdot\bm{1}_{p\nmid \Delta(\bm{c})}=0,
\quad
a_{\bm{c},3}(p^2)
\cdot\bm{1}_{p\nmid \Delta(\bm{c})}\ll p^{-1/2},
\quad
a_{\bm{c},3}(p^l)\ll_\eps p^{l\eps}.
\end{equation}
In particular, if $\bm{c}\in \mcal{S}_1$, then $\Phi^{\bm{c},3}(s)$ converges absolutely over $\Re(s)>1/3$.
\end{proposition}

\begin{proof}
Let $\bm{c}\in \mcal{S}_1$.
Let $\lambda^\natural_{\bm{c}}\defeq \lambda^\natural_{V_{\bm{c}}}$, as in \S\ref{SUBSEC:computing-local-averages}.
Suppose first that $p\mid \Delta(\bm{c})$.
Then $\Phi^{\map{G}}_p(\bm{c},s) = 1$ by \eqref{EQN:define-Phi^G}.
So by \eqref{EQN:define-Phi_3} and \eqref{geometric-Ramanujan-bound-via-WMC}, we have $a_{\bm{c},3}(p^l)\ll_\eps p^{l\eps}$.
So \eqref{INEQ:a_c,3(p^l)-estimate-goals} holds.

Now suppose $p\nmid \Delta(\bm{c})$.
Then $\Phi_p(\bm{c},s) = 1 + S^\natural_{\bm{c}}(p)p^{-s}$ by \eqref{EQN:S_c(p^l)-vanishing}, and
\begin{equation*}
S^\natural_{\bm{c}}(p)
= E^\natural_{\bm{c}}(p)-p^{-1/2}E^\natural_F(p)
= -\lambda^\natural_{\bm{c}}(p)-p^{-1/2}\lambda^\natural_V(p)
\end{equation*}
by \eqref{EQN:rewrite-S_c(p)-via-E_c}, \eqref{EQN:compute-lambda_V_c-p}, and \eqref{EQN:compute-lambda_V-p}.
In particular, $a_{\bm{c},3}(p^l)\ll_\eps p^{l\eps}$ by \eqref{EQN:define-Phi_3} and \eqref{geometric-Ramanujan-bound-via-WMC}.
Furthermore,
\begin{equation*}\begin{split}
\Phi_p(\bm{c},s) L_p(s,V_{\bm{c}})
&= (1 - \lambda^\natural_{\bm{c}}(p)p^{-s} - \lambda^\natural_V(p)p^{-1/2-s})
(1 + \lambda^\natural_{\bm{c}}(p)p^{-s} + \lambda^\natural_{\bm{c}}(p^2)p^{-2s} + O(p^{-3s})) \\
&= 1 - \lambda^\natural_V(p)p^{-1/2-s} + [\lambda^\natural_{\bm{c}}(p^2)-\lambda^\natural_{\bm{c}}(p)^2]p^{-2s} + O(p^{-1/2-2s}) + O(p^{-3s}).
\end{split}\end{equation*}
To get further cancellation,
we multiply $\Phi_p(\bm{c},s) L_p(s,V_{\bm{c}})$ by
\begin{equation*}
\begin{split}
L_p(1/2+s,V)
&= 1+\lambda^\natural_V(p)p^{-1/2-s}+O(p^{-1-2s}), \\
L_p(2s,V_{\bm{c}},\textstyle{\bigwedge^2})
&= 1+\lambda^\natural_{V_{\bm{c}},\bigwedge^2}(p)p^{-2s}+O(p^{-4s}),
\end{split}
\end{equation*}
to get (in view of $\lambda^\natural_{V_{\bm{c}},\bigwedge^2}(p) = \lambda^\natural_{\bm{c}}(p)^2-\lambda^\natural_{\bm{c}}(p^2)$ from \eqref{EQN:compute-exterior-square-lambda_V_c-p})
\begin{equation}
\label{EQN:final-second-order-cancellation-display}
\Phi_p(\bm{c},s)L_p(s,V_{\bm{c}})
L_p(1/2+s,V)
L_p(2s,V_{\bm{c}},\textstyle{\bigwedge^2})
= 1
+ O(p^{-1/2-2s})
+ O(p^{-3s}).
\end{equation}
By \eqref{EQN:define-Phi_3}, the left-hand side of \eqref{EQN:final-second-order-cancellation-display} is precisely the local factor $\Phi^{\bm{c},3}_p(s)$ of $\Phi^{\bm{c},3}(s)$.
Thus \eqref{EQN:final-second-order-cancellation-display} completes the proof of \eqref{INEQ:a_c,3(p^l)-estimate-goals}.
The convergence statement on $\Phi^{\bm{c},3}(s)$ follows from \eqref{INEQ:a_c,3(p^l)-estimate-goals}.
\end{proof}

\begin{corollary}
[$\Phi$3E]
\label{COR:average-bound-on-Phi_3}
Fix $A\in \RR_{>0}$.
Then uniformly over $Z, N\in \RR_{>0}$,
we have
\begin{equation}
\label{INEQ:Phi3E-goal}
\sum_{\bm{c}\in \mcal{S}_1\cap [-Z,Z]^m}
\biggl(\,\sum_{n\le N} \abs{a_{\bm{c},3}(n)}\biggr)^{\!A}
\ll_{A,\eps} Z^m N^{(1/3+\eps)A}.
\end{equation}
\end{corollary}

\begin{proof}
Recall $\mcal{N}^{\bm{c}}$, $\mcal{N}_{\bm{c}}$ from \eqref{EQN:define-moduli-sets-N^c,N_c}.
By multiplicativity (of $\abs{a_{\bm{c},3}}$) and positivity, we have
\begin{equation}
\label{INEQ:positive-upper-bound-factorization-for-|a_c,3|-sum}
    \sum_{n\le N} \abs{a_{\bm{c},3}(n)}
    \le \sum_{d\le N:\,d\in \mcal{N}_{\bm{c}}} \abs{a_{\bm{c},3}(d)} \sum_{r\le N:\,r\in \mcal{N}^{\bm{c}}} \abs{a_{\bm{c},3}(r)}.
\end{equation}
But for any $\eps>0$, the bounds in \eqref{INEQ:a_c,3(p^l)-estimate-goals} imply $a_{\bm{c},3}(d)\ll_\eps d^\eps$,
and that $\sum_{r\le N:\,r\in \mcal{N}^{\bm{c}}} \abs{a_{\bm{c},3}(r)}$ is 
$\le N^{1/3+\eps} \sum_{r\in \mcal{N}^{\bm{c}}} r^{-1/3-\eps} \abs{a_{\bm{c},3}(r)}
\ll_\eps N^{1/3+\eps}$
(cf.~\eqref{INEQ:easy-absolute-A_F,1-bound}).
So by \eqref{INEQ:positive-upper-bound-factorization-for-|a_c,3|-sum}, the left-hand side of \eqref{INEQ:Phi3E-goal} is
$\ll_{A,\eps} N^{(1/3+2\eps) A}
\sum_{\bm{c}\in \mcal{S}_1\cap [-Z,Z]^m}
\card{\mcal{N}_{\bm{c}} \cap [1,N]}^{A-1}$.
Now \eqref{INEQ:Phi3E-goal} follows from Lemma~\ref{LEM:N_c-small-divisor-moment-bound}.
\end{proof}

\subsection{From (R2'E) to (R2'E')}
\label{SUBSEC:proving-(R2'E')}

We now build on Conjecture~\ref{CNJ:(R2'E)}~(R2'E).

\begin{conjecture}
[R2'E']
\label{CNJ:(R2'E')}
Fix a function $D\in C^\infty_c(\RR_{>0})$.
Let $t\in \RR$.
Let $Z, N, \eps \in \RR_{>0}$ with $N\le Z^3$ and $\eps\le 1$.
Then for some real $A_6=A_6(F,\eps)>0$, we have
\begin{equation}
\label{INEQ:R2'E'-goal}
\sum_{\bm{c}\in \mcal{S}_1\cap [-Z,Z]^m}\,
\Bigl\lvert{
\sum_{n\in \mcal{N}^{\bm{c}}} D(n/N) \cdot n^{-it}S^\natural_{\bm{c}}(n)
}\Bigr\rvert^{2-\eps}
\ll_{F, D, \eps} (1+\abs{t})^{A_6} Z^m N^{(2-\eps)/2}.
\end{equation}
\end{conjecture}




\begin{proposition}
\label{PROP:R2'E-implies-R2'E'}
Assume Conjectures~\ref{CNJ:(HW2)} and~\ref{CNJ:(R2'E)}.
Then Conjecture~\ref{CNJ:(R2'E')} holds.
\end{proposition}

\begin{proof}

Let $\nu_2\in C^\infty_c(\RR_{>0})$ be as in \S\ref{SEC:separation-lemmas}, so that $\Supp{\nu_2}\belongs [1,2]$ and we have \eqref{EQN:nu_2-squared-integrates-to-1}.
In view of \eqref{EQN:define-Phi^G} and the factorization $\Phi^{\map{G}} = \Phi^{\bm{c},1} \Phi^{\bm{c},2} \Phi^{\bm{c},3}$, we may use Lemma~\ref{LEM:dyadic-partial-Mellin-summation} (with $k=3$ and $a(\bm{n}) = \prod_{1\le j\le 3} (n_j^{-it} a_{\bm{c},j}(n_j))$, and $f(\bm{r}) = D(r_1r_2r_3/N)$) to write
\begin{equation}
\label{EQN:smooth-dyadic-restriction-separation-for-R2'E'}
    \sum_{n\in \mcal{N}^{\bm{c}}} D(n/N) \cdot n^{-it}S^\natural_{\bm{c}}(n)
    = (2\pi)^{-3} \int_{\bm{N}\ge 1/2} d^\times\bm{N} \int_{\bm{t}\in \RR^3} d\bm{t}\,
    g_{\bm{N}}^\vee(i\bm{t})
    \prod_{1\le j\le 3} \Sigma_{10,\bm{N}}^{\bm{c},j}(\bm{t}),
\end{equation}
where $\bm{N}=(N_1,N_2,N_3)$,
where $g_{\bm{N}}(\bm{r})\defeq D(r_1r_2r_3/N) \prod_{1\le j\le 3} \nu_2(r_j/N_j)$,
and where
\begin{equation}
\label{EQN:define-smooth-dyadic-pieces-Sigma_10-for-R2'E'}
\Sigma_{10,\bm{N}}^{\bm{c},j}(\bm{t})=\Sigma_{10,\bm{N}}^{\bm{c},j}(t_1,t_2,t_3)
\defeq \sum_{n_j\ge 1} \nu_2(n_j/N_j) n_j^{-i(t_j+t)} a_{\bm{c},j}(n_j).
\end{equation}
(Note that $N$, $t$ are independent of $\bm{N}$, $\bm{t}$.)
For all $\bm{N}\ge 1/2$ and $b\ge 0$, Proposition~\ref{PROP:standard-general-Mellin-bound} gives
\begin{equation}
    \label{INEQ:Mellin-decay-for-g_N-for-R2'E'}
    g_{\bm{N}}^\vee(i\bm{t}) \ll_{b,D,\nu_2} (1+\norm{\bm{t}})^{-b}.
\end{equation}

Fix an integer $A\ge 1$ for which $\Supp{D}\belongs [A^{-1}, A]$.
If $\bm{r}\in \Supp{g_{\bm{N}}}$, then $N/A \le r_1r_2r_3 \le AN$ and $N_j\le r_j\le 2N_j$ for all $j$,
so $\bm{N}$ lies in the set
\begin{equation}
\label{EQN:dyadic-range-restrictions-for-R_10}
    \mscr{R}_{10}\defeq \set{\bm{N}\ge 1/2: N_1N_2N_3\in [N/8A, AN]}.
\end{equation}
Thus the equality \eqref{EQN:smooth-dyadic-restriction-separation-for-R2'E'} remains true if we restrict the integral over $\bm{N}\ge 1/2$ to the region $\mscr{R}_{10}$.
Now set $W_1(\bm{N}, \bm{t})\defeq N_1^\eta \cdot \abs{g_{\bm{N}}^\vee(i\bm{t})}$ and $W_2(\bm{N}, \bm{t})\defeq N_1^{-(\beta-1)\eta} \cdot \abs{g_{\bm{N}}^\vee(i\bm{t})}$, for a small constant $\eta>0$ to be chosen later.
Let $\mcal{I}_1\defeq \int_{\mscr{R}_{10} \times \RR^3} d^\times\bm{N}\, d\bm{t}\, W_1(\bm{N}, \bm{t})$.
Letting $\beta\defeq 2-\eps\in [1,2)$, and using H\"{o}lder over $(\log\bm{N}, \bm{t}) \belongs \RR^3 \times \RR^3$ (restricted to $\bm{N}\in \mscr{R}_{10}$), we obtain
\begin{equation}
\label{INEQ:first-Holder-for-R2'E'}
\Bigl\lvert{
\sum_{n\in \mcal{N}^{\bm{c}}} D(n/N) n^{-it}S^\natural_{\bm{c}}(n)
}\Bigr\rvert^\beta
\le \mcal{I}_1^{\beta-1} \cdot
\int_{\mscr{R}_{10} \times \RR^3} d^\times\bm{N}\, d\bm{t}\,
W_2(\bm{N}, \bm{t})
\prod_{1\le j\le 3} \abs{\Sigma_{10,\bm{N}}^{\bm{c},j}(\bm{t})}^\beta,
\end{equation}
since $\beta\ge 1$ and $W_1(\bm{N}, \bm{t})^{\beta-1} \cdot W_2(\bm{N}, \bm{t}) = \abs{g_{\bm{N}}^\vee(i\bm{t})}^\beta$.

By \eqref{INEQ:Mellin-decay-for-g_N-for-R2'E'} and \eqref{EQN:dyadic-range-restrictions-for-R_10}, we have
$\mcal{I}_1
\ll_{D,\nu_2} \int_{\mscr{R}_{10}} d^\times\bm{N}\, N_1^\eta
\ll_{A,\eta} N^\eta
\int_{N_2,N_3\ge 1/2} \frac{d^\times{N_2}\,d^\times{N_3}}{(N_2N_3)^\eta}
\ll_\eta N^\eta$.
Upon summing \eqref{INEQ:first-Holder-for-R2'E'} over $\bm{c}\in \mcal{S}_1\cap [-Z,Z]^m$, we thus find that the left-hand side of \eqref{INEQ:R2'E'-goal} is
\begin{equation}
\label{INEQ:simplified-Holder-and-Fubini-bound-for-R2'E'}
\ll_{D,\nu_2,\eta}
\int_{\mscr{R}_{10} \times \RR^3} d^\times\bm{N}\, d\bm{t}\,
(N/N_1)^{(\beta-1)\eta} \cdot \abs{g_{\bm{N}}^\vee(i\bm{t})}
\sum_{\bm{c}\in \mcal{S}_1\cap [-Z,Z]^m}
\prod_{1\le j\le 3} \abs{\Sigma_{10,\bm{N}}^{\bm{c},j}(\bm{t})}^\beta.
\end{equation}

Now let $(\gamma_1, \gamma_2, \gamma_3)\defeq (2, 4\beta/\eps, 4\beta/\eps)$.
Then $\sum_{1\le j\le 3} \beta/\gamma_j = 1$ (since $\beta=2-\eps$), so by H\"{o}lder over $\bm{c}$
(writing $\heartsuit_{10} = \prod_{1\le j\le 3} \Sigma_{10,\bm{N}}^{\bm{c},j}(\bm{t})$ and $\mscr{M}_j = \norm{\Sigma_{10,\bm{N}}^{\bm{c},j}(\bm{t})}_{\ell^{\gamma_j}_{\bm{c}}(\mcal{S}_1\cap [-Z,Z]^m)}^{\gamma_j}$ for brevity),
\begin{equation}
\label{INEQ:second-Holder-for-R2'E'}
    \sum_{\bm{c}\in \mcal{S}_1\cap [-Z,Z]^m} \abs{\heartsuit_{10}}^\beta
    = \norm{\heartsuit_{10}}_{\ell^\beta_{\bm{c}}(\mcal{S}_1\cap [-Z,Z]^m)}^\beta
    \le \prod_{1\le j\le 3} \mscr{M}_j^{\beta/\gamma_j}.
\end{equation}
We now bound the necessary $\ell^{\gamma_j}$-norms.
First,
$\mscr{M}_1
\ll_{\nu_2} (1+\abs{t_1+t})^{A_3}
(Z+N_1^{1/3})^m N_1^{\gamma_1/2}$,
by \eqref{EQN:define-smooth-dyadic-pieces-Sigma_10-for-R2'E'} and Conjecture~\ref{CNJ:(R2'E)} (with $Z+N_1^{1/3}$, $N_1$ in place of $Z$, $N$).
Second, by \eqref{EQN:define-smooth-dyadic-pieces-Sigma_10-for-R2'E'}, \eqref{EQN:relate-a_c,2-to-mu_c,2}, and \eqref{INEQ:general-weighted-wedge2E-goal} (with $A=\gamma_2$, with $f=\nu_2$, and with $Z+N_2^{1/3}$, $N_2^{1/2}$ in place of $Z$, $N$), we have
$\mscr{M}_2
\ll_{\gamma_2,\nu_2} (1+\abs{t_2+t})^{\gamma_2}
(Z+N_2^{1/3})^m(N_2^{1/2})^{(1-\eta_3(\gamma_2)) \gamma_2}$.
Third, by \eqref{EQN:define-smooth-dyadic-pieces-Sigma_10-for-R2'E'}, the bound $\nu_2(n_3/N_3) n_3^{-i(t_3+t)} \ll_{\nu_2} 1$, and Corollary~\ref{COR:average-bound-on-Phi_3}, we have
$\mscr{M}_3
\ll_{\gamma_3,\nu_2} Z^mN_3^{11\gamma_3/30}$.

Plugging into \eqref{INEQ:second-Holder-for-R2'E'}, and writing $1+\abs{t_j+t} \le (1+\norm{\bm{t}})(1+\abs{t})$,
we get that if $\bm{N}\in \mscr{R}_{10}$,
then the left-hand side of \eqref{INEQ:second-Holder-for-R2'E'} is
$\ll_{D,\nu_2,\beta} (1+\norm{\bm{t}})^{A_6} (1+\abs{t})^{A_6}$ times
\begin{equation*}
Z^{\sum_j m \beta/\gamma_j} N_1^{\beta/2} N_2^{(1-\eta_3(\gamma_2))\beta/2} N_3^{11\beta/30}
\asymp_{D,\beta} Z^m (N^{1/2} N_2^{-\eta_3(\gamma_2)/2} N_3^{-4/30})^\beta,
\end{equation*}
where $A_6 = (A_3 \beta/\gamma_1 + \gamma_2\beta/\gamma_2) = (A_3/2+1) \beta$.
Upon plugging this result into \eqref{INEQ:simplified-Holder-and-Fubini-bound-for-R2'E'}, we get that the left-hand side of \eqref{INEQ:R2'E'-goal} is at most $O_{D,\nu_2,\beta}(1)$ times
\begin{equation}
\label{EXPR:final-upper-bound-quantity-for-R2'E'}
    \int_{\mscr{R}_{10} \times \RR^3} d^\times\bm{N}\, d\bm{t}\,
    (N/N_1)^{(\beta-1)\eta} \cdot \abs{g_{\bm{N}}^\vee(i\bm{t})}
    \cdot (1+\norm{\bm{t}})^{A_6} \frac{(1+\abs{t})^{A_6} Z^m N^{\beta/2}}{(N_2^{\eta_3(\gamma_2)/2} N_3^{4/30})^\beta}.
\end{equation}

Finally, let $\eta = \min(\eta_3(\gamma_2)/2, 4/30)$.
Then for each $\bm{N}\in \mscr{R}_{10}$, we have $(N/N_1)^\eta \ll_{D,\beta} N_2^{\eta_3(\gamma_2)/2} N_3^{4/30}$.
Since $\beta-1\ge 0$, it follows that the quantity \eqref{EXPR:final-upper-bound-quantity-for-R2'E'} is
\begin{equation*}
\begin{split}
\ll_{D,\beta} \int_{\mscr{R}_{10} \times \RR^3} d^\times\bm{N}\, d\bm{t}\,
\abs{g_{\bm{N}}^\vee(i\bm{t})}
\cdot (1+\norm{\bm{t}})^{A_6} \frac{(1+\abs{t})^{A_6} Z^m N^{\beta/2}}{N_2^{\eta_3(\gamma_2)/2} N_3^{4/30}}
\ll_{D,\beta} (1+\abs{t})^{A_6} Z^m N^{\beta/2},
\end{split}
\end{equation*}
where we have used \eqref{INEQ:Mellin-decay-for-g_N-for-R2'E'} to integrate over $\bm{t}$, and \eqref{EQN:dyadic-range-restrictions-for-R_10} to integrate first over $N_1$ (given $N_2$, $N_3$) and then over $N_2,N_3\ge 1/2$.
Therefore, \eqref{INEQ:R2'E'-goal} holds (with $A_6 = (A_3/2+1) (2-\eps)$).
\end{proof}

\subsection{Handling variation of ``error factors'' for small fixed ``error moduli''}
\label{SUBSEC:handling-variation-of-error-factors}

We would like to build on Propositions~\ref{PROP:RA1o-implies-RA1o'E} and~\ref{PROP:RA1delta-implies-RA1delta'E}, but we must first improve our understanding of certain local factors.
For each $\bm{c}\in\mcal{S}_1$, recall $\Phi^{\bm{c},j}(s)$, $a_{\bm{c},j}(n)$ from Definition~\ref{DEFN:factor-Phi^G-into-Phi_1-Phi_2-Phi_3}, and let $a'_{\bm{c}}(n)$ denote the $n$th coefficient of the Dirichlet series
\begin{equation}
\label{EQN:equivalent-ways-to-factor-a'_c-Dirichlet-series}
\Phi^{\bm{c},2}\Phi^{\bm{c},3}
= \Phi^{\map{G}}/\Phi^{\bm{c},1}
= \Phi^{\map{G}}(\bm{c},s)
L(s,V_{\bm{c}})
L(1/2+s,V)\zeta(2s),
\end{equation}
so that for all $n\ge 1$, we have
\begin{equation}
\label{EQN:decompose-good-restricted-S_c(n)-via-a_c,1,a'_c}
    S^\natural_{\bm{c}}(n) \bm{1}_{n\in \mcal{N}^{\bm{c}}}
    = \sum_{n_1d=n} a_{\bm{c},1}(n_1) a'_{\bm{c}}(d).
\end{equation}
By Proposition~\ref{PROP:approx-Phi-past-1/2}, $a_{\bm{c},3}(n)\ll_\eps n^\eps$, and by \eqref{geometric-Ramanujan-bound-via-WMC} and Definition~\ref{DEFN:factor-Phi^G-into-Phi_1-Phi_2-Phi_3}, $a_{\bm{c},2}(n)\ll_\eps n^\eps$; so certainly
\begin{equation}
\label{INEQ:GRC-on-a'_c}
    a'_{\bm{c}}(d)\ll_\eps d^\eps.
\end{equation}
Moreover, if $d$ is small (or fixed), we would like $a'_{\bm{c}}(d)$ to not vary too wildly with $\bm{c}$.

Given $n\in \ZZ_{\ge 1}$, what data does $a'_{\bm{c}}(n)$ depend on?
Note that $L(1/2+s,V)\zeta(2s)$ is fixed (in terms of $F$).
So by \eqref{EQN:equivalent-ways-to-factor-a'_c-Dirichlet-series} and \eqref{EQN:define-Phi^G}, the coefficient $a'_{\bm{c}}(n)$ is determined by the residue class $\bm{c}+n\ZZ^m$ and the local factors $L_p(s, V_{\bm{c}})$ for $p\mid n$.
Therefore, if we define $l(p,\bm{c})$ as in Lemma~\ref{LEM:density-form-of-ineffective-Krasner-lemma}, then $a'_{\bm{c}}(n)$ is determined by the residue class $\bm{c}+r(n,\bm{c})\ZZ^m$, where
\begin{equation}
\label{EQN:define-r(n,c)-via-l(p,c)}
    r(n,\bm{c}) \defeq \prod_{p\mid n} \lcm(p^{v_p(n)}, p^{l(p,\bm{c})+1}).
\end{equation}

\begin{proposition}
\label{PROP:local-constancy-of-r(n,a)-exceptional-or-not-at-p}
Let $n, q\ge 1$.
Let $\bm{a}, \bm{b}\in \mcal{S}_1$ and suppose $\bm{a}\equiv \bm{b}\bmod{q}$.
Let $p\mid n$ be a prime.
Then $v_p(r(n,\bm{a})) \le v_p(q)$ if and only if $v_p(r(n,\bm{b})) \le v_p(q)$.
\end{proposition}

\begin{proof}
By symmetry, it suffices to prove the ``if'' direction.
So, say $v_p(q)\ge v_p(r(n,\bm{b}))$.
Then $v_p(q)\ge l(p,\bm{b})+1$, so $l(p,\bm{a})=l(p,\bm{b})$ by \eqref{EQN:local-constancy-of-l(p,b)}.
Thus $v_p(r(n,\bm{a})) = v_p(r(n,\bm{b})) \le v_p(q)$.
\end{proof}

In what follows, recall the notation $q_{\mcal{R}}$ from Definition~\ref{DEFN:residue-class-modulus-and-average-coefficient-terminology}.

\begin{definition}
\label{DEFN:define-set-E(n)-of-exceptional-residue-classes-R-to-modulus-n}
For every $n\ge 1$, let $\mscr{E}(n)$ be the set of residue classes $\mcal{R}\belongs \ZZ^m$ for which there exists a tuple $\bm{c}\in \mcal{S}_1\cap \mcal{R}$ with $r(n,\bm{c})\nmid q_{\mcal{R}}$.
\end{definition}

We now construct a partition $\mscr{P}(n,k) = \set{\mcal{R}}$ of $\ZZ^m$ into residue classes $\mcal{R}\belongs \ZZ^m$.

\begin{definition}
\label{DEFN:algorithmic-construction-of-partition-P(n,k)-exceptions-E(n,k)}
Fix integers $n, k \ge 1$.
We define $\mscr{P}(n,k)$ through a recursive decomposition process.
Let $\mscr{S}_0\defeq \set{\bm{a}+n\ZZ^m: 1\le \bm{a}\le n}$ denote the partition of $\ZZ^m$ into the $n^m$ residue classes modulo $n$.
For each $j\ge 0$, define $\mscr{S}_{j+1}$ in terms of $\mscr{S}_j$ as follows:
\begin{enumerate}
    \item If possible, choose a residue class $\mcal{R}\in \mscr{S}_j\cap \mscr{E}(n)$ with $q_{\mcal{R}}\leq n^{k-1}$.
    Otherwise, let $\mscr{S}_{j+1}\defeq \mscr{S}_j$, and skip step~(2).
    
    \item Write $\mcal{R} = \bm{a}+q\ZZ^m$, with $q\ge 1$.
    Choose $\bm{c}\in \mcal{S}_1\cap \mcal{R}$ with $r(n,\bm{c})\nmid q$.
    Choose a prime $p\mid n$ with $v_p(r(n,\bm{c})) > v_p(q)$.
    Create $\mscr{S}_{j+1}$ by replacing the element $\mcal{R}\in \mscr{S}_j$ with the $p^m$ lifted residue classes $(\bm{a}+q\bm{i})+pq\ZZ^m$ (with $1\le \bm{i}\le p$, say).
    Formally,
    \begin{equation*}
        \mscr{S}_{j+1}\defeq (\mscr{S}_j \setminus \set{\mcal{R}})
        \cup \set{(\bm{a}+q\bm{i})+pq\ZZ^m: 1\le \bm{i}\le p}.
    \end{equation*}
\end{enumerate}
Step~(2) can only occur finitely many times (because we require $q_{\mcal{R}}\le n^{k-1}$ in step~(1)).
Let $j_0\defeq \min{\set{j\ge 0: \mscr{S}_{j+1}=\mscr{S}_j}}$.
Let $\mscr{P}(n,k)\defeq \mscr{S}_{j_0}$.\footnote{The result may depend on the choices we make at each step, but all of our estimates based on $\mscr{P}(n,k)$ will apply uniformly over all possible outcomes.}
Let $\mscr{E}(n,k)\defeq \mscr{P}(n,k)\cap \mscr{E}(n)$.

\end{definition}

In Definition~\ref{DEFN:algorithmic-construction-of-partition-P(n,k)-exceptions-E(n,k)}, we allow the initial $\mscr{S}_0$ to branch into many different moduli.
If we did not do this, then to control $a'_{\bm{c}}(n)$ for $n\le M$ might require us to work with moduli exponentially large in $M$ (for some values of $n$), which would be fatal to our approach to Theorem~\ref{THM:level-3-power-saving} (though perhaps OK for Theorem~\ref{THM:level-2-almost-all-integers}).
We will eventually apply Propositions~\ref{PROP:RA1o-implies-RA1o'E} and~\ref{PROP:RA1delta-implies-RA1delta'E} in residue classes $\mcal{R}\in \mscr{P}(n,k)\setminus \mscr{E}(n)$, for some values of $n,k\ge 1$.
We first unravel the structure of $\mscr{P}(n,k)$, and provide some control on the exceptional set $\mscr{E}(n)$.


\begin{lemma}
\label{LEM:key-recursion-lemma}
Let $n,k\ge 1$.
Suppose $\mcal{R}\in \mscr{S}_j$ for some $j\ge 0$.
Suppose $q_{\mcal{R}} > n$.
Then $j\ge 1$.
Furthermore, there exist an index $i\in [0,j-1]$, a prime $p\mid n$, and a residue class $\mcal{R}'\in \mscr{S}_i\cap \mscr{E}(n)$ of modulus $q_{\mcal{R}'}\le n^{k-1}$ with $\mcal{R}\belongs \mcal{R}'$ and $q_{\mcal{R}}/q_{\mcal{R}'} = p$, such that
$v_p(q_{\mcal{R}})\le v_p(r(n,\bm{c}))$
holds for all $\bm{c}\in \mcal{S}_1\cap \mcal{R}'$.
\end{lemma}

\begin{proof}
Each element of $\mscr{S}_0$ has modulus $n$, so $\mcal{R}\notin \mscr{S}_0$; in particular, $j\ne 0$, so $j\ge 1$.
Choose $i\in [0,j-1]$ maximal with $\mcal{R}\notin \mscr{S}_i$.
Then by Definition~\ref{DEFN:algorithmic-construction-of-partition-P(n,k)-exceptions-E(n,k)}, there exists a residue class $\mcal{R}'\in \mscr{S}_i\cap \mscr{E}(n)$, a tuple $\bm{c}'\in \mcal{S}_1\cap \mcal{R}'$, and a prime $p\mid n$, with $q_{\mcal{R}'}\le n^{k-1}$ and $v_p(r(n,\bm{c}')) > v_p(q_{\mcal{R}'})$, such that $\mcal{R}\belongs \mcal{R}'$ and $q_{\mcal{R}}/q_{\mcal{R}'} = p$.
Now let $\bm{c}\in \mcal{S}_1\cap \mcal{R}'$.
By Proposition~\ref{PROP:local-constancy-of-r(n,a)-exceptional-or-not-at-p},
$v_p(r(n,\bm{c})) > v_p(q_{\mcal{R}'})$ (since $v_p(r(n,\bm{c}')) > v_p(q_{\mcal{R}'})$).
So $v_p(r(n,\bm{c})) \ge 1+v_p(q_{\mcal{R}'}) = v_p(q_{\mcal{R}})$.
\end{proof}

\begin{proposition}
\label{PROP:basic-properties-of-residue-classes-R-in-partition-P(n,k)}
Let $n,k\ge 1$.
Let $\mcal{R}\in \mscr{P}(n,k)$.
Then $q_{\mcal{R}}\le n^k$.
Furthermore, if $n\ge 2$ and $\mcal{R}\in \mscr{P}(n,k)\setminus \mscr{E}(n)$, then $\mcal{R}\belongs \mcal{S}_1$.
\end{proposition}

\begin{proof}
By Lemma~\ref{LEM:key-recursion-lemma}, we either have $q_{\mcal{R}}\le n$, or $q_{\mcal{R}} = pq_{\mcal{R}'}$ with $p\mid n$ and $q_{\mcal{R}'}\le n^{k-1}$.
So $q_{\mcal{R}}\le n^k$.
Now suppose $n\ge 2$ and $\mcal{R}\in \mscr{P}(n,k)\setminus \mscr{E}(n)$.
Choose $\bm{c}\in \mcal{S}_1\cap \mcal{R}$
(possible by \eqref{INEQ:dimension-growth-bound-on-S_0}).
Then $r(n,\bm{c})\mid q_{\mcal{R}}$, by the definition of $\mscr{E}(n)$.
Now choose $p\mid n$ (possible since $n\ge 2$).
Then $p^{l(p,\bm{c})+1} \mid q_{\mcal{R}}$, by \eqref{EQN:define-r(n,c)-via-l(p,c)}.
But by the definition of $l(p,\bm{c})$, we have $\Delta(\bm{c}')\ne 0$ for all $\bm{c}'\equiv \bm{c}\bmod{p^{l(p,\bm{c})+1}}$.
Thus $\mcal{R}\belongs \mcal{S}_1$ (since $\bm{c}\in \mcal{R}$).
\end{proof}



\begin{proposition}
\label{PROP:r(n,c)-large-if-R-exceptional}
Let $n,k\ge 2$.
If $\mcal{R}\in \mscr{E}(n,k)$, then $r(n,\bm{c})\ge n^{k-2}$ for all $\bm{c}\in \mcal{S}_1\cap \mcal{R}$.
\end{proposition}

\begin{proof}
Suppose $\mcal{R}\in \mscr{E}(n,k)$.
By Definition~\ref{DEFN:algorithmic-construction-of-partition-P(n,k)-exceptions-E(n,k)}, $\mscr{E}(n,k) = \mscr{S}_{j_0}\cap \mscr{E}(n)$.
So $\mcal{R}\in \mscr{S}_{j_0}\cap \mscr{E}(n)$, whence $q_{\mcal{R}}>n^{k-1}$ (or else we would have $\mscr{S}_{j_0+1}\ne \mscr{S}_{j_0}$ by the algorithm in Definition~\ref{DEFN:algorithmic-construction-of-partition-P(n,k)-exceptions-E(n,k)}).
By Lemma~\ref{LEM:key-recursion-lemma} (applied repeatedly), there exists a sequence of primes $p_1,\dots,p_s\mid n$, with $p_1\cdots p_s\mid q_{\mcal{R}}$ and $q_{\mcal{R}}/(p_1\cdots p_{s-1}) > n \ge q_{\mcal{R}}/(p_1\cdots p_s)$, such that
\begin{equation}
\label{INEQ:key-r(n,c)-divisibility-from-exceptional-recursion}
    v_{p_i}(q_{\mcal{R}}/(p_1\cdots p_{i-1}))\le v_{p_i}(r(n,\bm{c}))
\end{equation}
holds for all $i\in \set{1,2,\dots,s}$ and $\bm{c}\in \mcal{S}_1\cap \mcal{R}$.
If for each $p\mid p_1\cdots p_s$, we apply \eqref{INEQ:key-r(n,c)-divisibility-from-exceptional-recursion} with $i = \min \set{1\le u\le s: p_u = p}$, then we get $v_p(q_{\mcal{R}})\le v_p(r(n,\bm{c}))$.
Since $p_1\cdots p_s\mid q_{\mcal{R}}$, it follows that $p_1\cdots p_s\mid r(n,\bm{c})$, and thus $r(n,\bm{c})\ge p_1\cdots p_s\ge q_{\mcal{R}}/n > n^{k-2}$.
\end{proof}

Let $\mu(\mcal{R})\defeq q_{\mcal{R}}^{-m}$ be the \emph{density} of a residue class $\mcal{R}$ in $\ZZ^m$.

\begin{lemma}
[KL']
\label{LEM:(KL')}
Let $n\ge 2$.
Then $\lim_{k\to\infty} \sum_{\mcal{R}\in \mscr{E}(n,k)} \mu(\mcal{R}) = 0$.
\end{lemma}

\begin{proof}
Suppose $\mcal{R}\in \mscr{E}(n,k)$ and $\bm{c}\in \mcal{S}_1\cap \mcal{R}$, where $k\ge 3$.
Then by Proposition~\ref{PROP:r(n,c)-large-if-R-exceptional}, we have $r(n,\bm{c})\ge n^{k-2}$.
But by \eqref{EQN:define-r(n,c)-via-l(p,c)}, we have $r(n,\bm{c})\mid \prod_{p\mid n} p^{v_p(n)+l(p,\bm{c})} = n \prod_{p\mid n} p^{l(p,\bm{c})}$.
Therefore (since $n\ge 2$), there exists a prime $p\mid n$ with $l(p,\bm{c})\ge k-3$.
It follows that
\begin{equation*}
    \sum_{\mcal{R}\in \mscr{E}(n,k)}
    \EE_{\bm{c}\in [-Z,Z]^m}[\bm{1}_{\bm{c}\in \mcal{S}_1\cap \mcal{R}}]
    \le \sum_{p\mid n} \EE_{\bm{c}\in [-Z,Z]^m}[\bm{1}_{\bm{c}\in \mcal{S}_1} \bm{1}_{l(p,\bm{c})\ge k-3}]
\end{equation*}
for all reals $Z\ge 1$.
Taking $Z\to \infty$ (using
\eqref{INEQ:dimension-growth-bound-on-S_0}
on the left-hand side, and Lemma~\ref{LEM:density-form-of-ineffective-Krasner-lemma} on the right-hand side; cf.~the proof of Proposition~\ref{PROP:(LocAvSp)}), we get
\begin{equation}
\label{INEQ:lazy-exceptional-density-bound-via-p-adic-averages}
    \sum_{\mcal{R}\in \mscr{E}(n,k)} \mu(\mcal{R})
    \le \sum_{p\mid n} \EE_{\bm{c}\in \ZZ_p^m}[\bm{1}_{\Delta(\bm{c})\ne 0} \bm{1}_{l(p,\bm{c})\ge k-3}].
\end{equation}
But by Lemma~\ref{LEM:density-form-of-ineffective-Krasner-lemma}, the right-hand side of \eqref{INEQ:lazy-exceptional-density-bound-via-p-adic-averages} tends to $0$ as $k\to \infty$.
\end{proof}

\begin{lemma}
[EKL']
\label{LEM:(EKL')}
Assume Conjecture~\ref{CNJ:(EKL)}~(EKL).
Let $n\ge 2$ and $k\ge 3$.
Then $\sum_{\mcal{R}\in \mscr{E}(n,k)} \mu(\mcal{R}) \ll_{H,\eps} n^{k\eps} \cdot n^{-(k-3)/\deg{H}}$ (uniformly over $n$ and $k$).
\end{lemma}

\begin{proof}
Suppose $\mcal{R}\in \mscr{E}(n,k)$ and $\bm{c}\in \mcal{S}_1\cap \mcal{R}$.
As in the proof of Lemma~\ref{LEM:(KL')}, we have $r(n,\bm{c})\ge n^{k-2}$ and $r(n,\bm{c})\mid n\prod_{p\mid n} p^{l(p,\bm{c})}$.
By Remark~\ref{RMK:EKL-implication-on-l}, we have $\prod_{p\mid n} p^{l(p,\bm{c})} \mid H(\bm{c})$.
But $n\mid r(n,\bm{c})\mid n^\infty$ by \eqref{EQN:define-r(n,c)-via-l(p,c)}.
Hence $H(\bm{c})$ is divisible by the integer $q = r(n,\bm{c})/n\ge n^{k-3}$, where $q\mid n^\infty$.
Since every prime factor of $q$ is $\le n$, there exists $u\mid q$ with $n^{k-3}\le u < n^{k-2}$.
On the other hand, if $\bm{c}\in \mcal{S}_0\cap \mcal{R}$, then $H(\bm{c})=0$ (since $\Delta\mid H$), so $n^{k-3}\mid H(\bm{c})$.
So for every $\bm{c}\in \mcal{R}$ (if $\mcal{R}\in \mscr{E}(n,k)$), there exists $u\in [n^{k-3}, n^{k-2})$, with $u\mid n^\infty$, such that $u\mid H(\bm{c})$.
Thus
\begin{equation*}
    \sum_{\mcal{R}\in \mscr{E}(n,k)}
    \EE_{\bm{c}\in [-Z,Z]^m}[\bm{1}_{\bm{c}\in \mcal{R}}]
    \le \sum_{n^{k-3}\le u < n^{k-2}:\, u\mid n^\infty} \EE_{\bm{c}\in [-Z,Z]^m}[\bm{1}_{u\mid H(\bm{c})}]
\end{equation*}
for all reals $Z\ge 1$.
Taking $Z\to \infty$ (using \eqref{INEQ:multivarite-zero-density-mod-q} on the right-hand side), we get
\begin{equation*}
    \sum_{\mcal{R}\in \mscr{E}(n,k)} \mu(\mcal{R})
    \ll_H \sum_{n^{k-3}\le u < n^{k-2}:\, u\mid n^\infty} u^{-1/\deg{H}}
    \le \sum_{u\ge n^{k-3}:\, u\mid n^\infty} u^{-1/\deg{H}}.
\end{equation*}
But $\sum_{u\ge N:\,u\mid n^\infty} u^{-\beta}
\le N^{\eps-\beta} \sum_{u\ge 1:\,u\mid n^\infty} u^{-\eps}
\ll_\eps N^{\eps-\beta} n^\eps$, for all $N,\beta,\eps\in \RR_{>0}$ with $\eps<\beta$.
\end{proof}

\subsection{From (RA1'E) to (RA1'E')}
\label{SUBSEC:proving-(RA1'E')}

We now build on Propositions~\ref{PROP:RA1o-implies-RA1o'E} and~\ref{PROP:RA1delta-implies-RA1delta'E}.
Let $I\belongs \RR_{>0}$ be a compact set.
Let $\nu=\nu_{\bm{c}}(r)$ be a smooth function $\RR^m\times \RR \to \CC,\,(\bm{c}, r)\mapsto \nu_{\bm{c}}(r)$, supported on $[-1,1]^m\times I$.
Given $(\bm{a}, n_0)\in \ZZ^m\times \ZZ_{\ge 1}$ and reals $Z, N\ge 1$, let
\begin{equation}
\label{EQN:define-Sigma_11-for-RA1'E'}
\Sigma_{11}^{\bm{a},n_0}(\nu,Z,N)
\defeq \sum_{\bm{c}\in \mcal{S}_1:\, \bm{c}\equiv \bm{a}\bmod{n_0}}
\sum_{n\ge 1} \nu_{\bm{c}/Z}(n/N)
S^\natural_{\bm{c}}(n) \bm{1}_{n\in \mcal{N}^{\bm{c}}}
\end{equation}

\begin{conjecture}
[RA1$o$'E']
\label{CNJ:(RA1o'E')}
Let $M\in \RR_{\ge 1}$, and let $n_0\le M^{9/10}$ be a positive integer.
Let $Z, N \ge 2M$ be reals with $N\le Z^3$.
Then
$\sum_{1\le \bm{a}\le n_0} \abs{\Sigma_{11}^{\bm{a},n_0}(\nu,Z,N)}
\ll_{F,I} Z^m N^{1/2} \cdot (o_{F;M\to\infty}(1) + o_{F,M;Z\to\infty}(1)) \cdot \mcal{M}_{1, A_7}$
for some $A_7=A_7(F)>0$,
where $\mcal{M}_{1, k}$ is defined as in \eqref{EQN:define-norm-M_1,k(nu)}.
\end{conjecture}

The intermediate parameter $M$ here may seem strange, but it will ease our exposition.

\begin{proposition}
\label{PROP:RA1o'E-implies-RA1o'E'}
Assume Conjectures~\ref{CNJ:(HW2)}, \ref{CNJ:(R2'E)}, and~\ref{CNJ:(RA1o'E)}.
Then Conjecture~\ref{CNJ:(RA1o'E')} holds.
\end{proposition}

\begin{proof}
Plugging \eqref{EQN:decompose-good-restricted-S_c(n)-via-a_c,1,a'_c} into \eqref{EQN:define-Sigma_11-for-RA1'E'} reveals the equality
\begin{equation*}
\Sigma_{11}^{\bm{a},n_0}(\nu,Z,N)
= \sum_{\bm{c}\in \mcal{S}_1:\, \bm{c}\equiv \bm{a}\bmod{n_0}}
\sum_{n_1,d\ge 1} \nu_{\bm{c}/Z}(n_1d/N) a_{\bm{c},1}(n_1) a'_{\bm{c}}(d).
\end{equation*}
Fix a function $\upsilon_0\in C^\infty_c(\RR)$ with $\Supp{\upsilon_0}\belongs[-1,1]$ and $\upsilon_0\vert_{[-1/2,1/2]} = 1$.
Let $L\ge 1$ denote a real number to be chosen later.
We first analyze the piece
\begin{equation*}
\Sigma_{12}^{\bm{a},n_0}
\defeq \sum_{\bm{c}\in \mcal{S}_1:\, \bm{c}\equiv \bm{a}\bmod{n_0}}
\sum_{n_1,d\ge 1} (1-\upsilon_0(d/L)) \nu_{\bm{c}/Z}(n_1d/N) a_{\bm{c},1}(n_1) a'_{\bm{c}}(d)
\end{equation*}
of $\Sigma_{11}^{\bm{a},n_0}$.
For later reference, note that (since $\Supp{\upsilon_0}\belongs[-1,1]$)
\begin{align}
\Sigma_{11}^{\bm{a},n_0} - \Sigma_{12}^{\bm{a},n_0}
&= \sum_{d\le L} \upsilon_0(d/L) \Sigma_{13}^{\bm{a},n_0}(d),
\label{EQN:decompose-Sigma_11-for-RA1'E'-into-Sigma_12-pieces} \\
\textnormal{where}\quad
\Sigma_{13}^{\bm{a},n_0}(d)
&\defeq \sum_{\bm{c}\in \mcal{S}_1:\, \bm{c}\equiv \bm{a}\bmod{n_0}} a'_{\bm{c}}(d)
\sum_{n_1\ge 1} \nu_{\bm{c}/Z}(n_1d/N) a_{\bm{c},1}(n_1).
\label{EQN:define-Sigma_13(d)-for-RA1'E'}
\end{align}

We bound $\Sigma_{12}^{\bm{a},n_0}$ using the H\"{o}lder technique behind (the simplest case, $\eps=1$, of) Proposition~\ref{PROP:R2'E-implies-R2'E'}.
Since $a'_{\bm{c}}(d)$ is the $d$th coefficient of the Dirichlet series $\Phi^{\bm{c},2}(s) \Phi^{\bm{c},3}(s)$ (see \S\ref{SUBSEC:handling-variation-of-error-factors}), we may write $a'_{\bm{c}}$ in terms of $a_{\bm{c},2}$, $a_{\bm{c},3}$,
and then apply Lemma~\ref{LEM:dyadic-partial-Mellin-summation} (with $k=3$ and $a(\bm{n}) = \prod_{1\le j\le 3} a_{\bm{c},j}(n_j)$, and $f(\bm{r}) = (1-\upsilon_0(r_2r_3/L)) \nu_{\bm{c}/Z}(r_1r_2r_3/N)$), to get
\begin{equation}
\label{EQN:smooth-dyadic-restriction-separation-for-Sigma_12-in-RA1'E'}
    \Sigma_{12}^{\bm{a},n_0}
    = \sum_{\bm{c}\in \mcal{S}_1:\, \bm{c}\equiv \bm{a}\bmod{n_0}}
    (2\pi)^{-3} \int_{\bm{N}\in [1/2, \infty)^3} d^\times\bm{N} \int_{\bm{t}\in \RR^3} d\bm{t}\,
    g_{\bm{c},\bm{N}}^\vee(i\bm{t})
    \prod_{1\le j\le 3} \Sigma_{12,\bm{N}}^{\bm{c},j}(\bm{t})
\end{equation}
(cf.~\eqref{EQN:smooth-dyadic-restriction-separation-for-R2'E'}),
where $g_{\bm{c},\bm{N}}(\bm{r})\defeq (1-\upsilon_0(r_2r_3/L)) \nu_{\bm{c}/Z}(r_1r_2r_3/N) \prod_{1\le j\le 3} \nu_2(r_j/N_j)$ and
\begin{equation*}
\Sigma_{12,\bm{N}}^{\bm{c},j}(\bm{t})
\defeq \bm{1}_{\bm{c}\in [-Z,Z]^m} \sum_{n_j\ge 1} \nu_2(n_j/N_j) n_j^{-it_j} a_{\bm{c},j}(n_j).
\end{equation*}
For every integer $b\ge 0$, Proposition~\ref{PROP:standard-general-Mellin-bound} and \eqref{EQN:define-norm-M_1,k(nu)} imply (uniformly over $\bm{c}$, $\bm{N}$, $\bm{t}$)
\begin{equation}
    \label{INEQ:Mellin-decay-for-g_c,N-for-RA1'E'}
    g_{\bm{c},\bm{N}}^\vee(i\bm{t}) \ll_b \mcal{M}_{1,b}(\nu) (1+\norm{\bm{t}})^{-b}
\end{equation}
(where the implied constant may depend on $\upsilon_0$, $\nu_2$ as well as $b$).

Since $1-\upsilon_0$, $\nu_{\bm{c}/Z}$, $\nu_2$ are supported on $\RR\setminus [-1/2,1/2]$, $I$, $[1,2]$, respectively, we have $g_{\bm{c},\bm{N}}(\bm{r}) = 0$ unless $r_2r_3\ge L/2$, $r_1r_2r_3\in N\cdot I$, and $N_j\le r_j\le 2N_j$ for all $j$.
Fix an integer $A\ge 1$ satisfying $I\belongs [A^{-1}, A]$.
Then $g_{\bm{c},\bm{N}} = 0$ identically
unless $\bm{N}$ lies in the set
\begin{equation}
\label{EQN:dyadic-range-restrictions-for-R_12}
    \mscr{R}_{12}\defeq \set{\bm{N}\ge 1/2: N_2N_3\ge L/8,\;N_1N_2N_3\in [N/8A, AN]}
\end{equation}
(cf.~the region $\mscr{R}_{10}$ from \eqref{EQN:dyadic-range-restrictions-for-R_10}).
So \eqref{EQN:smooth-dyadic-restriction-separation-for-Sigma_12-in-RA1'E'} holds
even if we restrict $\bm{N}$ to $\mscr{R}_{12}$.

But if $\bm{N}\in \mscr{R}_{12}$ and $\bm{t}\in \RR^3$, then \eqref{INEQ:second-Holder-for-R2'E'} (with $\beta=1$, $(\gamma_1,\gamma_2,\gamma_3) = (2,4,4)$, $t=0$) and the subsequent arguments up to \eqref{EXPR:final-upper-bound-quantity-for-R2'E'} furnish (via Conjectures~\ref{CNJ:(HW2)} and~\ref{CNJ:(R2'E)}) the bound
\begin{equation}
\label{INEQ:main-Holder-bound-for-Sigma_12}
\int_{\bm{N}\in \mscr{R}_{12}} d^\times\bm{N}
\sum_{\bm{c}\in \mcal{S}_1} \prod_{1\le j\le 3} \abs{\Sigma_{12,\bm{N}}^{\bm{c},j}(\bm{t})}
\ll_I \int_{\bm{N}\in \mscr{R}_{12}} d^\times\bm{N}\,
(1+\norm{\bm{t}})^{A_6} \frac{Z^m N^{1/2}}{N_2^{\eta_3(4)/2} N_3^{4/30}}
\end{equation}
(where $A_6 = (A_3/2+1)\beta = A_3/2+1$, and $\eta_3(4)$ is as in Proposition~\ref{PROP:(bigwedge2E)}).
Upon taking absolute values in \eqref{EQN:smooth-dyadic-restriction-separation-for-Sigma_12-in-RA1'E'}
(after restricting $\bm{N}$ to $\mscr{R}_{12}$),
summing over $1\le \bm{a}\le n_0$, plugging in \eqref{INEQ:Mellin-decay-for-g_c,N-for-RA1'E'} (with $b = \ceil{A_6+4}$) and then \eqref{INEQ:main-Holder-bound-for-Sigma_12}, and integrating over $\bm{t}\in \RR^3$, we conclude that
\begin{equation}
\label{INEQ:final-Sigma_12-total-bound-under-HW2-and-R2'E}
    \sum_{1\le \bm{a}\le n_0} \abs{\Sigma_{12}^{\bm{a},n_0}}
    \ll_I \mcal{M}_{1,b}(\nu)
    \int_{\bm{N}\in \mscr{R}_{12}} d^\times\bm{N}\,
    \frac{Z^m N^{1/2}}{N_2^{\eta_3(4)/2} N_3^{4/30}}
    \ll_I \mcal{M}_{1,b}(\nu) \frac{Z^m N^{1/2}}{L^{\eta_4}},
\end{equation}
where $b = \ceil{A_6+4}$ and $\eta_4 = 0.9 \min(\eta_3(4)/2, 4/30)$.

We now turn to the sums $\Sigma_{13}^{\bm{a},n_0}(d)$, for $d\le L$.
We first treat $d=1$.
By \eqref{EQN:define-Sigma_13(d)-for-RA1'E'}, we have
$\Sigma_{13}^{\bm{a},n_0}(1)
= \sum_{\bm{c}\in \mcal{S}_1:\, \bm{c}\equiv \bm{a}\bmod{n_0}}
\sum_{n\ge 1} \nu_{\bm{c}/Z}(n/N) a_{\bm{c},1}(n)$.
Therefore, by the triangle inequality, $\sum_{1\le \bm{a}\le n_0} \abs{\Sigma_{13}^{\bm{a},n_0}(1)}$ is at most the sum of the quantity \eqref{EXPR:RA1'E-main-quantity} (with $\mscr{P} = \set{\bm{a}+n_0\ZZ^m: 1\le \bm{a}\le n_0}$) and the left-hand side of \eqref{INEQ:RA1'E-main-term-bound-goal} (with $\mcal{S}=\mcal{S}_1$).
So by Conjecture~\ref{CNJ:(RA1o'E)} (applicable since $Z, N\ge 2M\ge 2n_0 = 2Q(\mscr{P})$) and Lemma~\ref{LEM:basic-main-term-bound-for-RA1'E'} (with $\theta=m$ and $Q(\mscr{P})=n_0\le M$), the sum $\sum_{1\le \bm{a}\le n_0} \abs{\Sigma_{13}^{\bm{a},n_0}(1)}$ is at most $O_{I,\eps}(\mf{B}_0(N,\eps))$, where $\mf{B}_0(N,\eps)$ denotes the the expression\footnote{The replacement of $n_0$ with the weaker (larger) $M$ here may seem strange, but it will be convenient later.}
\begin{equation}
\label{EXPR:frakB_0-soft-bound-for-Sigma_13(1)-total}
    Z^m N^{1/2} o_{M;Z\to\infty}(1) \mcal{M}_{1, A_4}
    + Z^m N^{1/3+\eps} M^\eps \mcal{M}_{1, 0}.
\end{equation}
If $\eps$ is sufficiently small, then (since $M\le N$)
\begin{equation}
\label{INEQ:soft-prep-bound-on-frakB_0(N,eps)}
    \mf{B}_0(N,\eps) \le Z^m N^{1/2} (o_{M;Z\to\infty}(1) + N^{-1/6+2\eps}) \mcal{M}_{1, A_4}.
\end{equation}

Now suppose $2\le d\le L$.
Let $k=k(d)\ge 1$ denote an integer, with $n_0 d^{k+1} \le M$, to be chosen later.
Using Definition~\ref{DEFN:algorithmic-construction-of-partition-P(n,k)-exceptions-E(n,k)}, let
\begin{equation*}
\begin{split}
\Sigma_{13,0}^{\bm{a},n_0}(d)
&\defeq \sum_{\mcal{R}\in \mscr{E}(d,k)}
\sum_{\bm{c}\in \mcal{S}_1\cap \mcal{R}\cap (\bm{a}+n_0\ZZ^m)} a'_{\bm{c}}(d)
\sum_{n_1\ge 1} \nu_{\bm{c}/Z}(n_1d/N) a_{\bm{c},1}(n_1), \\
\Sigma_{13,1}^{\bm{a},n_0}(d)
&\defeq \sum_{\mcal{R}\in \mscr{P}(d,k)\setminus \mscr{E}(d)}
\sum_{\bm{c}\in \mcal{S}_1\cap \mcal{R}\cap (\bm{a}+n_0\ZZ^m)} a'_{\bm{c}}(d)
\sum_{n_1\ge 1} \nu_{\bm{c}/Z}(n_1d/N) a_{\bm{c},1}(n_1).
\end{split}
\end{equation*}
Then $\Sigma_{13}^{\bm{a},n_0}(d) = \Sigma_{13,0}^{\bm{a},n_0}(d) + \Sigma_{13,1}^{\bm{a},n_0}(d)$ by \eqref{EQN:define-Sigma_13(d)-for-RA1'E'} (since $\mscr{P}(d,k)$ is a partition of $\ZZ^m$).
By Proposition~\ref{PROP:basic-properties-of-residue-classes-R-in-partition-P(n,k)}, we have $q_{\mcal{R}}\le d^k$ for all $\mcal{R}\in \mscr{P}(d,k)$, so $Q(\mscr{P}(d,k))\le d^k$.

By the triangle inequality,
\begin{equation}
\label{INEQ:triangle-inequality-bound-for-Sigma_13,0(d)-total}
\sum_{1\le \bm{a}\le n_0} \abs{\Sigma_{13,0}^{\bm{a},n_0}(d)}
\le \sum_{\mcal{R}\in \mscr{E}(d,k)}
\sum_{\bm{c}\in \mcal{S}_1\cap \mcal{R}} \abs{a'_{\bm{c}}(d)}
\Bigl\lvert{
\sum_{n_1\ge 1} \nu_{\bm{c}/Z}(n_1d/N) a_{\bm{c},1}(n_1)
}\Bigr\rvert.
\end{equation}
By \eqref{INEQ:GRC-on-a'_c}, $a'_{\bm{c}}(d)\ll_\eps d^\eps$.
However, by Lemma~\ref{LEM:dyadic-partial-Mellin-summation} (with $k=1$ and $a(n_1) = a_{\bm{c},1}(n_1)$, and $f(r_1) = \nu_{\bm{c}/Z}(r_1d/N)$) and Proposition~\ref{PROP:standard-general-Mellin-bound}, we have (cf.~\eqref{EQN:smooth-dyadic-restriction-separation-for-Sigma_12-in-RA1'E'}, \eqref{INEQ:Mellin-decay-for-g_c,N-for-RA1'E'}, and \eqref{EQN:dyadic-range-restrictions-for-R_12})
\begin{equation*}
    \sum_{n_1\ge 1} \nu_{\bm{c}/Z}(n_1d/N) a_{\bm{c},1}(n_1)
    \ll_b \int_{1/2}^{AN/d} d^\times{N} \int_{t_1\in \RR} dt_1\,
    \mcal{M}_{1,b}(\nu) \frac{\abs{\Sigma_{12,(N_1,1,1)}^{\bm{c},1}(t_1,0,0)}}{(1+\abs{t_1})^b}
\end{equation*}
for all integers $b\ge 0$.
Plugging this and $a'_{\bm{c}}(d)\ll_\eps d^\eps$ into \eqref{INEQ:triangle-inequality-bound-for-Sigma_13,0(d)-total}, and then applying Cauchy--Schwarz over $\bm{c}\in \bigcup_{\mcal{R}\in \mscr{E}(d,k)} (\mcal{R}\cap [-Z,Z]^m)$, we get (by Lemma~\ref{LEM:(KL')} and Conjecture~\ref{CNJ:(R2'E)})
\begin{equation}
\label{INEQ:soft-Cauchy-bound-for-Sigma_13,0(d)-total-via-soft-KL'-and-R2'E}
\sum_{1\le \bm{a}\le n_0} \abs{\Sigma_{13,0}^{\bm{a},n_0}(d)}
\ll_{I,\eps} d^\eps \mcal{M}_{1,\ceil{A_3/2+2}}(\nu)
\cdot o_{d;k\to\infty}(Z^m)^{1/2} \cdot (Z^m N/d)^{1/2};
\end{equation}
cf.~the numerics in \eqref{INEQ:main-Holder-bound-for-Sigma_12} and \eqref{INEQ:final-Sigma_12-total-bound-under-HW2-and-R2'E}.
(Before applying Lemma~\ref{LEM:(KL')}, note that $Q(\mscr{E}(d,k))\le d^k\le M\le Z$, so $\card{\bigcup_{\mcal{R}\in \mscr{E}(d,k)} (\mcal{R}\cap [-Z,Z]^m)} \ll Z^m \sum_{\mcal{R}\in \mscr{E}(d,k)} \mu(\mcal{R})$.)

By Definition~\ref{DEFN:algorithmic-construction-of-partition-P(n,k)-exceptions-E(n,k)}, the quantity $a'_{\bm{c}}(d)$ is constant over $\mcal{S}_1\cap \mcal{R}$ for each $\mcal{R}\in \mscr{P}(d,k)\setminus \mscr{E}(d)$.
But $a'_{\bm{c}}(d)\ll_\eps d^\eps$ by \eqref{INEQ:GRC-on-a'_c}, so we get
\begin{equation*}
\Sigma_{13,1}^{\bm{a},n_0}(d)
\ll_\eps d^\eps \sum_{\mcal{R}\in \mscr{P}(d,k)\setminus \mscr{E}(d)}\,
\Bigl\lvert{
\sum_{\bm{c}\in \mcal{S}_1\cap \mcal{R}\cap (\bm{a}+n_0\ZZ^m)}
\sum_{n_1\ge 1} \nu_{\bm{c}/Z}(n_1d/N) a_{\bm{c},1}(n_1)
}\Bigr\rvert.
\end{equation*}
We may apply Conjecture~\ref{CNJ:(RA1o'E)} (with $N/d$ in place of $N$) and Lemma~\ref{LEM:basic-main-term-bound-for-RA1'E'} (with $\mcal{S}=\mcal{S}_1$ and $\theta=m$) with $\mscr{P} = \set{\mcal{R}\cap (\bm{a}+n_0\ZZ^m): \mcal{R}\in \mscr{P}(d,k),\; 1\le \bm{a}\le n_0}$, since $Q(\mscr{P})\le n_0 d^k\le M/d$ and $Z, N\ge 2M$ (so that $2M, Z, N/d\ge 2Q(\mscr{P})$ and $N/d\le N\le Z^3$).
By the triangle inequality applied to the right-hand side of the previous display, we then obtain the bound
\begin{equation}
\label{INEQ:Sigma_13,1(d)-total-bound-via-frakB_0(N/d)}
\sum_{1\le \bm{a}\le n_0} \abs{\Sigma_{13,1}^{\bm{a},n_0}(d)}
\ll_{I,\eps} d^\eps \mf{B}_0(N/d,\eps)
\quad \textnormal{(where $\mf{B}_0(N,\eps)$ denotes \eqref{EXPR:frakB_0-soft-bound-for-Sigma_13(1)-total} as before)}.
\end{equation}

Let $K\ge 1$ denote an integer to be chosen soon.
Assembling \eqref{EQN:decompose-Sigma_11-for-RA1'E'-into-Sigma_12-pieces}, \eqref{INEQ:final-Sigma_12-total-bound-under-HW2-and-R2'E}, and our work on $\Sigma_{13}^{\bm{a},n_0}(d)$ for $d\le L$ (see \eqref{EXPR:frakB_0-soft-bound-for-Sigma_13(1)-total} for $d=1$, and \eqref{INEQ:soft-Cauchy-bound-for-Sigma_13,0(d)-total-via-soft-KL'-and-R2'E}, \eqref{INEQ:Sigma_13,1(d)-total-bound-via-frakB_0(N/d)} for $2\le d\le L$), we get (by the triangle inequality) that the sum $\sum_{1\le \bm{a}\le n_0} \abs{\Sigma_{11}^{\bm{a},n_0}}$ is at most
\begin{equation}
\label{INEQ:near-final-soft-upper-bound-for-Sigma_11-total}
\ll_{I,\eps} \mcal{M}_{1,\ceil{A_3/2+5}}(\nu) (Z^m N^{1/2} L^{-\eta_4}
+ o_{L;K\to\infty}(Z^m N^{1/2}))
+ \sum_{d\le L} d^\eps \mf{B}_0(N/d,\eps),
\end{equation}
provided $k(d)\ge K$ and $n_0 d^{k(d)+1}\le M$ hold for all integers $d$ with $2\le d\le L$.
But upon replacing $N$ in \eqref{EXPR:frakB_0-soft-bound-for-Sigma_13(1)-total} with $N/d$, and summing over $d\le L$, we find that
\begin{equation}
\label{INEQ:sum-frakB_0-over-N/d}
\sum_{d\le L} d^\eps \mf{B}_0(N/d,\eps)
\ll L^{2/3} \mf{B}_0(N,\eps).
\end{equation}
It remains to carefully specify parameters.
Choose $K=K(L)\ge 1$
so that $o_{L;K\to\infty}(Z^m N^{1/2})\le L^{-\eta_4} Z^m N^{1/2}$.
Then let $k(d) = K(L)$ for all integers $d$ with $2\le d\le L$.
Let $L=L(M)\ge 1$
(in terms of $M$)
be the largest integer for which $L^{K(L)+1}\le M^{1/10}$; such an integer exists, because $1^{K(1)+1} \le M^{1/10}$ and $L^{K(L)+1} \ge L^2$.
Crucially, we have $\lim_{M\to \infty}{L(M)} = \infty$,
since the expression $L^{K(L)+1}$ is bounded on any finite set of integers $L$.

Since $n_0\le M^{9/10}$, we have $n_0 L^{K(L)+1} \le M$.
So the bound \eqref{INEQ:near-final-soft-upper-bound-for-Sigma_11-total} on $\sum_{1\le \bm{a}\le n_0} \abs{\Sigma_{11}^{\bm{a},n_0}}$ is valid, and therefore (by \eqref{INEQ:sum-frakB_0-over-N/d}) we get (letting $A_7 = \max(A_4, \ceil{A_3/2+5})$)
\begin{equation*}
\sum_{1\le \bm{a}\le n_0} \abs{\Sigma_{11}^{\bm{a},n_0}}
\ll_{I,\eps} \mcal{M}_{1,A_7}(\nu) Z^m N^{1/2} L(M)^{-\eta_4}
+ L(M)^{2/3} \mf{B}_0(N,\eps).
\end{equation*}
But $L(M)^2\le L(M)^{K(L(M))+1}\le M^{1/10}$, and $N\ge M$, so by \eqref{INEQ:soft-prep-bound-on-frakB_0(N,eps)}, we have
\begin{equation*}
L(M)^{2/3} \mf{B}_0(N,\eps)
\ll Z^m N^{1/2} (M^{1/30-1/6+2\eps} + o_{M;Z\to\infty}(1)) \mcal{M}_{1, A_7},
\end{equation*}
provided $\eps$ is sufficiently small.
Both $L(M)^{-\eta_4}$ and $M^{1/30-1/6+2\eps}$ are $o_{M\to\infty}(1)$.
\end{proof}

\begin{proposition}
[RA1$\delta$'E']
\label{PROP:RA1delta'E-implies-RA1delta'E'}
Assume Conjectures~\ref{CNJ:(HW2)}, \ref{CNJ:(RA1delta)}, and~\ref{CNJ:(EKL)}.
Suppose $Z, N \ge 2M$ and $N\le Z^3$.
Suppose $1\le M\le Z^{\eta_2}$
(with $\eta_2$ as in Proposition~\ref{PROP:RA1delta-implies-RA1delta'E}),
and let $n_0\le M^{1/2}$ be a positive integer.
Then
$\sum_{1\le \bm{a}\le n_0} \abs{\Sigma_{11}^{\bm{a},n_0}(\nu,Z,N)}
\ll_{F,I,\eps} Z^{m+\eps} N^{1/2} \cdot (M^{1/6\deg{H}} N^{-1/6} + M^{-\eta_5/4\deg{H}}) \cdot \mcal{M}_{1, A_8}$.
Here $\eta_5$, $A_8$ are positive reals depending only on $F$.
\end{proposition}

\begin{proof}
We adjust the proof of Proposition~\ref{PROP:RA1o'E-implies-RA1o'E'}.
Let $L\in [1, M^{1/8}]$ be a real to be chosen later.
For each integer $d$ with $2\le d\le L$, let $k(d)$ be the largest integer $k$ for which $d^{k+1}\le M^{1/2}$; then $k(d)\ge 3$ (since $d\le L \le M^{1/8}$), and $d^{k(d)+2}>M^{1/2}$ (by the maximality of $k(d)$).
Let
\begin{equation*}
\mf{B}_1(N,\eps) \defeq
Z^{m-\eta_2} N^{1/2} \mcal{M}_{1, A_5}
+ Z^m N^{1/3+\eps} M^\eps \mcal{M}_{1, 0}.
\end{equation*}
Using Proposition~\ref{PROP:RA1delta-implies-RA1delta'E}
in place of Conjecture~\ref{CNJ:(RA1o'E)}, we find that $\sum_{1\le \bm{a}\le n_0} \abs{\Sigma_{13}^{\bm{a},n_0}(1)} \ll_{I,\eps} \mf{B}_1(N,\eps)$ and (if $2\le d\le L$, then) $\sum_{1\le \bm{a}\le n_0} \abs{\Sigma_{13,1}^{\bm{a},n_0}(d)} \ll_{I,\eps} d^\eps \mf{B}_1(N/d,\eps)$; cf.~\eqref{EXPR:frakB_0-soft-bound-for-Sigma_13(1)-total} and \eqref{INEQ:Sigma_13,1(d)-total-bound-via-frakB_0(N/d)}.
Summing over $1\le d\le L$, and writing $\Sigma_{13,1}^{\bm{a},n_0}(1)\defeq \Sigma_{13}^{\bm{a},n_0}(1)$ for convenience, we obtain
\begin{equation}
\label{INEQ:total-hard-Sigma_13,1(d)-bound-over-d}
\sum_{1\le d\le L} \sum_{1\le \bm{a}\le n_0} \abs{\Sigma_{13,1}^{\bm{a},n_0}(d)}
\ll_{I,\eps} L^{2/3} \mf{B}_1(N,\eps)
\le L^{2/3} Z^{m+4\eps} N^{1/2} (Z^{-\eta_2} + N^{-1/6}) \mcal{M}_{1,A_5}.
\end{equation}

On the other hand, if we plug GRH (Proposition~\ref{PROP:HW2-consequences}(\ref{ITEM:1/L-bound-from-GRH})) into \eqref{INEQ:triangle-inequality-bound-for-Sigma_13,0(d)-total} (after partial summation over $n_1$, say, and recalling the definition of $\mcal{M}_{1,1}$ from \eqref{EQN:define-norm-M_1,k(nu)}) and then use Lemma~\ref{LEM:(EKL')} (applicable since $k(d)\ge 3$ and we assume Conjecture~\ref{CNJ:(EKL)}), then (if $2\le d\le L$) we get
\begin{equation}
\label{INEQ:hard-bound-for-Sigma_13,0(d)-total-via-EKL'-and-HW2}
\sum_{1\le \bm{a}\le n_0} \abs{\Sigma_{13,0}^{\bm{a},n_0}(d)}
\ll_{I,\eps} d^\eps \mcal{M}_{1,1}(\nu) Z^\eps (N/d)^{1/2+\eps}
\cdot Z^m d^{k(d)\eps} d^{-(k(d)-3)/\deg{H}};
\end{equation}
cf.~the use of Lemma~\ref{LEM:(KL')} towards \eqref{INEQ:soft-Cauchy-bound-for-Sigma_13,0(d)-total-via-soft-KL'-and-R2'E}.
On the right-hand side,
$d^\eps Z^\eps (N/d)^\eps d^{k(d)\eps} \le Z^{5\eps}$ (since $N\le Z^3$ and $d^{k(d)}\le M$), and
$d^{(k(d)-3)/\deg{H}} > M^{1/2\deg{H}} d^{-5/\deg{H}} \ge M^{1/2\deg{H}} d^{-5/12}$
(since $d^{k(d)+2} > M^{1/2}$ and $\deg{H} \ge 12$).
Thus, summing \eqref{INEQ:hard-bound-for-Sigma_13,0(d)-total-via-EKL'-and-HW2} over $2\le d\le L$ gives
\begin{equation}
\label{INEQ:total-hard-Sigma_13,0(d)-bound-over-d}
\sum_{2\le d\le L} \sum_{1\le \bm{a}\le n_0} \abs{\Sigma_{13,0}^{\bm{a},n_0}(d)}
\ll_{I,\eps} \mcal{M}_{1,1}(\nu) Z^{m+5\eps} N^{1/2} L^{11/12} \cdot M^{-1/2\deg{H}}.
\end{equation}

As for $\Sigma_{12}^{\bm{a},n_0}$, the bounds \eqref{INEQ:main-Holder-bound-for-Sigma_12} and \eqref{INEQ:final-Sigma_12-total-bound-under-HW2-and-R2'E} must be adjusted slightly, since we do not assume Conjecture~\ref{CNJ:(R2'E)}.
However, by Proposition~\ref{PROP:HW2-consequences}(\ref{ITEM:1/L-bound-from-GRH}) and partial summation, the bound \eqref{INEQ:R2'E-goal} in Conjecture~\ref{CNJ:(R2'E)} still holds up to a factor of $Z^\eps$ (for any $A_3>0$).
Therefore, \eqref{INEQ:main-Holder-bound-for-Sigma_12} and \eqref{INEQ:final-Sigma_12-total-bound-under-HW2-and-R2'E} still hold if we replace $Z^m N^{1/2}$ with $Z^{m+\eps} N^{1/2}$; so $\sum_{1\le \bm{a}\le n_0} \abs{\Sigma_{12}^{\bm{a},n_0}} \ll_I \mcal{M}_{1,\ceil{A_3/2+5}}(\nu) Z^m N^{1/2}/L^{\eta_4}$.

Let $A_8 = \max(A_5, \ceil{A_3/2+5})$.
Assembling \eqref{EQN:decompose-Sigma_11-for-RA1'E'-into-Sigma_12-pieces}, \eqref{INEQ:total-hard-Sigma_13,1(d)-bound-over-d}, \eqref{INEQ:total-hard-Sigma_13,0(d)-bound-over-d}, and our work on $\Sigma_{12}^{\bm{a},n_0}$, we get (by the triangle inequality) that the sum $\sum_{1\le \bm{a}\le n_0} \abs{\Sigma_{11}^{\bm{a},n_0}}$ is
\begin{equation*}
\ll_{I,\eps} \mcal{M}_{1,A_8}(\nu) Z^{m+5\eps} N^{1/2} (L^{2/3} (Z^{-\eta_2} + N^{-1/6}) + L^{11/12} M^{-1/2\deg{H}} + L^{-\eta_4}).
\end{equation*}
Let $L = M^{1/4\deg{H}}$ and $\eta_5 = \min(1,\eta_4)$ to finish.
\end{proof}

\section{New bounds on the integral factor}
\label{SEC:new-bounds-on-integral-J}

Recall $J_{\bm{c},X}(n)$ from \eqref{EQN:define-normalized-S-tilde-and-J}.
As we explained in \S\ref{SEC:background-on-discriminants-and-the-delta-method},
we need to go beyond the integral estimates from standard sources like \cites{duke1993bounds,heath1996new,heath1998circle,hooley2014octonary} (and the related estimates of \cites{hooley1986HasseWeil,hooley_greaves_harman_huxley_1997}), such as \eqref{INEQ:Hooley-best-Hessian-free-integral-bound-at-least-for-textbook-weights?}.
We will prove uniform bounds free of epsilons and logs, while also bringing discriminants into the picture
via the following consequence of \eqref{INEQ:Disc-grad-composite-divisibility-by-F}:
\begin{equation}
\label{INEQ:bound-Disc-grad-composite-in-terms-of-F}
\Delta(\grad{F}(\bm{x}))
\ll_F \abs{F(\bm{x})}\cdot \norm{\bm{x}}^{2\deg(\Delta)-\deg(F)}
\quad \textnormal{for $\bm{x}\in \RR^m$}.
\end{equation}
To give clean, general bounds, we assume \eqref{COND:clean-weight-condition-in-general}.
We prove the following on $J_{\bm{c},X}(n)$:

\begin{proposition}
\label{PROP:uniform-discriminating-integral-bound-with-derivatives}
Assume \eqref{COND:clean-weight-condition-in-general}.
Let $X,Z,n\in \RR_{>0}$.
Let $\bm{c}\in \RR^m$ with $\norm{\bm{c}}\le Z$.
Then
\begin{equation*}
\partial_{\log{n}}^j \partial_{\bm{c}}^{\bm{\alpha}} J_{\bm{c},X}(n)
\ll_{F,w,j,\bm{\alpha},b} \frac{(X/n)^{\abs{\bm{\alpha}}}
(1+X\norm{\bm{c}}/n)^{1-m/2}}
{(1+\norm{\bm{c}}/X^{1/2})^b
(1+X\norm{\Delta(\bm{c}/Z)\bm{c}}/n)^b}
\end{equation*}
for all integers $j, b\ge 0$ and multi-indices $\bm{\alpha}\ge 0$.
\end{proposition}

The fact that increasing $j$ is harmless can be interpreted as an instance of ``homogeneous dimensional analysis'' (and ultimately arises from the homogeneity of $F$, via a beautiful recursive structure due to \cite{heath1996new}; see \eqref{EQN:recursive-differentiate-generalized-J-by-log-n} below).
The factor $\Delta(\bm{c}/Z)\ll 1$ measures the ``degeneracy'' of the real hyperplane section $F(\bm{x})=\bm{c}\cdot\bm{x}=0$, and it arises in our proof for roughly the same reason that dual hypersurfaces arise in \cite{huang2020density}*{(5.4), (5.11)--(5.14)}.

Morally, in \eqref{EQN:delta-method}, Proposition~\ref{PROP:uniform-discriminating-integral-bound-with-derivatives} lets us ``imagine that there are sharp cutoffs'' $\norm{\bm{c}}\ll X^{1/2}$ and $n\gg X\norm{\Delta(\bm{c}/X^{1/2})\bm{c}}$.
Since by \eqref{INEQ:multivarite-near-zero-real-density} we typically have $\abs{\Delta(\bm{c}/X^{1/2})}\asymp 1$, one might thus expect (in view of Proposition~\ref{PROP:basic-integral-facts}, and our $\eps$-diagnosis at the end of \S\ref{SEC:background-on-discriminants-and-the-delta-method}) that $n\asymp X^{3/2}$ should be the ``dominant range'' on average, and there we have $J_{\bm{c},X}(n) \ll 1$.

To prove Proposition~\ref{PROP:uniform-discriminating-integral-bound-with-derivatives}, we must first dig into some technical aspects of \cites{duke1993bounds,heath1996new}'s $h$-function.
Let $\omega_{\textnormal{HB},0}(x)\defeq \bm{1}_{\abs{x}\leq 1}\cdot \exp(-(1-x^2)^{-1})\in C^\infty_c(\RR)$ and $c_{\textnormal{HB},0}\defeq \int_{x\in \RR} dx\,\omega_{\textnormal{HB},0}(x)$; note that $\omega_{\textnormal{HB},0}$ is supported on $[-1,1]$.
Following \cite{heath1996new}*{p.~165}, let
\begin{equation*}
    \omega_\textnormal{HB}(x)\defeq 4c_{\textnormal{HB},0}^{-1}
    \cdot \omega_{\textnormal{HB},0}(4x-3)\in C^\infty_c(\RR),
\end{equation*}
so that $\omega_\textnormal{HB}$ is supported on $[1/2,1]$.
For $x>0$ and $y\in\RR$, let
\begin{equation*}
    h(x,y) \defeq \sum_{j\geq 1} (xj)^{-1}[\omega_\textnormal{HB}(xj) - \omega_\textnormal{HB}(\abs{y}/(xj))].
\end{equation*}
By \eqref{EQN:define-I_c(n)} and \eqref{EQN:define-normalized-S-tilde-and-J}, and a change of variables from $\bm{x}$ to $\tilde{\bm{x}} = \bm{x}/X$, we have (since $Y^2 = X^3$)
\begin{equation}
\label{J_c,X-normalized-expression}
J_{\bm{c},X}(n)
= \int_{\tilde{\bm{x}}\in\RR^m}
d\tilde{\bm{x}}\,w(\tilde{\bm{x}})
h(n/Y, F(\tilde{\bm{x}}))
e(-X\bm{c}\cdot\tilde{\bm{x}}/n).
\end{equation}
The following shows that we may take $A_0 = \sup_{\bm{x}\in \Supp{w}} \max(1,2\abs{F(\bm{x})})$ in Proposition~\ref{PROP:basic-integral-facts}.

\begin{lemma}
[See \cite{heath1996new}*{Lemma~4}]
\label{LEM:basic-vanishing-of-h(x,y)}
If $\xi\in \RR$ and $r\ge \max(1,2\abs{\xi})$, then $h(r,\xi) = 0$.
\end{lemma}


Following \cite{heath1996new}, we now build a Fourier transform.
Fix $\upsilon_1\in C^\infty_c(\RR)$ such that $\upsilon_1(F(\bm{x})) \ge 1$ for all $\bm{x}\in \Supp{w}$.
Let $w_0(\bm{x})\defeq w(\bm{x})/\upsilon_1(F(\bm{x}))\in C^\infty_c(\RR)$.
For any $r>0$, let
\begin{equation*}
p_r(u)\defeq \int_{\xi\in\RR}d\xi\,\upsilon_1(\xi)h(r,\xi)e(-u\xi)
\end{equation*}
be the Fourier transform of $\upsilon_1(\xi)h(r,\xi)\in C^\infty_c(\RR)$.
Writing $w(\bm{x}) = w_0(\bm{x}) \upsilon_1(F(\bm{x}))$ and $\upsilon_1(\xi) h(r,\xi) = \int_{u\in \RR} du\, p_r(u) e(u \xi)$ (for $r = n/Y$ and $\xi = F(\bm{x})$) in \eqref{J_c,X-normalized-expression}, we get
\begin{equation}
\label{EQN:rewrite-Heath-Brown-by-Fourier}
J_{\bm{c},X}(n)
= \int_{u\in\RR}du\,p_r(u)
\int_{\tilde{\bm{x}}\in\RR^m}d\tilde{\bm{x}}\,w_0(\tilde{\bm{x}})e(uF(\tilde{\bm{x}})-\bm{v}\cdot\tilde{\bm{x}}),
\end{equation}
where $(r, \bm{v}) = (n/Y, X\bm{c}/n)$; cf.~\cite{heath1996new}*{Lemma~17}.

Let $\mf{q}_r(t)\defeq p_r(t/r)$, so that $p_r(u) = \mf{q}_r(ru)$.
It turns out (see Lemma~\ref{LEM:bounds-on-q_r(t)-and-its-t-and-r-derivatives}) that $\mf{q}_r$ behaves somewhat like a ``fixed'' Schwartz function independent of $r$.
Thus $u$ may be compared with $\beta Y^2$ in the classical Dirichlet arc theory, where $\abs{\beta}\leq 1/(nY) = 1/(rY^2)$.

\begin{lemma}
\label{LEM:bounds-on-q_r(t)-and-its-t-and-r-derivatives}
Let $r>0$ and $t\in \RR$.
Then $\partial_{\log{r}}^j \mf{q}_r^{(l)}(t) \ll_{l,j,k} (1+\abs{t})^{-k}$ for all $j,k,l\in \ZZ_{\ge 0}$.
\end{lemma}

\begin{proof}
We may assume $k=0$ if $\abs{t}\le 1$; this lets us treat all $t$ simultaneously.

Since $\mf{q}_r(t) = p_r(t/r) = \int_{\xi\in \RR} d\xi\, \upsilon_1(\xi) h(r,\xi) e(-t\xi/r)$, we have
\begin{equation*}
\mf{q}_r^{(l)}(t)
= \int_{\xi\in \RR} d\xi\,
\upsilon_1(\xi) h(r,\xi) \cdot (-2\pi i\xi/r)^l e(-t\xi/r).
\end{equation*}
Since $\partial_{\log{r}}((\xi/r)^l) = -l (\xi/r)^l$ and $\partial_{\log{r}}(e(-t\xi/r)) = (2\pi it\xi/r) e(-t\xi/r)$, we find by induction (and the Leibniz rule) that for some constants $c^{l,j}_{a,b}\in \CC$, we have
\begin{equation*}
\partial_{\log{r}}^j \mf{q}_r^{(l)}(t)
= \sum_{a,b\ge 0:\,a+b\le j} c^{l,j}_{a,b} \int_{\xi\in \RR} d\xi\,
\upsilon_1(\xi) \cdot (\partial_{\log{r}}^a h(r,\xi)) \cdot (\xi/r)^l (t\xi/r)^b e(-t\xi/r).
\end{equation*}
Integrating by parts $k$ times in $\xi$ (repeatedly integrating $e(-t\xi/r)$ and differentiating the complementary factor), and then taking absolute values, we get
\begin{equation*}
\partial_{\log{r}}^j \mf{q}_r^{(l)}(t)
\ll_{l,j,k,\upsilon_1} \sum_{\substack{a,b,\alpha,\beta\ge 0: \\ a+b\le j,\; \alpha+\beta\le k,\; \beta\le l+b}} \int_{\xi\in \Supp{\upsilon_1}} d\xi\,
(\partial_\xi^\alpha \partial_{\log{r}}^a h(r,\xi)) \cdot \frac{\abs{\xi/r}^l \abs{t\xi/r}^b}{\abs{\xi}^\beta} \frac{1}{\abs{t/r}^k}.
\end{equation*}

By Lemma~\ref{LEM:basic-vanishing-of-h(x,y)}, we may assume $r\ll_{\upsilon_1} 1$.
Then $\partial_r^c \partial_\xi^\alpha h(r,\xi) \ll_{c,\alpha,A} r^{-1-c-\alpha} \min(1, (r/\abs{\xi})^A)$ for $c,\alpha,A\ge 0$,
by \cite{heath1996new}*{Lemma~5}.
Thus $\partial_{\log{r}}^a \partial_\xi^\alpha h(r,\xi) \ll_a \sum_{0\le c\le a} r^c\cdot \abs{\partial_r^c \partial_\xi^\alpha h(r,\xi)} \ll_{a,\alpha,A} r^{-1-\alpha} \min(1,(r/\abs{\xi})^A)$, which when inserted in the previous display gives
\begin{equation*}
\partial_{\log{r}}^j \mf{q}_r^{(l)}(t)
\ll_{l,j,k,A} \sum_{\substack{a,b,\alpha,\beta\ge 0: \\ a+b\le j,\; \alpha+\beta\le k,\; \beta\le l+b}} \int_{\xi\in \RR} d\xi\,
\frac{r^{-1-\alpha-l-b+k} \abs{\xi}^{l+b-\beta} \abs{t}^{b-k}}{(1+\abs{\xi/r})^A}.
\end{equation*}
Let $A = l+j+2$; then each integral here over $\xi$ is $\ll_{l,j,k} r^{-1-\alpha-l-b+k} r^{1+l+b-\beta} \abs{t}^{b-k}$ (by Lemma~\ref{LEM:general-dyadic-sum-split-into-2-geometric-series}), and thus $\ll_{l,j,k} r^{k-\alpha-\beta} \abs{t}^{b-k} \ll \abs{t}^{b-k}$ (since $k-\alpha-\beta\ge 0$ and $r\ll 1$).
If $\abs{t}\le 1$, we are done (since $k=0$).
If $\abs{t}>1$, we may replace $k$ with $k+j$ to finish (since $b\le j$).
\end{proof}

We now put \eqref{EQN:rewrite-Heath-Brown-by-Fourier} in a broader framework.
For any Schwartz functions $q\maps \RR \to \CC$ and $\phi\in C^\infty_c(\RR^m)\otimes \CC$ with $\Supp{\phi}\belongs \Supp{w}$, let
\begin{equation}
\label{EQN:define-generalized-integral-J_r,v,q,phi}
\mscr{J}_{r,\bm{v}}(q,\phi) \defeq
\int_{u\in \RR} du\,q(ru)
\int_{\bm{x}\in \RR^m} d\bm{x}\,\phi(\bm{x})e(uF(\bm{x})-\bm{v}\cdot\bm{x}).
\end{equation}
Then by \eqref{EQN:rewrite-Heath-Brown-by-Fourier}, we have $J_{\bm{c},X}(n) = \mscr{J}_{r,\bm{v}}(\mf{q}_r,w_0)$, where $(r, \bm{v}) = (n/Y, X\bm{c}/n)$.

\begin{proposition}
\label{PROP:uniform-integral-estimate-for-clean-weights}
Assume \eqref{COND:clean-weight-condition-in-general}.
Let $q$, $\phi$ be as above (with $q$, $\phi$ Schwartz and $\Supp\phi\belongs\Supp{w}$).
Then for all $(k,\bm{v})\in\ZZ_{\geq0}\times\RR^m$
and positive reals $r\le A_0$ and $M\ge \norm{\bm{v}}$,
we have
\begin{equation*}
\mscr{J}_{r,\bm{v}}(q,\phi) \ll_{F,w,q,\phi,k}
(1+\norm{\bm{v}})^{1-m/2}
\cdot (1+\norm{r\bm{v}})^{-k} (1+\norm{\Delta(\bm{v}/M)\bm{v}})^{-k}.
\end{equation*}
\end{proposition}

Before proving Proposition~\ref{PROP:uniform-integral-estimate-for-clean-weights},
we first explain why it implies the desired Proposition~\ref{PROP:uniform-discriminating-integral-bound-with-derivatives}.
In order to handle $\partial_{\log{n}}^{\geq1}$,
we need a recursion, \eqref{EQN:recursive-differentiate-generalized-J-by-log-n}, originally observed to first order by \cite{heath1996new}.
Without such a recursion,
we might suffer for small moduli $n$,
as in \cite{hooley1986HasseWeil}*{\S9}'s analysis of ``junior arcs''
(repaired in \cite{hooley_greaves_harman_huxley_1997} for some purposes,
by a clever averaging argument).

\begin{lemma}
\label{LEM:uniform-recursive-derivative-structure}
Let $(r, \bm{v}) = (n/Y, X\bm{c}/n)$.
Let $q_1(t)\defeq t\cdot q(t)$ and let $\phi_1(\bm{x})\defeq \map{div}(\phi_1(\bm{x})\bm{x})=\bm{x}\cdot \grad{\phi_1}(\bm{x}) = \sum_{1\le j\le m} x_j\cdot \partial_{x_j}{\phi_1}$.
Let $\phi_{2,j}(\bm{x})\defeq x_j\cdot \phi(\bm{x})$.
Then
\begin{align}
\partial_{\log{n}} \mscr{J}_{r,\bm{v}}(q,\phi)
&= (m-3) \mscr{J}_{r,\bm{v}}(q,\phi)
- 2\mscr{J}_{r,\bm{v}}(q_1,\phi)
+ \mscr{J}_{r,\bm{v}}(q,\phi_1),
\label{EQN:recursive-differentiate-generalized-J-by-log-n} \\
\partial_{c_j} \mscr{J}_{r,\bm{v}}(q,\phi)
&= (-2\pi i X/n) \mscr{J}_{r,\bm{v}}(q,\phi_{2,j}).
\label{EQN:recursive-differentiate-generalized-J-by-c_j}
\end{align}
\end{lemma}

\begin{proof}
The formula \eqref{EQN:recursive-differentiate-generalized-J-by-c_j} immediately follows upon differentiating \eqref{EQN:define-generalized-integral-J_r,v,q,phi} by $c_j$.
It is possible to prove \eqref{EQN:recursive-differentiate-generalized-J-by-log-n} by a clever integration by parts (cf.~\cite{heath1996new}*{p.~182, proof of Lemma~14}).
We give a slightly shorter argument.
Write $u = t/r^3$ and $\bm{x} = r\bm{y}$ in \eqref{EQN:define-generalized-integral-J_r,v,q,phi} to get
\begin{equation}
\label{EQN:convenient-formula-for-J_r,v-with-phase-independent-of-n}
\mscr{J}_{r,\bm{v}}(q,\phi)
= r^{m-3} \int_{\RR} dt\,q(t/r^2)
\int_{\RR^m} d\bm{y}\,\phi(r\bm{y})e(tF(\bm{y})-r\bm{v}\cdot\bm{y}).
\end{equation}
Differentiating both sides of \eqref{EQN:convenient-formula-for-J_r,v-with-phase-independent-of-n} by $\log{n}$, using $\partial_{\log{n}}{r} = r$ and the fact that $r\bm{v} = X\bm{c}/Y$ is independent of $\log{n}$,
we get \eqref{EQN:recursive-differentiate-generalized-J-by-log-n},
by \eqref{EQN:convenient-formula-for-J_r,v-with-phase-independent-of-n}
applied to each of $(q,\phi)$, $(q_1,\phi)$, $(q,\phi_1)$.
\end{proof}

\begin{proof}
[Proof of Proposition~\ref{PROP:uniform-discriminating-integral-bound-with-derivatives}, assuming Proposition~\ref{PROP:uniform-integral-estimate-for-clean-weights}]
Recall that by \eqref{EQN:rewrite-Heath-Brown-by-Fourier}, we have $J_{\bm{c},X}(n) = \mscr{J}_{r,\bm{v}}(\mf{q}_r,w_0)$, where $(r, \bm{v}) = (n/Y, X\bm{c}/n)$.
By Lemmas~\ref{LEM:bounds-on-q_r(t)-and-its-t-and-r-derivatives} and~\ref{LEM:uniform-recursive-derivative-structure} (first using the chain rule, \eqref{EQN:recursive-differentiate-generalized-J-by-log-n}, and Lemma~\ref{LEM:bounds-on-q_r(t)-and-its-t-and-r-derivatives} when differentiating by $\log{n}$, and then using \eqref{EQN:recursive-differentiate-generalized-J-by-c_j} when differentiating $\bm{c}$),
we may thus write $\partial_{\bm{c}}^{\bm{\alpha}} \partial_{\log{n}}^j J_{\bm{c},X}(n)$ as a finite linear combination of integrals $\mscr{J}_{r,\bm{v}}(q,\phi)$,
with coefficients $O_{m,j,\bm{\alpha}}((X/n)^{\abs{\bm{\alpha}}})$,
running over a set of at most $4^j$ pairs $(q,\phi)$.
(For example, for $j=1$ and $\bm{\alpha}=\bm{0}$, we would use the chain rule and \eqref{EQN:recursive-differentiate-generalized-J-by-log-n} to write $\partial_{\log{n}} \mscr{J}_{r,\bm{v}}(\mf{q}_r,w_0)$ as
\begin{equation*}
\mscr{J}_{r,\bm{v}}(\partial_{\log{r}}{\mf{q}_r},w_0)
+ (m-3) \mscr{J}_{r,\bm{v}}(\mf{q}_r,w_0)
- 2\mscr{J}_{r,\bm{v}}(\mf{q}_{r,1},w_0)
+ \mscr{J}_{r,\bm{v}}(q,w_{0,1}),
\end{equation*}
where $\mf{q}_{r,1}(t) = t\cdot \mf{q}_r(t)$ and $w_{0,1}(\bm{x})\defeq \map{div}(w_0(\bm{x})\bm{x})$.)
Proposition~\ref{PROP:uniform-discriminating-integral-bound-with-derivatives} then immediately follows from Proposition~\ref{PROP:uniform-integral-estimate-for-clean-weights} (applied to each individual $\mscr{J}_{r,\bm{v}}(q,\phi)$).
\end{proof}

To prove Proposition~\ref{PROP:uniform-integral-estimate-for-clean-weights}, we need the following lemma; the basic principle is familiar (see e.g.~\cite{hormander1990analysis}*{Theorem~7.7.1}) but the treatment of \cite{heath1996new} is ideal for us.

\begin{lemma}
[Non-stationary phase]
\label{LEM:uniform-non-stationary-phase}
Let $\lambda, A\in \RR_{>0}$ and $d, k\in \ZZ_{\ge 1}$.
Suppose $a\in C^\infty_c(\RR^d)\otimes \CC$ is supported on $[-A,A]^d$, with $\sum_{\abs{\bm{\alpha}}\le k} \abs{\partial_{\bm{x}}^{\bm{\alpha}}{a}(\bm{x})} \le A$ for all $\bm{x}\in \RR^d$.
Suppose $f\maps \RR^d\to \RR$ is smooth, with $\norm{\grad{f}(\bm{x})}\ge \lambda/A$ and $\sum_{2\le \abs{\bm{\alpha}}\le k+1} \abs{\partial_{\bm{x}}^{\bm{\alpha}}{f}(\bm{x})} \le A\lambda$ for all $\bm{x}\in \Supp{a}$.
Then
\begin{equation*}
    \int_{\bm{x}\in \RR^d} d\bm{x}\, a(\bm{x}) e(f(\bm{x})) \ll_{d,A,k} \lambda^{-k}.
\end{equation*}
\end{lemma}

\begin{proof}
See \cite{heath1996new}*{Lemma~10 and its proof (repeated integration by parts)}; see \cite{heath1996new}*{\S2} for the definition of $\mscr{C}(S)$ (which allows for the required uniformity over weight functions).
Note that we do not require any explicit upper bound on $\sup_{\bm{x}\in \Supp{a}}{\abs{f(\bm{x})}}$, or any control on the shape of the compact set $\Supp{a}\belongs [-A,A]^d$.
\end{proof}

\begin{proof}
[Proof of Proposition~\ref{PROP:uniform-integral-estimate-for-clean-weights}]
Certainly $\bm{0}\notin\Supp{w}$ by \eqref{COND:clean-weight-condition-in-general}.
Since $V$ is smooth, we thus have $\norm{\grad{F}(\bm{x})} > 0$ for all $\bm{x}\in \Supp{w}$.
Since $\Supp{w}$ is compact, there thus exists $A_9 = A_9(F,w) \ge 2$ such that for all $\bm{x}\in \Supp{w}$, we have
$\norm{\bm{x}}, \norm{\grad{F}(\bm{x})} \in [A_9^{-1}, A_9]$.
Let $\mscr{W}\defeq \set{a\bm{x}: (a, \bm{x})\in [1/3A_9, 3A_9] \times \Supp{w}}$ be the union of all dilates of $\Supp{w}$ by a scale factor $a\in [1/3A_9, 3A_9]$.
The set $\mscr{W}$ is compact, being the continuous image of a product of compact sets.
Also,
since the right-hand side of \eqref{COND:clean-weight-condition-in-general} is invariant under scaling, we have
\begin{equation}
\label{COND:extended-clean-Hessian-free-condition-to-scaled-region-W}
\mscr{W}\belongs \set{\bm{w}\in \RR^m: \det(\map{Hess}{F}(\bm{w}))\neq 0}.
\end{equation}
Let $A_{10}=A_{10}(F,w)\ge 2$ be a constant such that for all $\bm{w}\in \mscr{W}$, we have
\begin{equation}
\label{INEQ:A_10-bounds-on-w,gradFw}
    A_{10}^{-1} \le \norm{\bm{w}} \le A_{10},
    \quad A_{10}^{-1} \le \norm{\grad{F}(\bm{w})} \le A_{10}.
\end{equation}

Our plan is to first consider the $\bm{x}$-integral in $\mscr{J}_{r,\bm{v}}(q,\phi)$ (see \eqref{EQN:define-generalized-integral-J_r,v,q,phi}) for a single value of $u$ at a time, and then integrate over $u$.
Given $u$, $\bm{v}$, let
$\psi_0(\bm{x})\defeq uF(\bm{x})-\bm{v}\cdot\bm{x}$ and
$\mcal{J}_{u,\bm{v}}(\phi)\defeq \int_{\bm{x}\in \RR^m} d\bm{x}\,\phi(\bm{x})e(\psi_0(\bm{x}))$.
Taking absolute values in $\mcal{J}$ gives the trivial bound
\begin{equation}
\label{INEQ:trivial-J-u,v-bound}
    \mcal{J}_{u,\bm{v}}(\phi)\ll_\phi 1.
\end{equation}
On the other hand, $\grad{\psi_0}(\bm{x}) = u\grad{F}(\bm{x})-\bm{v}$.
In particular, if $\abs{u}\ge 2A_9\norm{\bm{v}}$ or $\abs{u}\le \norm{\bm{v}}/2A_9$,
then we have $\norm{\grad{\psi_0}(\bm{x})}\ge \max(\abs{u/2A_9}, \norm{\bm{v}}/2)$ and $\partial_{\bm{x}}^{\bm{\alpha}}{\psi_0}(\bm{x}) \ll_{F,w,\bm{\alpha}} \abs{u}+\norm{\bm{v}}$ for all $\bm{x}\in \Supp{w}$, and thus $\mcal{J}_{u,\bm{v}}(\phi)\ll_{\phi,k} (\abs{u}+\norm{\bm{v}})^{-k}$ by Lemma~\ref{LEM:uniform-non-stationary-phase} (provided $\abs{u}+\norm{\bm{v}}>0$).
This, together with \eqref{INEQ:trivial-J-u,v-bound} and the bound $q(ru) \ll_q 1$, gives
\begin{align}
\int_{\abs{u}\ge 2A_9\norm{\bm{v}}} du\, \abs{q(ru) \cdot \mcal{J}_{u,\bm{v}}(\phi)}
&\ll_{q,\phi,b} \int_{\abs{u}\ge 2A_9\norm{\bm{v}}} \frac{du}{(1+\abs{u})^{b+1}}
\ll_b (1+\norm{\bm{v}})^{-b},
\label{INEQ:integrated-absolute-large-u-bound} \\
\int_{\abs{u}\le \norm{\bm{v}}/2A_9} du\, \abs{q(ru) \cdot \mcal{J}_{u,\bm{v}}(\phi)}
&\ll_{q,\phi,b} \int_{\abs{u}\le \norm{\bm{v}}/2A_9} \frac{du}{(1+\norm{\bm{v}})^{b+1}}
\ll_b (1+\norm{\bm{v}})^{-b},
\label{INEQ:integrated-absolute-small-u-bound}
\end{align}
for all integers $b\ge 1$.
Fix $w_1\in C^\infty_c(\RR)$ with
\begin{equation*}
w_1\vert_{[-2A_9, 2A_9]\setminus [-1/2A_9, 1/2A_9]} = 1,
\quad \Supp{w_1}\belongs [-3A_9, 3A_9]\setminus [-1/3A_9, 1/3A_9].
\end{equation*}
Then by \eqref{INEQ:integrated-absolute-large-u-bound}, \eqref{INEQ:integrated-absolute-small-u-bound}, and the triangle inequality, we have (for integers $b\ge 1$)
\begin{equation}
\label{INEQ:reduce-to-u-roughly-proportional-to-v}
\mscr{J}_{r,\bm{v}}(q,\phi) - \int_{u\in \RR} du\, q(ru) w_1(u/\norm{\bm{v}}) \mcal{J}_{u,\bm{v}}(\phi)
\ll_{q,\phi,b} (1+\norm{\bm{v}})^{-b}.
\end{equation}

It remains to handle $\abs{u}/\norm{\bm{v}}\in [1/3A_9, 3A_9]$.
Suppose first that $\norm{\bm{v}}\leq 1$.
Plugging \eqref{INEQ:trivial-J-u,v-bound} and the bound $q(ru) \ll_q 1$ into \eqref{INEQ:reduce-to-u-roughly-proportional-to-v}, we get
$\mscr{J}_{r,\bm{v}}(q,\phi) \ll_{q,\phi} (1+\norm{\bm{v}})^{-1} + \int_{\abs{u}\le 3A_9} du\,1 \ll 1$.
This bound on $\mscr{J}_{r,\bm{v}}(q,\phi)$ fits in Proposition~\ref{PROP:uniform-integral-estimate-for-clean-weights}, since $\norm{r\bm{v}}\le A_0$ and $\norm{\Delta(\bm{v}/M)\bm{v}}\ll_{F} 1$.

For the rest of the proof,
suppose $\norm{\bm{v}}\ge 1$,
let $\tilde{u}\defeq u/\norm{\bm{v}}$ and $\tilde{\bm{v}}\defeq \bm{v}/\norm{\bm{v}}$,
let $\bm{z}\defeq \abs{\tilde{u}}^{1/2} \bm{x}$,
let $\sgn(u)\defeq u/\abs{u}\in \set{-1,1}$, and let
\begin{equation*}
\psi_1(\bm{z})
\defeq \abs{\tilde{u}}^{1/2} \psi_0(\bm{x})/\norm{\bm{v}}
= \abs{\tilde{u}}^{1/2} \psi_0(\bm{z}/\abs{\tilde{u}}^{1/2})/\norm{\bm{v}}
= \sgn(u) F(\bm{z}) - \tilde{\bm{v}}\cdot \bm{z}.
\end{equation*}
(The normalization by $\abs{\tilde{u}}^{1/2}$ makes $\psi_1$ nearly independent of $u$;
this is crucial later.)
Note that $\grad{\psi_1}(\bm{z}) = \sgn(u)\grad{F}(\bm{z}) - \tilde{\bm{v}} = \tilde{u}\grad{F}(\bm{x}) - \tilde{\bm{v}}$.
Now consider an individual $\tilde{u}\in \RR$ with
\begin{equation}
\label{INEQ:key-range-of-u-in-J-analysis}
1/3A_9\le \abs{\tilde{u}}\le 3A_9.
\end{equation}
For all $\bm{x}\in \Supp{w}$, we have $\bm{z}\in \mscr{W}$, so $\norm{\bm{z}}\in [1/A_{10}, A_{10}]$ and $\norm{\grad{\psi_1}(\bm{z})}\le \norm{\grad{F}(\bm{z})} + 1 \le 2A_{10}$.
Furthermore, the condition \eqref{COND:extended-clean-Hessian-free-condition-to-scaled-region-W} implies
\begin{equation}
\label{INEQ:Hessian-free-volume-bound-on-psi_1-singular-level-sets}
\vol{\set{\bm{z}\in \mscr{W}: \norm{\grad{\psi_1}(\bm{z})}\le m\lambda}} \ll_{F,w} \lambda^m,
\end{equation}
uniformly over $\lambda\in \RR_{>0}$, by calculus (see \cite{heath1996new}*{Lemma~21 and its proof}).
We will split $\mcal{J}_{u,\bm{v}}(\phi)$ into pieces according to the size of $\norm{\grad{\psi_1}(\bm{z})}$; cf.~\cite{huang2020density}*{$\mscr{K}_i$ on p.~2061}.
To avoid the need for explicit stationary phase expansions (which can be messy after the leading term),
we will also use a subdivision process inspired by \cite{heath1996new}*{\S8}.

Recall the weights $\nu_0$, $\nu_1$ from \S\ref{SUBSEC:conventions}.
Let $\ol{\nu}_0\defeq 1-\nu_0$.
For reals $\lambda>0$, let
\begin{align}
\mcal{J}_{0,u,\bm{v}}
&\defeq \int_{\RR^m} \frac{d\bm{z}}{\abs{\tilde{u}}^{m/2}}\,
\nu_0{\left(\frac{\grad{\psi_1}(\bm{z})}{\norm{\bm{v}}^{-1/2}}\right)}
\phi(\bm{x}) e{\left(\frac{\norm{\bm{v}} \psi_1(\bm{z})}{\abs{\tilde{u}}^{1/2}}\right)},
\label{EQN:define-grad-localized-x-integral-J_0} \\
\mcal{J}_{1,u,\bm{v},\lambda}
&\defeq \int_{\RR^m} \frac{d\bm{z}}{\abs{\tilde{u}}^{m/2}}\,
\ol{\nu}_0{\left(\frac{\grad{\psi_1}(\bm{z})}{\norm{\bm{v}}^{-1/2}}\right)}
\nu_1{\left(\frac{\grad{\psi_1}(\bm{z})}{\lambda}\right)}
\phi(\bm{x}) e{\left(\frac{\norm{\bm{v}} \psi_1(\bm{z})}{\abs{\tilde{u}}^{1/2}}\right)},
\label{EQN:define-grad-localized-x-integral-J_1}
\end{align}
where $\phi(\bm{x}) = \phi(\bm{z}/\abs{\tilde{u}}^{1/2})$.
By the definitions of $\nu_0$, $\nu_1$,
we have $\mcal{J}_{1,u,\bm{v},\lambda} = 0$ unless there exists $\bm{z}\in \mscr{W}$ (corresponding to $\bm{x}\in \Supp{w}$) satisfying $\norm{\grad{\psi_1}(\bm{z})}>\norm{\bm{v}}^{-1/2}/2$ and $\lambda\le \norm{\grad{\psi_1}(\bm{z})}\le m\lambda$.
Therefore, $\mcal{J}_{1,u,\bm{v},\lambda} = 0$ unless $\norm{\bm{v}}^{-1/2}/2m < \lambda\le 2A_{10}$.
Since $\int_{\lambda>0} d^\times{\lambda}\,\nu_1(\bm{t}/\lambda) = 1$ for all $\bm{t}\in \RR^m\setminus \set{\bm{0}}$, we obtain (under \eqref{INEQ:key-range-of-u-in-J-analysis}) the decomposition
\begin{equation}
\label{INEQ:decomposition-of-x-integral-into-nice-pieces}
\mcal{J}_{u,\bm{v}}(\phi) - \mcal{J}_{0,u,\bm{v}}
= \int_{\norm{\bm{v}}^{-1/2}/2m}^{2A_{10}} d^\times{\lambda}\,
\mcal{J}_{1,u,\bm{v},\lambda}.
\end{equation}
Fix functions $\omega_1,\dots,\omega_m\in C^\infty_c(\RR^m)$, all supported on $\Supp{\nu_1}$, such that $\omega_1+\dots+\omega_m = \nu_1$ and we have $1/2\le \abs{t_j}\le 2m$ for all $\bm{t}\in \Supp{\omega_j}$.
Let $\mcal{J}_{1,j,u,\bm{v},\lambda}$ denote $\mcal{J}_{1,u,\bm{v},\lambda}$ with $\omega_j$ in place of $\nu_1$; clearly $\mcal{J}_{1,u,\bm{v},\lambda} = \sum_{1\le j\le m} \mcal{J}_{1,j,u,\bm{v},\lambda}$.

We proceed based on the following rough idea:
if $\norm{\grad{\psi_1(\bm{z})}}\gg \norm{\bm{v}}^{-1/2}$ with a large constant, integration by parts over $\bm{z}$ is useful;
and if $\norm{\grad{\psi_1(\bm{z})}}\ll \abs{\Delta(\tilde{\bm{v}})}^{1/2}$ with a small constant, integration by parts over $u$ is useful (if $\norm{\Delta(\tilde{\bm{v}})\bm{v}}\gg 1$ with any constant).

Let $\lambda\in [4\norm{\bm{v}}^{-1/2}, 2A_{10}]$.
Since $\lambda\ge 4\norm{\bm{v}}^{-1/2}$, the supports of $\nu_0(\grad{\psi_1}(\bm{z})/\norm{\bm{v}}^{-1/2})$ and $\nu_1(\grad{\psi_1}(\bm{z})/\lambda)$ are disjoint, so $\ol{\nu}_0(\grad{\psi_1}(\bm{z})/\norm{\bm{v}}^{-1/2}) \omega_j(\grad{\psi_1}(\bm{z})/\lambda) = \omega_j(\grad{\psi_1}(\bm{z})/\lambda)$.
Thus
\begin{equation}
\label{EQN:J_1,j,u,v,phi,y,lambda-change-of-variables-large-lambda-conceptual}
\abs{\tilde{u}}^{m/2} \mcal{J}_{1,j,u,\bm{v},\lambda}
= \int_{\RR^m} d\bm{z}\, w_{2,j,0}(\bm{z}) e(\norm{\bm{v}} \psi_1(\bm{z})/\abs{\tilde{u}}^{1/2}),
\end{equation}
where $w_{2,j,0}=w_{2,j,0,u,\bm{v},\lambda}(\bm{z})\defeq
\omega_j(\grad{\psi_1}(\bm{z})/\lambda) \phi(\bm{z}/\abs{\tilde{u}}^{1/2})$.

Since $\Supp{\omega_j}\belongs \set{\bm{t}\in \RR^m: 1/2\le \abs{t_j}\le 2m}\cap \Supp{\nu_1}$, we have
\begin{equation}
\label{INEQ:good-bounds-on-first-derivatives-on-Supp-w_2,j,0}
\norm{\partial_{z_j}{\psi_1}(\bm{z})} \ge \lambda/2,
\quad \norm{\grad{\psi_1}(\bm{z})}\le m\lambda,
\end{equation}
for all $\bm{z}\in \Supp{w_{2,j,0}}$.
For $k\ge 1$, recursively define $w_{2,j,k}=w_{2,j,k,u,\bm{v},\lambda}(\bm{z})$ (a sequence of smooth functions $\RR^m\to \CC$ supported on $\Supp{w_{2,j,0}}$) by setting
\begin{equation}
\label{EQN:recursively-define-normalized-IBP-weights-w,2,j,k}
w_{2,j,k}\defeq -(\partial_{z_j}{\psi_1})^k
\cdot \partial_{z_j}(w_{2,j,k-1}/(\partial_{z_j}{\psi_1})^k)
\end{equation}
Integrating by parts $k$ times in \eqref{EQN:J_1,j,u,v,phi,y,lambda-change-of-variables-large-lambda-conceptual}, we get (cf.~\cite{heath1996new}*{(5.4), from the proof of Lemma~10})
\begin{equation}
\label{INEQ:J_1,j,u,v,phi,y,lambda-integrated-by-parts-k-times}
\abs{\tilde{u}}^{m/2} \mcal{J}_{1,j,u,\bm{v},\lambda}
= \int_{\RR^m} d\bm{z}\, \frac{w_{2,j,k,u,\bm{v},\lambda}(\bm{z})}{(\norm{\bm{v}} \partial_{z_j}{\psi_1}(\bm{z})/\abs{\tilde{u}}^{1/2})^k}
\cdot e(\norm{\bm{v}} \psi_1(\bm{z})/\abs{\tilde{u}}^{1/2}).
\end{equation}

We claim that for each $k\ge 0$, there exists a smooth
function $w_{3,j,k}\maps \RR^{m(k+2) + 3} \to \CC$ of the form
$w_{3,j,k}=w_{3,j,k}(\bm{a}_{1,0},\dots,\bm{a}_{1,k},\bm{a}_{2,0},a_{2,1},a_3,a_4)$,
defined in terms of $F$, $j$, $k$ (independently of $u$, $\bm{v}$, $\lambda$), such that for all $u$, $\bm{v}$, $\lambda$ currently under consideration (namely, satisfying $\norm{\bm{v}}\ge 1$, \eqref{INEQ:key-range-of-u-in-J-analysis}, and $\lambda\in [4\norm{\bm{v}}^{-1/2}, 2A_{10}]$), and for all $\bm{z}\in \Supp{w_{2,j,0}}$, we have
\begin{equation}
\label{EQN:key-uniformity-structure-of-w_2,j,k,u,v,lambda}
    \lambda^k \cdot w_{2,j,k}(\bm{z})
    = w_{3,j,k}{\left(\frac{\grad{\psi_1}(\bm{z})}{\lambda},
    \grad{\partial_{z_j}^1{\psi_1}}(\bm{z}),\dots,\grad{\partial_{z_j}^k{\psi_1}}(\bm{z}),
    \frac{\bm{z}}{\abs{\tilde{u}}^{1/2}}, \frac{1}{\abs{\tilde{u}}^{1/2}},
    \frac{\lambda}{\partial_{z_j}{\psi_1}}, \lambda\right)}.
\end{equation}
This claim can be easily proven by induction on $k\ge 0$, where for $k=0$ we take
\begin{equation*}
    w_{3,j,0}(\bm{a}_{1,0},\bm{a}_{2,0},a_{2,1},a_3,a_4)\defeq \omega_j(\bm{a}_{1,0}) \phi(\bm{a}_{2,0}),
\end{equation*}
and for $k\ge 1$ we use the product and chain rules in \eqref{EQN:recursively-define-normalized-IBP-weights-w,2,j,k} (to compute $\lambda^k w_{2,j,k} = \lambda\cdot \lambda^{k-1}w_{2,j,k}$ in terms of $\lambda^{k-1} w_{2,j,k-1}$ and its $z_j$-derivatives), and observe that (for $l\ge 1$)
\begin{gather*}
\lambda\cdot \partial_{z_j}{\left(\frac{\psi_1}{\lambda}\right)}
= \partial_{z_j}^1{\psi_1},
\quad
\lambda\cdot \partial_{z_j}(\partial_{z_j}^l) = \lambda\cdot \partial_{z_j}^{l+1},
\quad
\lambda\cdot \partial_{z_j}(\bm{z}) = \lambda, \\
\lambda\cdot \partial_{z_j}{\left(\frac{\lambda}{\partial_{z_j}{\psi_1}}\right)}
= -(\partial_{z_j}^2{\psi_1})\cdot \left(\frac{\lambda}{\partial_{z_j}{\psi_1}}\right)^2,
\quad
(\partial_{z_j}{\psi_1})^k\cdot \partial_{z_j}{\left(\frac{\lambda}{(\partial_{z_j}{\psi_1})^k}\right)}
= (\partial_{z_j}^2{\psi_1})\cdot \frac{-k\lambda}{\partial_{z_j}{\psi_1}}.
\end{gather*}
(This proof uses the smoothness of $\psi_1$, and the nonvanishing of $\partial_{z_j}{\psi_1}$ on $\Supp{w_{2,j,0}}$.)

By \eqref{INEQ:A_10-bounds-on-w,gradFw}, we have $\partial^{\bm{\alpha}}{\psi_1}(\bm{z}) \ll_{\bm{\alpha}} 1$ whenever $\abs{\bm{\alpha}}\ge 2$ and $\bm{z}\in \mscr{W}$.
Inserting this and the bounds \eqref{INEQ:good-bounds-on-first-derivatives-on-Supp-w_2,j,0}, \eqref{INEQ:key-range-of-u-in-J-analysis}, \eqref{INEQ:A_10-bounds-on-w,gradFw}, and $\lambda\ll 1$ into \eqref{EQN:key-uniformity-structure-of-w_2,j,k,u,v,lambda}, we find that $\lambda^k w_{2,j,k}(\bm{z}) \ll_{w_{3,j,k}} 1$ for all $\bm{z}\in \RR^m$.
Thus \eqref{INEQ:J_1,j,u,v,phi,y,lambda-integrated-by-parts-k-times} and \eqref{INEQ:good-bounds-on-first-derivatives-on-Supp-w_2,j,0} immediately give (for $\lambda\in [4\norm{\bm{v}}^{-1/2}, 2A_{10}]$)
\begin{equation}
\label{INEQ:basic-x-IBP-J_1,j,u,v,phi,y,lambda-bound}
\abs{\tilde{u}}^{m/2} \mcal{J}_{1,j,u,\bm{v},\lambda}
\ll_k \int_{\bm{z}\in \mscr{W}} d\bm{z}\,
\frac{\bm{1}_{\norm{\grad{\psi_1}(\bm{z})}\le m\lambda}}
{(\lambda^2 \norm{\bm{v}}/\abs{\tilde{u}}^{1/2})^k},
\end{equation}
under \eqref{INEQ:key-range-of-u-in-J-analysis}.
Our derivation of \eqref{INEQ:basic-x-IBP-J_1,j,u,v,phi,y,lambda-bound} has essentially followed the proof of Lemma~\ref{LEM:uniform-non-stationary-phase} (with a small but important twist: we use $\Supp{w_{2,j,k}}\belongs \Supp{w_{2,j,0}}$ to get the factor $\bm{1}_{\norm{\grad{\psi_1}(\bm{z})}\le m\lambda}$), but in some key ranges we need to go further.
We need to decide when to integrate by parts over $\tilde{u}\in \RR$; this will be informed by the next two paragraphs.

By \eqref{COND:extended-clean-Hessian-free-condition-to-scaled-region-W} and the inverse function theorem, we know that for each $\bm{z}\in \mscr{W}$, there exists an open neighborhood $U\belongs \RR^m$ of $\bm{z}$ such that the gradient map $\grad{F}\maps \RR^m\to \RR^m$ maps $U$ diffeomorphically onto $\grad{F}(U)$.
But $\grad{F}(\mscr{W})$ is compact (since $\mscr{W}$ is compact), so there exists a constant $\eta_6=\eta_6(F,w)>0$ such that for any $\bm{z}\in \mscr{W}$ and $\bm{b}\in [-1,1]^m$ with $\norm{\grad{F}(\bm{z})-\bm{b}}\le \eta_6$, there exists $\bm{s}\in (\grad{F})^{-1}(\bm{b})$ with $\norm{\bm{z}-\bm{s}} \le \eta_6^{-1} \norm{\grad{F}(\bm{z})-\bm{b}}$.

We apply this as follows.
Let $\bm{b} = \sgn(u) \tilde{\bm{v}}$ and $\bm{z}\in \mscr{W}$, and suppose $\norm{\grad{\psi_1(\bm{z})}}\le \eta_6$.
Then
there exists $\bm{s}\in \RR^m$ with $\grad{\psi_1}(\bm{s})=\bm{0}$ and $\norm{\bm{z}-\bm{s}}\le \eta_6^{-1}\norm{\grad{\psi_1(\bm{z})}}$.
Since $\grad{\psi_1}(\bm{s})=\bm{0}$,
Taylor expansion of $\psi_1$ at $\bm{s}$ gives $\psi_1(\bm{z}) = \psi_1(\bm{s}) + O_{F,w}(\norm{\bm{z}-\bm{s}}^2)$.
Furthermore, $\grad{\psi_1}(\bm{s})=\bm{0}$ implies $\bm{b}\cdot \bm{s} = \grad{F}(\bm{s})\cdot \bm{s} = 3F(\bm{s})$, so $\abs{\psi_1(\bm{s})} = \abs{2F(\bm{s})}$.
But $\Delta(\bm{b})=\Delta(\grad{F}(\bm{s})) \ll_{F,w} \abs{F(\bm{s})}$ by \eqref{INEQ:bound-Disc-grad-composite-in-terms-of-F}.
Combining the above, we get
\begin{equation}
\label{INEQ:lower-bound-on-phase-psi_1-given-eta_6-bound-on-gradient}
\abs{\psi_1(\bm{z})}\ge \eta_7 \abs{\Delta(\tilde{\bm{v}})} - \eta_7^{-1} \norm{\grad{\psi_1(\bm{z})}}^2
\quad\textnormal{(for some $\eta_7=\eta_7(F,w)>0$)}.
\end{equation}
Let $A_{11}\defeq 1+\sup_{\bm{a}\in [-1,1]^m}{\abs{\Delta(\bm{a})}}$, and let $\eta_8\defeq m^{-1} \min(\eta_6/A_{11}^{1/2}, \eta_7/2)$.

Our remaining analysis breaks into two cases: $\norm{\Delta(\tilde{\bm{v}})\bm{v}}\ge (4/\eta_8)^2$ and $\norm{\Delta(\tilde{\bm{v}})\bm{v}} < (4/\eta_8)^2$.

Suppose first that $\norm{\Delta(\tilde{\bm{v}})\bm{v}} < (4/\eta_8)^2$.
Consider a $\tilde{u}\in \RR$ satisfying \eqref{INEQ:key-range-of-u-in-J-analysis}.
By \eqref{INEQ:Hessian-free-volume-bound-on-psi_1-singular-level-sets}, the right-hand side of \eqref{INEQ:basic-x-IBP-J_1,j,u,v,phi,y,lambda-bound} is $\ll_k \lambda^m (\lambda^2\norm{\bm{v}})^{-k} = \norm{\bm{v}}^{-m/2} (\lambda\norm{\bm{v}}^{1/2})^{m-2k}$.
Choosing $k = \ceil{m/2+1}$ and inserting \eqref{INEQ:basic-x-IBP-J_1,j,u,v,phi,y,lambda-bound} into \eqref{INEQ:decomposition-of-x-integral-into-nice-pieces}, we then get
\begin{equation*}
\mcal{J}_{u,\bm{v}}(\phi) - \mcal{J}_{0,u,\bm{v}}
- \int_{\norm{\bm{v}}^{-1/2}/2m}^{4\norm{\bm{v}}^{-1/2}} d^\times{\lambda}\, \mcal{J}_{1,u,\bm{v},\lambda}
\ll_m \norm{\bm{v}}^{-m/2},
\end{equation*}
since $\int_{4\norm{\bm{v}}^{-1/2}}^{\infty} d^\times{\lambda}\, (\lambda\norm{\bm{v}}^{1/2})^{m-2k} = \int_{4}^{\infty} d^\times{a}\, a^{m-2k} \ll 1$ (via the substitution $a = \lambda\norm{\bm{v}}^{1/2}$).
But if we simply take absolute values in $\mcal{J}_{0,u,\bm{v}}$ (see \eqref{EQN:define-grad-localized-x-integral-J_0}) and in $\mcal{J}_{1,u,\bm{v},\lambda}$ (see \eqref{EQN:define-grad-localized-x-integral-J_1}), and then apply \eqref{INEQ:Hessian-free-volume-bound-on-psi_1-singular-level-sets}, we get $\mcal{J}_{0,u,\bm{v}} \ll \norm{\bm{v}}^{-m/2}$ and $\mcal{J}_{1,u,\bm{v},\lambda} \ll \norm{\bm{v}}^{-m/2}$ for $\lambda\in [\norm{\bm{v}}^{-1/2}/2m, 4\norm{\bm{v}}^{-1/2}]$.
Integrating over $\lambda$ in the previous display, we conclude that $\mcal{J}_{u,\bm{v}}(\phi) \ll_m \norm{\bm{v}}^{-m/2}$ (under \eqref{INEQ:key-range-of-u-in-J-analysis}).
Hence by \eqref{INEQ:reduce-to-u-roughly-proportional-to-v} (after writing $u=\norm{\bm{v}}\tilde{u}$) we have
\begin{equation*}
    \mscr{J}_{r,\bm{v}}(q,\phi) \ll_b (1+\norm{\bm{v}})^{-b} + \norm{\bm{v}} \int_{1/3A_9\le \abs{\tilde{u}}\le 3A_9} d\tilde{u}\, \abs{q(r\norm{\bm{v}}\tilde{u})} \cdot \norm{\bm{v}}^{-m/2}.
\end{equation*}
But $q(r\norm{\bm{v}}\tilde{u}) \ll_b (1+r\norm{\bm{v}})^{-b}$ (under \eqref{INEQ:key-range-of-u-in-J-analysis}).
Thus $\mscr{J}_{r,\bm{v}}(q,\phi) \ll_b \norm{\bm{v}}^{1-m/2} (1+r\norm{\bm{v}})^{-b}$, which suffices for Proposition~\ref{PROP:uniform-integral-estimate-for-clean-weights} (since $\norm{\bm{v}}\ge 1$ and $\norm{\Delta(\bm{v}/M)\bm{v}}\le \norm{\Delta(\tilde{\bm{v}})\bm{v}} \ll 1$).

For the rest of the proof, assume $\norm{\Delta(\tilde{\bm{v}})\bm{v}} \ge (4/\eta_8)^2$, so that $4\norm{\bm{v}}^{-1/2}\le \eta_8 \abs{\Delta(\tilde{\bm{v}})}^{1/2}$.
We will integrate by parts over $\tilde{u}$ in the integrals
\begin{equation*}
\mscr{J}_0\defeq \int_{\RR} d\tilde{u}\, q(r\norm{\bm{v}}\tilde{u}) w_1(\tilde{u}) \mcal{J}_{0,u,\bm{v}},
\quad \mscr{J}_{1,\lambda} \defeq \int_{\RR} d\tilde{u}\, q(r\norm{\bm{v}}\tilde{u}) w_1(\tilde{u}) \mcal{J}_{1,u,\bm{v},\lambda}
\end{equation*}
for $\lambda\le \eta_8 \abs{\Delta(\tilde{\bm{v}})}^{1/2}$.
It is crucial to work in terms of $\bm{z}$ rather than $\bm{x} = \bm{z}/\abs{\tilde{u}}^{1/2}$.
Before proceeding, note that inserting \eqref{INEQ:decomposition-of-x-integral-into-nice-pieces} into \eqref{INEQ:reduce-to-u-roughly-proportional-to-v} (and writing $u=\norm{\bm{v}}\tilde{u}$) gives
\begin{equation}
\label{INEQ:final-decomposition-of-main-integral-J_r,v-over-u,x}
\mscr{J}_{r,\bm{v}}
- \norm{\bm{v}}\cdot \mscr{J}_0 - \norm{\bm{v}}\cdot \int_{\norm{\bm{v}}^{-1/2}/2m}^{2A_{10}} d^\times{\lambda}\, \mscr{J}_{1,\lambda}
\ll_b (1+\norm{\bm{v}})^{-b}.
\end{equation}

Let $\bm{z}\in \mscr{W}$, and suppose $\norm{\grad{\psi_1}(\bm{z})}\le m\eta_8 \abs{\Delta(\tilde{\bm{v}})}^{1/2}$.
By the definition of $\eta_8$, we then have $\norm{\grad{\psi_1}(\bm{z})}\le \eta_6$ and $\norm{\grad{\psi_1}(\bm{z})} \le \eta_7\abs{\Delta(\tilde{\bm{v}})}^{1/2}/2$, so by \eqref{INEQ:lower-bound-on-phase-psi_1-given-eta_6-bound-on-gradient}, we have
\begin{equation}
\label{INEQ:final-simplified-lower-bound-on-phase-psi_1}
    \abs{\psi_1(\bm{z})}\ge 3\eta_7\abs{\Delta(\tilde{\bm{v}})}/4.
\end{equation}
For convenience, let $\mscr{W}(\lambda)\defeq \set{\bm{z}\in \mscr{W}: \norm{\grad{\psi_1}(\bm{z})}\le m\lambda}$ for each $\lambda>0$.

We first bound $\mscr{J}_0$.
For each $\tilde{u}\in \Supp{w_1}$, the condition \eqref{INEQ:key-range-of-u-in-J-analysis} holds, so if we plug in the definition of $\mcal{J}_{0,u,\bm{v}}$
(see \eqref{EQN:define-grad-localized-x-integral-J_0}),
and then switch the order of $u$, $\bm{z}$, we may rewrite $\mscr{J}_0$ as
\begin{equation*}
\int_{\mscr{W}(\norm{\bm{v}}^{-1/2})} d\bm{z}\,
\nu_0{\left(\frac{\grad{\psi_1}(\bm{z})}{\norm{\bm{v}}^{-1/2}}\right)}
\int_{\RR} \frac{d\tilde{u}}{\abs{\tilde{u}}^{m/2}}\,
q(r\norm{\bm{v}}\tilde{u}) w_1(\tilde{u})
\phi{\left(\frac{\bm{z}}{\abs{\tilde{u}}^{1/2}}\right)}
e{\left(\frac{\norm{\bm{v}} \psi_1(\bm{z})}{\abs{\tilde{u}}^{1/2}}\right)}.
\end{equation*}
Here $\psi_1$ is independent of $\abs{\tilde{u}}$ (when $\sgn(u)$ is fixed).
Therefore, for each $\bm{z}\in \mscr{W}(\norm{\bm{v}}^{-1/2})$, the inner integral (over $\tilde{u}$) is $\ll_k (1+r\norm{\bm{v}})^{-k} \norm{\Delta(\tilde{\bm{v}})\bm{v}}^{-k}$ by \eqref{INEQ:final-simplified-lower-bound-on-phase-psi_1} and Lemma~\ref{LEM:uniform-non-stationary-phase}, because $\abs{\tilde{u}}\asymp 1$ for $\tilde{u}\in \Supp{w_1}$ and we have $\partial_{\tilde{u}}^l q(r\norm{\bm{v}}\tilde{u}) = (r\norm{\bm{v}})^l q^{(l)}(r\norm{\bm{v}}\tilde{u}) \ll_{l,k} (1+r\norm{\bm{v}})^{-k}$ for all $l,k\ge 0$ (since $q$ is Schwartz).
Applying this inner integral estimate for each $\bm{z}\in \mscr{W}(\norm{\bm{v}}^{-1/2})$, and then using \eqref{INEQ:Hessian-free-volume-bound-on-psi_1-singular-level-sets}, we get
$\mscr{J}_0 \ll_k \norm{\bm{v}}^{-m/2} (1+r\norm{\bm{v}})^{-k} \norm{\Delta(\tilde{\bm{v}})\bm{v}}^{-k}$.

One can similarly prove $\mscr{J}_{1,\lambda} \ll_k \norm{\bm{v}}^{-m/2} (1+r\norm{\bm{v}})^{-k}$ for $\lambda\in [\norm{\bm{v}}^{-1/2}/2m, 4\norm{\bm{v}}^{-1/2}]$.

Now suppose $4\norm{\bm{v}}^{-1/2}\le \lambda\le \min(\eta_8\abs{\Delta(\tilde{\bm{v}})}^{1/2}, 2A_{10})$.
Let $\mscr{J}_{1,j,\lambda}$ denote $\mscr{J}_{1,\lambda}$ with $\mcal{J}_{1,j,u,\bm{v},\lambda}$ in place of $\mcal{J}_{1,u,\bm{v},\lambda}$.
Then $\mscr{J}_{1,\lambda} = \sum_{1\le j\le m} \mscr{J}_{1,j,\lambda}$.
By \eqref{INEQ:J_1,j,u,v,phi,y,lambda-integrated-by-parts-k-times}, we have (for all integers $l\ge 0$)
\begin{equation*}
\mscr{J}_{1,j,\lambda} = \int_{\bm{z}\in \mscr{W}(\lambda)} d\bm{z}
\int_{\RR} \frac{d\tilde{u}}{\abs{\tilde{u}}^{m/2}}\,
q(r\norm{\bm{v}}\tilde{u}) w_1(\tilde{u})
\frac{w_{2,j,l,u,\bm{v},\lambda}(\bm{z})
e(\norm{\bm{v}} \psi_1(\bm{z})/\abs{\tilde{u}}^{1/2})}
{(\norm{\bm{v}} \partial_{z_j}{\psi_1}(\bm{z})/\abs{\tilde{u}}^{1/2})^l}.
\end{equation*}
Since $\psi_1$ is independent of $\abs{\tilde{u}}$, and we have $\norm{\bm{z}}\ll 1$ and $\abs{\tilde{u}}\asymp 1$ for all $(\bm{z}, \tilde{u})\in \mscr{W}\times \Supp{w_1}$, we find by \eqref{INEQ:good-bounds-on-first-derivatives-on-Supp-w_2,j,0}, \eqref{EQN:key-uniformity-structure-of-w_2,j,k,u,v,lambda}, \eqref{INEQ:final-simplified-lower-bound-on-phase-psi_1}, and Lemma~\ref{LEM:uniform-non-stationary-phase} (applied to the inner integral over $\tilde{u}$, for each $\bm{z}\in \mscr{W}(\lambda)$ for which there exists $u\in \RR$ with $w_1(\tilde{u}) w_{2,j,l,u,\bm{v},\lambda}(\bm{z})\ne 0$) that
\begin{equation*}
\mscr{J}_{1,j,\lambda} \ll_{k,l} \int_{\bm{z}\in \mscr{W}(\lambda)} d\bm{z}\,
\frac{(1+r\norm{\bm{v}})^{-k} \norm{\Delta(\tilde{\bm{v}})\bm{v}}^{-k}}
{(\lambda^2\norm{\bm{v}})^l}
\ll \frac{\lambda^m (1+r\norm{\bm{v}})^{-k} \norm{\Delta(\tilde{\bm{v}})\bm{v}}^{-k}}
{(\lambda^2\norm{\bm{v}})^l}
\end{equation*}
for all integers $k,l\ge 0$, where in the final step we use \eqref{INEQ:Hessian-free-volume-bound-on-psi_1-singular-level-sets}.

Finally, suppose $\eta_8\abs{\Delta(\tilde{\bm{v}})}^{1/2}\le \lambda\le 2A_{10}$.
Then \eqref{INEQ:basic-x-IBP-J_1,j,u,v,phi,y,lambda-bound} and \eqref{INEQ:Hessian-free-volume-bound-on-psi_1-singular-level-sets} directly give (for $k,l\ge 0$)
\begin{equation*}
\mscr{J}_{1,\lambda} \ll_l \int_{\RR} d\tilde{u}\,
\frac{\abs{q(r\norm{\bm{v}}\tilde{u}) w_1(\tilde{u})}\cdot \lambda^m}
{(\lambda^2\norm{\bm{v}})^l}
\ll_k \frac{(1+r\norm{\bm{v}})^{-k} \norm{\bm{v}}^{-m/2}}
{(\lambda\norm{\bm{v}}^{1/2})^{2l-m}}.
\end{equation*}

Inserting our work from the last four paragraphs (ignoring the last one if $\eta_8\abs{\Delta(\tilde{\bm{v}})}^{1/2} > 2A_{10}$) into the left-hand side of \eqref{INEQ:final-decomposition-of-main-integral-J_r,v-over-u,x}, and applying \eqref{INEQ:final-decomposition-of-main-integral-J_r,v-over-u,x} with $b=\ceil{m/2+2k}$, we get the bound
\begin{equation*}
\frac{\mscr{J}_{r,\bm{v}}}{\norm{\bm{v}}} \ll_{k,l} \frac{1
+ \int_{\lambda\ge 4\norm{\bm{v}}^{-1/2}} d^\times{\lambda}\, (\lambda\norm{\bm{v}}^{1/2})^{m-2l}}
{\norm{\bm{v}}^{m/2} (1+r\norm{\bm{v}})^k \norm{\Delta(\tilde{\bm{v}})\bm{v}}^k}
+ \frac{\int_{\lambda\ge \eta_8\abs{\Delta(\tilde{\bm{v}})}^{1/2}} d^\times{\lambda}\, (\lambda\norm{\bm{v}}^{1/2})^{m-2l+2k}}
{\norm{\bm{v}}^{m/2} (1+r\norm{\bm{v}})^k \norm{\Delta(\tilde{\bm{v}})\bm{v}}^k}.
\end{equation*}
Taking $l = \ceil{m/2+k+1}$, and evaluating both integrals over $\lambda$ (the first being $\ll 1$ since $4\norm{\bm{v}}^{-1/2}\cdot \norm{\bm{v}}^{1/2} = 4 \gg 1$, and the second being $\ll 1$ since $\eta_8\norm{\Delta(\tilde{\bm{v}})\bm{v}}^{1/2} \ge 4\gg 1$), we get $\mscr{J}_{r,\bm{v}} \ll_k \norm{\bm{v}}^{1-m/2} (1+r\norm{\bm{v}})^{-k} \norm{\Delta(\tilde{\bm{v}})\bm{v}}^{-k}$, which suffices for Proposition~\ref{PROP:uniform-integral-estimate-for-clean-weights}.
\end{proof}

\section{New bounds on bad exponential sums}
\label{SEC:new-bounds-on-bad-sums-S}

We first provide some new, general, \emph{vanishing and boundedness criteria} for $S_{\bm{c}}(p^l)$,
which when combined with classical estimates from \cites{hooley1986HasseWeil,heath1998circle}
will (under Conjecture~\ref{CNJ:(SFSCp)}) allow us to break a critical $\eps$-barrier behind \eqref{INEQ:near-optimal-diagonal-GRH-bound-over-smooth-locus-S_1}.
Recall $S_{\bm{c}}(n)$, $S^\natural_{\bm{c}}(n)$ from \eqref{EQN:define-S_c(n)}, \eqref{EQN:define-normalized-S-tilde-and-J}.

\begin{lemma}
\label{LEM:bad-sum-vanishing-and-boundedness-criteria}
Let $\bm{c}\in \mcal{S}_1$, and let $p$ be a prime.
\begin{enumerate}
    \item If $p\nmid \gcd(\Delta(\bm{c}), \partial_{c_1}{\Delta}(\bm{c}),\dots,\partial_{c_m}{\Delta}(\bm{c}))$,
    then $S^\natural_{\bm{c}}(p) \ll_m 1$.
    
    \item If $v_p(\Delta(\bm{c}))\le 1$,
    then $S^\natural_{\bm{c}}(p^2) \ll_F 1$.
    
    \item We have $S_{\bm{c}}(p^l)=0$ for all integers $l\ge 2+v_p(\Delta(\bm{c}))$.
\end{enumerate}
\end{lemma}

\begin{proof}
(1):
If $p\mid 6$, then $\abs{S_{\bm{c}}(p)}\le p^{1+m}$ trivially, so $\abs{S^\natural_{\bm{c}}(p)}\le p^{(1+m)/2} \ll_m 1$.
Now suppose $p\nmid 6 \gcd(\Delta(\bm{c}), \partial_{c_1}{\Delta}(\bm{c}),\dots,\partial_{c_m}{\Delta}(\bm{c}))$.
Then by \eqref{EQN:define-Delta-as-product-of-disc}, we have
\begin{equation*}
p\nmid \disc(F),
\quad p\nmid 6,
\quad p\nmid \gcd(\disc(F,\bm{c}),
\partial_{c_1}{\disc(F,\bm{c})},\dots,\partial_{c_m}{\disc(F,\bm{c})}).
\end{equation*}
Therefore, by \eqref{EQN:define-normalized-point-count-errors-E_F,E_c} and \cite{wang2023dichotomous}*{Theorem~1.1, and Proposition~2.7(2)$\Rightarrow$(1) with $(n,r,d)=(m-1,2,3)$}, we have $\abs{E^\natural_{\bm{c}}(p)} \le 72\cdot 9^m$.
Yet $E^\natural_F(p) \ll_m 1$, by \eqref{EQN:define-normalized-point-count-errors-E_F,E_c} and the Weil conjectures (since $p\nmid \disc(F)$).
Plugging these two estimates into \eqref{EQN:rewrite-S_c(p)-via-E_c}, we get $S^\natural_{\bm{c}}(p) \ll_m 1$.

(If $m=6$ and $F$ is diagonal, one can also improve (1) to a ``codimension-three'' statement, by using \cite{wang2023dichotomous}*{Theorem~1.3} in place of \cite{wang2023dichotomous}*{Theorem~1.1}.)

We now turn to (2)--(3).
For any vector $\bm{u}\in \ZZ^m$, let $v_p(\bm{u})\defeq v_p(\gcd(u_1,\dots,u_m))$.
Write $\bm{c} = p^g \tilde{\bm{c}}$, where $g=v_p(\bm{c})<\infty$ and $\tilde{\bm{c}}\in \mcal{S}_1$.
For integers $u,d\ge 0$, let
\begin{equation*}
\mscr{B}^{(d)}(\bm{c};p,u)\defeq \bigcup_{\lambda\in \ZZ:\, p\nmid \lambda} \set{\bm{x}\in \ZZ^m: p\nmid \bm{x},\; p^u\mid F(\bm{x}),\; p^u\mid \tilde{\bm{c}}\cdot \bm{x},\; p^d\mid \grad{F}(\bm{x}) - \lambda\bm{c}}.
\end{equation*}
For any $u,d\ge 0$ and $\bm{x}\in \mscr{B}^{(d)}(\bm{c};p,u)$, Corollary~\ref{COR:Delta(c)-in-double-saturation} (with $\tilde{\bm{c}}$ in place of $\bm{c}$) implies
\begin{equation}
\label{INEQ:key-gradient-descent-type-divisibility}
\gcd(p^u, p^{2d}) \mid \Delta(\tilde{\bm{c}}).
\end{equation}

For any set $A\belongs \ZZ^m$ that can be written as a finite union of residue classes $\mcal{R}\belongs \ZZ^m$, let $\mu(A)$ be the density of $A$ in $\ZZ^m$.
Now let $l\ge 2+2g$ be an integer, and let $d=\floor{l/2}\ge 1+g$; then by \cite{wang2023_isolating_special_solutions}*{Proposition~7.4}, we have
\begin{equation}
\label{EQN:compute-S_c-via-mscrB-singular-sets}
p^{-lm}\varphi(p^l)S'_{\bm{c}}(p^l)
= p^{2l+g} \mu(\mscr{B}^{(d)}(\bm{c};p,l))
- p^{2l-2+g} \mu(\mscr{B}^{(d)}(\bm{c};p,l-1)),
\end{equation}
where $S'_{\bm{c}}(p^l)$ (the ``restriction to $p\nmid \bm{x}$'' of $S_{\bm{c}}(p^l)$) is defined as in \cite{wang2023_isolating_special_solutions}*{(7.5)}.

(2):
Suppose $v_p(\Delta(\bm{c}))\le 1$.
Then $p\nmid \bm{c}$ (since $\deg{\Delta}\ge 2$), so $g=0$ and $\tilde{\bm{c}}=\bm{c}$.
In particular, $S_{\bm{c}}(p^2) = S'_{\bm{c}}(p^2)$ by \cite{wang2023_isolating_special_solutions}*{Lemma~7.2}.
By \eqref{EQN:compute-S_c-via-mscrB-singular-sets} (with $l=2$ and $d=\floor{l/2}=1$), it remains to analyze $\mscr{B}^{(1)}(\bm{c};p,u)$ for $u\in \set{1,2}$.
By \eqref{INEQ:key-gradient-descent-type-divisibility}, we have $\mscr{B}^{(1)}(\bm{c};p,2) = \emptyset$ (since $p^2\nmid \Delta(\tilde{\bm{c}})$).
We now analyze $\mscr{B}^{(1)}(\bm{c};p,1)$.

For each $\bm{a}\in \mscr{B}^{(1)}(\bm{c};p,1)$, the variety $F(\bm{x})=\bm{c}\cdot\bm{x}=0$ in $\PP^{m-1}_{\FF_p}$ is singular at the point $[a_1:\cdots:a_m]\in \PP^{m-1}(\FF_p)$.
Since $p\nmid \bm{c}$, this variety is isomorphic to a cubic hypersurface $f(y_1,\dots,y_{m-1})=0$ in $\PP^{m-2}_{\FF_p}$, and has a singular locus of dimension $\le 0$ if $p\nmid \disc(F)$ (by \cite{wang2023dichotomous}*{Theorem~2.3, due to Zak}).
So if $p\nmid 3\disc(F)$, then $p^m \mu(\mscr{B}^{(1)}(\bm{c};p,1))\le (p-1)\cdot 2^{m-1}$ by B\'{e}zout's theorem (applied to the system $\grad{f}=0$ in $\PP^{m-2}_{\FF_p}$; note that $3f = \bm{y}\cdot \grad{f}$).

But if $p\mid 3\disc(F)$, then $p^m \mu(\mscr{B}^{(1)}(\bm{c};p,1))\le p^m - 1$ trivially.
So in every case, $p^m \mu(\mscr{B}^{(1)}(\bm{c};p,1)) \ll_F p-1$.
Now by \eqref{EQN:define-normalized-S-tilde-and-J} and \eqref{EQN:compute-S_c-via-mscrB-singular-sets}, we have
\begin{equation*}
S^\natural_{\bm{c}}(p^2)
= \frac{S_{\bm{c}}(p^2)}{p^{1+m}}
= \frac{S'_{\bm{c}}(p^2)}{p^{1+m}}
= \frac{p^{m-1}}{\varphi(p^2)} \cdot \left(p^4\cdot 0 - p^2 \cdot \mu(\mscr{B}^{(1)}(\bm{c};p,1))\right)
\ll_F \frac{p(p-1)}{\varphi(p^2)} = 1.
\end{equation*}

(3):
Take a counterexample $(\bm{c}, l)\in \mcal{S}_1\times \ZZ_{\ge 0}$ with $l$ minimal.
Then $l\ge 2+v_p(\Delta(\bm{c}))$ and $S_{\bm{c}}(p^l)\ne 0$.
So $l\ge 2+2g$, because $\deg{\Delta}\ge 2$.
Furthermore, $d=\floor{l/2}\ge (l-1)/2$.
By \eqref{INEQ:key-gradient-descent-type-divisibility}, we have $\mscr{B}^{(d)}(\bm{c};p,u) = \emptyset$ for all $u\ge l-1$, since $2d\ge l-1$ and $p^{l-1}\nmid \Delta(\bm{c})$.
So \eqref{EQN:compute-S_c-via-mscrB-singular-sets} gives $S'_{\bm{c}}(p^l) = 0$, whence $S_{\bm{c}}(p^l)\ne S'_{\bm{c}}(p^l)$.
If $g=0$, this contradicts \cite{wang2023_isolating_special_solutions}*{Lemma~7.2}.
If $g=1$, then $l\ge 4$, so $S_{\bm{c}}(p^l)-S'_{\bm{c}}(p^l) = 0$ by \cite{wang2023_isolating_special_solutions}*{Lemma~7.2(2)}; again, a contradiction.
Now suppose $g\ge 2$; then $p^2\mid \bm{c}$ and $l-3\ge 2+v_p(\Delta(\bm{c}/p^2))$ (since $\deg{\Delta}\ge 3/2$), so $S_{\bm{c}/p^2}(p^{l-3}) = 0$ by the minimality hypothesis.
But \cite{wang2023_isolating_special_solutions}*{Lemma~7.2(2)} then gives $S_{\bm{c}}(p^l)-S'_{\bm{c}}(p^l) = 0$; another contradiction.
Therefore, no counterexample $(\bm{c}, l)\in \mcal{S}_1\times \ZZ_{\ge 0}$ to (3) in fact exists.

For an alternative, more algorithmic and computational approach to (2)--(3) (at least when $F$ is diagonal), see \cite{wang2022thesis}*{\S7.2} and \cite{wang2021_HLH_vs_RMT}*{Appendix~D}.
\end{proof}

\begin{remark}
\label{RMK:failure-of-prime-boundedness-criteria-for-quadratic-F}
If $F$ were quadratic (rather than cubic), then Lemma~\ref{LEM:bad-sum-vanishing-and-boundedness-criteria}(1) would be false whenever $2\mid m$.
See \cite{wang2022thesis}*{Remark~7.2.4} or \cite{wang2023dichotomous}*{sentence after Theorem~1.1}.
\end{remark}

Let $\mcal{N}_{\le}(t) \defeq \set{n\geq 1: p\mid n\Rightarrow v_p(n)\le t}$
and $\mcal{N}_{\ge}(t) \defeq \set{n\geq 1: p\mid n\Rightarrow v_p(n)\ge t}$
for each integer $t\ge 1$.
For any integers $N,t\ge 1$, we have (see e.g.~\cite{bateman1958theorem})
\begin{equation}
\label{INEQ:count-square-full-and-cube-full}
\card{\set{N\le n<2N: n\in \mcal{N}_{\ge}(t)}}\ll_t N^{1/t}.
\end{equation}

Recall $\mcal{N}_{\bm{c}}$ from \eqref{EQN:define-moduli-sets-N^c,N_c}.
Lemma~\ref{LEM:bad-sum-vanishing-and-boundedness-criteria}(1) has a useful unconditional consequence:

\begin{proposition}
\label{PROP:unconditional-Ekedahl-sieve-application-B3Gp}
Let $A\in \RR_{\ge 1}$.
Uniformly over reals $Z, N > 0$ with $N\leq Z^3$,
we have
\begin{equation}
\label{INEQ:B3Gp-goal}
\sum_{\bm{c}\in \mcal{S}_1\cap [-Z,Z]^m}
\biggl(\,
\sum_{n\in \mcal{N}_{\bm{c}}\cap \mcal{N}_{\le}(1)\cap [N,2N)} n^{-1/2} \abs{S^\natural_{\bm{c}}(n)}
\biggr)^{\!2A}
\ll_{F,A,\eps} Z^m \min(N, Z)^{-1+\eps}.
\end{equation}
\end{proposition}

\begin{proof}

The case $N<1$ is trivial,
so assume $N\ge 1$.
Let $\mcal{D} = \set{1,2,4,8,\ldots}$.
Let
\begin{equation}
\label{EQN:define-N_c,0-doubly-singular-moduli}
\mcal{N}_{\bm{c},0} \defeq \set{n\in \mcal{N}_{\le}(1): p\mid n\Rightarrow p\mid \gcd(\Delta(\bm{c}), \partial_{c_1}{\Delta}(\bm{c}),\dots,\partial_{c_m}{\Delta}(\bm{c}))}
\end{equation}
for each $\bm{c}\in \ZZ^m$.
For each $Q\in \mcal{D}$,
let $\mcal{S}_3(Q)$ be the set of tuples $\bm{c}\in \mcal{S}_1$ for which $\mcal{N}_{\bm{c},0}\cap [Q,2Q) \ne \emptyset$.
Note that $1\in \mcal{N}_{\bm{c},0}$, so $\mcal{S}_3(1) = \mcal{S}_1$.

It is known (e.g.~by \cite{heath1983cubic}*{Lemma~11}) that $S^\natural_{\bm{c}}(p) \ll_m p^{1/2}$ for all $\bm{c}\in \ZZ^m$ and primes $p$.
So for all $\bm{c}\in \mcal{S}_1$ and $n\in \mcal{N}_{\le}(1)$, Lemma~\ref{LEM:bad-sum-vanishing-and-boundedness-criteria}(1) implies
\begin{equation*}
S^\natural_{\bm{c}}(n)
\ll_{m,\eps} n^\eps \cdot \max_{q\mid n:\, q\in \mcal{N}_{\bm{c},0}}{q^{1/2}}
\ll n^\eps \cdot \left(\max{\set{Q\in \mcal{D}\cap [1,n]: \bm{c}\in \mcal{S}_3(Q)}}\right)^{1/2}.
\end{equation*}
Therefore, the left-hand side of \eqref{INEQ:B3Gp-goal} is at most $O_{m,A,\eps}(1)$ times the quantity
\begin{equation}
\label{EQN:switch-order-preparation-for-Ekedahl-bound}
\sum_{\substack{\bm{c}\in \mcal{S}_1: \\ \norm{\bm{c}}\le Z}}
\card{\mcal{N}_{\bm{c}}\cap [N,2N)}^{2A}
\sum_{\substack{Q\in \mcal{D}\cap [1,2N): \\ \bm{c}\in \mcal{S}_3(Q)}}
\frac{N^\eps Q^A}{N^A}
= \sum_{Q\in \mcal{D}\cap [1,2N)} \frac{N^\eps Q^A}{N^A}
\sum_{\substack{\bm{c}\in \mcal{S}_3(Q): \\ \norm{\bm{c}}\le Z}}
\card{\mcal{N}_{\bm{c}}\cap [N,2N)}^{2A}.
\end{equation}

Let $\eps\in (0,1)$, and consider an individual $Q\in \mcal{D}$.
Lemma~\ref{LEM:N_c-small-divisor-moment-bound}, and H\"{o}lder over $\bm{c}$, give
\begin{equation}
\label{INEQ:over-extended-N_c-Holder-from-S_3-to-S_1}
\sum_{\bm{c}\in \mcal{S}_3(Q)\cap [-Z,Z]^m} \card{\mcal{N}_{\bm{c}}\cap [N,2N)}^{2A}
\ll_{A,\eps} \card{\mcal{S}_3(Q)\cap [-Z,Z]^m}^{1-\eps} \cdot (Z^m N^\eps)^\eps.
\end{equation}
On the other hand, the scheme $\Delta = \partial_{c_1}{\Delta} = \dots = \partial_{c_m}{\Delta} = 0$ in $\Aff^m_\QQ$ has dimension $\le m-2$ (since $\Delta=0$ is generically smooth, due to the absolute irreducibility of $\Delta$).
So by a quantitative form of \cite{ekedahl1991infinite}'s geometric sieve (see \cite{bhargava2014geometric}*{Theorem~3.3} for the case of prime moduli, which extends to square-free moduli as in \cite{bhargava2021galois}*{\S5, Case~III}), we have
\begin{equation}
\label{INEQ:Ekedahl-bound-on-S_3(Q)}
\card{\mcal{S}_3(Q)\cap [-Z,Z]^m} \ll_{\Delta,\eps} Z^m Q^{-1+\eps} + Z^{m-1+\eps};
\end{equation}
this bound follows from Lang--Weil if $Q\le Z$ (since $\sum_{n\in \mcal{N}_{\le}(1)\cap [Q, 2Q)} O_\eps(\frac{Z^m}{n^{2-\eps}}) \ll_\eps Z^m Q^{-1+\eps}$; cf.~\cite{bhargava2014geometric}*{(16)}), and from elimination theory if $Q>Z$ (cf.~\cite{bhargava2014geometric}*{(17)} and \cite{bhargava2021galois}*{the three paragraphs after Proposition~33}).

Upon inserting \eqref{INEQ:Ekedahl-bound-on-S_3(Q)} into \eqref{INEQ:over-extended-N_c-Holder-from-S_3-to-S_1}, we find that the right-hand side of \eqref{EQN:switch-order-preparation-for-Ekedahl-bound} is
\begin{equation*}
\ll_{A,\eps} \sum_{Q\in \mcal{D}\cap [1,2N)} \frac{N^\eps Q^A}{N^A} \cdot Z^m N^{\eps^2}(Q^{-(1-\eps)^2} + Z^{-(1-\eps)^2})
\ll_{A,\eps} N^{2\eps} \cdot Z^m (N^{-(1-\eps)^2} + Z^{-(1-\eps)^2}),
\end{equation*}
since $A\ge 1$.
Since $N^\eps\le \min(N, Z)^{3\eps}$, the bound \eqref{INEQ:B3Gp-goal} follows immediately.
\end{proof}

We now build on Conjecture~\ref{CNJ:(SFSCp)}.

\begin{conjecture}
[SFSC$_{q,3}$]
\label{CNJ:(SFSCq)}
There exists a real $\eta_9=\eta_9(\Delta)>0$ such that
\begin{equation}
\label{INEQ:SFSCq-goal}
\#\set{\bm{c} \in [-Z,Z]^m
:\exists\;\textnormal{$q \in \mcal{N}_{\ge}(2)\cap [Q,2Q)$ with $q \mid \Delta(\bm{c})$}}
\ll_{\Delta} Z^m Q^{-\eta_9}
\end{equation}
holds uniformly over reals $Z, Q \ge 1$ with $Q\leq Z^3$.
\end{conjecture}

\begin{proposition}
\label{PROP:SFSCp-implies-SFSCq}
Assume Conjecture~\ref{CNJ:(SFSCp)}.
Then Conjecture~\ref{CNJ:(SFSCq)} holds.
\end{proposition}

\begin{proof}
Suppose $Z, Q \ge 1$ are reals with $Q\leq Z^3$.

\emph{Case~1: $Q\le Z$.}
Recall the notation $N(P;q)$ from \S\ref{SEC:local-control-on-polynomials-and-L-functions}.
Let $\delta = (\deg{\Delta})^{-1}$ and $\eps\in (0, \delta)$.
The left-hand side of \eqref{INEQ:SFSCq-goal} is at most
\begin{equation}
\label{INEQ:reduce-unlikely-complete-sum-bound-to-Euler-product}
\sum_{q \in \mcal{N}_{\ge}(2)\cap [Q,2Q)} N(\Delta;q)\cdot O{\left(\frac{Z^m}{q^m}\right)}
\ll_\eps \frac{Z^m}{Q^{\delta/2-\eps}}
\sum_{q \in \mcal{N}_{\ge}(2)} q^{\delta/2-\eps}\cdot \frac{N(\Delta;q)}{q^m}.
\end{equation}
Let $p$ be a prime and $l\ge 2$ an integer.
By \eqref{INEQ:multivarite-zero-density-mod-q}, we have $N(\Delta;p^l)\ll_\Delta p^{l(m-\delta)}$.
On the other hand, by Hensel's lemma (over the smooth locus of $\Delta=0$) and Lang--Weil (over the singular locus of $\Delta=0$), we have $N(\Delta;p^2)\ll_\Delta p^{2m-2}$, and thus $N(\Delta;p^l)\le p^{(l-2)m} N(\Delta;p^2) \ll_\Delta p^{lm-2}$.
So by the Chinese remainder theorem, the right-hand side of \eqref{INEQ:reduce-unlikely-complete-sum-bound-to-Euler-product} is
\begin{equation*}
\le \frac{Z^m}{Q^{\delta/2-\eps}} \prod_p
\biggl(1+\sum_{l\ge 2} p^{(\delta/2-\eps)l} \cdot O_\Delta(p^{-\max(\delta l, 2)})\biggr)
\le \frac{Z^m}{Q^{\delta/2-\eps}}
\biggl(1+\sum_{l\ge 2} O_\Delta(p^{-1-\eps l})\biggr)
\end{equation*}
(since $(\delta/2-\eps)l - \max(\delta l, 2) \le -1 - \eps l$).
Since $\eps>0$, it follows that the right-hand side of \eqref{INEQ:reduce-unlikely-complete-sum-bound-to-Euler-product} is $\ll_{\Delta,\eps} Z^m Q^{-\delta/2+\eps}$.
Hence the left-hand side of \eqref{INEQ:SFSCq-goal} is $\ll_{\Delta,\eps} Z^m Q^{-\delta/2+\eps}$.

\emph{Case~2: $Q>Z$.}
Suppose $q\in \mcal{N}_{\ge}(2)\cap [Q, 2Q)$.
Let $s$ be the largest square divisor of $q$; then $s\ge q^{2/3}\ge Q^{2/3}>Z^{2/3}$.
The integer $s^{1/2}$ thus lies in $(Z^{1/4}, 2Q^{1/2})$, and hence either has a prime factor $p>Z^{1/4}$ (so that $p\in [P, 2P)$ for some $P\in [Z^{1/4}, Q^{1/2}]$), or an integer factor $d\in (Z^{1/4}, Z^{1/2}]$ (so that $d^2\in [D, 2D)$ for some $D\in [Z^{1/2}, Z]$).
So by Conjecture~\ref{CNJ:(SFSCp)} (with $P\in [Z^{1/4}, Q^{1/2}]$) and Case~1 (with $D\in [Z^{1/2}, Z]$ in place of $Q$), the left-hand side of \eqref{INEQ:SFSCq-goal} is
\begin{equation*}
\ll_{\Delta,\eps} Z^m (Z^{1/4})^{-\eta_0} + Z^m (Z^{1/2})^{-\delta/2+\eps}.
\end{equation*}

Since $Z\ge Q^{1/3}$, it follows that \eqref{INEQ:SFSCq-goal} holds with $\eta_9 = \frac{1}{12} \min(\eta_0, \frac{9}{10} (\deg{\Delta})^{-1})$.
\end{proof}

For all $W\in \ZZ[\bm{c}]$ and $Z,N,A\in \RR_{>0}$, let
\begin{equation}
\label{EQN:define-Sigma_14-general-B3G-type-sum}
\Sigma_{14}^{W,A}(Z,N) \defeq 
\sum_{\bm{c}\in \mcal{S}_1\cap [-Z,Z]^m:\, W(\bm{c})\ne 0}
\biggl(\,
\sum_{n\in \mcal{N}_{\bm{c}}\cap [N,2N)} n^{-1/2} \abs{S^\natural_{\bm{c}}(n)}
\biggr)^{\!A}.
\end{equation}
Close analogs of the moment $\Sigma_{14}^{1,A}(Z,N)$ have been considered for $A=1$ classically (e.g.~in \cites{heath1983cubic,hooley1986HasseWeil,heath1998circle}), and for $A=2$ in \cite{wang2023_large_sieve_diagonal_cubic_forms}*{Propositions~4.12 and~A.1}.

\begin{conjecture}
[B2]
\label{CNJ:(B2)}
Let $Z, N\in \RR_{\ge 1}$ with $N\leq Z^3$.
Then $\Sigma_{14}^{1,2}(Z,N) \ll_{F,\eps} Z^{m+\eps}$.
\end{conjecture}

\begin{proposition}
\label{PROP:B2-for-diagonal-F}
Suppose $F$ is diagonal.
Then Conjecture~\ref{CNJ:(B2)} holds.
Also, for $W=c_1\cdots c_m$ and $1\le N\le Z^3$, we have
$\Sigma_{14}^{1,2}(Z,N) - \Sigma_{14}^{W,2}(Z,N) \ll_{F,\eps} Z^{m-1+\eps} N^{1/3}$.
\end{proposition}

Proposition~\ref{PROP:B2-for-diagonal-F} is essentially due to \cite{wang2023_large_sieve_diagonal_cubic_forms}; but we need to go one step further.

\begin{conjecture}
[B3G]
\label{CNJ:(B3G)}
Let $1\le A\le 2$.
There exists a nonzero polynomial $W=W_{F,A}\in \ZZ[\bm{c}]$, and a real $\eta_{10}=\eta_{10}(F,A)>0$, such that
if $1\le N\le Z^3$, then $\Sigma_{14}^{W,A}(Z,N) \ll_{F,A} Z^m N^{-\eta_{10}}$.
\end{conjecture}

One can extend Conjecture~\ref{CNJ:(B3G)} to $A=2+\delta$, almost for free (see Proposition~\ref{PROP:B3G+}), but it is not clear to us what the limit is.
Also, for $1\le A\le 2$, one might hope for $\eta_{10}(F,A)\approx A/2$ to be admissible, by comparing with Proposition~\ref{PROP:unconditional-Ekedahl-sieve-application-B3Gp} or \cite{sarnak1991bounds}*{Conjecture~1}.

\begin{proposition}
\label{PROP:SFSCp-implies-B3G}
Suppose $F$ is diagonal.
Assume Conjecture~\ref{CNJ:(SFSCp)}.
Then Conjecture~\ref{CNJ:(B3G)} holds with $W=c_1\cdots c_m$ (for all $A\in [1,2]$).
\end{proposition}

\begin{remark}
If for each $m'\in \set{2,3,\dots,m}$, one assumed an analog of Conjecture~\ref{CNJ:(SFSCp)} with $m'$ in place of $m$,
then one could prove Proposition~\ref{PROP:SFSCp-implies-B3G} with $W=1$.
See \cite{wang2022thesis}*{Conjecture~7.3.7~(B3) and Remark~7.3.13} and \cite{wang2021_HLH_vs_RMT}*{Lemma~7.43} for details.
\end{remark}

We need a classical pointwise bound on $S^\natural_{\bm{c}}(n)$.
For integers $c\neq 0$, let $\map{sq}(c)\defeq \prod_{p^2\mid c} p^{v_p(c)}$ and $\map{cub}(c)\defeq \prod_{p^3\mid c} p^{v_p(c)}$.
Also let $\map{sq}(0)\defeq 0$ and $\map{cub}(0)\defeq 0$.

\begin{proposition}
[\cites{hooley1986HasseWeil,heath1998circle};
see e.g.~\cite{wang2023_large_sieve_diagonal_cubic_forms}*{Proposition~4.9}]
\label{PROP:diagonal-pointwise-bound-on-S_c(n)}
Assume $F$ is diagonal.
Let $\eps>0$.
There exists $A_{12}(F,\eps)\in \RR_{\ge 1}$ such that if $\bm{c}\in \ZZ^m$ and $n\in \ZZ_{\ge 1}$, then
\begin{equation*}
n^{-1/2}\abs{S^\natural_{\bm{c}}(n)}
\le A_{12}(F,\eps) n^\eps
\prod_{1\leq j\leq m} {\gcd}{\left(\map{cub}(n)^2, \gcd(\map{cub}(n),\map{sq}(c_j))^3\right)}^{1/12}.
\end{equation*}
\end{proposition}

One could replace the $n^\eps$ in Proposition~\ref{PROP:diagonal-pointwise-bound-on-S_c(n)} with $O(1)^{\omega(n)}$.
This would be important if one wanted to try to prove a softer version of Proposition~\ref{PROP:SFSCp-implies-B3G} (conditional on a softer version of Conjecture~\ref{CNJ:(SFSCp)}).
We do not explore this direction in the present paper, since it would seem to involve complicated divisor-type sums with square-full parts and discriminant divisibility conditions (which might require complicated parameterizations to handle).

Proposition~\ref{PROP:diagonal-pointwise-bound-on-S_c(n)} is an explicit stratification result for $S_{\bm{c}}(n)$, based on $v_p(c_j)$ for $p\mid n$.
It would be very nice to have a usable (perhaps less explicit) replacement for Proposition~\ref{PROP:diagonal-pointwise-bound-on-S_c(n)} when $F$ is no longer diagonal.
Work such as \cites{denef1984rationality,pas1989uniform}
could conceivably help.

We now prove Propositions~\ref{PROP:B2-for-diagonal-F} and~\ref{PROP:SFSCp-implies-B3G}.
For convenience, let
$\mcal{S}_4\defeq \set{\bm{c}\in \mcal{S}_1: c_1\cdots c_m\ne 0}$.

\begin{proof}
[Proof of Proposition~\ref{PROP:B2-for-diagonal-F}]
Let $W = c_1\cdots c_m$.
Suppose $Z, N \in \RR_{\ge 1}$ with $N\leq Z^3$.
For each $\bm{c}\in \mcal{S}_4$, Proposition~\ref{PROP:diagonal-pointwise-bound-on-S_c(n)} implies $n^{-1/2} S^\sharp_{\bm{c}}(n) \ll_{F,\eps} n^\eps \prod_{1\le j\le m} \map{sq}(c_j)^{1/4}$.
Plugging this into $\Sigma_{14}^{W,2}(Z,N)$ (see \eqref{EQN:define-Sigma_14-general-B3G-type-sum}), and using Lemma~\ref{LEM:count-R-N_c-infty-divisors} to bound $\card{\mcal{N}_{\bm{c}}\cap [N,2N)}$, we get $\Sigma_{14}^{W,2}(Z,N) \ll_{F,\eps} Z^\eps \prod_{1\le j\le m} \sum_{1\le \abs{c}\le Z} \map{sq}(c_j)^{1/2}$.
However, by \eqref{INEQ:count-square-full-and-cube-full}, we have
\begin{equation}
\label{INEQ:sum-sq(c)^(1/2)}
\sum_{1\le \abs{c}\le Z} \map{sq}(c)^{1/2} \ll \sum_{e\in \mcal{N}_{\ge}(2)} e^{1/2}\cdot (Z/e) \ll Z \log(1+Z).
\end{equation}
Therefore, $\Sigma_{14}^{W,2}(Z,N) \ll_{F,\eps} Z^{m+\eps}$.
The statement $\Sigma_{14}^{1,2}(Z,N) - \Sigma_{14}^{W,2}(Z,N) \ll_{F,\eps} Z^{m-1+\eps} N^{1/3}$ holds by a similar calculation with Proposition~\ref{PROP:diagonal-pointwise-bound-on-S_c(n)}, Lemma~\ref{LEM:count-R-N_c-infty-divisors}, and \eqref{INEQ:sum-sq(c)^(1/2)},
ending with
\begin{equation*}
\map{cub}(n)^{(m-t)/3} \prod_{1\le i\le t} \sum_{1\le \abs{c_i}\le Z} \map{sq}(c_i)^{1/2}
\ll_{t,\eps} N^{(m-t)/3} Z^{t+\eps} \le N^{1/3} Z^{m-1+\eps}
\end{equation*}
for $n\in [N,2N)$, $t\in [1,m-1]$;
cf.~\cite{wang2023_large_sieve_diagonal_cubic_forms}*{(4.11) in the proof of Proposition~4.12}.
Since $N\le Z^3$, we obtain $\Sigma_{14}^{1,2}(Z,N) \ll_{F,\eps} Z^{m+\eps}$, proving Conjecture~\ref{CNJ:(B2)}.
\end{proof}

\begin{proof}
[Proof of Proposition~\ref{PROP:SFSCp-implies-B3G}]
Suppose $1\le N\le Z^3$.
For each $\bm{c}\in \ZZ^m$ and $n\in \ZZ_{\ge 1}$, let
\begin{equation*}
S^\sharp_{\bm{c}}(n)\defeq \prod_{p\mid n} \max\left(1, \frac{\abs{S^\natural_{\bm{c}}(p^{v_p(n)})}}{p^{v_p(n)/2}}\right).
\end{equation*}
Let $P\ge 1$ be a real parameter to be specified later.
Let
\begin{equation*}
\Sigma_{15} \defeq \sum_{\substack{\bm{c}\in \mcal{S}_1: \\ \norm{\bm{c}}\le Z}}\,
\biggl(\,
\sum_{\substack{n\in \mcal{N}_{\bm{c}}\cap [N,2N): \\ S^\sharp_{\bm{c}}(n) < P}} \frac{\abs{S^\natural_{\bm{c}}(n)}}{n^{1/2}}
\biggr)^{\!2},
\quad
\Sigma_{16} \defeq \sum_{\substack{\bm{c}\in \mcal{S}_4: \\ \norm{\bm{c}}\le Z}}\,
\biggl(\,
\sum_{\substack{n\in \mcal{N}_{\bm{c}}\cap [N,2N): \\ S^\sharp_{\bm{c}}(n) \ge P}} \frac{\abs{S^\natural_{\bm{c}}(n)}}{n^{1/2}}
\biggr)^{\!2}.
\end{equation*}
Let $W = c_1\cdots c_m$; then $\Sigma_{14}^{W,2}(Z,N) \le 2(\Sigma_{15}+\Sigma_{16})$, since $(a+b)^2 \le 2(a^2+b^2)$ for $a,b\in \RR$.
We will bound $\Sigma_{15}$ conditionally (using Conjecture~\ref{CNJ:(SFSCp)}), and $\Sigma_{16}$ unconditionally (using the diagonality of $F$).
We will lose factors of $Z^\eps$ at first, and then remove $Z^\eps$ later.

We first handle $\Sigma_{15}$.
Given $\bm{c}\in \mcal{S}_1\cap [-Z,Z]^m$, let
\begin{equation*}
\begin{split}
\mcal{N}_{\bm{c},1} &\defeq \set{n\in \mcal{N}_{\ge}(2):
p\mid n\Rightarrow \lcm(p^2,p^{v_p(n)-1})\nmid \Delta(\bm{c})}, \\
\mcal{N}_{\bm{c},2} &\defeq \set{n\in \mcal{N}_{\ge}(2):
p\mid n\Rightarrow \lcm(p^2,p^{v_p(n)-1})\mid \Delta(\bm{c})}.
\end{split}
\end{equation*}
Any integer $n\ge 1$ can be written uniquely as $qn_1n_2$, where $q$, $n_1$, $n_2$ are pairwise
coprime integers satisfying $q\in \mcal{N}_{\le}(1)$ and $n_j\in \mcal{N}_{\bm{c},j}$ (for $j\in \set{1,2}$).
We then have
\begin{equation*}
n^{-1/2} \abs{S^\natural_{\bm{c}}(n)}
\le \min_{d\mid n}(S^\sharp_{\bm{c}}(n/d)\cdot d^{-1/2} \abs{S^\natural_{\bm{c}}(d)})
\le S^\sharp_{\bm{c}}(n) \min(1, n_1^{-1/2} \abs{S^\natural_{\bm{c}}(n_1)},
q^{-1/2} \abs{S^\natural_{\bm{c}}(q)})
\end{equation*}
by the definition of $S^\sharp$.
Here $S^\natural_{\bm{c}}(n_1) \ll_F 1$ by Lemma~\ref{LEM:bad-sum-vanishing-and-boundedness-criteria}(2)--(3).
Also, the integer $\prod_{p\mid n_2} \lcm(p^2, p^{v_p(n_2)-1})\in \mcal{N}_{\ge 2}\cap [n_2^{2/3}, n_2]$ divides $\Delta(\bm{c})$.
Since $\card{\mcal{N}_{\bm{c}}\cap [1, 2N)} \ll_\eps (ZN)^\eps$ (by Lemma~\ref{LEM:count-R-N_c-infty-divisors}) and $\min(q, n_1, n_2)\ge n^{1/3}$, we conclude that
\begin{equation*}
\Sigma_{15} \ll_{F,\eps} P^2 Z^\eps
\sum_{\bm{c}\in \mcal{S}_1\cap [-Z,Z]^m}
\biggl(\, \bm{1}_{\bm{c}\in \mcal{S}_5} + \frac{1}{(N^{1/3})^{1/2}}
+ \sum_{q\in \mcal{N}_{\bm{c}}\cap \mcal{N}_{\le}(1)\cap [N^{1/3}, 2N)}
\frac{\abs{S^\natural_{\bm{c}}(q)}}{q^{1/2}}
\biggr)^{\!2},
\end{equation*}
where $\mcal{S}_5$ denotes the set of tuples $\bm{c}\in \mcal{S}_1$ for which there exists $d\in \mcal{N}_{\ge 2}\cap [(N^{1/3})^{2/3}, 2N)$ with $d\mid \Delta(\bm{c})$.
Since $(r+s+t)^2 \le 3(r^2+s^2+t^2)$ (for $r,s,t\in \RR$), the sum over $\bm{c}$ on the right above is $\ll_{F,\eps} Z^m (N^{2/9})^{-\eta_9} + Z^m N^{-1/3} + Z^m N^\eps \min(N^{1/3}, Z)^{-1+\eps}$, by Propositions~\ref{PROP:SFSCp-implies-SFSCq} and~\ref{PROP:unconditional-Ekedahl-sieve-application-B3Gp}.
Since $N\le Z^3$, it follows that
\begin{equation}
\label{INEQ:final-Sigma_15-bound-via-Ekedahl-and-SFSC}
\Sigma_{15} \ll_{F,\eps} P^2 Z^{m+\eps} (N^{-2\eta_9/9} + N^{-1/3}).
\end{equation}

We next handle $\Sigma_{16}$, by introducing a new Ekedahl-type idea.
For each $\bm{c}\in \mcal{S}_1$, let
\begin{equation*}
\mcal{N}_{\bm{c},3}(P) \defeq \set{n\in \mcal{N}_{\bm{c}}: S^\natural_{\bm{c}}(n)\ne 0,\; S^\sharp_{\bm{c}}(n)\ge P}.
\end{equation*}
Let $\eps\in (0, 1)$ (to be specified); suppose $P^{1/10}\ge A_{12}(F,\eps) \cdot (2N)^\eps$.
Now consider an individual $\bm{c}\in \mcal{S}_4\cap [-Z,Z]^m$ and $n\in \mcal{N}_{\bm{c},3}(P)\cap [N, 2N)$.
Here $S^\natural_{\bm{c}}(n)\ne 0$, so $\prod_{p\mid \map{cub}(n)} p^{v_p(n)-1} \mid \Delta(\bm{c})$ by Lemma~\ref{LEM:bad-sum-vanishing-and-boundedness-criteria}(3).
And $S^\sharp_{\bm{c}}(n)\ge P$, so
\begin{equation*}
P^{9/10} \le S^\sharp_{\bm{c}}(n)/P^{1/10} \le S^\sharp_{\bm{c}}(n)/A_{12}(F,\eps)n^\eps
\le \prod_{1\le j\le m} \gcd(\map{cub}(n), \map{sq}(c_j))^{1/4},
\end{equation*}
by Proposition~\ref{PROP:diagonal-pointwise-bound-on-S_c(n)}.
Thus there exists $j$ with $d\defeq \gcd(\map{cub}(n), \map{sq}(c_j))\ge P^{3.6/m}$.
Clearly $d\in \mcal{N}_{\ge}(2)$.
Also, $d\mid c_j$, so $d\le Z$ (since $W(\bm{c})\ne 0$ implies $c_j\ne 0$); and $d\mid \map{cub}(n)\mid \prod_{p\mid \map{cub}(n)} p^{2(v_p(n)-1)} \mid \Delta(\bm{c})^2$.
Letting $\mcal{S}_6(j,d)\defeq \set{\bm{c}\in \mcal{S}_4: d\mid \gcd(c_j, \Delta(\bm{c})^2)}$, it follows (by taking $n\in \mcal{N}_{\bm{c},3}(P)\cap [N,2N)$ with $n^{-1/2} \abs{S^\natural_{\bm{c}}(n)}$ maximal, if such an $n$ exists) that
\begin{equation*}
\Sigma_{16} \le \sum_{1\le j\le m}
\sum_{\substack{d\in \mcal{N}_{\ge}(2): \\ P^{3.6/m}\le d\le Z}}\,
\sum_{\substack{\bm{c}\in \mcal{S}_6(j,d): \\ \norm{\bm{c}}\le Z}}\, \card{\mcal{N}_{\bm{c}}\cap [N, 2N)}^2
\cdot d^{1/2} \prod_{\substack{1\le i\le m: \\ i\ne j}} \map{sq}(c_i)^{1/2}.
\end{equation*}
Using Lemma~\ref{LEM:count-R-N_c-infty-divisors} to bound $\card{\mcal{N}_{\bm{c}}\cap [N, 2N)}$, we get
\begin{equation}
\label{INEQ:Sigma_16-Ekedahl-type-preparation-bound}
\Sigma_{16} \ll_{F,\eps} Z^\eps \sum_{1\le j\le m}
\sum_{\substack{d\in \mcal{N}_{\ge}(2): \\ P^{3.6/m}\le d\le Z}}\,
\sum_{\substack{\bm{c}\in \mcal{S}_6(j,d): \\ \norm{\bm{c}}\le Z}}
d^{1/2} \prod_{\substack{1\le i\le m: \\ i\ne j}} \map{sq}(c_i)^{1/2}.
\end{equation}

For notational simplicity, suppose $j=m$.
For $k\in \set{0,1}$, let
\begin{equation*}
\mcal{S}_{6,k}(m)\defeq \set{(c_1,\dots,c_{m-1})\in (\ZZ\setminus \set{0})^{m-1}:
\bm{1}_{\Delta(c_1,\dots,c_{m-1},0)\ne 0} = k}.
\end{equation*}
For each $(c_1,\dots,c_{m-1})\in \mcal{S}_{6,1}(m)\cap [-Z,Z]^{m-1}$, we have
\begin{equation*}
\sum_{\substack{d\in \mcal{N}_{\ge}(2): \\ P^{3.6/m}\le d\le Z}}\,
\sum_{\substack{1\le \abs{c_m}\le Z: \\ (c_1,\dots,c_m)\in \mcal{S}_6(m,d)}} d^{1/2}
\le \sum_{\substack{d\mid \Delta(c_1,\dots,c_{m-1},0)^2: \\ P^{3.6/m}\le d\le Z}}\,
\sum_{\substack{1\le \abs{c_m}\le Z: \\ d\mid c_m}} d^{1/2}
\ll_{F,\eps} Z^\eps \cdot Z/(P^{3.6/m})^{1/2},
\end{equation*}
by the divisor bound for $\Delta(c_1,\dots,c_{m-1},0)^2\ne 0$.
Therefore, the contribution to the right-hand side of \eqref{INEQ:Sigma_16-Ekedahl-type-preparation-bound} from $\bm{c}$ with $(c_1,\dots,c_{m-1})\in \mcal{S}_{6,1}(m)$ (or the analogous condition if $j\ne m$) is
\begin{equation*}
\ll_{F,\eps} \frac{Z^{1+2\eps}}{P^{1.8/m}}
\prod_{1\le i\le m-1} \sum_{1\le \abs{c_i}\le Z} \map{sq}(c_i)^{1/2}
\ll_{m,\eps} \frac{Z^{m+3\eps}}{P^{1.8/m}},
\end{equation*}
where we use \eqref{INEQ:sum-sq(c)^(1/2)} in the final step.

On the other hand, for any $c_1,\dots,c_{m-2}\in \ZZ$, there are at most $\deg{\Delta}$ integers $c_{m-1}$ for which $(c_1,\dots,c_{m-1})\in \mcal{S}_{6,0}(m)$.
(This is because for diagonal $F$,
the $c_i^{\deg{\Delta}}$ coefficient of $\Delta$ is nonzero for all $i$;
see e.g.~\cite{heath1998circle}*{(4.2)}.)
Therefore, the contribution to the right-hand side of \eqref{INEQ:Sigma_16-Ekedahl-type-preparation-bound} from $\bm{c}$ with $(c_1,\dots,c_{m-1})\in \mcal{S}_{6,0}(m)$ (or the analogous condition if $j\ne m$) is
\begin{equation*}
\ll_{F,\eps} Z^\eps \sum_{\substack{d\in \mcal{N}_{\ge}(2): \\ P^{3.6/m}\le d\le Z}}\,
\sum_{\substack{1\le \abs{c_m}\le Z: \\ d\mid c_m}} d^{1/2}
\max_{1\le \abs{c_{m-1}}\le Z}{\map{sq}(c_{m-1})^{1/2}}
\prod_{1\le i\le m-2} \sum_{1\le \abs{c_i}\le Z} \map{sq}(c_i)^{1/2},
\end{equation*}
and thus (by \eqref{INEQ:sum-sq(c)^(1/2)}) $\ll_{F,\eps} Z^{2\eps} \cdot Z \cdot Z^{1/2} \cdot Z^{m-2} = Z^{m-1/2+2\eps}$.

The previous two paragraphs imply that the right-hand side of \eqref{INEQ:Sigma_16-Ekedahl-type-preparation-bound} is $\ll_{F,\eps} Z^{m+3\eps} (P^{-1.8/m} + Z^{-1/2})$.
This, combined with \eqref{INEQ:final-Sigma_15-bound-via-Ekedahl-and-SFSC}, gives
\begin{equation*}
\Sigma_{14}^{W,2}(Z,N) \le 2(\Sigma_{15}+\Sigma_{16}) \ll_{F,\eps}
P^2 Z^{m+\eps} (N^{-2\eta_9/9} + N^{-1/3}) + Z^{m+3\eps} (P^{-1.8/m} + Z^{-1/2}).
\end{equation*}
Suppose $\eps \le 2\min(1,\eta_9)/270$, and let $P = A_{12}(F,\eps)^{10} 2^{10\eps} N^{2\min(1,\eta_9)/27}\ge 1$; then we get
\begin{equation*}
\Sigma_{14}^{W,2}(Z,N) \ll_{F,\eps} Z^{m+3\eps} (N^{-3.6\min(1,\eta_9)/27m} + Z^{-1/2})
\ll Z^{m+3\eps} N^{-\min(1,\eta_9)/9m},
\end{equation*}
since $m\ge 4$ and $Z\ge N^{1/3}$.
It follows (upon redefining $\eps$, now that we can forget about $P$) that $\Sigma_{14}^{W,2}(Z,N) \ll_{F,\eps} Z^{m+\eps} N^{-\min(1,\eta_9)/9m}$ for all $\eps>0$.

It remains to replace $Z^\eps$ with $N^\eps$.
This is trivial unless $N$ is very small relative to $Z$.
Expanding the square in \eqref{EQN:define-Sigma_14-general-B3G-type-sum}, and noting that $S_{\bm{c}}(n)$ depends only on $\bm{c}\bmod{n}$, gives
\begin{equation*}
\Sigma_{14}^{W,2}(Z,N)
\le \sum_{n_1,n_2\in [N,2N)} \left(1+\frac{2Z}{n_1n_2}\right)^m
\sum_{\substack{1\le \bm{a}\leq n_1n_2: \\ \rad(n_1n_2)\mid \Delta(\bm{a})}}
(n_1n_2)^{-1/2} \abs{S^\natural_{\bm{a}}(n_1)S^\natural_{\bm{a}}(n_2)}
\end{equation*}
for all $Z, N \ge 1$.
Yet for any $n_1,n_2\ge 1$ and $\bm{a}\in [1,n_1n_2]^m$, there exists $\bm{c}\in \mcal{S}_4\cap [1,Bn_1n_2]^m$ with $\bm{c}\equiv \bm{a}\bmod{n_1n_2}$, where $B = 1+\deg(c_1\cdots c_m \Delta(\bm{c}))$.
So
\begin{equation*}
\Sigma_{14}^{W,2}(Z,N)
\ge \sum_{n_1,n_2\in [N,2N)}
\sum_{\substack{1\le \bm{a}\leq n_1n_2: \\ \rad(n_1n_2)\mid \Delta(\bm{a})}}
(n_1n_2)^{-1/2} \abs{S^\natural_{\bm{a}}(n_1)S^\natural_{\bm{a}}(n_2)}
\end{equation*}
for all $Z\ge 4BN^2$.
Hence for all $Z\ge N^2$ we have
\begin{equation*}
\Sigma_{14}^{W,2}(Z,N)\le (3Z/N^2)^m \cdot \Sigma_{14}^{W,2}(4BN^2,N)
\ll_{F,\eps} (Z/N^2)^m (4BN^2)^{m+\eps} N^{-\min(1,\eta_9)/9m},
\end{equation*}
since $1\le N\le (4BN^2)^3$.
For all $Z\ge N^{1/3}$ (including $Z<N^2$, where $Z^\eps<N^{2\eps}$), then,
\begin{equation*}
\Sigma_{14}^{W,2}(Z,N) \ll_{F,\eps} Z^m N^{\eps - \min(1,\eta_9)/9m}.
\end{equation*}
Therefore, Conjecture~\ref{CNJ:(B3G)} holds with $W=c_1\cdots c_m$ and $\eta_{10}(F,2r) = r\min(1,\eta_9)/10m$ (for $r=1$ at first, and then for $r\in [1/2, 1]$ by H\"older over $\bm{c}$).
\end{proof}

\begin{remark}
Handling $\map{sq}(c_i)^{1/2}$ takes care, because $\sum_{1\le \abs{c}\le Z} \map{sq}(c)^{1/2}$ remains roughly the same size even if we restrict $\abs{c}$ to $\mcal{N}_{\ge}(2)$ (a rather sparse set).
One could take $W=1$ in Proposition~\ref{PROP:SFSCp-implies-B3G} if we were only interested in $\Sigma_{14}^{W,A}(Z,N)$ for $A\in [1,2)$, or for $N\le Z^{3-\delta}$.
\end{remark}

To end \S\ref{SEC:new-bounds-on-bad-sums-S},
we prove a result clarifying the nature of Conjecture~\ref{CNJ:(SFSCp)} for diagonal $F$:

\begin{proposition}
\label{PROP:SFSC-for-m=4-implies-SFSC-for-diagonal}
Assume Conjecture~\ref{CNJ:(SFSCp)} holds when $(m, F) = (4, x_1^3+x_2^3+x_3^3+x_4^3)$.
Then Conjecture~\ref{CNJ:(SFSCp)} holds whenever $F$ is diagonal and $m\ge 4$.
\end{proposition}

\begin{proof}
(This result is not used in the rest of the paper, so we confine ourselves to a sketch.)

Assume $1\le P\le Z^{3/2}$.
Recall $\mcal{N}_{\bm{c},0}$
from \eqref{EQN:define-N_c,0-doubly-singular-moduli}.
By \cite{bhargava2014geometric}*{Theorem~3.3}, we have
\begin{equation}
\label{INEQ:Ekedahl-bound-singular-locus-in-SFSCp}
\#\set{\bm{c} \in [-Z,Z]^m
:\exists\;\textnormal{$p \in [P,2P)\cap \mcal{N}_{\bm{c},0}$}}
\ll_{\Delta} Z^m P^{-\Theta(1)}.
\end{equation}
By \eqref{INEQ:Ekedahl-bound-singular-locus-in-SFSCp}, we see that for any given $F$, Conjecture~\ref{CNJ:(SFSCp)} is equivalent to the statement
\begin{equation}
\label{INEQ:SFSCp-equivalent-1}
\#\set{\bm{c} \in [-Z,Z]^m
:\exists\;\textnormal{$p \in [P,2P)\setminus \mcal{N}_{\bm{c},0}$ with $p^2 \mid \Delta(\bm{c})$}}
\ll Z^m P^{-\Theta(1)}.
\end{equation}
Let $\mcal{U}_1, \mcal{U}_2\belongs \PP^{m-1}_\ZZ$ be the smooth loci of the hypersurfaces $F=0$, $\Delta=0$, respectively, over $\ZZ$.
The gradient $\grad{F}$ defines a Gauss map $[\grad{F}]\maps \mcal{U}_1\to \PP^{m-1}_\ZZ$; let $\mcal{U}_3\belongs \mcal{U}_1$ be the inverse image of $\mcal{U}_2$ under this map.
The map $[\grad{F}]\maps \mcal{U}_3\to \mcal{U}_2$ is an isomorphism over $\QQ$ (by the biduality theorem), and thus an isomorphism over $\ZZ$ (since $\mcal{U}_2$, $\mcal{U}_3$ are flat over $\ZZ$).
In particular, $\mcal{U}_2(\ZZ/p^2\ZZ) = [\grad{F}](\mcal{U}_3(\ZZ/p^2\ZZ))$.
Furthermore, it is known that $\mcal{U}_3$ lies in the open subscheme $\det(\map{Hess}{F}(\bm{x}))\ne 0$ of $\PP^{m-1}_\ZZ$.
These two facts, plus \eqref{INEQ:Ekedahl-bound-singular-locus-in-SFSCp}, imply that \eqref{INEQ:SFSCp-equivalent-1} is equivalent to
\begin{equation}
\label{INEQ:SFSCp-equivalent-2}
\begin{split}
\#\{\bm{c} \in [-Z,Z]^m
:\exists\;&p \in [P,2P),\; \bm{x}\in \ZZ^m,\; \lambda\in [1, p-1]\textnormal{ with }\\
&p\nmid \det(\map{Hess}{F}(\bm{x})),\;
p^2\mid F(\bm{x}), \grad{F}(\bm{x}) - \lambda\bm{c}\}
\ll Z^m P^{-\Theta(1)}.
\end{split}
\end{equation}

We would like to convert each ``exists'' into a sum.
Each $p$ on the left-hand side of \eqref{INEQ:SFSCp-equivalent-2} satisfies $p\mid \Delta(\bm{c})$, so there are at most finitely many possibilities for $p$ if $\bm{c}\in \mcal{S}_1$.
Also, $\card{\mcal{S}_0\cap [-Z,Z]^m}\cdot P\ll_m Z^{m-1/2}$ by \eqref{INEQ:dimension-growth-bound-on-S_0}.
Therefore, the statement \eqref{INEQ:SFSCp-equivalent-2} is equivalent to
\begin{equation}
\label{INEQ:SFSCp-equivalent-3}
\sum_{\bm{c}\in [-Z,Z]^m} \sum_{p\in [P, 2P)} \sum_{1\le \bm{x}\le p^2} \bm{1}_{p\nmid \det(\map{Hess}{F}(\bm{x}))} \bm{1}_{p^2\mid F(\bm{x})} \EE_{1\le \lambda\le p-1}[\bm{1}_{p^2 \mid \grad{F}(\bm{x}) - \lambda\bm{c}}]
\ll Z^m P^{-\Theta(1)}.
\end{equation}
Now suppose $F=F_1x_1^3+\dots+F_mx_m^3$ (where $F_j\in \ZZ\setminus \set{0}$), and for each $G, \lambda\in \ZZ$ let
\begin{equation*}
T_{3/2, G, \lambda}(a, p^2) \defeq \sum_{\substack{c\in [-Z,Z],\; 1\le x\le p^2: \\
p\nmid x,\; p^2\mid 3Gx^2-\lambda c}} e_{p^2}(aGx^3) \in \RR.
\end{equation*}
Then the left-hand side of \eqref{INEQ:SFSCp-equivalent-3} equals $1/(p-1)$ times
\begin{equation}
\label{INEQ:key-Holder-type-positive-bound-for-diagonal-SFSCp-equivalent}
\sum_{\substack{p\in [P, 2P), \\ 1\le a\le p^2,\; 1\le \lambda\le p-1}}
\prod_{1\le j\le m} T_{3/2, F_j, \lambda}(a, p^2)
\ll_F Z^{m-4} \sum_{1\le j\le 4}
\sum_{\substack{p\in [P, 2P), \\ 1\le a\le p^2,\; 1\le \lambda\le p-1}} T_{3/2, F_j, \lambda}(a, p^2)^4,
\end{equation}
since $T_{3/2, G, \lambda}(a, p^2) \ll_G Z$ trivially (by considering the cases $p\nmid 3G$ and $p\mid 3G$ separately) and $\prod_{1\le j\le 4} \abs{T_{3/2, F_j, \lambda}(a, p^2)} \le \sum_{1\le j\le 4} T_{3/2, F_j, \lambda}(a, p^2)^4$ (since $T_{3/2, G, \lambda}\in \RR$).

Since Conjecture~\ref{CNJ:(SFSCp)} holds by assumption when $F = x_1^3+x_2^3+x_3^3+x_4^3$, it also holds when $F = F_j \cdot (x_1^3+x_2^3+x_3^3+x_4^3)$ for any fixed $j$ (since scaling $\Delta$ does not affect the truth of Conjecture~\ref{CNJ:(SFSCp)}).
So by \eqref{INEQ:key-Holder-type-positive-bound-for-diagonal-SFSCp-equivalent}, the statement \eqref{INEQ:SFSCp-equivalent-3} holds for all diagonal $F$ (assuming $m\ge 4$).
Hence \eqref{INEQ:SFSCp-equivalent-2}, \eqref{INEQ:SFSCp-equivalent-1}, and Conjecture~\ref{CNJ:(SFSCp)} hold too.
\end{proof}

\section{Delta endgame}
\label{SEC:delta-endgame}

Throughout \S\ref{SEC:delta-endgame}, assume $2\mid m$.
Recall $\Sigma(X,\mcal{S})$, $\Sigma^\natural(X,\mcal{S})$ from \eqref{EQN:define-sum-Sigma(X,S)}.
Explicitly, we have
\begin{equation}
\label{EQN:Fubini-Sigma(X,S)-for-endgame}
\Sigma^\natural(X,\mcal{S})
= \sum_{n\ge 1} \sum_{\bm{c}\in \mcal{S}}
n^{(1-m)/2} S^\natural_{\bm{c}}(n) J_{\bm{c},X}(n)
= \sum_{\bm{c}\in \mcal{S}} \sum_{n\ge 1}
n^{(1-m)/2} S^\natural_{\bm{c}}(n) J_{\bm{c},X}(n)
\end{equation}
(by Proposition~\ref{PROP:basic-integral-facts} and Fubini).
We are finally prepared to analyze these sums for $\mcal{S}\belongs \mcal{S}_1$.



\subsection{Delta decomposition}

Consider an individual $\bm{c}\in \mcal{S}_1$.
In view of \eqref{EQN:define-Phi^G}, \eqref{EQN:define-Phi^B} and the factorization $\Phi = \Phi^{\map{G}} \Phi^{\map{B}}$, we may use Lemma~\ref{LEM:dyadic-partial-Mellin-summation} with $k=2$,
\begin{equation*}
a(\bm{n})=a(n_0,n_1)
= (n_0^{-it_0} S^\natural_{\bm{c}}(n_0)\bm{1}_{n_0\in \mcal{N}_{\bm{c}}})
\cdot (n_1^{-it_1} S^\natural_{\bm{c}}(n_1) \bm{1}_{n_1\in \mcal{N}^{\bm{c}}}),
\end{equation*}
and $f(\bm{r})=f(r_0,r_1) = (r_0r_1)^{(1-m)/2} J_{\bm{c},X}(r_0r_1)$, to write
\begin{equation}
\label{EQN:smooth-dyadic-restriction-separation-for-endgame}
    \sum_{n\ge 1} n^{(1-m)/2} S^\natural_{\bm{c}}(n) J_{\bm{c},X}(n)
    = (2\pi)^{-2} \int_{\bm{N}\ge 1/2} d^\times\bm{N} \int_{\bm{t}\in \RR^2} d\bm{t}\,
    g_{\bm{c},X,\bm{N}}^\vee(i\bm{t})
    \prod_{0\le j\le 1} \Sigma_{17,\bm{N}}^{\bm{c},j}(\bm{t}),
\end{equation}
where $\bm{N}=(N_0,N_1)$ runs over $[1/2, \infty)^2$, where $\bm{t}=(t_0,t_1)\in \RR^2$, and where
\begin{align}
g_{\bm{c},X,\bm{N}}(\bm{r})
&\defeq (r_0r_1)^{(1-m)/2} J_{\bm{c},X}(r_0r_1)
\cdot \nu_2(r_0/N_0) \nu_2(r_1/N_1),
\label{EQN:define-final-dyadic-modulus-weight-for-endgame} \\
\Sigma_{17,\bm{N}}^{\bm{c},0}(\bm{t}) &\defeq \sum_{n_0\in \mcal{N}_{\bm{c}}}
\nu_2(n_0/N_0) n_0^{-it_0} S^\natural_{\bm{c}}(n_0),
\quad
\Sigma_{17,\bm{N}}^{\bm{c},1}(\bm{t}) \defeq \sum_{n_1\in \mcal{N}^{\bm{c}}}
\nu_2(n_1/N_1) n_1^{-it_1} S^\natural_{\bm{c}}(n_1).
\label{EQN:define-smooth-dyadic-pieces-Sigma_17-for-endgame}
\end{align}



Since $J_{\bm{c},X}(r)$ is supported on $r\le A_0Y$ (by Proposition~\ref{PROP:basic-integral-facts}) and $\nu_2$ is supported on $[1,2]$, we have $g_{\bm{c},X,\bm{N}}(\bm{r}) = 0$ unless $r_0r_1\le A_0Y$ and $N_j\le r_j\le 2N_j$ for all $j$.
Thus $g_{\bm{c},X,\bm{N}} = 0$ identically unless $N_0N_1\le A_0Y$.
So $g_{\bm{c},X,\bm{N}}^\vee(i\bm{t}) = 0$ unless $\bm{N}$ lies in the set
\begin{equation}
\label{EQN:dyadic-range-restrictions-for-R_17}
    \mscr{R}_{17}(X)\defeq \set{\bm{N}\ge 1/2: N_0N_1\le A_0Y}
\end{equation}
(cf.~the region $\mscr{R}_{10}$ from \eqref{EQN:dyadic-range-restrictions-for-R_10}).
So \eqref{EQN:smooth-dyadic-restriction-separation-for-endgame} holds even if we restrict $\bm{N}$ to $\mscr{R}_{17}(X)$.


Given $\bm{N}=(N_0,N_1)$, let $N\defeq N_0N_1$.
For integers $C\ge 1$ and reals $\lambda>0$, let
\begin{equation*}
\mcal{S}_1(C,\lambda)\defeq \set{\bm{c}\in \mcal{S}_1: \norm{\bm{c}}\in [C,2C),\; \abs{\Delta(\bm{c}/C)}\in (\lambda/2,\lambda]}.
\end{equation*}
Assume \eqref{COND:clean-weight-condition-in-general}.
The definition \eqref{EQN:define-final-dyadic-modulus-weight-for-endgame} and Propositions~\ref{PROP:standard-general-Mellin-bound} and~\ref{PROP:uniform-discriminating-integral-bound-with-derivatives} imply, for $\bm{c}\in \mcal{S}_1(C,\lambda)$ and $b\in \ZZ_{\ge 0}$ (and multi-indices $\bm{\alpha}\ge 0$), that
\begin{equation}
\label{INEQ:Mellin-decay-for-g_c,X,N-for-endgame}
N^{1/2} \partial_{\bm{c}}^{\bm{\alpha}}{g_{\bm{c},X,\bm{N}}^\vee(i\bm{t})}
\ll_{F,w,\nu_2,b}
\frac{(X/N)^{\abs{\bm{\alpha}}} (N+XC)^{1-m/2}}{(1+\norm{\bm{t}})^b (1+C/X^{1/2})^b (1+XC\lambda/N)^b}.
\end{equation}

For each $\bm{c}\in \mcal{S}_1$, we have $\bm{c}\in \ZZ^m$ and $\Delta(\bm{c})\ne 0$, so $\norm{\bm{c}}\ge 1$ and $1\le \abs{\Delta(\bm{c})}\ll \norm{\bm{c}}^{\deg{\Delta}}$.
Let $\mcal{D} = \set{1,2,4,8,\ldots}$ and $\mcal{D}_2 = \mcal{D}\cup \set{1/d: d\in \mcal{D}}$; then we deduce that
\begin{equation}
\label{INEQ:decompose-S_1-dyadically}
\mcal{S}_1 \belongs \bigcup_{C\in \mcal{D},\; \lambda\in \mcal{D}_2:\, 1/C^{\deg{\Delta}}\le \lambda\le A_{13}} \mcal{S}_1(C,\lambda)
\end{equation}
for some $A_{13}=A_{13}(F)\ge 1$.
Also, the volume bound \eqref{INEQ:multivarite-near-zero-real-density} implies (for all $C\ge 1$ and $\lambda\le A_{13}$)
\begin{equation}
\label{INEQ:lambda-decay-of-S_1(C,lambda)}
\card{\mcal{S}_1(C,\lambda)}
\ll_F C^m (\lambda + C^{-1})^{1/\deg{\Delta}},
\end{equation}
because $\abs{\Delta(\bm{z})}\ll_F \abs{\Delta(\bm{y})} + C^{-1}$ for all $\bm{y}\in [-2,2]^m$ and $\bm{z}\in \bm{y} + [-C^{-1}, C^{-1}]^m$.

\subsection{Sharp delta bounds}
\label{SUBSEC:sharp-upper-bounds-via-R2'}

In \cite{wang2023_large_sieve_diagonal_cubic_forms},
we used Cauchy--Schwarz on $1/L$, $\Phi L$ over $\bm{c}$ to conditionally prove \eqref{INEQ:near-optimal-diagonal-GRH-bound-over-smooth-locus-S_1} (for diagonal $F$ with $2\mid m$) under a large-sieve hypothesis,
based on the first-order approximation $\Phi\approx 1/L$ (see \eqref{EQN:define-Psi_1,Psi_2}).
Now, in \S\ref{SUBSEC:sharp-upper-bounds-via-R2'},
we will use H\"{o}lder in a similar spirit to prove Theorem~\ref{THM:level-1-positive-density}.
This relies crucially on several new features, including the more precise $\Phi$-data captured in $\Sigma_{17,\bm{N}}^{\bm{c},1}(\bm{t})$ (compared to $1/L$ in \cite{wang2023_large_sieve_diagonal_cubic_forms}).

For the next four results, let $C\in \mcal{D}$, let $\bm{N}\in \mscr{R}_{17}(X)$, let $\bm{t}\in \RR^2$, and let $\mcal{S}\belongs \mcal{S}_1\cap [-C,C]^m$.
Roughly speaking, \eqref{INEQ:X^eps-loss-very-small-S-range-correlation-bound-goal} will be useful when $C\le X^{1/2-\delta}$;
and otherwise, \eqref{INEQ:non-dominant-N_0-range-correlation-bound-goal} will be useful for small $N_0$ (relative to $XC/N$), and \eqref{INEQ:N_0-dominant-correlation-bound-goal} useful for larger $N_0$.

\begin{lemma}
Assume Conjecture~\ref{CNJ:(R2'E')}.
Then for some $A_{14}=A_{14}(F,\eps)>0$, we have
\begin{equation}
\label{INEQ:non-dominant-N_0-range-correlation-bound-goal}
\norm{\Sigma_{17,\bm{N}}^{\bm{c},0}(\bm{t})
\cdot \Sigma_{17,\bm{N}}^{\bm{c},1}(\bm{t})}_{\ell^1_{\bm{c}}(\mcal{S})}
\ll_{F,\eps} (1+\abs{t_1})^{A_{14}}
\cdot \card{\mcal{S}}^{1/2-\eps} (C^m + N_1^{m/3})^{1/2+\eps} N_0^{m/2+\eps} N^{1/2}.
\end{equation}
\end{lemma}

\begin{proof}
Plugging the trivial bound $\abs{S^\natural_{\bm{c}}(n_0)} \le n_0^{(1+m)/2}$ into \eqref{EQN:define-smooth-dyadic-pieces-Sigma_17-for-endgame} gives
\begin{equation}
\label{INEQ:trivial-bound-for-Sigma_17}
\Sigma_{17,\bm{N}}^{\bm{c},0}(\bm{t}) \ll_{\nu_2} N_0^{(1+m)/2}\cdot \card{\mcal{N}_{\bm{c}}\cap [N_0, 2N_0)}
\end{equation}
for all $\bm{c}\in \mcal{S}$.
So by Lemma~\ref{LEM:N_c-small-divisor-moment-bound} (with $A = (2-\eps)/\eps$), Conjecture~\ref{CNJ:(R2'E')} (with $C+N_1^{1/3}$, $N_1$ in place of $Z$, $N$), and H\"{o}lder over $\bm{c}\in \mcal{S}$, the left-hand of \eqref{INEQ:non-dominant-N_0-range-correlation-bound-goal} is
\begin{equation*}
\ll_{F,\eps} N_0^{(1+m)/2} \cdot (C^m N_0^\eps)^{\eps/(2-\eps)}
\cdot \card{\mcal{S}}^{(1-2\eps)/(2-\eps)}
\cdot [(1+\abs{t_1})^{A_6(F,\eps)} (C^m+N_1^{m/3}) N_1^{(2-\eps)/2}]^{1/(2-\eps)}
\end{equation*}
for all $\eps\in (0,\frac12)$, since $1 = \frac{\eps}{2-\eps} + \frac{1-2\eps}{2-\eps} + \frac{1}{2-\eps}$.
Writing $N_0^{(1+m)/2} N_1^{1/2} = N_0^{m/2} N^{1/2}$ gives \eqref{INEQ:non-dominant-N_0-range-correlation-bound-goal},
since $\frac{1-2\eps}{2-\eps} \ge \frac12 - O(\eps)$ and $\card{\mcal{S}}\le C^m\le C^m+N_1^{m/3}$.
\end{proof}

\begin{lemma}
Assume Conjectures~\ref{CNJ:(HW2)} and~\ref{CNJ:(B2)}.
Then
\begin{equation}
\label{INEQ:X^eps-loss-very-small-S-range-correlation-bound-goal}
\norm{\Sigma_{17,\bm{N}}^{\bm{c},0}(\bm{t})
\cdot \Sigma_{17,\bm{N}}^{\bm{c},1}(\bm{t})}_{\ell^1_{\bm{c}}(\mcal{S})}
\ll_{F,\eps} (1+\abs{t_1})^\eps
\cdot \card{\mcal{S}}^{1/2} (C^m + N_0^{m/3})^{1/2+\eps} N^{1/2+\eps}.
\end{equation}
\end{lemma}

\begin{proof}
By Proposition~\ref{PROP:HW2-consequences}(\ref{ITEM:1/L-bound-from-GRH}) and partial summation, $\abs{\Sigma_{17,\bm{N}}^{\bm{c},1}(\bm{t})} \ll_{\nu_2,\eps} C^\eps (1+\abs{t_1})^\eps N_1^{1/2+\eps}$.
By Cauchy--Schwarz and \eqref{EQN:define-Sigma_14-general-B3G-type-sum}, the left-hand side of \eqref{INEQ:X^eps-loss-very-small-S-range-correlation-bound-goal} is $\ll_{F,\eps} [N_0\cdot \Sigma_{14}^{1,2}(C,N_0)]^{1/2}
\cdot \card{\mcal{S}}^{1/2} \cdot C^\eps (1+\abs{t_1})^\eps N_1^{1/2+\eps}$.
Now use Conjecture~\ref{CNJ:(B2)} (with $C+N_0^{1/3}$, $N_0$ in place of $Z$, $N$).
\end{proof}

\begin{proposition}
\label{PROP:B3G+}
Assume Conjecture~\ref{CNJ:(B3G)}
with some $W$.
Then for any real $\delta\ge 0$, we have
\begin{equation}
\label{INEQ:B3G+-goal}
\sum_{\bm{c}\in \mcal{S}:\, W(\bm{c})\ne 0}
\abs{
N_0^{-1/2} \cdot \Sigma_{17,\bm{N}}^{\bm{c},0}(\bm{t})
}^{2+\delta}
\ll_{F,\delta}
(C^m+N_0^{m/3}) N_0^{m\delta/2 - 9\eta_{10}(F,2)/10}.
\end{equation}
\end{proposition}

\begin{proof}
We may assume $\mcal{S}\belongs \set{W\ne 0}$.
By H\"{o}lder over $\mcal{S}$, the left-hand side of \eqref{INEQ:B3G+-goal} is
\begin{equation*}
\le \norm{
(N_0^{-1/2} \cdot \Sigma_{17,\bm{N}}^{\bm{c},0}(\bm{t}))^{\delta+\eps}
}_{\ell^{2/\eps}_{\bm{c}}(\mcal{S})}
\cdot \norm{
(N_0^{-1/2} \cdot \Sigma_{17,\bm{N}}^{\bm{c},0}(\bm{t}))^{2-\eps}
}_{\ell^{2/(2-\eps)}_{\bm{c}}(\mcal{S})},
\end{equation*}
for any $\eps\in [0,\delta]$, since $\frac{\eps}{2} + \frac{2-\eps}{2} = 1$.
The first factor here is $\ll_{\delta,\eps} (N_0^{m/2})^{\delta+\eps} (C^m N_0^\eps)^{\eps/2}$ by \eqref{INEQ:trivial-bound-for-Sigma_17} and Lemma~\ref{LEM:N_c-small-divisor-moment-bound};
the second factor is $\ll_\eps [(C+N_0^{1/3})^m N_0^{-\eta_{10}(F,2)}]^{(2-\eps)/2}$ by Conjecture~\ref{CNJ:(B3G)}.
Taking $\eps$ sufficiently small in terms of $m$, $\eta_{10}(F,2)$ gives \eqref{INEQ:B3G+-goal}.
\end{proof}

\begin{lemma}
Fix $\xi\in \RR_{\ge 0}$.
Assume Conjectures~\ref{CNJ:(HW2)}, \ref{CNJ:(B2)}, and~\ref{CNJ:(B3G)}.
Also assume Conjecture~\ref{CNJ:(R2'E')} if $\xi=0$.
Assume $N \ll C^6$.
Then for some $A_{15}=A_{15}(F,\eps)>0$, we have
\begin{equation}
\label{INEQ:N_0-dominant-correlation-bound-goal}
\norm{\Sigma_{17,\bm{N}}^{\bm{c},0}(\bm{t})
\cdot \Sigma_{17,\bm{N}}^{\bm{c},1}(\bm{t})}_{\ell^1_{\bm{c}}(\mcal{S})}
\ll_{F,\eps} (1+\abs{t_1})^{A_{15}}
\cdot \frac{(CN_1)^\xi (C^m)^{1/2-\eps}
(C^m+N^{m/3})^{1/2+\eps} N^{1/2}}{\min(C^{9/20}, N_0^{4\eta_{10}(F,2)/5})}.
\end{equation}
\end{lemma}

\begin{proof}
\emph{Case~1: $\mcal{S}\belongs \set{W=0}$.}
Then $\card{\mcal{S}} \ll_W C^{m-1}$ (see e.g.~\cite{bhargava2014geometric}*{Lemma~3.1}).
Inserting this into \eqref{INEQ:X^eps-loss-very-small-S-range-correlation-bound-goal}, we get \eqref{INEQ:N_0-dominant-correlation-bound-goal} upon writing $C^{(m-1)/2} N^\eps \ll_\eps C^{(m-1)/2+6\eps}$ (for a small $\eps$).

\emph{Case~2: $\mcal{S}\belongs \set{W\ne 0}$.}
Since $1 = \frac{1-\eps}{2-\eps} + \frac{1}{2-\eps}$,
we may use \eqref{INEQ:B3G+-goal}
(with $\delta=\frac{\eps}{1-\eps}$),
Conjecture~\ref{CNJ:(R2'E')} (if $\xi=0$) or GRH (if $\xi>0$, using \eqref{INEQ:second-Holder-for-R2'E'} and Proposition~\ref{PROP:HW2-consequences}(\ref{ITEM:1/L-bound-from-GRH}) to prove Conjecture~\ref{CNJ:(R2'E')} up to a factor of $Z^\eps$), and H\"{o}lder over $\bm{c}\in \mcal{S}$ to bound the left-hand side of \eqref{INEQ:N_0-dominant-correlation-bound-goal} by $O_{F,\eps}(N^{1/2} (CN_1)^\xi)
\cdot [(C^m+N_0^{m/3}) N_0^{m\delta/2 - 9\eta_{10}(F,2)/10}]^{(1-\eps)/(2-\eps)}
\cdot [(1+\abs{t_1})^{A_6(F,\eps)} (C^m+N_1^{m/3})]^{1/(2-\eps)}$.
Since $C^m+N_j^{m/3} \ll C^m$ for some $j\in \set{0,1}$, we get \eqref{INEQ:N_0-dominant-correlation-bound-goal} if $\eps$ is small.

\emph{Case~3: The general case.}
Decompose $\mcal{S}$ as $(\mcal{S}\cap \set{W=0})\cup (\mcal{S}\setminus \set{W=0})$.
\end{proof}

\begin{theorem}
\label{THM:axiomatic-level-1-endgame-theorem}
Fix $\xi\in \RR_{\ge 0}$.
Assume \eqref{COND:clean-weight-condition-in-general} and $2\mid m$.
Assume Conjectures~\ref{CNJ:(HW2)}, \ref{CNJ:(B2)}, and~\ref{CNJ:(B3G)}.
Also assume Conjecture~\ref{CNJ:(R2')} if $\xi=0$.
Then $\Sigma^\natural(X,\mcal{S}_1)\ll_{F,w} X^{(6-m)/4+\xi}$; in fact,
\begin{equation}
\label{INEQ:absolute-log-free-bound-in-delta-method}
\sum_{\bm{c}\in \mcal{S}_1}\,
\Bigl\lvert{
\sum_{n\ge 1} n^{(1-m)/2} S^\natural_{\bm{c}}(n) J_{\bm{c},X}(n)
}\Bigr\rvert
\ll_{F,w} X^{(6-m)/4+\xi}.
\end{equation}
Moreover, for any $X\in \RR_{\ge 1}$ and $P_0, P \in \RR_{>0}$, we have (for some constant $\eta_{11}=\eta_{11}(F)>0$)
\begin{equation}
\label{INEQ:Sigma_gen-bound-with-P_0,P-parameters-integrated-absolutely-over-N}
\int_{\bm{N}\in \mscr{R}_{17,P_0,P}(X)} d^\times\bm{N} \int_{\bm{t}\in \RR^2} d\bm{t}
\sum_{\bm{c}\in \mcal{S}_1}
\abs{g_{\bm{c},X,\bm{N}}^\vee(i\bm{t})}
\prod_{0\le j\le 1} \abs{\Sigma_{17,\bm{N}}^{\bm{c},j}(\bm{t})}
\ll_{F,w}
\frac{X^{(6-m)/4+\xi}}{\min(P_0 P, X)^{\eta_{11}}},
\end{equation}
where $\mscr{R}_{17,P_0,P}(X)\defeq
\set{\bm{N}\in \mscr{R}_{17}(X): N_0\geq P_0/2,\; N\leq A_0Y/P}$.
\end{theorem}

\begin{proof}
By \eqref{EQN:smooth-dyadic-restriction-separation-for-endgame} and \eqref{EQN:dyadic-range-restrictions-for-R_17}, the bound \eqref{INEQ:Sigma_gen-bound-with-P_0,P-parameters-integrated-absolutely-over-N} (with $P_0=1$ and $P=1$) implies \eqref{INEQ:absolute-log-free-bound-in-delta-method}.
And by \eqref{EQN:Fubini-Sigma(X,S)-for-endgame}, the bound \eqref{INEQ:absolute-log-free-bound-in-delta-method} implies $\Sigma^\natural(X,\mcal{S}_1)\ll_{F,w} X^{(6-m)/4+\xi}$.
It remains to prove \eqref{INEQ:Sigma_gen-bound-with-P_0,P-parameters-integrated-absolutely-over-N}; for this, we may assume $P_0\ge 1$ and $P\ge 1$ (since if $P_0<1$ or $P<1$, we may increase $P_0$ or $P$ while keeping the set $\mscr{R}_{17,P_0,P}(X)$ the same).
To begin,
we decompose $\mcal{S}_1$ using \eqref{INEQ:decompose-S_1-dyadically},
and then plug in \eqref{INEQ:Mellin-decay-for-g_c,X,N-for-endgame} (with $\bm{\alpha}=\bm{0}$);
this bounds the left-hand side of \eqref{INEQ:Sigma_gen-bound-with-P_0,P-parameters-integrated-absolutely-over-N} by $O_{F,w,b}(1)$ times
\begin{equation}
\int_{\mscr{R}_{17,P_0,P}(X)} d^\times\bm{N} \int_{\RR^2} d\bm{t}
\sum_{C\in \mcal{D},\; \lambda\in \mcal{D}_2(C)}
f(\bm{N},\bm{t},C,\lambda),
\end{equation}
where we write (for convenience) $\mcal{D}_2(C) = \set{\lambda\in \mcal{D}_2:
1/C^{\deg{\Delta}}\le \lambda\le A_{13}}$ and
\begin{equation}
\label{EQN:define-key-dyadic-upper-bound-endgame-quantity-f}
f(\bm{N},\bm{t},C,\lambda) = \frac{\norm{\Sigma_{17,\bm{N}}^{\bm{c},0}(\bm{t})
\cdot \Sigma_{17,\bm{N}}^{\bm{c},1}(\bm{t})}_{\ell^1_{\bm{c}}(\mcal{S}_1(C,\lambda))}
\cdot N^{-1/2} (N+XC)^{1-m/2}}{(1+\norm{\bm{t}})^b (1+C/X^{1/2})^b (1+XC\lambda/N)^b}.
\end{equation}

Next, note that if $\xi=0$, then Conjecture~\ref{CNJ:(R2'E')} holds by Propositions~\ref{PROP:R2'-implies-R2'E} and~\ref{PROP:R2'E-implies-R2'E'}.
So \eqref{INEQ:non-dominant-N_0-range-correlation-bound-goal}, \eqref{INEQ:X^eps-loss-very-small-S-range-correlation-bound-goal}, \eqref{INEQ:N_0-dominant-correlation-bound-goal} are all at our disposal, no matter what $\xi$ is.
Let $\theta = \frac{1}{10m}$ and $\eps = \frac{\theta}{100m}$, say.
Suppose $\bm{N}\in \mscr{R}_{17}(X)$ and $b \ge 3+\max(A_{14}(F,\eps), \eps, A_{15}(F,\eps))$.
Let $C\in \mcal{D}$ and $\lambda\in \mcal{D}_2(C)$.

Plugging \eqref{INEQ:lambda-decay-of-S_1(C,lambda)} into \eqref{INEQ:non-dominant-N_0-range-correlation-bound-goal} (for $\mcal{S} = \mcal{S}_1(C,\lambda)$) and integrating \eqref{EQN:define-key-dyadic-upper-bound-endgame-quantity-f} over $\bm{t}$ gives
\begin{equation*}
\int_{\RR^2} d\bm{t}\, f(\bm{N},\bm{t},C,\lambda)
\ll \frac{[C^m (\lambda+C^{-1})^{1/\deg{\Delta}}]^{1/2-\eps}
(C^m + N_1^{m/3})^{1/2+\eps} N_0^{m/2+\eps}}
{(N+XC)^{m/2-1} (1+C/X^{1/2})^b (1+XC\lambda/N)^b}.
\end{equation*}
Let $\alpha = (1/2-\eps)/\deg{\Delta}$.
Summing over $\lambda\in \mcal{D}_2(C)$ (using Lemma~\ref{LEM:general-dyadic-sum-split-into-2-geometric-series} with $r=\lambda$, $\tau=\frac{XC}{N}$, $a=0$, and $q\in \set{\alpha, 0}$, noting that $0\le \alpha<b$ and $(\lambda+C^{-1})^\alpha\ll \lambda^\alpha+C^{-\alpha}$), we get
\begin{equation*}
\sum_{\lambda\in \mcal{D}_2(C)}
\int_{\RR^2} d\bm{t}\, f
\ll \frac{C^{(1/2-\eps)m} [\min(\frac{N}{XC}, 1)^\alpha + C^{-(1-\eps)\alpha}]
(C^m + N_1^{m/3})^{1/2+\eps} N_0^{m/2+\eps}}
{(N+XC)^{m/2-1} (1+C/X^{1/2})^b}.
\end{equation*}
Writing $\min(\frac{N}{XC}, 1)\le \frac{N}{XC}$ and $(N+XC)^{m/2-1}\ge (XC)^{m/2-1}$,
and then summing over $C\in \mcal{D}$ using Lemma~\ref{LEM:general-dyadic-sum-split-into-2-geometric-series} (with $r=C$ and $\tau=X^{-1/2}$,
noting that $(\frac12-\eps)m > \alpha + (\frac{m}{2}-1)$),
we get
\begin{equation*}
\sum_{C\in \mcal{D},\; \lambda\in \mcal{D}_2(C)}
\int_{\RR^2} d\bm{t}\, f
\ll \frac{[(\frac{N}{Y})^\alpha + X^{-(1-\eps)\alpha/2}]
(X^{m/2} + X^{(1/2-\eps)m/2} N_1^{(1/2+\eps)m/3}) N_0^{m/2+\eps}}
{Y^{m/2-1}},
\end{equation*}
since $Y = X\cdot X^{1/2}$.
But $N_1\ll X^{3/2}$ (by \eqref{EQN:dyadic-range-restrictions-for-R_17}),
so the right-hand side is
$\ll [(\tfrac{N}{Y})^{\alpha} + X^{-(1-\eps)\alpha/2}]
\cdot N_0^{m/2+\eps} \cdot X^{(6-m)/4}$.
Letting $\mscr{R}_{18}\defeq \set{\bm{N}\in \mscr{R}_{17,P_0,P}(X): N_0^{m+\eps}\le \min(\tfrac{Y}{N}, X^{(1-\eps)/2})^{\alpha/2}}$, and writing $N_0^{m/2+\eps}\le N_0^{-m/2}\cdot \min(\tfrac{Y}{N}, X^{(1-\eps)/2})^{\alpha/2}$ for each $\bm{N}\in \mscr{R}_{18}$, we get
\begin{equation*}
\int_{\mscr{R}_{18}} d^\times\bm{N}
\sum_{C\in \mcal{D},\; \lambda\in \mcal{D}_2(C)}
\int_{\RR^2} d\bm{t}\, f
\ll \int_{1/4}^{A_0Y/P} d^\times{N}\,
[(\tfrac{N}{Y})^{\alpha/2} + X^{-(1-\eps)\alpha/4}] P_0^{-m/2} X^{(6-m)/4},
\end{equation*}
by integrating over $N_0\ge P_0/2$ for each fixed valued of $N_0N_1=N$.
Thus
\begin{equation}
\label{INEQ:final-non-dominant-N_0-endgame-bound}
\int_{\mscr{R}_{18}} d^\times\bm{N}
\sum_{C\in \mcal{D},\; \lambda\in \mcal{D}_2(C)}
\int_{\RR^2} d\bm{t}\, f
\ll [P^{-\alpha/2} + X^{-(1-2\eps)\alpha/4}] P_0^{-m/2} X^{(6-m)/4}.
\end{equation}

On the other hand, discarding the factor $(1+XC\lambda/N)^b$ in \eqref{EQN:define-key-dyadic-upper-bound-endgame-quantity-f} (using $1+XC\lambda/N\ge 1$),
plugging $\card{\mcal{S}}\ll C^m$ into \eqref{INEQ:X^eps-loss-very-small-S-range-correlation-bound-goal} (for $\mcal{S} = \bigcup_{\lambda\in \mcal{D}_2(C)} \mcal{S}_1(C,\lambda)$),
and integrating over $\bm{t}$ gives
\begin{equation*}
\sum_{\lambda\in \mcal{D}_2(C)}
\int_{\RR^2} d\bm{t}\, f(\bm{N},\bm{t},C,\lambda)
\ll \frac{C^{m/2} (C^m + N_0^{m/3})^{1/2+\eps} N^\eps}
{(N+XC)^{m/2-1} (1+C/X^{1/2})^b}.
\end{equation*}
Writing $(N+XC)^{m/2-1}\ge (XC)^{m/2-1}$,
and then summing over $1\le C\le X^{1/2-\theta}$, we get
\begin{equation*}
\sum_{C\in \mcal{D}:\, 1\le C\le X^{1/2-\theta}}
\sum_{\lambda\in \mcal{D}_2(C)}
\int_{\RR^2} d\bm{t}\, f
\ll \frac{(X^{1/2-\theta})^{m/2} ((X^{1/2-\theta})^m + N_0^{m/3})^{1/2+\eps} N^\eps}
{(X^{3/2-\theta})^{m/2-1}}.
\end{equation*}
But $N_0^{m/3}\ll X^{m/2}$ by \eqref{EQN:dyadic-range-restrictions-for-R_17};
integrating the previous display over $\mscr{R}_{17}(X)$ thus gives
\begin{equation}
\label{INEQ:final-small-C-endgame-bound}
\int_{\mscr{R}_{17}(X)} d^\times\bm{N}
\sum_{C\in \mcal{D}:\, 1\le C\le X^{1/2-\theta}}
\sum_{\lambda\in \mcal{D}_2(C)}
\int_{\RR^2} d\bm{t}\, f
\ll \frac{X^{(6-m)/4} X^{(m/2+2)\eps}}{X^\theta} \le \frac{X^{(6-m)/4}}{X^{9\theta/10}}.
\end{equation}

Now suppose $C>X^{1/2-\theta}$; then $C>X^{1/2-1/10}=X^{2/5} \gg N^{1/6}$.
Discarding $(1+XC\lambda/N)^b$ in \eqref{EQN:define-key-dyadic-upper-bound-endgame-quantity-f},
taking $\mcal{S} = \bigcup_{\lambda\in \mcal{D}_2(C)} \mcal{S}_1(C,\lambda)$ in \eqref{INEQ:N_0-dominant-correlation-bound-goal},
and integrating over $\bm{t}$ gives
\begin{equation*}
\sum_{\lambda\in \mcal{D}_2(C)}
\int_{\RR^2} d\bm{t}\, f(\bm{N},\bm{t},C,\lambda)
\ll \frac{(CN_1)^\xi \cdot C^{(1/2-\eps)m} (C^m+N^{m/3})^{1/2+\eps} \cdot (N+XC)^{1-m/2}}
{\min(C^{9/20}, N_0^{4\eta_{10}(F,2)/5}) (1+C/X^{1/2})^b}.
\end{equation*}
Writing $(N+XC)^{1-m/2}\le (XC)^{1-m/2}$, and summing over $C$ using Lemma~\ref{LEM:general-dyadic-sum-split-into-2-geometric-series},
gives
\begin{equation*}
\sum_{C\in \mcal{D}:\, C>X^{1/2-\theta}}
\sum_{\lambda\in \mcal{D}_2(C)}
\int_{\RR^2} d\bm{t}\, f
\ll \frac{(X^{1/2} N_1)^\xi \cdot X^{(1/2-\eps)m/2} (X^{m/2}+N^{m/3})^{1/2+\eps} \cdot Y^{1-m/2}}
{\min(X^{9/50}, N_0^{4\eta_{10}(F,2)/5})}.
\end{equation*}
Here $N_1,N\ll X^{3/2}$, so the right-hand side is
$\ll X^{2\xi} \cdot X^{(6-m)/4}
\cdot (X^{-9/50} + N_0^{-4\eta_{10}(F,2)/5})$.
Let $\mscr{R}_{19}\defeq \mscr{R}_{17,P_0,P}(X)\setminus \mscr{R}_{18}$.
Each $\bm{N}\in \mscr{R}_{19}$ satisfies $N_0^{m+\eps} > \min(\tfrac{Y}{N}, X^{(1-\eps)/2})^{\alpha/2}$, and thus $N_0^{-\eta_{10}(F,2)} < \min(\tfrac{Y}{N}, X^{(1-\eps)/2})^{-\beta} \le (\tfrac{N}{Y})^\beta + X^{-(1-\eps)\beta/2}$, where $\beta = \frac{\eta_{10}(F,2) \alpha}{2(m+\eps)}$.
Thus the integral $\int_{\mscr{R}_{19}} d^\times\bm{N}$ of (the left-hand side of) the previous display is
\begin{equation*}
\ll \int_{\mscr{R}_{17,P_0,P}(X)} d^\times\bm{N}\,
X^{2\xi} \cdot X^{(6-m)/4}
\cdot (X^{-9/50} + N_0^{-2\eta_{10}(F,2)/5} [(\tfrac{N}{Y})^{2\beta/5} + X^{-(1-\eps)\beta/5}]),
\end{equation*}
which is $\ll (X^{\eps-9/50} + P_0^{-2\eta_{10}(F,2)/5}[P^{-2\beta/5} + X^{-(1-2\eps)\beta/5}]) \cdot X^{(6-m)/4+2\xi}$ (by integrating first over $N_0\ge P_0/2$ when $N$ is fixed, and then integrating over $N\le A_0Y/P$).
This, when combined with \eqref{INEQ:final-non-dominant-N_0-endgame-bound} and \eqref{INEQ:final-small-C-endgame-bound}, establishes \eqref{INEQ:Sigma_gen-bound-with-P_0,P-parameters-integrated-absolutely-over-N} with $\eta_{11} = \min(\frac{\alpha}{2}, \frac{0.9\alpha}{4}, \frac{9\theta}{10}, \frac{17}{100}, \frac{2\eta_{10}(F,2)}{5}, \frac{2\beta}{5}, \frac{0.9\beta}{5})$ (where $\theta = \frac{1}{10m}$, $\alpha \ge \frac{0.4}{\deg{\Delta}}$, and $\beta \ge \frac{\eta_{10}(F,2)}{10m}$), after replacing $\xi$ with $\xi/2$ if $\xi>0$.
\end{proof}

\begin{remark}
Our use of H\"{o}lder above is fairly uniform (each of \eqref{INEQ:non-dominant-N_0-range-correlation-bound-goal}, \eqref{INEQ:X^eps-loss-very-small-S-range-correlation-bound-goal}, \eqref{INEQ:N_0-dominant-correlation-bound-goal} being based on ``approximately Cauchy--Schwarz''), but the input required could maybe be slightly relaxed by applying H\"{o}lder with more varied exponents.
For instance, if $C\approx X^{1/2}$ and $N_0\ge X^\delta$, one could work with $\Sigma_{17,\bm{N}}^{\bm{c},0}(\bm{t})$ in $\ell^1$ and $\Sigma_{17,\bm{N}}^{\bm{c},1}(\bm{t})$ in $\ell^\infty$ (the $N_0^\delta$ saving from the former drowning out the $X^\eps$ loss from the latter); cf.~\eqref{INEQ:N_0-dominant-correlation-bound-goal}.
And if $N_0\le X^\delta$, one might hope to work with $\Sigma_{17,\bm{N}}^{\bm{c},1}(\bm{t})$ in $\ell^1$ over $\bm{c}$ in a residue class to modulus $n_0\asymp N_0$ (with some nontrivial archimedean restrictions on $\bm{c}$), and then work with $\Sigma_{17,\bm{N}}^{\bm{c},0}(\bm{t})$ in $\ell^1$ afterwards.
\end{remark}

\begin{proof}
[Proof of Theorem~\ref{THM:level-1-positive-density}]
Let $\mcal{D}=\set{1,2,4,8,\ldots}$.
Let $T_{0,X}(\theta) \defeq \sum_{\abs{x}\leq X}e(\theta x^3)$ for $\theta\in \RR$; then
\begin{equation}
\label{EQN:Fermat-Weyl-sum-integral-circle-method}
N_F(X) = \int_{[0,1]} d\theta\, \abs{T_{0,X}(\theta)}^6
= \norm{T_{0,X}(\theta)^6}_{L^1_\theta([0,1])}.
\end{equation}
Now let $T_{1,X}(\theta) \defeq \sum_{\abs{x}\in (X/2,X]} e(\theta x^3)$.
Let $\norm{f(\lambda)}_{\ell^q_\lambda(\mcal{D})} \defeq (\sum_{\lambda\in \mcal{D}} \abs{f(\lambda)}^q)^{1/q}$ for $q\ge 1$; then
\begin{equation}
\label{INEQ:Holder-to-bound-T_0-via-T_1}
T_{0,X}(\theta)
= 1 + \sum_{\lambda\in \mcal{D}} T_{1,X/\lambda}(\theta)
\ll 1
+ \norm{(X/\lambda)^\rho}_{\ell^{6/5}_\lambda(\mcal{D})}
\cdot \norm{(X/\lambda)^{-\rho}T_{1,X/\lambda}(\theta)}_{\ell^6_\lambda(\mcal{D})}
\end{equation}
(by H\"{o}lder),
where $\rho=1/100$, say.
Since $\norm{(X/\lambda)^\rho}_{\ell^{6/5}_\lambda(\mcal{D})}\ll X^\rho$, we deduce
(upon inserting \eqref{INEQ:Holder-to-bound-T_0-via-T_1} into \eqref{EQN:Fermat-Weyl-sum-integral-circle-method}) that if $\norm{T_{1,X/\lambda}(\theta)^6}_{L^1_\theta([0,1])}\ll (X/\lambda)^3$ holds, then (since $6\rho<3$)
\begin{equation}
\label{INEQ:bound-6th-moment-of-T_0-via-that-of-T_1}
N_F(X)
\ll 1
+ X^{6\rho}
\cdot \norm{(X/\lambda)^{-6\rho}T_{1,X/\lambda}(\theta)^6}_{L^1_{(\lambda,\theta)}(\mcal{D}\times [0,1])}
\ll X^3.
\end{equation}

Now let $w_\star(\bm{x})\defeq\prod_{1\le j\le 6} \bm{1}_{\abs{x_j}\in(1/2, 1]}$,
and choose $w\in C^\infty_c((\RR\setminus\set{0})^6)$ with $w\geq w_\star$.
Then
\begin{equation}
\label{INEQ:bound-6th-moment-of-T_1-via-smooth-clean-weight}
\norm{T_{1,X}(\theta)^6}_{L^1_\theta([0,1])}=N_{F,w_\star}(X)\leq N_{F,w}(X).
\end{equation}
But $w$ satisfies \eqref{COND:clean-weight-condition-in-diagonal-case} (and thus \eqref{COND:clean-weight-condition-in-general}).
And we are assuming (for Theorem~\ref{THM:level-1-positive-density}) Conjectures~\ref{CNJ:(HW2)}, \ref{CNJ:(R2')}, and~\ref{CNJ:(SFSCp)}; in particular, Conjectures~\ref{CNJ:(B2)} and~\ref{CNJ:(B3G)} hold by Propositions~\ref{PROP:B2-for-diagonal-F} and~\ref{PROP:SFSCp-implies-B3G}, respectively.
So Theorem~\ref{THM:axiomatic-level-1-endgame-theorem} (with $\xi=0$) gives $\Sigma^\natural(X,\mcal{S}_1)\ll 1$.
But $\Sigma^\natural(X,\mcal{S}_0)\ll 1$ by Theorem~\ref{THM:Sigma(X,S_0)-theorem-black-box}.
So by \eqref{EQN:delta-method} (and the definitions \eqref{EQN:define-sum-Sigma(X,S)}, \eqref{EQN:define-S_0,S_1-for-Delta-vanishing-and-nonzero-loci}), we have
\begin{equation*}
N_{F,w}(X)/X^3 = O_A(Y^{-A}) + \Sigma^\natural(X,\mcal{S}_0)
+ \Sigma^\natural(X,\mcal{S}_1) \ll 1.
\end{equation*}
Hence $N_F(X)\ll X^3$ by \eqref{INEQ:bound-6th-moment-of-T_1-via-smooth-clean-weight}, \eqref{INEQ:bound-6th-moment-of-T_0-via-that-of-T_1}.
Given \eqref{INEQ:main-level-1-goal-eps-free-point-count}, a standard Cauchy--Schwarz argument then leads to the desired application to $F_0(S^3)$ (producing a ``positive lower density'').
\end{proof}

\subsection{Delta cancellation}
\label{SUBSEC:cancellation-via-RA1}

To prove Theorems~\ref{THM:level-2-almost-all-integers} and~\ref{THM:level-3-power-saving}, we will complement \eqref{INEQ:Sigma_gen-bound-with-P_0,P-parameters-integrated-absolutely-over-N} using Propositions~\ref{PROP:RA1o'E-implies-RA1o'E'} and~\ref{PROP:RA1delta'E-implies-RA1delta'E'}, identifying $\bm{c}$-cancellation in some pieces of $\Sigma^\natural(X,\mcal{S}_1)$.
Since $J_{\bm{c},X}(n)$ is not compactly supported in $\bm{c}$,
we begin with a decomposition resembling \eqref{INEQ:decomposition-of-x-integral-into-nice-pieces}.

Recall $\nu_0$, $\nu_1$ from \S\ref{SUBSEC:conventions}.
Let $\ol{\nu}_0\defeq 1-\nu_0$.
For all $\kappa\in \RR_{>0}$, let (cf.~\eqref{EQN:define-grad-localized-x-integral-J_0}, \eqref{EQN:define-grad-localized-x-integral-J_1})
\begin{align}
\Sigma_{18,0}(X,\bm{N},\bm{t}) &\defeq
\sum_{\bm{c}\in \mcal{S}_1}
\nu_0(\bm{c}/X^{1/2})
g_{\bm{c},X,\bm{N}}^\vee(i\bm{t})
\prod_{0\le j\le 1} \Sigma_{17,\bm{N}}^{\bm{c},j}(\bm{t}),
\label{EQN:define-Sigma_18,0-endgame-piece} \\
\Sigma_{18,1,\kappa}(X,\bm{N},\bm{t}) &\defeq
\sum_{\bm{c}\in \mcal{S}_1}
\ol{\nu}_0(\bm{c}/X^{1/2}) \nu_1(\bm{c}/\kappa)
g_{\bm{c},X,\bm{N}}^\vee(i\bm{t})
\prod_{0\le j\le 1} \Sigma_{17,\bm{N}}^{\bm{c},j}(\bm{t}).
\label{EQN:define-Sigma_18,1-endgame-piece}
\end{align}
We have $\Sigma_{18,1,\kappa}(X,\bm{N},\bm{t}) = 0$ unless there exists $\bm{c}\in \mcal{S}_1$ satisfying $X^{1/2}/2 < \norm{\bm{c}}\le m\kappa$.
Using $\int_{\kappa>0} d^\times{\kappa}\,\nu_1(\bm{c}/\kappa) = 1$ (valid for $\bm{c}\in \RR^m\setminus \set{\bm{0}}$), we may thus write (cf.~\eqref{INEQ:decomposition-of-x-integral-into-nice-pieces})
\begin{equation}
\label{INEQ:decomposition-of-N-t-endgame-piece-into-Sigma_18-pieces}
\sum_{\bm{c}\in \mcal{S}_1}
g_{\bm{c},X,\bm{N}}^\vee(i\bm{t})
\prod_{0\le j\le 1} \Sigma_{17,\bm{N}}^{\bm{c},j}(\bm{t})
= \Sigma_{18,0}(X,\bm{N},\bm{t})
+ \int_{X^{1/2}/2m}^{\infty} d^\times{\kappa}\, \Sigma_{18,1,\kappa}(X,\bm{N},\bm{t}).
\end{equation}

Recall $\mscr{R}_{17}(X)$ from \eqref{EQN:dyadic-range-restrictions-for-R_17}.
Let $\mscr{R}_{17}^{P_0,P}(X)\defeq
\set{\bm{N}\in \mscr{R}_{17}(X): N_0 < P_0/2,\; N > A_0Y/P}$
for reals $P_0,P>0$.
In terms of $\mscr{R}_{17,P_0,P}(X)$ from Theorem~\ref{THM:axiomatic-level-1-endgame-theorem}, we have
\begin{equation}
\label{INEQ:decompose-region-R_17-for-endgame}
\mscr{R}_{17}^{P_0,P}(X)
\belongs \mscr{R}_{17}(X)
\belongs \mscr{R}_{17}^{P_0,P}(X)
\cup \mscr{R}_{17,1,P}(X) \cup \mscr{R}_{17,P_0,1}(X).
\end{equation}
By \eqref{EQN:Fubini-Sigma(X,S)-for-endgame}, \eqref{EQN:smooth-dyadic-restriction-separation-for-endgame}, \eqref{EQN:dyadic-range-restrictions-for-R_17}, and \eqref{INEQ:decompose-region-R_17-for-endgame}, the bound \eqref{INEQ:Sigma_gen-bound-with-P_0,P-parameters-integrated-absolutely-over-N} (when applicable) implies
\begin{equation}
\label{INEQ:approximate-Sigma(X,S_1)-by-small-N_0,large-N-range}
\Sigma^\natural(X,\mcal{S}_1) - \int_{\mscr{R}_{17}^{P_0,P}(X)} d^\times\bm{N}
\int_{\RR^2} d\bm{t}
\sum_{\bm{c}\in \mcal{S}_1}
g_{\bm{c},X,\bm{N}}^\vee(i\bm{t})
\prod_{0\le j\le 1} \Sigma_{17,\bm{N}}^{\bm{c},j}(\bm{t})
\ll \frac{X^{(6-m)/4+\xi}}{\min(P_0, P, X)^{\eta_{11}}}.
\end{equation}
This leads to the following results.

\begin{theorem}
\label{THM:axiomatic-level-2-endgame-theorem}
Assume \eqref{COND:clean-weight-condition-in-general} and $2\mid m$.
Assume Conjectures~\ref{CNJ:(HW2)}, \ref{CNJ:(R2')}, \ref{CNJ:(RA1o)},
\ref{CNJ:(B2)}, and~\ref{CNJ:(B3G)}.
Then for $X\ge 1$, we have
$\Sigma^\natural(X,\mcal{S}_1) = o_{F,w;X\to\infty}(X^{(6-m)/4})$.
\end{theorem}

\begin{proof}
Before proceeding, note that Conjecture~\ref{CNJ:(RA1o'E')} holds by Propositions~\ref{PROP:R2'-implies-R2'E}, \ref{PROP:RA1o-implies-RA1o'E}, and~\ref{PROP:RA1o'E-implies-RA1o'E'}, since we assume Conjectures~\ref{CNJ:(HW2)}, \ref{CNJ:(R2')}, and~\ref{CNJ:(RA1o)}.
Let $M$ be a real number to be chosen later.

Let $\bm{N}\in \mscr{R}_{17}^{P_0,P}(X)$.
For each real $Z\ge 2X^{1/2}$, let (in the context of \eqref{EQN:define-Sigma_11-for-RA1'E'})
\begin{equation*}
\nu_{18,0}=\nu_{18,0,X,\bm{N},\bm{t},Z}
= \nu_0(Z\bm{c}/X^{1/2}) g_{Z\bm{c},X,\bm{N}}^\vee(i\bm{t}) \nu_2(r)
\end{equation*}
(noting that $\Supp{\nu_0}\belongs [-2,2]^m$, so $\nu_{18,0}$ is supported on $[-1,1]^m \times [1,2]$).
Then after plugging \eqref{EQN:define-smooth-dyadic-pieces-Sigma_17-for-endgame} into \eqref{EQN:define-Sigma_18,0-endgame-piece} and decomposing $\mcal{S}_1$ into residue classes modulo $n_0$, we get
\begin{equation*}
\Sigma_{18,0}(X,\bm{N},\bm{t})
= \sum_{n_0\ge 1} \nu_2(n_0/N_0) n_0^{-it_0}
\sum_{1\le \bm{a}\le n_0:\, n_0\in \mcal{N}_{\bm{a}}}
S^\natural_{\bm{a}}(n_0) \cdot \Sigma_{11}^{\bm{a},n_0}(\nu_{18,0},Z,N_1)
\end{equation*}
(in terms of $\Sigma_{11}^{\bm{a},n_0}$ from \eqref{EQN:define-Sigma_11-for-RA1'E'}).
By Conjecture~\ref{CNJ:(RA1o'E')} (with $\nu = \nu_{18,0}$, and with $2X^{1/2}+N_1^{1/3}$, $N_1$ in place of $Z$, $N$) and the trivial bound $\abs{S^\natural_{\bm{a}}(n_0)} \le n_0^{(1+m)/2}$, we get
\begin{equation*}
\Sigma_{18,0}(X,\bm{N},\bm{t}) \ll N_0^{(3+m)/2} X^{m/2} N_1^{1/2}
(o_{M\to\infty}(1) + o_{M;X\to\infty}(1)) \cdot \mcal{M}_{1,A_7}(\nu_{18,0}),
\end{equation*}
provided $\min(2X^{1/2}, N_1)\ge 2M\ge 2(2N_0)^{10/9}$.
Additionally, by \eqref{INEQ:Mellin-decay-for-g_c,X,N-for-endgame} (with $Z\bm{c} = (2X^{1/2}+N_1^{1/3})\bm{c}$ in place of $\bm{c}$) and the definition \eqref{EQN:define-norm-M_1,k(nu)} of $\mcal{M}_{1,k}$ (with $\nu=\nu_{18,0}$, $k=A_7$), we have
\begin{equation*}
\mcal{M}_{1,A_7}(N^{1/2} \nu_{18,0}) \ll_b \frac{(X^{1/2}/X^{1/2}+X\cdot X^{1/2}/N)^1 \cdot (N+0)^{1-m/2}}{(1+\norm{\bm{t}})^b (1+0)^b (1+0)^b}
\ll \frac{(Y/N)\cdot N^{1-m/2}}{(1+\norm{\bm{t}})^b}.
\end{equation*}

Now, for each $\kappa\ge X^{1/2}/2m$ and $Z\ge m\kappa$ (noting that $\Supp{\nu_1}\belongs [-m,m]^m$), let
\begin{equation*}
\nu_{18,1,\kappa}=\nu_{18,1,\kappa,X,\bm{N},\bm{t},Z}
= \ol{\nu}_0(Z\bm{c}/X^{1/2}) \nu_1(Z\bm{c}/\kappa)
g_{Z\bm{c},X,\bm{N}}^\vee(i\bm{t}) \nu_2(r).
\end{equation*}
Applying Conjecture~\ref{CNJ:(RA1o'E')} with $m\kappa+N_1^{1/3}$, $N_1$ in place of $Z$, $N$, we get
\begin{equation*}
\Sigma_{18,1,\kappa}(X,\bm{N},\bm{t}) \ll N_0^{(3+m)/2} \kappa^m N_1^{1/2}
(o_{M\to\infty}(1) + o_{M;X\to\infty}(1)) \cdot \mcal{M}_{1,A_7}(\nu_{18,1,\kappa}),
\end{equation*}
provided $\min(X^{1/2}/2, N_1)\ge 2M\ge (2N_0)^{10/9}$.
This time, since $\bm{0}\notin \Supp{\nu_1}$, the bound \eqref{INEQ:Mellin-decay-for-g_c,X,N-for-endgame} (with $(m\kappa+N_1^{1/3})\bm{c}$ in place of $\bm{c}$) and \eqref{EQN:define-norm-M_1,k(nu)} (with $\nu=\nu_{18,1,\kappa}$, $k=A_7$) give
\begin{equation*}
\mcal{M}_{1,A_7}(N^{1/2} \nu_{18,1,\kappa}) \ll_b \frac{(\kappa/X^{1/2}+X\kappa/N)^1 \cdot (N+X\kappa)^{1-m/2}}{(1+\norm{\bm{t}})^b (1+\kappa/X^{1/2})^b (1+0)^b}
\ll \frac{(Y/N)\cdot (X\kappa)^{1-m/2}}{(1+\norm{\bm{t}})^b (\kappa/X^{1/2})^{b-1}}.
\end{equation*}

Inserting our bounds on $\Sigma_{18,0}$, $\Sigma_{18,1,\kappa}$ into \eqref{INEQ:decomposition-of-N-t-endgame-piece-into-Sigma_18-pieces},
assuming $b\ge 3+m/2$,
we find that if $\min(X^{1/2}/2, N_1)\ge 2M\ge (2N_0)^{10/9}$, then the left-hand side of \eqref{INEQ:decomposition-of-N-t-endgame-piece-into-Sigma_18-pieces} is
\begin{equation}
\label{INEQ:soft-bound-on-Sigma_18-total}
\ll N_0^{(2+m)/2} X^{m/2}
(o_{M\to\infty}(1) + o_{M;X\to\infty}(1))
\cdot (Y/N)\cdot N^{1-m/2} / (1+\norm{\bm{t}})^b.
\end{equation}
Since $N_0<P_0/2$ and $N_1=N/N_0>A_0Y/P_0P$, we conclude (upon integrating over $\bm{t}\in \RR^2$, then over $N_0<P_0/2$ with $N_0N_1=N$ fixed, and finally over $A_0Y/P<N\le A_0Y$) that
\begin{equation*}
\int_{\mscr{R}_{17}^{P_0,P}(X)} d^\times\bm{N}
\int_{\RR^2} d\bm{t}
\sum_{\bm{c}\in \mcal{S}_1}
g_{\bm{c},X,\bm{N}}^\vee(i\bm{t})
\prod_{0\le j\le 1} \Sigma_{17,\bm{N}}^{\bm{c},j}(\bm{t})
\ll \frac{P_0^{(2+m)/2} X^{m/2} o_{M\to\infty}(1) \cdot Y}{(Y/P)^{m/2}},
\end{equation*}
provided $X\gg_{M,P_0,P} 1$ and $M\gg_{P_0} 1$ hold with large enough implied constants.
Therefore, there exist functions $f_3, f_4\maps \RR_{>0}\to \ZZ_{\ge 1}$ such that if $M\ge f_3(P_0+P)$ and $X\ge f_4(M)$, then the left-hand side of the previous display is $o_{P_0+P\to \infty}(Y^{1-m/2} X^{m/2})$.
It follows from \eqref{INEQ:approximate-Sigma(X,S_1)-by-small-N_0,large-N-range} (with $\xi=0$) that if $X\ge f_4(f_3(P_0+P))$, then
$\Sigma^\natural(X,\mcal{S}_1)
= o_{\min(P_0, P, X)\to \infty}(X^{(6-m)/4})$.
Finally, suppose $X\ge f_4(f_3(2))$, let $P$ be the largest integer in $[1,X]$ with $f_4(f_3(2P))\le X$ (so that $P\to \infty$ as $X\to \infty$), and take $P_0=P$, to get $\Sigma^\natural(X,\mcal{S}_1) = o_{X\to \infty}(X^{(6-m)/4})$.
\end{proof}

\begin{proof}
[Proof of Theorem~\ref{THM:level-2-almost-all-integers}]

Unconditionally,
by \eqref{EQN:define-HLH-error-E_w(X)}, \eqref{EQN:delta-method}, and \eqref{EQN:disc-locus-S_0-evaluation}, we have
\begin{equation}
\label{EQN:approximate-E_w(X)-by-Sigma(X,S_1)}
E_{F,w}(X)/X^3 = O(X^{2.75+\eps})/X^3 + \Sigma^\natural(X,\mcal{S}_1).
\end{equation}
Now assume \eqref{COND:clean-weight-condition-in-diagonal-case} (so \eqref{COND:clean-weight-condition-in-general} holds).
Then Theorem~\ref{THM:axiomatic-level-2-endgame-theorem} gives $\Sigma^\natural(X,\mcal{S}_1) = o_{X\to \infty}(X^{(6-m)/4})$,
since Conjectures~\ref{CNJ:(B2)} and~\ref{CNJ:(B3G)} hold by Propositions~\ref{PROP:B2-for-diagonal-F} and~\ref{PROP:SFSCp-implies-B3G}, respectively (since we assume Conjecture~\ref{CNJ:(SFSCp)}).
Upon plugging this bound into \eqref{EQN:approximate-E_w(X)-by-Sigma(X,S_1)}, we get \eqref{EQN:soft-HLH-general-homogeneous-weight}.

The Hasse principle for $V$ follows from \eqref{EQN:soft-HLH-general-homogeneous-weight},
upon choosing $w$ with $\sigma_{\infty,F,w}>0$ (possible since $V(\RR)$ contains a point with $x_1\cdots x_6\ne 0$).
Now suppose $F=x_1^3+\dots+x_6^3$.
Then by \cite{wang2023prime}*{Theorem~1.1} (or \cite{wang2022thesis}*{Theorem~2.1.8}), we find (from \eqref{EQN:soft-HLH-general-homogeneous-weight}) that
$F_0(\ZZ^3)$ has density $1$ in $\set{a\in \ZZ: a\not\equiv \pm 4 \bmod{9}}$.
(See \cite{wang2023prime} for details on this last deduction, which is based on \cite{diaconu2019admissible}.
Diaconu assumes an analog of \eqref{EQN:soft-HLH-general-homogeneous-weight} over rather quantitatively deformed regions $R^\ast$,
whereas we work with fixed weights $w$.
It would be very interesting to see if there could be any
miraculous cancellation or symmetries in the analog of $J_{\bm{c},X}(n)$ over $R^\ast$, but at the moment it seems easier to handle minimally deformed regions.)
\end{proof}

\begin{proof}
[Proof of Corollary~\ref{COR:remove-Hessian-assumption-for-free-in-soft-Fermat-case}]
(Here we drop the assumption \eqref{COND:clean-weight-condition-in-general}.)
Let $\rho>0$ be a parameter tending to $0$ slowly as $X\to \infty$.
Use \eqref{INEQ:main-level-1-goal-eps-free-point-count} and H\"{o}lder's inequality to upper bound the contribution to $N_{F,w}(X)$ from points $\bm{x}\in \ZZ^m$ with $\min(\abs{x_1},\dots,\abs{x_m})\le \rho \cdot X$.
Then use Theorem~\ref{THM:level-2-almost-all-integers} to estimate the remaining contribution to $N_{F,w}(X)$.
This gives \eqref{EQN:soft-HLH-general-homogeneous-weight} for arbitrary $w\in C^\infty_c(\RR^m)$.
Hooley's conjecture follows upon taking a suitable sequence of weights $w$.
\end{proof}

\begin{theorem}
\label{THM:axiomatic-level-3-endgame-theorem}
Assume \eqref{COND:clean-weight-condition-in-general} and $2\mid m$.
Assume Conjectures~\ref{CNJ:(HW2)}, \ref{CNJ:(RA1delta)}, \ref{CNJ:(EKL)},
\ref{CNJ:(B2)}, and~\ref{CNJ:(B3G)}.
Then for some constant $\eta_{12}=\eta_{12}(F)>0$, we have
$\Sigma^\natural(X,\mcal{S}_1) \ll_{F,w} X^{(6-m)/4-\eta_{12}}$.
\end{theorem}

\begin{proof}
Proposition~\ref{PROP:RA1delta'E-implies-RA1delta'E'} applies, since we assume Conjectures~\ref{CNJ:(HW2)}, \ref{CNJ:(RA1delta)}, \ref{CNJ:(EKL)}.
We now mimic the proof of Theorem~\ref{THM:axiomatic-level-2-endgame-theorem}, while replacing each use of Conjecture~\ref{CNJ:(RA1o'E')} with Proposition~\ref{PROP:RA1delta'E-implies-RA1delta'E'}.

Let $\bm{N}\in \mscr{R}_{17}^{P_0,P}(X)$.
Proposition~\ref{PROP:RA1delta'E-implies-RA1delta'E'} implies that for any real $M\in [1, X^{\eta_2/2}]$ satisfying $\min(2X^{1/2}, N_1)\ge 2M\ge 2(2N_0)^2$, we have
\begin{equation*}
\Sigma_{18,0}(X,\bm{N},\bm{t}) \ll_\eps N_0^{(3+m)/2} X^{m/2+\eps} N_1^{1/2}
(M^{1/6\deg{H}} N_1^{-1/6} + M^{-\eta_5/4\deg{H}}) \mcal{M}_{1,A_8}(\nu_{18,0}).
\end{equation*}
Furthermore, for each $\kappa\ge X^{1/2}/2m$, Proposition~\ref{PROP:RA1delta'E-implies-RA1delta'E'} implies that
\begin{equation*}
\Sigma_{18,1,\kappa}(X,\bm{N},\bm{t}) \ll_\eps N_0^{(3+m)/2} \kappa^{m+\eps} N_1^{1/2}
(M^{1/6\deg{H}} N_1^{-1/6} + M^{-\eta_5/4\deg{H}}) \mcal{M}_{1,A_8}(\nu_{18,1,\kappa}).
\end{equation*}
for any real $M\in [1, (X^{1/2}/2)^{\eta_2}]$ satisfying $\min(X^{1/2}/2, N_1)\ge 2M\ge 2(2N_0)^2$.

Now let $M = X^{\min(1,\eta_2)/2}/(4+2^{\eta_2})$, and suppose
\begin{equation}
\label{INEQ:condition-for-P_0,P-small-hard-endgame-in-terms-of-M}
(A_0Y/P_0P)^{1/2}\ge 2M\ge 2P_0^2,
\end{equation}
so that $N_1\gg M^2$ and (thus) $M^{1/6\deg{H}} N_1^{-1/6} + M^{-\eta_5/4\deg{H}} \ll M^{-\min(1,\eta_5)/6\deg{H}}$.
Arguing as we did for \eqref{INEQ:soft-bound-on-Sigma_18-total},
we find (assuming $b\ge 3+m/2$) that the left-hand side of \eqref{INEQ:decomposition-of-N-t-endgame-piece-into-Sigma_18-pieces} is
\begin{equation*}
\ll_\eps N_0^{(2+m)/2} X^{m/2+\eps}
(M^{-\min(1,\eta_5)/6\deg{H}})
\cdot (Y/N)\cdot N^{1-m/2}/(1+\norm{\bm{t}})^b.
\end{equation*}
So (upon integrating over $\bm{t}$, then over $N_0$, and finally over $N$)
\begin{equation*}
\int_{\mscr{R}_{17}^{P_0,P}(X)} d^\times\bm{N}
\int_{\RR^2} d\bm{t}
\sum_{\bm{c}\in \mcal{S}_1}
g_{\bm{c},X,\bm{N}}^\vee(i\bm{t})
\prod_{0\le j\le 1} \Sigma_{17,\bm{N}}^{\bm{c},j}(\bm{t})
\ll_\eps \frac{P_0^{(2+m)/2} X^{m/2+\eps}\cdot Y}{M^{\min(1,\eta_5)/6\deg{H}}\cdot (Y/P)^{m/2}}.
\end{equation*}
Now let $\gamma = \min(1,\eta_5)/12\deg{H}\in (0, 1/12]$ and $P_0 = M^{\gamma/(2+m)}$, $P = A_0 M^{\gamma/m}/4$.
Then \eqref{INEQ:condition-for-P_0,P-small-hard-endgame-in-terms-of-M} holds
(since $\gamma\le 1/2$ and $Y\ge M^3$),
and the right-hand side of the previous display is $\ll M^{\gamma/2+\gamma/2-2\gamma} X^{m/2+\eps} Y^{1-m/2} = M^{-\gamma} X^{(6-m)/4+\eps}$.
So by \eqref{INEQ:approximate-Sigma(X,S_1)-by-small-N_0,large-N-range} (with $\xi=\eps>0$) we have
\begin{equation*}
\Sigma^\natural(X,\mcal{S}_1) \ll_\eps X^{(6-m)/4+\eps} (M^{-\gamma} + P_0^{-\eta_{11}})
\ll X^{(6-m)/4-\eta_{12}},
\end{equation*}
where $\eta_{12} = \frac{9}{10}\cdot \min(1, \eta_{11}/(2+m))\cdot \gamma \cdot \min(1,\eta_2)/2$, say.
\end{proof}

\begin{proof}
[Proof of Theorem~\ref{THM:level-3-power-saving}]
Proceed as in the proof of Theorem~\ref{THM:axiomatic-level-2-endgame-theorem}, but use Theorem~\ref{THM:axiomatic-level-3-endgame-theorem} instead of Theorem~\ref{THM:axiomatic-level-2-endgame-theorem}, to get (from \eqref{EQN:approximate-E_w(X)-by-Sigma(X,S_1)}) the bound $E_{F,w}(X)/X^3 \ll_\eps X^{-0.25+\eps} + X^{-\eta_{12}}$.
\end{proof}

    
    



\section*{Acknowledgements}

Many thanks to Amit Ghosh and Peter Sarnak for suggesting the main problem addressed here.
I also thank my advisor, Peter Sarnak, for his guidance and encouragement.{\let\thefootnote\relax\footnote{This work was partially supported by NSF grant DMS-1802211.}}
I thank Calvin Deng and Yotam Hendel for sharing references for \eqref{INEQ:univariate-zero-density-mod-q} and \eqref{INEQ:univariate-near-zero-real-density}, and Yotam for discussions related to \S\ref{SEC:new-bounds-on-bad-sums-S}.
For early comments, I thank
Manjul Bhargava,
Andy Booker,
Tim Browning,
Brian Conrey,
Simona Diaconu,
Bill Duke,
Roger Heath-Brown,
Nick Katz,
Will Sawin,
Bob Vaughan,
and Trevor Wooley.
Thanks also to everyone listed in \cite{wang2022thesis},
and to people from seminars, conferences, and other events.
For more recent interactions, I thank
Louis-Pierre Arguin,
Emma Bailey,
Paul Bourgade,
Alex Gamburd,
Jayce Getz,
Jeff Hoffstein,
Junehyuk Jung,
Victor Kolyvagin,
Valeriya Kovaleva,
Lillian Pierce,
Kannan Soundararajan,
Yuri Tschinkel,
Katy Woo,
Max Xu,
Liyang Yang,
and Peter Zenz.
Finally, I thank my family for their exceptional support.

\bibliographystyle{amsxport}
\bibliography{final.bib}

\begin{bibdiv}
\begin{biblist}

\bib{bettin2020averages}{article}{
      author={Bettin, Sandro},
      author={Conrey, J.~B.},
       title={Averages of long {D}irichlet polynomials},
        date={2020-06-22},
     journal={Preprint},
      eprint={arXiv:2002.09466v2},
}

\bib{bui2023negative}{article}{
      author={Bui, H.~M.},
      author={Florea, A.},
       title={Negative moments of the {R}iemann zeta-function},
        date={2023-02-14},
     journal={Preprint},
      eprint={arXiv:2302.07226v1},
}

\bib{bui2021ratios}{article}{
      author={Bui, H.~M.},
      author={Florea, A.},
      author={Keating, J.~P.},
       title={The {R}atios {C}onjecture and upper bounds for negative moments
  of {$L$}-functions over function fields},
        date={2021-09-21},
     journal={Preprint},
      eprint={arXiv:2109.10396v1},
}

\bib{bateman1958theorem}{article}{
      author={Bateman, Paul~T.},
      author={Grosswald, Emil},
       title={On a theorem of {E}rd\"{o}s and {S}zekeres},
        date={1958},
        ISSN={0019-2082},
     journal={Illinois J. Math.},
      volume={2},
       pages={88\ndash 98},
         url={http://projecteuclid.org/euclid.ijm/1255380836},
      review={\MR{95804}},
}

\bib{bhargava2014geometric}{article}{
      author={Bhargava, Manjul},
       title={The geometric sieve and the density of squarefree values of
  invariant polynomials},
        date={2014-01-31},
     journal={Preprint},
      eprint={arXiv:1402.0031v1},
}

\bib{bhargava2021galois}{article}{
      author={Bhargava, Manjul},
       title={Galois groups of random integer polynomials and {van der
  Waerden's Conjecture}},
        date={2022-10-03},
     journal={Preprint},
      eprint={arXiv:2111.06507v2},
}

\bib{buse2014discriminant}{article}{
      author={Bus\'{e}, Laurent},
      author={Jouanolou, Jean-Pierre},
       title={On the discriminant scheme of homogeneous polynomials},
        date={2014},
        ISSN={1661-8270},
     journal={Math. Comput. Sci.},
      volume={8},
      number={2},
       pages={175\ndash 234},
         url={https://doi.org/10.1007/s11786-014-0188-7},
      review={\MR{3224627}},
}

\bib{bombieri2006riemann}{incollection}{
      author={Bombieri, E.},
       title={The {R}iemann hypothesis},
        date={2006},
   booktitle={The millennium prize problems},
   publisher={Clay Math. Inst., Cambridge, MA},
       pages={107\ndash 124},
      review={\MR{2238277}},
}

\bib{booker2021question}{article}{
      author={Booker, Andrew~R.},
      author={Sutherland, Andrew~V.},
       title={On a question of {M}ordell},
        date={2021},
        ISSN={0027-8424},
     journal={Proc. Natl. Acad. Sci. USA},
      volume={118},
      number={11},
       pages={Paper No. 2022377118, 11},
         url={https://doi.org/10.1073/pnas.2022377118},
      review={\MR{4279690}},
}

\bib{castryck2020dimension}{article}{
      author={Castryck, Wouter},
      author={Cluckers, Raf},
      author={Dittmann, Philip},
      author={Nguyen, Kien~Huu},
       title={The dimension growth conjecture, polynomial in the degree and
  without logarithmic factors},
        date={2020},
        ISSN={1937-0652},
     journal={Algebra Number Theory},
      volume={14},
      number={8},
       pages={2261\ndash 2294},
         url={https://doi.org/10.2140/ant.2020.14.2261},
      review={\MR{4172708}},
}

\bib{vcech2022ratios}{article}{
      author={{\v{C}}ech, Martin},
       title={The {R}atios conjecture for real {D}irichlet characters and
  multiple {D}irichlet series},
        date={2022-01-05},
     journal={Preprint},
      eprint={arXiv:2110.04409v2},
}

\bib{conrey2005integral}{article}{
      author={Conrey, J.~B.},
      author={Farmer, D.~W.},
      author={Keating, J.~P.},
      author={Rubinstein, M.~O.},
      author={Snaith, N.~C.},
       title={Integral moments of {$L$}-functions},
        date={2005},
        ISSN={0024-6115},
     journal={Proc. London Math. Soc. (3)},
      volume={91},
      number={1},
       pages={33\ndash 104},
         url={https://doi.org/10.1112/S0024611504015175},
      review={\MR{2149530}},
}

\bib{conrey2008autocorrelation}{article}{
      author={Conrey, J.~B.},
      author={Farmer, D.~W.},
      author={Zirnbauer, M.~R.},
       title={Autocorrelation of ratios of {$L$}-functions},
        date={2008},
        ISSN={1931-4523},
     journal={Commun. Number Theory Phys.},
      volume={2},
      number={3},
       pages={593\ndash 636},
         url={https://doi.org/10.4310/CNTP.2008.v2.n3.a4},
      review={\MR{2482944}},
}

\bib{conrey2007applications}{article}{
      author={Conrey, J.~B.},
      author={Snaith, N.~C.},
       title={Applications of the {$L$}-functions ratios conjectures},
        date={2007},
        ISSN={0024-6115},
     journal={Proc. Lond. Math. Soc. (3)},
      volume={94},
      number={3},
       pages={594\ndash 646},
         url={https://doi.org/10.1112/plms/pdl021},
      review={\MR{2325314}},
}

\bib{colliot2012groupe}{article}{
      author={Colliot-Th\'{e}l\`ene, Jean-Louis},
      author={Wittenberg, Olivier},
       title={Groupe de {B}rauer et points entiers de deux familles de surfaces
  cubiques affines},
        date={2012},
        ISSN={0002-9327},
     journal={Amer. J. Math.},
      volume={134},
      number={5},
       pages={1303\ndash 1327},
         url={https://doi.org/10.1353/ajm.2012.0036},
      review={\MR{2975237}},
}

\bib{davenport1939waring}{article}{
      author={Davenport, H.},
       title={On {W}aring's problem for cubes},
        date={1939},
        ISSN={0001-5962},
     journal={Acta Math.},
      volume={71},
       pages={123\ndash 143},
         url={https://doi.org/10.1007/BF02547752},
      review={\MR{0000026}},
}

\bib{deligne1969constantes}{incollection}{
      author={Deligne, Pierre},
       title={Les constantes des \'{e}quations fonctionnelles},
        date={1970},
   booktitle={S\'{e}minaire {D}elange-{P}isot-{P}oitou: 1969/70, {T}h\'{e}orie
  des {N}ombres, {F}asc. 2},
   publisher={Secr\'{e}tariat math\'{e}matique, Paris},
       pages={Exp. 19 bis, 13},
      review={\MR{0280450}},
}

\bib{denef1984rationality}{article}{
      author={Denef, J.},
       title={The rationality of the {P}oincar\'{e} series associated to the
  {$p$}-adic points on a variety},
        date={1984},
        ISSN={0020-9910},
     journal={Invent. Math.},
      volume={77},
      number={1},
       pages={1\ndash 23},
         url={https://doi.org/10.1007/BF01389133},
      review={\MR{751129}},
}

\bib{duke1993bounds}{article}{
      author={Duke, W.},
      author={Friedlander, J.~B.},
      author={Iwaniec, H.},
       title={Bounds for automorphic {$L$}-functions},
        date={1993},
        ISSN={0020-9910},
     journal={Invent. Math.},
      volume={112},
      number={1},
       pages={1\ndash 8},
         url={https://doi.org/10.1007/BF01232422},
      review={\MR{1207474}},
}

\bib{diaconu2003multiple}{article}{
      author={Diaconu, Adrian},
      author={Goldfeld, Dorian},
      author={Hoffstein, Jeffrey},
       title={Multiple {D}irichlet series and moments of zeta and
  {$L$}-functions},
        date={2003},
        ISSN={0010-437X},
     journal={Compositio Math.},
      volume={139},
      number={3},
       pages={297\ndash 360},
         url={https://doi.org/10.1023/B:COMP.0000018137.38458.68},
      review={\MR{2041614}},
}

\bib{deshouillers2006density}{incollection}{
      author={Deshouillers, Jean-Marc},
      author={Hennecart, Fran\c{c}ois},
      author={Landreau, Bernard},
       title={On the density of sums of three cubes},
        date={2006},
   booktitle={Algorithmic number theory},
      series={Lecture Notes in Comput. Sci.},
      volume={4076},
   publisher={Springer, Berlin},
       pages={141\ndash 155},
         url={https://doi.org/10.1007/11792086_11},
      review={\MR{2282921}},
}

\bib{diaconu2019admissible}{thesis}{
      author={Diaconu, Simona},
       title={On admissible integers of cubic forms},
        type={Senior Thesis, Princeton University},
        date={2019},
        note={URL: \url{http://arks.princeton.edu/ark:/88435/dsp0112579w10h}},
}

\bib{ding2020variance}{article}{
      author={Ding, Pengyong},
       title={On a variance associated with the distribution of $r_3(n)$ in
  arithmetic progressions},
        date={2020-10-09},
     journal={Preprint},
      eprint={arXiv:2010.04319v1},
}

\bib{eisenbud1995commutative}{book}{
      author={Eisenbud, David},
       title={Commutative algebra},
      series={Graduate Texts in Mathematics},
   publisher={Springer-Verlag, New York},
        date={1995},
      volume={150},
        ISBN={0-387-94268-8; 0-387-94269-6},
         url={https://doi.org/10.1007/978-1-4612-5350-1},
        note={With a view toward algebraic geometry},
      review={\MR{1322960}},
}

\bib{ekedahl1991infinite}{article}{
      author={Ekedahl, Torsten},
       title={An infinite version of the {C}hinese remainder theorem},
        date={1991},
        ISSN={0010-258X},
     journal={Comment. Math. Univ. St. Paul.},
      volume={40},
      number={1},
       pages={53\ndash 59},
      review={\MR{1104780}},
}

\bib{florea2021negative}{article}{
      author={Florea, A.},
       title={Negative moments of {$L$}-functions with small shifts over
  function fields},
        date={2022-11-28},
     journal={Preprint},
      eprint={arXiv:2111.10477v2},
}

\bib{franke1989rational}{article}{
      author={Franke, Jens},
      author={Manin, Yuri~I.},
      author={Tschinkel, Yuri},
       title={Rational points of bounded height on {F}ano varieties},
        date={1989},
        ISSN={0020-9910},
     journal={Invent. Math.},
      volume={95},
      number={2},
       pages={421\ndash 435},
         url={https://doi.org/10.1007/BF01393904},
      review={\MR{974910}},
}

\bib{farmer2019analytic}{article}{
      author={Farmer, David~W.},
      author={Pitale, Ameya},
      author={Ryan, Nathan~C.},
      author={Schmidt, Ralf},
       title={Analytic {$L$}-functions: definitions, theorems, and
  connections},
        date={2019},
        ISSN={0273-0979},
     journal={Bull. Amer. Math. Soc. (N.S.)},
      volume={56},
      number={2},
       pages={261\ndash 280},
         url={https://doi.org/10.1090/bull/1646},
      review={\MR{3923345}},
}

\bib{ganzburg2001polynomial}{article}{
      author={Ganzburg, M.~I.},
       title={Polynomial inequalities on measurable sets and their
  applications},
        date={2001},
        ISSN={0176-4276},
     journal={Constr. Approx.},
      volume={17},
      number={2},
       pages={275\ndash 306},
         url={https://doi.org/10.1007/s003650010020},
      review={\MR{1814358}},
}

\bib{glas2022question}{article}{
      author={Glas, Jakob},
      author={Hochfilzer, Leonhard},
       title={On a question of {D}avenport and diagonal cubic forms over
  {$\FF_q(t)$}},
        date={2022-08-10},
     journal={Preprint},
      eprint={arXiv:2208.05422v1},
}

\bib{gorodetsky2021magic}{article}{
      author={Gorodetsky, Ofir},
       title={Magic squares and the symmetric group},
        date={2022-07-11},
     journal={Preprint},
      eprint={arXiv:2102.11966v3},
}

\bib{granville1998abc}{article}{
      author={Granville, Andrew},
       title={{$ABC$} allows us to count squarefrees},
        date={1998},
        ISSN={1073-7928},
     journal={Internat. Math. Res. Notices},
      number={19},
       pages={991\ndash 1009},
         url={https://doi.org/10.1155/S1073792898000592},
      review={\MR{1654759}},
}

\bib{ghosh2017integral}{article}{
      author={Ghosh, Amit},
      author={Sarnak, Peter},
       title={Integral points on {M}arkoff type cubic surfaces},
        date={2022},
        ISSN={0020-9910},
     journal={Invent. Math.},
      volume={229},
      number={2},
       pages={689\ndash 749},
         url={https://doi.org/10.1007/s00222-022-01114-z},
      review={\MR{4448994}},
}

\bib{hormander1990analysis}{book}{
      author={H\"{o}rmander, Lars},
       title={The analysis of linear partial differential operators. {I}},
     edition={Second},
      series={Grundlehren der mathematischen Wissenschaften [Fundamental
  Principles of Mathematical Sciences]},
   publisher={Springer-Verlag, Berlin},
        date={1990},
      volume={256},
        ISBN={3-540-52345-6},
         url={https://doi.org/10.1007/978-3-642-61497-2},
        note={Distribution theory and Fourier analysis},
      review={\MR{1065993}},
}

\bib{harper2013sharp}{article}{
      author={Harper, Adam~J.},
       title={Sharp conditional bounds for moments of the {R}iemann zeta
  function},
        date={2013-05-20},
     journal={Preprint},
      eprint={arXiv:1305.4618v1},
}

\bib{harper2023typical}{article}{
      author={Harper, Adam~J.},
       title={The typical size of character and zeta sums is $o(\sqrt{x})$},
        date={2023-01-11},
     journal={Preprint},
      eprint={arXiv:2301.04390v1},
}

\bib{heath1983cubic}{article}{
      author={Heath-Brown, D.~R.},
       title={Cubic forms in ten variables},
        date={1983},
        ISSN={0024-6115},
     journal={Proc. London Math. Soc. (3)},
      volume={47},
      number={2},
       pages={225\ndash 257},
         url={https://doi.org/10.1112/plms/s3-47.2.225},
      review={\MR{703978}},
}

\bib{heath1992density}{article}{
      author={Heath-Brown, D.~R.},
       title={The density of zeros of forms for which weak approximation
  fails},
        date={1992},
        ISSN={0025-5718},
     journal={Math. Comp.},
      volume={59},
      number={200},
       pages={613\ndash 623},
         url={https://doi.org/10.2307/2153078},
      review={\MR{1146835}},
}

\bib{heath1996new}{article}{
      author={Heath-Brown, D.~R.},
       title={A new form of the circle method, and its application to quadratic
  forms},
        date={1996},
        ISSN={0075-4102},
     journal={J. Reine Angew. Math.},
      volume={481},
       pages={149\ndash 206},
         url={https://doi.org/10.1515/crll.1996.481.149},
      review={\MR{1421949}},
}

\bib{heath1998circle}{article}{
      author={Heath-Brown, D.~R.},
       title={The circle method and diagonal cubic forms},
        date={1998},
        ISSN={1364-503X},
     journal={R. Soc. Lond. Philos. Trans. Ser. A Math. Phys. Eng. Sci.},
      volume={356},
      number={1738},
       pages={673\ndash 699},
         url={https://doi.org/10.1098/rsta.1998.0181},
      review={\MR{1620820}},
}

\bib{hooley2014octonary}{article}{
      author={Hooley, Christopher},
       title={On octonary cubic forms},
        date={2014},
        ISSN={0024-6115},
     journal={Proc. Lond. Math. Soc. (3)},
      volume={109},
      number={1},
       pages={241\ndash 281},
         url={https://doi.org/10.1112/plms/pdt066},
      review={\MR{3237742}},
}

\bib{hooley1986some}{article}{
      author={Hooley, Christopher},
       title={On some topics connected with {W}aring's problem},
        date={1986},
        ISSN={0075-4102},
     journal={J. Reine Angew. Math.},
      volume={369},
       pages={110\ndash 153},
         url={https://doi.org/10.1515/crll.1986.369.110},
      review={\MR{850631}},
}

\bib{hooley1986HasseWeil}{article}{
      author={Hooley, Christopher},
       title={On {W}aring's problem},
        date={1986},
        ISSN={0001-5962},
     journal={Acta Math.},
      volume={157},
      number={1-2},
       pages={49\ndash 97},
         url={https://doi.org/10.1007/BF02392591},
      review={\MR{857679}},
}

\bib{hooley1988nonary}{article}{
      author={Hooley, Christopher},
       title={On nonary cubic forms},
        date={1988},
        ISSN={0075-4102},
     journal={J. Reine Angew. Math.},
      volume={386},
       pages={32\ndash 98},
         url={https://doi.org/10.1515/crll.1988.386.32},
      review={\MR{936992}},
}

\bib{hooley_greaves_harman_huxley_1997}{incollection}{
      author={Hooley, Christopher},
       title={On {H}ypothesis {$K^\ast$} in {W}aring's problem},
        date={1997},
   booktitle={Sieve methods, exponential sums, and their applications in number
  theory ({C}ardiff, 1995)},
      series={London Math. Soc. Lecture Note Ser.},
      volume={237},
   publisher={Cambridge Univ. Press, Cambridge},
       pages={175\ndash 185},
         url={https://doi.org/10.1017/CBO9780511526091.013},
      review={\MR{1635754}},
}

\bib{huang2020density}{article}{
      author={Huang, Jing-Jing},
       title={The density of rational points near hypersurfaces},
        date={2020},
        ISSN={0012-7094},
     journal={Duke Math. J.},
      volume={169},
      number={11},
       pages={2045\ndash 2077},
         url={https://doi.org/10.1215/00127094-2020-0004},
      review={\MR{4132580}},
}

\bib{iwaniec2004analytic}{book}{
      author={Iwaniec, Henryk},
      author={Kowalski, Emmanuel},
       title={Analytic number theory},
      series={American Mathematical Society Colloquium Publications},
   publisher={American Mathematical Society, Providence, RI},
        date={2004},
      volume={53},
        ISBN={0-8218-3633-1},
         url={https://doi.org/10.1090/coll/053},
      review={\MR{2061214}},
}

\bib{kisin1999local}{article}{
      author={Kisin, Mark},
       title={Local constancy in {$p$}-adic families of {G}alois
  representations},
        date={1999},
        ISSN={0025-5874},
     journal={Math. Z.},
      volume={230},
      number={3},
       pages={569\ndash 593},
         url={https://doi.org/10.1007/PL00004706},
      review={\MR{1680032}},
}

\bib{kloosterman1926representation}{article}{
      author={Kloosterman, H.~D.},
       title={On the representation of numbers in the form
  {$ax^2+by^2+cz^2+dt^2$}},
        date={1926},
        ISSN={0001-5962},
     journal={Acta Math.},
      volume={49},
      number={3-4},
       pages={407\ndash 464},
         url={https://doi.org/10.1007/BF02564120},
      review={\MR{1555249}},
}

\bib{konjagin1979number}{article}{
      author={Konjagin, V.~S.},
       title={The number of solutions of congruences of the {$n$}th degree with
  one unknown},
        date={1979},
        ISSN={0368-8666},
     journal={Mat. Sb. (N.S.)},
      volume={109(151)},
      number={2},
       pages={171\ndash 187, 327},
      review={\MR{542556}},
}

\bib{laskar2017local}{article}{
      author={Laskar, Abhijit},
       title={On local monodromy filtrations attached to motives},
        date={2017},
        ISSN={1793-0421},
     journal={Int. J. Number Theory},
      volume={13},
      number={3},
       pages={801\ndash 817},
         url={https://doi.org/10.1142/S1793042117500427},
      review={\MR{3606955}},
}

\bib{li2022moments}{article}{
      author={Li, Xiannan},
       title={Moments of quadratic twists of modular {$L$}-functions},
        date={2022-08-21},
     journal={Preprint},
      eprint={arXiv:2208.07343v2},
}

\bib{miller2004one}{article}{
      author={Miller, Steven~J.},
       title={One- and two-level densities for rational families of elliptic
  curves: evidence for the underlying group symmetries},
        date={2004},
        ISSN={0010-437X},
     journal={Compos. Math.},
      volume={140},
      number={4},
       pages={952\ndash 992},
         url={https://doi.org/10.1112/S0010437X04000582},
      review={\MR{2059225}},
}

\bib{mordell1953integer}{article}{
      author={Mordell, L.~J.},
       title={On the integer solutions of the equation {$x^2+y^2+z^2+2xyz=n$}},
        date={1953},
        ISSN={0024-6107},
     journal={J. London Math. Soc.},
      volume={28},
       pages={500\ndash 510},
         url={https://doi.org/10.1112/jlms/s1-28.4.500},
      review={\MR{56619}},
}

\bib{ng2004distribution}{article}{
      author={Ng, Nathan},
       title={The distribution of the summatory function of the {M}\"{o}bius
  function},
        date={2004},
        ISSN={0024-6115},
     journal={Proc. London Math. Soc. (3)},
      volume={89},
      number={2},
       pages={361\ndash 389},
         url={https://doi.org/10.1112/S0024611504014741},
      review={\MR{2078705}},
}

\bib{pas1989uniform}{article}{
      author={Pas, Johan},
       title={Uniform {$p$}-adic cell decomposition and local zeta functions},
        date={1989},
        ISSN={0075-4102},
     journal={J. Reine Angew. Math.},
      volume={399},
       pages={137\ndash 172},
         url={https://doi.org/10.1515/crll.1989.399.137},
      review={\MR{1004136}},
}

\bib{peyre1995hauteurs}{article}{
      author={Peyre, Emmanuel},
       title={Hauteurs et mesures de {T}amagawa sur les vari\'{e}t\'{e}s de
  {F}ano},
        date={1995},
        ISSN={0012-7094},
     journal={Duke Math. J.},
      volume={79},
      number={1},
       pages={101\ndash 218},
         url={https://doi.org/10.1215/S0012-7094-95-07904-6},
      review={\MR{1340296}},
}

\bib{poonen2003squarefree}{article}{
      author={Poonen, Bjorn},
       title={Squarefree values of multivariable polynomials},
        date={2003},
        ISSN={0012-7094},
     journal={Duke Math. J.},
      volume={118},
      number={2},
       pages={353\ndash 373},
         url={https://doi.org/10.1215/S0012-7094-03-11826-8},
      review={\MR{1980998}},
}

\bib{saito2003weight}{article}{
      author={Saito, Takeshi},
       title={Weight spectral sequences and independence of {$l$}},
        date={2003},
        ISSN={1474-7480},
     journal={J. Inst. Math. Jussieu},
      volume={2},
      number={4},
       pages={583\ndash 634},
         url={https://doi.org/10.1017/S1474748003000173},
      review={\MR{2006800}},
}

\bib{scholze2012perfectoid}{article}{
      author={Scholze, Peter},
       title={Perfectoid spaces},
        date={2012},
        ISSN={0073-8301},
     journal={Publ. Math. Inst. Hautes \'{E}tudes Sci.},
      volume={116},
       pages={245\ndash 313},
         url={https://doi.org/10.1007/s10240-012-0042-x},
      review={\MR{3090258}},
}

\bib{swinnerton2001solubility}{article}{
      author={Swinnerton-Dyer, Peter},
       title={The solubility of diagonal cubic surfaces},
        date={2001},
        ISSN={0012-9593},
     journal={Ann. Sci. \'{E}cole Norm. Sup. (4)},
      volume={34},
      number={6},
       pages={891\ndash 912},
         url={https://doi.org/10.1016/S0012-9593(01)01080-1},
      review={\MR{1872424}},
}

\bib{segre1943note}{article}{
      author={Segre, B.},
       title={A note on arithmetical properties of cubic surfaces},
        date={1943},
        ISSN={0024-6107},
     journal={J. London Math. Soc},
      volume={18},
       pages={24\ndash 31},
         url={https://doi.org/10.1112/jlms/s1-18.1.24},
      review={\MR{0009471}},
}

\bib{serre1969facteurs}{incollection}{
      author={Serre, Jean-Pierre},
       title={Facteurs locaux des fonctions z\^{e}ta des variet\'{e}s
  alg\'{e}briques (d\'{e}finitions et conjectures)},
        date={1970},
   booktitle={S\'{e}minaire {D}elange-{P}isot-{P}oitou. 11e ann\'{e}e: 1969/70.
  {T}h\'{e}orie des nombres. {F}asc. 1: {E}xpos\'{e}s 1 \`a 15; {F}asc. 2:
  {E}xpos\'{e}s 16 \`a 24},
   publisher={Secr\'{e}tariat Math., Paris},
       pages={15},
      review={\MR{3618526}},
}

\bib{serre1981quelques}{article}{
      author={Serre, Jean-Pierre},
       title={Quelques applications du th\'{e}or\`eme de densit\'{e} de
  {C}hebotarev},
        date={1981},
        ISSN={0073-8301},
     journal={Inst. Hautes \'{E}tudes Sci. Publ. Math.},
      number={54},
       pages={123\ndash 201},
         url={http://archive.numdam.org/article/PMIHES_1981__54__123_0.pdf},
      review={\MR{644559}},
}

\bib{soundararajan2009moments}{article}{
      author={Soundararajan, Kannan},
       title={Moments of the {R}iemann zeta function},
        date={2009},
        ISSN={0003-486X},
     journal={Ann. of Math. (2)},
      volume={170},
      number={2},
       pages={981\ndash 993},
         url={https://doi.org/10.4007/annals.2009.170.981},
      review={\MR{2552116}},
}

\bib{sarnak2016families}{incollection}{
      author={Sarnak, Peter},
      author={Shin, Sug~Woo},
      author={Templier, Nicolas},
       title={Families of {$L$}-functions and their symmetry},
        date={2016},
   booktitle={Families of automorphic forms and the trace formula},
      series={Simons Symp.},
   publisher={Springer, [Cham]},
       pages={531\ndash 578},
      review={\MR{3675175}},
}

\bib{shin2014fields}{article}{
      author={Shin, Sug~Woo},
      author={Templier, Nicolas},
       title={On fields of rationality for automorphic representations},
        date={2014},
        ISSN={0010-437X},
     journal={Compos. Math.},
      volume={150},
      number={12},
       pages={2003\ndash 2053},
         url={https://doi.org/10.1112/S0010437X14007428},
      review={\MR{3292292}},
}

\bib{sarnak1991bounds}{article}{
      author={Sarnak, Peter},
      author={Xue, Xiao~Xi},
       title={Bounds for multiplicities of automorphic representations},
        date={1991},
        ISSN={0012-7094},
     journal={Duke Math. J.},
      volume={64},
      number={1},
       pages={207\ndash 227},
         url={https://doi.org/10.1215/S0012-7094-91-06410-0},
      review={\MR{1131400}},
}

\bib{taylor2004galois}{article}{
      author={Taylor, Richard},
       title={Galois representations},
        date={2004},
        ISSN={0240-2963},
     journal={Ann. Fac. Sci. Toulouse Math. (6)},
      volume={13},
      number={1},
       pages={73\ndash 119},
         url={http://afst.cedram.org/item?id=AFST_2004_6_13_1_73_0},
      review={\MR{2060030}},
}

\bib{tian2017hasse}{article}{
      author={Tian, Zhiyu},
       title={Hasse principle for three classes of varieties over global
  function fields},
        date={2017},
        ISSN={0012-7094},
     journal={Duke Math. J.},
      volume={166},
      number={17},
       pages={3349\ndash 3424},
         url={https://doi.org/10.1215/00127094-2017-0034},
      review={\MR{3724220}},
}

\bib{vaughan2020some}{article}{
      author={Vaughan, R.~C.},
       title={On some questions of partitio numerorum: Tres cubi},
        date={2020},
     journal={Glasgow Mathematical Journal},
       pages={1\ndash 22},
}

\bib{vaughan1995certain}{article}{
      author={Vaughan, R.~C.},
      author={Wooley, T.~D.},
       title={On a certain nonary cubic form and related equations},
        date={1995},
        ISSN={0012-7094},
     journal={Duke Math. J.},
      volume={80},
      number={3},
       pages={669\ndash 735},
         url={https://doi.org/10.1215/S0012-7094-95-08023-5},
      review={\MR{1370112}},
}

\bib{wang2021_HLH_vs_RMT}{article}{
      author={Wang, Victor~Y.},
       title={Approaching cubic {Diophantine} statistics via mean-value
  {$L$-function} conjectures of {Random Matrix Theory} type},
        date={2021-08-07},
     journal={Preprint},
      eprint={arXiv:2108.03398v1},
}

\bib{wang2022thesis}{thesis}{
      author={Wang, Victor~Y.},
       title={Families and dichotomies in the circle method},
        type={Ph.D. Thesis, Princeton University},
        date={2022},
        note={URL: \url{http://arks.princeton.edu/ark:/88435/dsp01rf55zb86g}},
}

\bib{wang2023_large_sieve_diagonal_cubic_forms}{article}{
      author={Wang, Victor~Y.},
       title={Diagonal cubic forms and the large sieve},
        date={2023-04-18},
     journal={Preprint},
      eprint={arXiv:2108.03395v2},
}

\bib{wang2023dichotomous}{article}{
      author={Wang, Victor~Y.},
       title={Dichotomous point counts over finite fields},
        date={2023},
     journal={J. Number Theory, to appear},
}

\bib{wang2023prime}{article}{
      author={Wang, Victor~Y.},
       title={Prime {Hasse} principles via {Diophantine} second moments},
        date={2023-04-18},
     journal={Preprint},
      eprint={arXiv:2304.08674v1},
}

\bib{wang2023_isolating_special_solutions}{article}{
      author={Wang, Victor~Y.},
       title={Special cubic zeros and the dual variety},
        date={2023-04-18},
     journal={Preprint},
      eprint={arXiv:2108.03396v2},
}

\bib{wooley2000sums}{article}{
      author={Wooley, T.~D.},
       title={Sums of three cubes},
        date={2000},
        ISSN={0025-5793},
     journal={Mathematika},
      volume={47},
      number={1-2},
       pages={53\ndash 61 (2002)},
         url={https://doi.org/10.1112/S0025579300015710},
      review={\MR{1924487}},
}

\bib{wooley2015sums}{article}{
      author={Wooley, T.~D.},
       title={Sums of three cubes, {II}},
        date={2015},
        ISSN={0065-1036},
     journal={Acta Arith.},
      volume={170},
      number={1},
       pages={73\ndash 100},
         url={https://doi.org/10.4064/aa170-1-6},
      review={\MR{3373831}},
}

\bib{wooley1995breaking}{article}{
      author={Wooley, T.~D.},
       title={Breaking classical convexity in {W}aring's problem: sums of cubes
  and quasi-diagonal behaviour},
        date={1995},
        ISSN={0020-9910},
     journal={Invent. Math.},
      volume={122},
      number={3},
       pages={421\ndash 451},
         url={https://doi.org/10.1007/BF01231451},
      review={\MR{1359599}},
}

\end{biblist}
\end{bibdiv}
\end{document}